\documentclass{conm-p-l}

\usepackage{amssymb}
\usepackage{amsfonts}
\usepackage{amscd}
\usepackage[mathscr]{eucal}
\usepackage{verbatim}
\usepackage{epsfig}

\newtheorem{thespecial}{Theorem}
\newtheorem{thm}{Theorem}[section]

\newtheorem{prop}[thm]{Proposition}

\newtheorem{sth}{Theorem}[thm]
\newtheorem{slem}[sth]{Lemma}
\newtheorem{sprop}[sth]{Proposition}
\newtheorem{scor}[sth]{Corollary}

\theoremstyle{definition}
 \newtheorem{parag}[thm]{}
 \newtheorem{sparag}[sth]{}
 \newtheorem{defi}[thm]{Definition}
 \newtheorem{sdef}[sth]{Definition}
 \newtheorem{srem}[sth]{Remark}
 \newtheorem{srems}[sth]{Remarks}
 
 \newtheorem{exams}[sth]{Examples}

\newcounter{numb}[section]
\newcounter{t}

\newcommand{\Cref}[1]{Corollary~\textup{\ref{#1}}}
\newcommand{\Dref}[1]{Definition~\textup{\ref{#1}}}
\newcommand{\Eref}[1]{Example~\textup{\ref{#1}}}
\newcommand{\Lref}[1]{Lemma~\textup{\ref{#1}}}
\newcommand{\Pref}[1]{Proposition~\textup{\ref{#1}}}
\newcommand{\Rref}[1]{Remark~\textup{\ref{#1}}}
\newcommand{\Sref}[1]{Section~\textup{\ref{#1}}}
\newcommand{\Tref}[1]{Theorem~\textup{\ref{#1}}}

\def\bilap#1{\hbox to 0pt{\hss#1\hss}}
 \def\Rarrow#1{\bilap{\hbox to#1{\rightarrowfill}}}
 \def\Larrow#1{\bilap{\hbox to#1{\leftarrowfill}}}
\def\Equals#1{\bilap
                  {\hbox{\rule[3.5pt]{#1} {.5pt}}
                   \kern-#1
                   \hbox{\rule[1pt]{#1}{.5pt}}%
                 }}

\def\UnderElement#1#2#3#4{\vbox to 0pt{
\hbox{$
\llap{$\scriptstyle#1$}
\left#2\vbox to #3{}\right.
\rlap{$\scriptstyle#4$}
     $}
\vss}}
\newcommand{\EQAL}[1]%
{\,\begin{picture}(#1,0)%
\put(0,3){\line(1,0){#1}}%
\put(0,1){\line(1,0){#1}}%
\end{picture}\,}%

\newcommand{\vlto}[1]%
{\,\begin{picture}(#1,3)%
\put(0,2){\vector(1,0){#1}}%
\end{picture}\,}%

\newcommand{\vllarrow}[1]%
{\,\begin{picture}(#1,3)%
\put(#1,2){\vector(-1,0){#1}}%
\end{picture}\,}%

\newcommand{\dirlm}[1]%
  {
     {\lim\hskip-1.58em\lower.65ex
       \hbox{$
                {}_{\stackrel{\lower1ex\hbox
                                        {$\scriptstyle -\!\!\!\longrightarrow$}
                                      }{\vbox to0pt{\vss\vskip1ex
                                            \hbox{$\scriptstyle{}^{#1}$}\vss}}
                   }
            $}
     }
\:}

\newcommand{\subdirlm}[1]%
  {
     {\lim\hskip-1.5em\lower.6ex
       \hbox{$
                   {}_{\stackrel{\lower1ex\hbox
                                           {$\scriptstyle\longrightarrow$}
                                }{ ^{#1} }
                      }
             $}
     }
\:}

\newcommand{\inlm}[1]%
   {
      {\lim\hskip-1.58em\lower.65ex
        \hbox{$
                 {}_{\stackrel{\lower1ex\hbox
                                        {$\scriptstyle \longleftarrow\!\!\<-$}
                              }{\vbox to0pt{\vss\vskip.6ex
                                            \hbox{$\scriptstyle{}^{#1}$}\vss}}
                    }
             $}
      }
\:}

\def\halfsize#1{{\scalebox{.6}{#1}}}
\def\hz#1{{\hbox to 0pt{#1}}}
\def\Iso{\vbox to 0pt{\vss\hbox{$\widetilde{\phantom{nn}}$}\vskip-7pt}}
\def\Is{\vbox to 0pt{\vss\hbox{$\widetilde{\phantom{n}}$}\vskip-8.5pt}}
\def\>{\mspace {1mu}}
\def\<{\mspace{-1mu}}
\def\({{\textup(}}
\def\){{\textup)}}
\def\bigl#1{{\textup{\begin{large}#1\end{large}}}}
\def\bigr#1{{\textup{\begin{large}#1\end{large}}}}
\newcommand{\X}{{\mathscr X}}

\newcommand{\Y}{{\mathscr Y}}
\newcommand{\Z}{{\mathscr Z}}
\newcommand{\V}{{\mathscr V}}
\newcommand{\U}{{\mathscr U}}
\newcommand{\W}{{\mathscr W}}
\newcommand{\I}{{\mathscr I}}
\newcommand{\J}{{\mathscr J}}

\newcommand{\Spec}{{\mathrm {Spec}}}
\newcommand{\Spf}{{\mathrm {Spf}}}
\newcommand{\Supp}{{\mathrm {Supp}}}
\newcommand{\A}{{\mathcal A}}

\newcommand{\C}{{\mathcal C}}
\newcommand{\cD}{{\mathcal D}}
\newcommand{\cDc}{\cD_{\<\<\mathrm c}}
\newcommand{\cDt}{\cD_{\<\mathrm t}}
\newcommand{\E}{{\mathcal E}}
\newcommand{\F}{{\mathcal F}}
\newcommand{\G}{{\mathcal G}}
\renewcommand{\H}{{\mathcal H}}
\newcommand{\Hr}{{\mathrm H}}
\newcommand{\Hp}[2]{\Hr_{#1}^{\<\prime\>#2}}
\newcommand{\M}{{\mathcal M}}
\newcommand{\N}{{\mathcal N}}
\newcommand{\cH}{{\mathcal H}}
\newcommand{\cL}{{\mathcal L}}
\newcommand{\cI}{{\mathcal I}}
\newcommand{\cJ}{{\mathcal J}}
\newcommand{\cK}{{\mathcal K}}
\newcommand{\cO}{{\mathcal O}}

\newcommand{\cR}{{\mathcal R}}
\newcommand{\cRc}{{\mathcal R_{\mathrm c}}}
\newcommand{\cRt}{{\mathcal R_{\>\mathrm t}}}

\newcommand{\cV}{{\mathcal V}}

\newcommand{\D}{{\mathbf D}}
\newcommand{\K}{{\mathbf K}}
\newcommand{\vc}{{\vec{\mathrm{c}}}}
\newcommand{\Kvc}{\K_{\vc}}
\newcommand{\Dqc}{\D_{\mkern-1.5mu\mathrm {qc}}}
\newcommand{\wDqc}{ \widetilde
         {\vbox to6.5pt{\vss\hbox{$\mathbf D$}}}
   _{\mkern-1.5mu\mathrm {qc}} }
\newcommand{\wDqcp}{\wDqc^{\lower.5ex\hbox{$\scriptstyle+$}}}
\newcommand{\Dvc}{\D_{\<\vc}}
\newcommand{\Dqct}{\D_{\mkern-1.5mu\mathrm{qct}}}
\newcommand{\DqcZ}{\D_{\mkern-1.5mu{\mathrm {qc}}Z}}
\newcommand{\Dc}{\D_{\mkern-1.5mu\mathrm c}}
\newcommand{\Qt}{Q^{\mathrm t}}

\newcommand{\qc}{{\mathrm{qc}}}

\newcommand{\bj}{{\boldsymbol j}}
\newcommand{\bi}{{\boldsymbol i}}

\newcommand{\bh}{{\boldsymbol h}}
\newcommand{\R}{{\mathbf R}}
\newcommand{\Rfs}{{\mathbf R f_{\!*}}}
\newcommand{\bL}{{\mathbf L}}
\newcommand{\Hom}{{\mathrm {Hom}}}

\newcommand{\Homb}{{\mathrm {Hom}}^{\bullet}}
\newcommand{\Ac}{\A_{\mathrm c}}
\newcommand{\Aqc}{\A_{\qc}}
\newcommand{\Avc}{\A_{\vec {\mathrm c}}}
\newcommand{\Aqct}{\A_{\mathrm {qct}}\<}
\newcommand{\AqcZ}{\A_{{\qc}Z}}
\newcommand{\At}{\A_{\mathrm t}\<}
\newcommand{\Dt}{\D_{\mathrm t}\<}
\newcommand{\ft}{f_{\mathrm t}^\times}
\newcommand{\fs}{f^!}
\newcommand{\gt}{g_{\mathrm t}^{\<\times}}
\newcommand{\gs}{g^!}
\newcommand{\pd}{{\phantom{.}}}
\newcommand{\ush}[1]{{#1^{\textup{\texttt\#}}}}
\newcommand{\Ext}{\E{xt}}
\newcommand{\sHom}{\cH{om}}
\newcommand{\sHomb}{\cH{om}^{\bullet}}
\newcommand{\iGp}[1]{{\varGamma_{\<\!#1}'}}
\newcommand{\iG}[1]{{\varGamma_{\<\!#1}^{\phantom\prime}}}

\newcommand{\set}{\!:=}
\newcommand{\lra}{\longrightarrow}
\newcommand{\iso}%
{{\mkern8mu\longrightarrow \mkern-25.5mu{}^\sim\mkern17mu}}
\newcommand{\osi}%
{{\mkern8mu\longleftarrow \mkern-24.5mu{}^\sim\mkern16mu}}
\newcommand{\Otimes}{\underset
  {\vbox to 0pt {\vskip-1ex\hbox{$\scriptscriptstyle=$}\vss}}
    {\otimes}\vadjust{\kern.4pt}} 
\newcommand{\BG}{\boldsymbol{\varGamma}}
\newcommand{\LL}{{\mathbf L}\Lambda}
\newcommand{\BL}{{\boldsymbol\Lambda}}
\newcommand{\BLc}{{\boldsymbol\Lambda{\!^\vc}}}
\newcommand{\smcirc}%
  {{\raise.15ex\hbox to.7em{$\hss \scriptstyle\circ\hss$}}} 

\begin{document}

\author[Leovigildo Alonso, Ana Jerem\'{\i}as,
and Joseph Lipman]
{\smash{\ \ \  \  \, Leovigildo Alonso Tarr\'{\i}o, Ana Jerem\'{\i}as,
L\'opez \and Joseph Lipman}}

\address{Universidade de Santiago de Compostela\\
               Facultade de Matem\'aticas\\
               E-15771  Santiago de Compostela, SPAIN}
\email{lalonso@zmat.usc.es}
\urladdr{http://web.usc.es/\~{}lalonso/}

\address{Universidade de Santiago de Compostela\\
               Facultade de Matem\'aticas\\
               E-15771  Santiago de Compostela, SPAIN}
\email{jeremias@zmat.usc.es}

\address{Dept.\ of Mathematics, Purdue University\\
              W. Lafayette IN 47907, USA}
\email {lipman@math.purdue.edu}
\urladdr{www.math.purdue.edu/\~{}lipman/}

 \def\copyrightyear{1999}
\def\copyrightholder{American Mathematical Society}

\setcounter{page}{3}

\title[Duality and Flat Base Change on Formal Schemes]%
{Duality and Flat Base Change on Formal Schemes}


\thanks{First two authors partially supported by
            Xunta de Galicia research project XUGA20701A96
            and Spain's DGES grant PB97-0530. They also
            thank the Mathematics Department of Purdue University for its hospitality,
            help and support. .}

\thanks{Third author partially
           supported by the National Security Agency.\vadjust{\kern1.5pt}}



\begin{abstract}
We give several related versions of global Grothendieck Duality for
unbounded complexes on noetherian formal schemes. The proofs, based on
a non-trivial adaptation of Deligne's method for the special case of
ordinary schemes, are reasonably self-contained, modulo the Special
Adjoint Functor \hbox{Theorem.}  An~alternative approach, inspired by
Neeman and based on recent results about ``Brown Representability," is
indicated as well. A section on applications and examples illustrates
how our results synthesize a number of different duality-related
topics (local duality, formal duality, residue theorems, dualizing
complexes,\mbox{\dots\!\!).}

A flat-base-change theorem for pseudo\kern.6pt-proper maps
leads in particular to sheafified versions of duality for bounded-below complexes
with quasi-coherent homology.  Thanks to Greenlees-May duality, the results take a
specially nice form for proper maps and bounded-below complexes with coherent
homology.
\end{abstract}

\maketitle
\tableofcontents

\section{Preliminaries and main theorems.}
\label{S:prelim}

First we need some notation and terminology.  Let $X$ be a ringed
space,\index{ringed space} i.e., a topological space together with a sheaf of
commutative rings $\cO_{\<\<X}$.%
\index{ $\R$@$\cO_{\<\<\<X}$ (structure sheaf of ringed space $X$)} 
Let $\A(X)$\index{ $\A$ (module category)} be the
\hbox{category} of
$\cO_{\<\<X}$-modules, and $\Aqc\<(X)$\index{ $\A$ (module category)!$\Aqc$}
 (resp.\ $\Ac(X)$,\index{ $\A$ (module category)!$\Ac$} resp.~
$\Avc(X)$)\index{ $\A$ (module category)!$\Ac$@$\Avc$} the full subcategory of~
$\A(X)$ whose objects are the quasi-coherent (resp.\;coherent,
resp.\;\smash{$\dirlm{}\!\!$}'s\index{lim@\smash{$\subdirlm{}\mkern-4mu$}} of
coherent)
$\cO_{\<\<X}$-modules.%
\footnote{%
``\smash{$\subdirlm{}\!\!$}" always  denotes 
a direct limit over a small ordered index set 
in which any two elements have an upper bound. More general 
direct limits will be referred to as \emph{colimits}.%
}
Let~ $\K(X)$\index{ $\K$ (homotopy category)} be the homotopy category of
$\A(X)$-complexes, and let~$\D(X)$\index{ $\D$ (derived category)} be the
corresponding derived category, obtained from~$\K(X)$ by adjoining\- an inverse
for every quasi-isomorphism (=\:homotopy class of maps of complexes
inducing homology isomorphisms).\index{quasi-isomorphism}

\penalty-1000

For any full subcategory\-
$\A_{\scriptscriptstyle{\ldots}}(X)$ of~$\A(X)$, denote by
$\D_{\scriptscriptstyle{\ldots}}(X)$ the full subcategory of~$\>\D(X)$
whose objects are those complexes whose homology sheaves all lie in~
$\A_{\scriptscriptstyle{\ldots}}(X)$, and by
$\D^+_{\raise.3ex\hbox{$\scriptscriptstyle{\ldots}$}}(X)$
\index{ $\D$ (derived
category)!$\D^\pm_{\raise.3ex\hbox{$\scriptscriptstyle{\ldots}$}}\>$}
(resp.~$\D^-_{\raise.3ex\hbox{$\scriptscriptstyle{\ldots}$}}(X)$) the full
subcategory of~$\D_{\scriptscriptstyle{\ldots}}(X)$ whose objects are those
complexes~$\F\in\D_{\scriptscriptstyle{\ldots}}(X)$ such that the homology
$H^m(\F\>)$ vanishes for all $m\ll0$ (resp.~$m\gg0$).

The full subcategory $\A_{\scriptscriptstyle{\ldots}}(X)$ of $\A(X)$
is \emph{plump}\index{plump} if it contains 0 and for every exact sequence
$\M_1\to\M_2\to\M\to\M_3\to \M_4$ in~$\A(X)$ with $\M_1$, $\M_2$,
$\M_3$ and~ $\M_4$ in~$\A_{\scriptscriptstyle{\ldots}}(X)$, $\M$ is
in~$\A_{\scriptscriptstyle{\ldots}}(X)$ too.  If
$\A_{\scriptscriptstyle{\ldots}}(X)$ is plump then it is abelian, and has a
derived category $\D(\A_{\scriptscriptstyle{\ldots}}(X))$.
 For~example, $\Ac(X)$ is plump
\cite[p.\,113,~(5.3.5)]{GD}.  If $\X$ is a locally noetherian formal\index{formal
scheme} scheme,%
\footnote
{Basic properties of formal schemes can be found in \cite[Chap.\,1, \S10]{GD}.%
} 
then $\Avc(\X)\subset\Aqc(\X)$ (\Cref{C:vec-c is qc})---with equality
when $\X$  is an ordinary scheme, i.e., when $\cO_\X$ has discrete topology
\cite[p.\,319, (6.9.9)]{GD}---and both of these are plump subcategories
of~$\A(\X)$, see \Pref{(3.2.2)}.

\smallskip
Let $\K_1$, $\K_2$ be triangulated categories\index{triangulated category} with
respective translation functors $T_1\>,\,T_2$ \cite[p.~20]{H1}.  A (covariant)
\emph{$\Delta$-functor}\index{Delta fun@$\Delta$-functors (on
triangulated categories)}  is a pair $(F, \Theta)$
consisting of an additive functor $F\colon \K_1\to\K_2$ together with an
isomorphism of~functors $ \Theta:FT_1\iso T_2F $ such that for every triangle $
A\stackrel{u}{\longrightarrow}B\stackrel{v}{\longrightarrow}
C\stackrel{w}{\longrightarrow} T_1A $ in $\K_1\mspace{.6mu}$, the
diagram
$$
FA\xrightarrow{\ Fu\ } FB\xrightarrow{\ Fv\ } FC
  \xrightarrow{\Theta\<\smcirc\< Fw\>} T_2FA
$$
is a triangle in $\K_2\mspace{.6mu}$. Explicit reference to
$\Theta$ is often suppressed---but one should keep it in mind.
(For example, if
$\A_{\scriptscriptstyle{\ldots}}(X)\subset\A(X)$ is plump, then each
of $\D_{\scriptscriptstyle{\ldots}}(X)$ and~
$\D^\pm_{\raise.3ex\hbox{$\scriptscriptstyle{\ldots}$}}(X)$ carries a
unique triangulation for which the translation is the restriction of
that on~$\D(X)$ and such that inclusion into $\D(X)$ together with
$\Theta\!:=$identity is a $\Delta$-functor; in other words, they are
all \emph{triangulated subcategories}\index{triangulated category!triangulated
subcategory} of~$\D(X)$. See e.g.,
\Pref{P:Rhom} for the usefulness of this remark.) 
Compositions\index{Delta fun@$\Delta$-functors (on
triangulated categories)!composition of} of
$\Delta$-functors, and morphisms between\index{Delta fun@$\Delta$-functors (on
triangulated categories)!morphism between}
$\Delta$-functors, are defined in the natural way.%
\footnote{%
See also \cite[\S0,\,\S1]{De} for the multivariate case,
where signs come into play---and $\Delta$-functors are called
``exact functors."%
} 
A $\Delta$-functor $(G,\Psi)\colon\K_2\to\K_1$ is a \emph{right
$\Delta$-adjoint}\index{Delta adj@$\Delta$-adjoint} of~$(F,\Theta)$ if $G$~is a
right adjoint of~$F$ and the resulting functorial map $FG\to \mathbf 1$ (or
equivalently,
$\mathbf 1\to GF$) is a morphism of $\Delta$-functors.

We use $\R$\index{ $\R$ (right-derived functor)}\index{derived functor} to denote
right-derived functors, constructed e.g., via K-injective\index{K-injective
resolution} resolutions (which exist for all $\A(X)$-complexes  \cite[p.\,138,
Thm.~4.5]{Sp}).%
\footnote
{A complex $F$ in an abelian category~$\A$ is K-injective if for each exact
$\A$-complex~$G$ the abelian-group complex $\Homb_\A(G,F)$ is again
exact. In particular, any bounded-below complex of injectives is K-injective.
If every $\A$-complex~$E$ admits a K-injective resolution
$E\to I(E)$ (i.e., a quasi-isomorphism into a K-injective complex~$I(E)$), then
every additive functor~$\Gamma\colon\A\to\A'$ ($\A'$~abelian) has a
right-derived functor~$\R\Gamma\colon\D(\A)\to\D(\A')$ which satisfies
$\R\Gamma(E)=\Gamma(I(E))$. For~example, 
$\R\Homb_\A(E_1,E_2)=\Homb_\A(E_1,I(E_2))$.%
} 
For a map $f\colon X\to Y$\index{ringed
space!ringed-space map}
of ringed spaces (i.e., a continuous map $f\colon X\to Y$ together with a
ring-homomorphism $\cO_Y\to f_{\!*}\cO_{\<\<X}$),
$\bL f^*$%
\index{ $\mathbf L$@$\bL$ (left-derived functor)} 
denotes the left-derived functor of
$f^*\<$, constructed via K-flat resolutions
\cite[p.\,147, 6.7]{Sp}. Each derived functor in this paper comes equipped,
implicitly, with a $\Theta$ making it into a $\Delta$-functor (modulo
obvious modifications for contravariance),
cf.~\cite[Example~(2.2.4)]{Derived categories}.%
\footnote
{We do not know, for instance, whether $\bL f^*$---which is
defined only up to isomorphism---can always be chosen so as
to commute with translation, i.e., so that $\Theta={}$Identity will
do.%
} 
Conscientious readers may verify that such morphisms
between derived functors as occur in this paper are in fact
morphisms of $\Delta$-functors.

\begin{parag}\setcounter{sth}{0} 
Our \textbf{first main result,} 
global Grothendieck Duality\index{Grothendieck Duality!global}
for a map \mbox{$f\colon\X\to\Y$} of quasi-compact formal schemes with $\X$
noetherian, is that, $\D(\Avc(\X))$ being the derived category of~$\Avc(\X)$
and $\bj \colon\D(\Avc(\X))\to\D(\X)$ being\index{ $\iG{\<\cJ\>}$@$\bj$}
the natural functor, \emph{the
$\Delta$-functor\/~$\Rfs\<\<\smcirc\<\bj$ has a right\/ 
$\Delta$-adjoint.} \vspace{1pt}

A more elaborate---but readily shown equivalent---statement is:
\begin{thespecial}
\label{Th1}
Let\/ $f\colon\<\X \to \Y$ be a map of quasi-compact formal schemes, with $\X$
noetherian, and let\/~$\bj\colon\<\D(\Avc(\X))
\to\D(\X)$ be the natural functor.  Then there exists a\/
$\Delta$-functor%
\index{ $\iG$@$f^{{}^{\>\ldots}}$ (right adjoint of
$\R f_{\<\<*}\cdots$)!$f^\times\<\<$} 
$f^{\times}\<\colon\D(\Y)
\to\D\left({\Avc(\X)}\right)$\vadjust{\kern.3pt} together with a
morphism of\/ $\Delta$-functors 
$\tau:\R f_{\!*}\>\>\bj\> f^{\times}\to{\bf 1}$\index{ {}$\tau$ (trace map)}
such that for all\/
$\G\in\D(\Avc(\X))$ and\/ $\F\in\D(\Y),$\ the composed map $($in the derived
category of abelian groups\/$)$
\begin{align*}
\R\Homb_{\Avc(\X)\!}(\G,\>f^\times\<\<\F\>)
&\xrightarrow{\<\<\mathrm{natural}\,}
   \R\Homb_{\A(\Y)\!}(\R f_{\!*}\>\G,\R f_{\!*} f^\times\<\<\F\>)\\
&\xrightarrow{\;\>\mathrm{via}\ \tau\ }
  \R\Homb_{\A(\Y)\!}(\R f_{\!*}\>\G,\F\>)
\end{align*}
is an isomorphism.
\end{thespecial}

Here we think of the $\Avc(\X)$-complexes $\G$ and $f^\times\<\<\F$
as objects in
both $\D(\Avc(\X))$ and~$\D(\X)$. But as far as we know, the natural
map $\Hom_{\D(\Avc(\X))}\to\Hom_{\D(\X)}$ need not always be an
isomorphism. It \emph{is} when $\X$ is \emph{properly algebraic,}\index{properly
algebraic} i.e., the $J$-adic completion of a proper $B$-scheme with $B$ a
noetherian ring and $J$ a $B$-ideal: then $\bj$ induces an
equivalence of categories $\D(\Avc(\X))\to\Dvc(\X)$, see
\Cref{corollary}.  So for properly algebraic $\X$, we can replace
$\D(\Avc(\X))$ in~Theorem~1 by~$\Dvc(\X)$, 
and let~$\G$ be any $\A(\X)$-complex with $\Avc(\X)$-homology.

We prove \Tref{Th1} (=\:\Tref{prop-duality}) in \S\ref{sec-th-duality},
adapting the argument of Deligne\index{Deligne, Pierre} in \cite[Appendix]{H1} (see
also
\cite[\S1.1.12]{De}) to the category~$\Avc(\X)$, which  presents itself as an
appropriate generalization to formal schemes of the category of quasi-coherent
sheaves on an ordinary noetherian scheme.  For this adaptation what is needed,
mainly, is the plumpness of
$\Avc(\X)$ in
$\A(\X)$, a non-obvious fact mentioned above.    
In addition, we need some facts on
``boundedness" of certain derived functors in order to extend the argument to
unbounded complexes. (See section~\ref{SS:bounded}, which makes use of
techniques from~\cite{Sp}.)%
\footnote{A  $\Delta$-functor $\phi$  is \emph{bounded above} if there is an
integer $b$ such that for any $n$ and any complex~$\E$ such that $H^i\E=0$ for all
$i\le n$ it holds that $H^j\<(\phi\E)=0$ for all $j<n+b$.
\emph{Bounded below} and \emph{bounded} (above and below)
are defined analogously. Boundedness (way-outness)\index{boundedness
(way-outness) of $\Delta$-functors} is what makes the very useful ``way-out
Lemma"
\cite[p.\,68, 7.1]{H1} applicable.}

In Deligne's approach\index{Deligne, Pierre} the
``Special Adjoint Functor Theorem"\index{Special Adjoint Functor Theorem}  is
used to get right adjoints for certain  functors on $\Aqc(X)$, and then these right
adjoints are applied to injective resolutions of complexes\dots There is now a
neater approach to duality on a quasi-compact separated ordinary scheme~$X\<$,
due to Neeman\index{Neeman, Amnon}
\cite{N1}, in which ``Brown Representability"\index{Brown Representability}  
shows directly that a $\Delta$-functor~$F$ on~$\D(\Aqc(X))$ 
has  a right adjoint if and only if
$F$ commutes with coproducts.   Both approaches need a small
set of category-generators: coherent sheaves for $\Aqc(X)$ in Deligne's,
and perfect complexes for $\D(\Aqc(X))$ in Neeman's.  Lack of knowledge about
perfect\- complexes over formal schemes discouraged us from pursuing Neeman's
strategy. Recently however (after this paper was essentially written), Franke 
showed in~\cite{BR} that Brown Representability\index{Brown Representability}
holds for the derived category of an arbitrary Grothen\-dieck category~$\A$.%
\footnote{So does the closely-related
existence of K-injective resolutions for all $\A$-complexes. (See also
\cite[\S5]{AJS}.)} %
Consequently \Tref{Th1} also follows from the fact that
$\Avc(\X)$ is a Grothendieck category (straightforward to see once we know 
it---by plumpness in~$\A(\X)$---to be abelian) together with the fact that
$\Rfs\<\smcirc\bj$ commutes with coproducts (\Pref{P:coprod}).
\end{parag}

\medskip
\begin{parag}\label{Gamma'} 
Two other, probably more useful, generalizations---from ordinary
schemes to formal schemes---of global Grothendieck Duality are stated
below in \Tref{Th2} and treated in detail in
\S\ref{S:t-duality}. To describe them, and related results, we need
some preliminaries about \emph{torsion functors}.

\smallskip
\begin{sparag}\label{Gamma'1}
Once again let $(X, \cO_{\<\<X})$ be a ringed space. For any
$\cO_{\<\<X}$-ideal~$\cJ\<$, set
$$
\iG{\<\cJ\>}{\M} := \dirlm{n>0\,\,\>} 
{\cH}om_{\cO_{\<\<X}\!}(\cO_{\<\<X}/\cJ^{n},\, \M)\qquad\bigl(\M\in\A(X)\bigr),
$$
and regard~$\iG{\<\cJ\>}$ as a subfunctor of the identity functor on
$\cO_{\<\<X}$-modules.  If $\N\subset \M$ then 
$\iG{\<\cJ\>}{\N}=\iG{\<\cJ\>}{\M}\cap\N\>$; and it follows formally that 
the functor $\iG{\<\cJ\>}$ is  idempotent
($\iG{\<\cJ\>}\iG{\<\cJ\>}\M=\iG{\<\cJ\>}\M$) 
and left exact \cite[p.\,138,
Proposition~1.7\kern.5pt]{Stenstrom}.\vadjust{\kern.7pt}

Set $\A_\cJ(X)\!:=\iG{\<\cJ\>}(\A(X))$,\index{ $\A$ (module
category)!$\Ac$@$\A_\cJ$} the full subcategory of~$\A(X)$ whose objects are the
\emph{$\cJ\!$-torsion sheaves,}\index{torsion sheaf} i.e., the
$\cO_{\<\<X}$-modules~$\M$ such that
$\iG{\<\cJ\>}{\M}=\M$.  Since
$\iG{\<\cJ\>}$ is an idempotent subfunctor of the identity functor,
therefore it is right-adjoint to the inclusion
$i=i_\cJ\colon\A_\cJ(X)\hookrightarrow\A(X)$. Moreover, $\A_\cJ(X)$ is
closed under $\A(X)$-colimits: if $F$ is any functor into $\A_\cJ(X)$
such that $iF\/$ has a colimit $\M\in\A(X)$, then, since $i$ and
$\iG{\<\cJ}$ are adjoint, the corresponding functorial map from $iF$
to the constant functor with value~$\M$ factors via a functorial map
from $iF$ to the constant functor with value $\iG{\<\cJ\>}{\M}$, and
from the definition of colimits it follows that the monomorphism
$\iG{\<\cJ\>}{\M}\hookrightarrow\M$ has a right inverse,\vspace{.6pt}
so that it is an isomorphism, and thus $\M\in\A_\cJ(X)$. In
particular, if the domain of a functor~$G$ into $\A_\cJ(X)$ is a small
category, then $iG$ does have a colimit, which is also a colimit of
$G$; and so $\A_\cJ(X)$ has small colimits, i.e., it is
small-cocomplete.

Submodules and quotient modules of~$\cJ\!$-torsion sheaves are 
$\cJ\!$-torsion sheaves.
If $\cJ$ is  \emph{finitely-generated} (locally) and if $\N\subset\M$ are 
$\cO_{\<\<X}$-modules such that $\N$ and~$\M/\N$ are $\cJ\!$-torsion sheaves
then $\M$  is a $\cJ\!$-torsion sheaf too; and hence 
$\A_\cJ(\<X)$~is plump in~$\A(\<X)$.%
\footnote{Thus the subcategory $\A_{\cJ}(X)$ is a \emph{hereditary torsion
class} in
$\A(X)$,  in the sense of Dickson, see \cite[pp.\;139--141]{Stenstrom}. } 
In this case, the stalk
of~$\iG{\<\cJ\>}\M$ at $x\in X$ is%
\index{ $\iG{\raise.3ex\hbox{$\scriptscriptstyle{\ldots}$}}$ (torsion functor)} 
$$
(\iG{\<\cJ\>}\M)_x=\dirlm{n>0\,\,\>}
\mathrm{Hom}_{\cO_{\!X\!,\>x}}\<(\cO_{\!X\!,\>x}/\cJ_x^{n},\, \M_x).
$$ 

Let $X$ be a locally noetherian scheme and $Z\subset X$ a closed
subset, the support of~$\cO_{\<\<X}/\cJ$ for some quasi-coherent
$\cO_{\<\<X}$-ideal $\cJ\<$. The functor $\iGp{Z}\set\iG{\<\cJ\>}$
\index{ $\iG{\raise.3ex\hbox{$\scriptscriptstyle{\ldots}$}}$ 
(torsion functor)!$\varGamma'_{\raise.3ex\hbox{$\scriptscriptstyle{\ldots}$}}$} 
does not depend on the quasi-coherent ideal~$\cJ$ determining $Z$. It is a
subfunctor of the left-exact functor
$\iG{Z}^{\phantom{.}}$ which associates to each $\cO_{\<\<X}$-module~$\M$ its
subsheaf of sections supported in~$Z$. If $\M$ is quasi-coherent, then
\mbox{$\iGp{Z}(\M) =\iG{Z}(\M)$}.  

\pagebreak[3]
More generally, for any complex $\E\in\Dqc\<(X)$,  the $\D(X)$-map
$\R\iGp{Z}\E\to\R\iG{Z}\E$ induced by the inclusion
$\iGp{Z}\hookrightarrow \iG Z$ is an isomorphism \cite[p.\,25, Corollary
(3.2.4)]{AJL}; so for such~$\E$ we usually identify  $\R\iGp{Z}\E$ and
$\R\iG{Z}\E$.  

Set $\A_Z(X)\set \A_\cJ(X)$,\index{ $\A$ (module category)!$\A_Z$} the plump
subcategory of~$\A(X)$ whose objects are the \mbox{\emph{$Z$-torsion sheaves,}}
that is,  the $\cO_{\<\<X}$-modules~$\M$ such that \hbox{$\iGp{Z}{\M}=\M\>$;} and
set
\hbox{$\AqcZ(X)\!:=\Aqc(X)\cap\A_Z(X)$,}\index{ $\A$ (module category)!$\AqcZ$}
the plump subcategory of~$\A(X)$ whose objects are the quasi-coherent
$\cO_{\<\<X}$-modules  supported in~$Z$.
\enlargethispage{-\baselineskip}

\smallskip
\pagebreak[3]

For a locally noetherian formal scheme $\X$ with ideal of
definition~$\J$,  set $\iGp{\X}\!:=\iG{\J}\mspace{-.5mu}$,%
\index{ $\iG{\raise.3ex\hbox{$\scriptscriptstyle{\ldots}$}}$ 
(torsion functor)!$\varGamma'_{\raise.3ex\hbox{$\scriptscriptstyle{\ldots}$}}$} a
left-exact functor depending only on the sheaf of topological rings
$\cO_{\X}\>$, not on the choice of~$\J$---for $\M\in\A(\X)$,
$\iGp{\X}\M\subset\M$ is  the submodule whose sections are those
of~$\M$ annihilated locally by an open ideal. Say that $\M$ is a
\emph{torsion~sheaf}\index{torsion sheaf} if \mbox{$\iGp{\X}\M=\M$}.
Let $\At(\X)\set\A_\J(\X)$\index{ $\A$ (module category)!$\At\>$} be  the plump
subcategory of~$\A(\X)$ whose objects\- are all the torsion sheaves; and set
$\Aqct(\X)\!:=\Aqc\<(\X)\cap\At(\X)$,\index{ $\A$ (module category)!$\Aqct$} 
the full (in fact plump, see
\Cref{qct=plump}) subcategory of~$\A(\X)$ whose objects are the quasi-coherent
torsion sheaves.  It holds that $\Aqct(\X)\subset\Avc(\X)$, see \Cref{Gamma'+qc}.
If $\X$ is an ordinary locally noetherian scheme (i.e., $\J=0$), then
$\At(\X)=\A(\X)$ and~$\Aqct(\X)=\Aqc\<(\X)=\Avc(\X)$.
\end{sparag}

\begin{sparag}\label{maptypes}
For any map $f\colon\X \to \Y$ of locally noetherian formal schemes
there are ideals of definition $\I\subset \cO_{\Y}$ and
$\J\subset\cO_{\X\>}$ such that $\I\cO_{\X}\subset \J$ \cite[p.\,416,
(10.6.10)]{GD}; and correspondingly there is a map of ordinary schemes
($=\:$formal schemes having~(0) as ideal of definition)
$(\X,\cO_\X/\J)\to(\Y,\cO_\Y/\I)$ \cite[p.\,410, (10.5.6)]{GD}. We~say
that $f$~is \emph{separated}%
\index{formal-scheme map!separated}
(resp.~\emph{affine,}%
\index{formal-scheme map!affine}
resp.~\emph{pseudo\kern.6pt-proper,}%
\index{formal-scheme map!pseudo\kern.6pt-proper}
resp.~\emph{pseudo\kern.6pt-finite,}%
\index{formal-scheme map!pseudo\kern.6pt-finite}
resp.~\emph{of pseudo\kern.6pt-finite type}\kern-1pt)%
\index{formal-scheme map!of pseudo\kern.6pt-finite type}
if for some---and hence
any---such~$\I,\J$ the corresponding scheme-map is separated (resp.~affine,
resp.~proper, resp.~finite,  resp.~of finite type),
see \cite[\S\S10.15--10.16, p.\,444\:\emph{ff.}]{GD}, 
keeping in mind \cite[p.\,416, (10.6.10)(ii)]{GD}.%
\footnote
{In \cite[Definition 1.14]{Ye}, pseudo-finite-type maps
are called ``maps of formally finite type." 
The proof  of Prop.\,1.4 in \cite{Ye} (with $A'=A$)  yields
the following characterization of pseudo-finite-type maps of affine formal
schemes (cf.~\cite[p.\,439,  Prop.\,(10.13.1)]{GD}):
The map\/ $f\colon \Spf(B) \to \Spf(A)$ corresponding to a continuous
homomorphism~$h\colon A\to B$ of noetherian adic rings is of
pseudo-finite type $\Leftrightarrow$
for any ideal of definition\/ $I$ of\/ $A,$  there
exists\vadjust{\kern.8pt} an\/
$A$-algebra of finite type\/~$A'\<,$ an\/ $A'\<$-ideal\/~ $I' \supset IA',$
and an\/ $A$-algebra homomorphism\/
\mbox{$A'\to B$} inducing an adic surjective map\/ 
$\widehat {A'}\twoheadrightarrow B$
 where\/ $\widehat{A'}$ is the\/ $I'$-adic
completion of\/~$A'$.} 
Any affine map is separated. Any pseudo\kern.6pt-proper map is separated and of
pseudo\kern.6pt-finite type. The map~$f$ is pseudo\kern.6pt-finite
$\Leftrightarrow$ it is  pseudo\kern.6pt-proper and affine
$\Leftrightarrow$ it is  pseudo\kern.6pt-proper and has finite fibers
\cite[p.\,136, (4.4.2)]{EGA}.

 We say that $f$ is \emph{adic}%
\index{formal-scheme map!adic}
if for some---and hence any---ideal of definition 
$\I\subset\cO_\Y\>$, $\I\cO_{\X}$~is an ideal of definition of~$\X$
\cite[p.\,436, (10.12.1)]{GD}.  We say that $f$ is \emph{proper}%
\index{formal-scheme map!proper}
(resp.~\emph{finite,}%
\index{formal-scheme map!finite}  
resp.~\emph{of finite type}\kern-1pt)%
\index{formal-scheme map!of finite type} 
if $f$ is pseudo\kern.6pt-proper (resp.~pseudo\kern.6pt-finite,  resp.~of
pseudo\kern.6pt-finite type) and adic, see \cite[p.\,119, (3.4.1)]{EGA}, \cite[p.\,148,
(4.8.11)]{EGA} and \cite[p.\,440, (10.13.3)]{GD}.
\end{sparag}

\begin{sparag} Here is our \textbf{second main result}, 
Torsion Duality\index{Grothendieck Duality!Torsion (global)} for
formal schemes. (See \Tref{T:qct-duality} and
\Cref{C:f*gam-duality} for more elaborate statements.)  In
the assertion, $\wDqc(\X)\set\R\iGp\X{}^{-1}(\Dqct(\X))$%
\index{ $\D$ (derived category)!z@${ \widetilde
{\vbox to5pt{\vss\hbox{$\mathbf D$}}}_{\mkern-1.5mu\mathrm {qc}} }$} 
is the least $\Delta$-subcategory of~$\D(\X)$ containing both
$\Dqc(\X)$\vspace{.6pt} and $\R\iGp\X{}^{-1}(0)$ (\Dref{D:Dtilde},
Remarks~\ref{R:Dtilde}, (1) and~(2)).\vspace{-1.3pt} For example, when
$\X$ is an ordinary scheme then $\wDqc(\X)=\Dqc(\X)$.

\begin{thespecial}
\label{Th2}
Let $f\colon \X\to\Y$ be a map of noetherian formal schemes. Assume either
that\/ $f$ is separated or that\/ $\X$ has finite Krull dimension, or else
restrict~to bounded-below complexes.

\smallskip
\textup{(a)} The restriction of\/ $\Rfs\colon\D(\X)\to\D(\Y)$ takes\/ $\Dqct(\X)$
to\/
$\Dqct(\Y),$ and it has a right\/ $\Delta$-adjoint\/ 
$
\ft\colon\D(\Y)\to\Dqct(\X).\index{ $\iG$@$f^{{}^{\>\ldots}}$ (right adjoint of
$\R f_{\<\<*}\cdots$)!$\ft\<\<$}
$

\smallskip
\textup{(b)} The restriction of\/ $\Rfs\R\iGp\X$ takes\/
$\wDqc(\X)$ to\/
$\Dqct(\Y)\subset\wDqc(\Y),$ and it has a~right\/ $\Delta$-adjoint\/ 
$
\ush f\colon\D(\Y)\to\wDqc(\X).\index{ $\iG$@$f^{{}^{\>\ldots}}$ (right adjoint of
$\R f_{\<\<*}\cdots$)!$\ush f$}
$
\end{thespecial}

\penalty -1000
\begin{srems}\label{R:Th2}
(1) The ``homology localization" functor
$$
\BL^{}_\X(-)\set\R\sHomb(\R\iGp\X\cO_\X^{}\>,-)\index{ $\mathbf {La}$@$\BL$
(homology localization)}
$$
is right-adjoint to $\R\iGp\X$, and $\BL_\X^{-1}(0)=\R\iGp\X{}^{-1}(0)$
(Remarks~\ref{R:Gamma-Lambda}). The $\Delta$-functors
$\ush f$ and~$\ft$ are connected thus (Corollaries~\ref{C:f*gam-duality} 
and~\ref{C:identities}(a)):
$$
\ush f=\BL^{}_\X\ft\<, \qquad
\ft=\R\iGp\X\ush f.
$$

\smallskip
(2) In the footnote
on page \pageref{C:completion-proper} it is indicated that
$\R\iGp\X{}^{-1}(0)$ admits a ``Bousfield colocalization"\index{Bousfield
colocalization} in
$\D(\X)$, with associated ``cohomology colocalization" functor $\R\iGp\X\>$; and
in \Rref{R:Gamma-Lambda}(3), \Tref{Th2} is interpreted as
duality\- with coefficients in the corresponding quotient 
$\wDqc(\X)/\R\iGp\X{}^{-1}(0)\cong
\Dqc(\X)/\bigl(\Dqc(\X)\cap\R\iGp\X{}^{-1}(0)\bigr)$.

\smallskip
(3) The proof of \Tref{Th2} is similar to that of
\Tref{Th1}, at least when the formal scheme~$\X$ is separated
(i.e., the unique formal-scheme map $\X\to\Spec(\mathbb Z)$ is
separated) or finite-dimensional, in which case there is an
\emph{equivalence of categories} $\D(\Aqct(\X))\to \Dqct(\X)$
(\Pref{1!}).  (As mentioned before, we know the
corresponding result with ``$\mspace{2mu}\vec{\mathrm c}\:$" in place
of ``qct" only for \emph{properly algebraic} formal schemes.)  In
addition, replacing separatedness of~$\X$ by separatedness of~$f$
takes a technical pasting argument.

\smallskip
(4) For an ordinary scheme~$X$ (having $(0)$ as ideal of definition),
$\iGp X$ is just the identity functor of $\A(X)$, and
$\Dqct(X)=\Dqc(X)$. In this case, Theorems~\ref{Th1}~and~\ref{Th2}
both reduce to the usual global (non-sheafified) version of
Grothendieck Duality. In \S\ref{S:apps} we will describe how
\Tref{Th2} generalizes and ties together various strands in the
literature on local, formal, and global duality.  In particular, the
behavior of \Tref{Th2} vis-\`a-vis variable $f$ gives
compatibility of local and global duality, at least on an
abstract level---i.e., without the involvement of differentials,
residues,~etc. (See \Cref{C:kappa-f^times-tors}.)
\end{srems}
\end{sparag}
\end{parag}

\begin{parag}\label{culminate}  As in the classic paper \cite{f!} of Verdier,%
\index{Verdier, Jean-Louis} the \textbf{culminating results} devolve from
flat-base-change isomorphisms, established here for the functors $\ft$ and $\ush
f$ of
\mbox{\Tref{Th2},} with $f$
\emph{pseudo\kern.6pt-proper}---in which case we denote 
$\ft$ by $f^!\<$.\index{ $\iG$@$f^{{}^{\>\ldots}}$ (right adjoint of
$\R f_{\<\<*}\cdots$)!$\mathstrut \fs\<$} 

\begin{thespecial}\label{Th3}\index{base-change isomorphism}

Let\/ $\X,$ $\Y$ and\/ $\U$ be noetherian formal schemes, let\/
$f\colon\X\to\Y$ be a pseudo\kern.6pt-proper map, and let\/ $u\colon
\U\to\Y$ be flat, so that in the natural diagram
$$
\begin{CD}
\X\times_\Y\U=:\>@.\V@>v>>\X \\
@. @VgVV @VVfV \\
@.\U@>>\vbox to 0pt{\vskip-1ex\hbox{$\scriptstyle u$}\vss}>\Y
\end{CD}
$$
the formal scheme\/ $\V$ is noetherian, $g$ is pseudo\kern.6pt-proper,
and $v$ is flat \textup(\Pref{P:basechange}\kern.5pt\textup). 

\pagebreak[3]
Then for all\/
$\F\in\wDqcp(\Y)\set\wDqc(\Y)\cap\D^+(\Y)$ the base-change map\/~$\beta_\F$
of~\Dref{D:basechange} is an isomorphism
$$
\beta_\F\colon\R\iGp\V\>v^*\<\fs\F \iso 
\gs\>\>\R\iGp\U u^*\<\F\ 
\underset{\textup{\ref{C:identities}(b)}}\cong
\gs\<u^*\<\F.
$$
In particular, if\/ $u$ is \emph{adic} then we have a functorial
isomorphism\/ $v^*\<\<\fs\F \iso \gs u^*\<\F\<.$
\end{thespecial}

This theorem is proved in \S7 (\Tref{T:basechange}).  The functor
$\R\iGp\V$ has a right adjoint~$\BL_{\<\V}$, see \eqref{adj}. \Tref{Th3}
leads quickly to the corresponding result for~$\ush f$ (see
\Tref{T:sharp-basechange} and \Cref{C:coh-basechange}):

\penalty -1000
\begin{thespecial}\label{Th4}\index{base-change isomorphism}
Under the preceding conditions, let 
$$
\ush{\beta_{\<\<\F}}\colon v^*\<\<\ush f\<\F\to\ush g
u^*\<\F\qquad\bigl(\F\in\wDqcp(\Y)\bigr)
$$
be the map adjoint to the natural composition
$$
\R g_*\R\iGp\V\>v^*\<\<\ush f\<\F
\underset{\textup{Thm.\,\ref{Th3}}}\iso
\R g_*g^!u^*\<\F\to u^*\<\F.
$$
Then the map $\BL_{\<\V}(\ush{\beta_{\<\<\F}})$ is an \emph{isomorphism}
$$
\BL_{\<\V}(\ush{\beta_{\<\<\F}})\colon\BL_{\<\V}\> v^*\<\<\ush f\<\F\iso
\BL_{\<\V}\>\ush g u^*\<\F
\underset{\textup{\ref{C:identities}(a)}}\cong
\ush g u^*\<\F.
$$
Moreover, if\/ $u$ is an open immersion, or if\/ $\F\in\Dc^+(\Y),$
then\/ $\ush{\beta_{\<\<\F}}$ itself is an isomorphism.
\end{thespecial}

The special case of Theorems~\ref{Th3} and~\ref{Th4} when $u$ is an open
immersion is equivalent to what may be properly referred to as Grothendieck
Duality (unqualified by the prefix ``global"), 
namely the following \emph{sheafified} version of \Tref{Th2} (see
\Tref{T:sheafify}):

\begin{thespecial}\label{Th5}\index{Grothendieck Duality!Torsion
(sheafified)} Let\/ $\X$ and\/ $\Y$ be noetherian formal schemes 
and let\/ $f\colon\X\to\Y$ be
a pseudo\kern.6pt-proper map. Then the following natural compositions are
\emph{isomorphisms:}
\begin{align*}
\Rfs\R\sHomb_{\>\X}\<(\G\<,\>\ush f\<\F\>)
&\to\<
\R\sHomb_{\>\Y}\<(\Rfs\R\iGp\X\>\G\<,\>\Rfs\R\iGp\X\ush f\<\F\>) \\
&\to
\R\sHomb_{\>\Y}\<(\Rfs\R\iGp\X\>\G\<,\>\F\>)\vspace{-4pt}
\quad\ 
\bigl(\G\in\wDqc(\X),\;\F\>\in\wDqcp(\Y)\bigr);
\end{align*}
\vspace{-6pt}

\noindent
$
\quad\Rfs\R\sHomb_{\>\X}\<(\G\<,\>f^!\F\>)
\to
\R\sHomb_{\>\Y}\<(\Rfs\G\<,\>\Rfs f^!\F\>)
\to
\R\sHomb_{\>\Y}\<(\Rfs\G\<,\>\F\>)\vspace{3pt}
$

\rightline{$\bigl(\G\in\Dqct(\X),\;\F\>\in\wDqcp(\Y)\bigr).$}
\end{thespecial}

Finally, if $f$ is \emph{proper} and $\F\in\Dc^+(\Y)$, then
$\ush f\colon\Dc^+(\Y)\to\Dc^+(\X)$ \emph{is right-adjoint to
$\Rfs\colon\Dc^+(\X)\to\Dc^+(\Y)$, and
 $\R\iGp\X$ in
\Tref{Th5} can be deleted,} see \Tref{T:properdual}. 

\smallskip
In this---and several other results about complexes with coherent
homology---an essential ingredient is \Pref{formal-GM}, deduced
here from Greenlees-May duality for ordinary affine schemes, see \cite{AJL}:

\smallskip
\emph{%
Let\/ $\X$ be a locally noetherian formal scheme, and let\/ $\E\in\D(\X)$. 
Then for all\/ $\F\in\Dc(\X)$ the natural map\/ $\R\iGp{\X}\E\to \E$ induces an
isomorphism%
}
$$
\R\sHomb(\E, \>\F\>) \iso \R\sHomb(\R\iGp{\X}\E,\>\F\>).
$$

\end{parag}

\medskip
In closing this introductory section, we wish to express our appreciation for
illuminating interchanges with Amnon Neeman\index{Neeman, Amnon} and Amnon
Yekutieli\index{Yekutieli, Amnon}.

\section{Applications and examples.}
\label{S:apps}

  Again,  \Tref{Th2} generalizes global Grothendieck
Duality on ordinary schemes.  This section illustrates further how
\Tref{Th2} provides a common home for a number of different
duality-related results (local duality, formal duality, residue theorems,
dualizing complexes\dots\!\!).
For a quick example, see \Rref{R:d-vein}.

\Sref{bf (a)} reviews several forms of local duality.
In section~\ref{sheafify} we  sheafify these
results, and connect them to \Tref{Th2}. 
In particular, \Pref{(2.2)} is an abstract version
of the Local Duality theorem of~\cite[p.\,73, Theorem 3.4]{Integration}; and
\Tref{T:pf-duality} (Pseudo-finite Duality)\index{Duality!Pseudo-finite}
globalizes it to formal schemes.\looseness=-1 

\Sref{residue thm} relates Theorems~\ref{Th1}
and~\ref{Th2}  to the central ``Residue Theorems"\index{Residue theorems}
in~\cite{Asterisque} and~\cite{HS} 
(but does not subsume those results).

\Sref{bf (d)} indicates how both the Formal Duality\index{Duality!Formal}
theorem of \cite[p.\,48, Proposition~(5.2)]{De-Rham-cohomology} and the
Local-Global Duality\index{Duality!Local-Global} theorem in
\cite[p.\,188]{Desingularization} can be deduced  from \Tref{Th2}.

\Sref{bf (e)}, building on work of Yekutieli \cite[\S5]{Ye},\index{Yekutieli, Amnon}
treats \emph{dualizing complexes} on formal schemes, and their associated
dualizing functors. For a pseudo\kern.6pt-proper map~$f\<$, the functor~$\ush f$
of
\Tref{Th2} lifts dualizing complexes to dualizing complexes (\Pref{P:twisted
inverse}). For any map $f\colon\X\to\Y$ of noetherian formal schemes,
there is natural isomorphism
$$
\R\sHomb_\X(\bL f^*\F\<,\ush f\G)\iso\ush f\R\sHomb_\Y(\F\<,\G)
\qquad \bigl(\F\in\Dc^-(\Y),\;\G\in\D^+(\Y)\bigr),
$$
(\Pref{P:Hom!}). For pseudo\kern.6pt-proper
$f\<$,  if $\Y$ has a dualizing complex~$\cR$, so that $\ush f\cR$ is
a dualizing complex on~$\X$, and if $\cD^\Y\set\R\sHomb(-,\cR)$ and
$\cD^\X$ are the corresponding dualizing functors, one deduces a natural
isomorphism (well-known for ordinary schemes)
$$
\ush f\<\E\cong\cD^\X\bL f^*\cD^\Y\E\qquad\bigl(\E\in\Dc^+(\Y)\bigr),
$$
see \Pref{P:Dual!}.

There are corresponding results for~$\ft\<$ as well.

\smallskip

\def\GG#1{\Gamma_{\!\!#1}^{\phantom{.}}}
\def\fij{\ush{\varphi_{\<\!J}}}
\begin{parag}
\label{bf (a)}
\renewcommand{\theequation}{\theparag.\arabic{numb}}
(Local Duality.)\index{Duality!Local|(} All rings are commutative, unless otherwise
specified.

Let $\varphi\colon R\to S$ be a ring homomorphism with
$S$ noetherian,  let $J$ be an $S$-ideal, and let $\GG {J}$ be
the functor taking any $S$-module to its submodule of elements which are
annihilated by some power of~$J$. Let $E$ and $E'$ be complexes
in~$\D(S)$, the derived category of $S$-modules, and let $F\in\D(R)$.
With \smash{$\Otimes$}\vadjust{\kern.7pt} denoting derived tensor product
in~$\D(S)$ (defined via K-flat resolutions \cite[p.\,147, Proposition 6.5]{Sp}), there
is a natural isomorphism
\smash{$E \Otimes \R \GG {\<J} E' \iso
 \R \GG {J} (E \Otimes E')$},\vadjust{\kern.7pt} 
see~e.g., \cite[p.\,20, Corollary(3.1.2)]{AJL}. Also, viewing
$\R\Hom_R^\bullet(E'\<, F)$  as a functor from 
$\D(S)^{\rm{op}}\<\times\D(R)$ to~$\D(S)$, one has a canonical
$\D(S)$-isomorphism 
$$
\R\Homb_R( E \Otimes E'\<, F )\iso
 \R\Homb_S\bigl(E, \R\Hom_R^\bullet(E'\<, F)\bigr),\
$$
see \cite [p.\,147; 6.6]{Sp}. Thus, with
 $\fij\colon\D(R)\to \D(S)$ the functor given by
$$
\fij(-)\!:=\R\Homb_R(\R\GG {\mspace{-.5mu}J} S,-)
\cong
\R\Homb_S\bigl(\R\GG {\mspace{-.5mu}J}S,\R\Homb_R(S,-)\bigr),
$$
there is a composed isomorphism
$$
\R\Homb_S(E,\fij\<F) 
 \iso 
  \R\Homb_R(E\Otimes \R\GG {\mspace{-.5mu}J} S,  F )
\iso  
  \R\Homb_R (\R\GG {\mspace{-.5mu}J} E, F).
$$
Application of homology  $H^0$ yields the (rather trivial)
{\it local duality isomorphism}
\stepcounter{numb}
\begin{equation}
\label{(2.1)}
\Hom_{\D(S)}\<\<(E,\fij\<F)\iso 
\Hom_{\D(R)}\<\<(\R\GG {\mspace{-.5mu}J}E, F).
\end{equation}

``Non-trivial" versions of~\eqref{(2.1)} include more information
about $\fij$.
For example, Greenlees-May duality\index{Greenlees-May Duality} \cite[p.\,4,
$(0.3)_{\text{aff}}$]{AJL} gives a canonical isomorphism
\stepcounter{numb}
\begin{equation}\label{2.1.1}
\fij\<F \cong\LL_J\R\Hom_R^\bullet(S,F),
\end{equation}
where $\Lambda_J$ is the {\it $J\<$-adic completion functor,} and ${\bf L}$ denotes
``\kern.5pt left-derived.'' In particular, if $R$ is noetherian,
$S$ is a finite $R$-module, and 
$F\in\Dc(R)$ (i.e., each homology module of~$F$ is finitely generated),
then as in~\cite[p.\,6, Proposition (0.4.1)]{AJL},
\stepcounter{numb}
\begin{equation}
\fij\<F=\R\Hom_R^\bullet(S,F)\otimes_S{\hat S}
\qquad ({\hat S}=J \mbox{-adic completion of}\, S).
\label{2.1.2}
\end{equation}
More particularly, for  $S=R$ and $\varphi=\text{id}$ (the identity map)
we get
$$
\ush{\text{id}_{\<\<J}}F= F\otimes_R {\hat R}
\qquad \bigl(F\in\Dc(R)\bigr).
$$
Hence, {\it classical local duality} \cite[p.~278 (modulo Matlis dualization)]{H1} is
just~(\ref{(2.1)}) when $R$ is local, $\varphi=\text{id}$, 
$J$ is the maximal ideal of~$R$, and $F$ is a normalized dualizing complex---so
that, as in  ~\Cref{C:Hom-Rgamma}, and by~\cite[p.\,276,
Proposition~6.1]{H1},
$$
\Hom_{\D(R)}\<\<(\R\GG {\mspace{-.5mu}J}E, F)=
\Hom_{\D(R)}\<\<(\R\GG {\mspace{-.5mu}J}E, \R\GG {\mspace{-.5mu}J}F)=
\Hom_{\D(R)}\<\<(\R\GG {\mspace{-.5mu}J}E, I)
$$
where $I$ is an $R$-injective hull of the residue field~$R/\<J$.
(See also \Lref{L:dualizing}.)\looseness=-1

\smallskip
For another example, let $S=R[[{\bf t}]]$ 
where $\mathbf t\!:= (t_1,\dots,t_d)$ 
is a sequence of variables, and set $J\!:= \mathbf tS$. The
standard calculation (via Koszul complexes) gives an isomorphism
$\R\Gamma_ {\!\! J}S\cong\nu[-d\>]$ 
where $\nu$ is the free $R$-submodule of the 
localization~$S_{t_1\dots t_d}$ generated by 
those monomials~$t_1^{n_1}\!\dots t_d^{n_d}$ with all exponents~$n_i<0$, the
$S$-module structure being induced by that of 
$S_{t_1\<\dots\> t_d}/S\supset \nu\>$. The \emph{relative canonical module}
$\omega_{R[[\mathbf t]]/R}\!:=\Hom_R(\nu,R)$  is  a \emph{free, rank one,
$S$-module.}  There result,  for  finitely-generated \mbox{$R$-modules}~$F$,
functorial isomorphisms
\stepcounter{numb}
\begin{equation}\label{2.1.3}
\ush {\varphi_{\mathbf tR[[\mathbf t]]}}F\cong
\Hom_R(\nu[-d\>],F)\cong
\omega_{R[[\mathbf t]]/R}^{\phantom{.}}[d\>]\otimes_R F\cong
R[[\mathbf t]]\otimes_RF[d\>];
\end{equation}
and when $R$ is noetherian, the usual way-out argument 
\cite[p.\,69, (ii)]{H1} yields the same for any $F\in\Dc^+\<(\<R)$.

\smallskip

Next, we give a commutative-algebra analogue of \Tref{Th2} in
\S\ref{S:prelim}, in the form of a ``torsion" variant of the duality
isomorphism~\eqref{(2.1)}. \Pref{P:affine} will clarify the relation
between the algebraic and formal-scheme contexts.

With $\varphi\colon R\to S$ and $J$ an $S$-ideal as before, let
$\A_J(S)$ be the category of \mbox{$J$-torsion} $S$-modules, i.e.,
\mbox{$S$-modules}~$M$ such that 
$$
M=\GG{\<J}M:=\{\,m\in M\mid J^nm=0\textup{ for some }n>0\,\}.
$$
The derived category of $\A_J(S)$ is equivalent to the full
subcategory~$\D_{\!J\<}(S)$  of~$\D(S)$ with objects those
$S$-complexes~$E$  whose homology lies in~$\A_J(S)$, (or equivalently, such that
the natural map \smash{$\R\GG{\<J}E\to E$} is an isomorphism),
and the functor~$\R\GG{\<J}$ is right-adjoint to~the inclusion
$\D_{\!J\<}(S)\hookrightarrow\D(S)$  (cf.~\Pref{Gamma'(qc)} and its
proof).  
\goodbreak
\noindent Hence the functor~$\varphi_{\!J}^\times\colon\D(R)\to \D_{\!J\<}(S)$
defined by
$$
\varphi_{\!J}^\times(-)\!:=\R\GG{\<J}\R\Homb_R(S, -)
\cong \smash{\R\GG{\<J}S\Otimes \R\Homb_R(S, -)}
$$
is  right-adjoint to the natural composition
$\D_{\!J\<}(S)\hookrightarrow\D(S)\to\D(R)$: in fact, for
$E\in\D_{\!J\<}(S)$ and $F\in\D(R)$ there are natural isomorphisms
\stepcounter{numb}
\begin{equation}\label{2.2.1}
\R\Homb_S(E,\varphi_{\!J}^\times\<F)\iso
\R\Homb_S(E,\R\Homb_R(S, F))\iso
\R\Homb_R(E,F).
\end{equation}

Here is another interpretation of $\varphi_{\!J}^\times\<F$. For 
$S$-modules~$A$ and $R$-modules~$B$~set
$$
\Hom_{R,J}(A, B)\!:=\GG{\<J}\Hom_R(A,B),
$$
the $S$-module of $R$-homomorphisms $\alpha$ vanishing on $J^n\<A$
for some $n$ (depending on~$\alpha$), i.e.,  \emph{continuous} when $A$ is
$J$-adically topologized and $B$ is discrete. If $E$~is a K-flat
$S$-complex and $F$ is a K-injective $R$-complex, then
$\Homb_R(E,F)$ is  a K-injective $S$-complex; and it follows for all
$E\in\D(S)$ and $F\in\D(R)$ that
$$
\R\Hom_{R,J}^\bullet(E,F)\cong\R\GG{\<J}\R\Hom_R^\bullet(E,F).
$$
Thus,\vspace{-1pt}
$$
\varphi_{\!J}^\times\<F=\R\Hom_{R,J}^\bullet(S,F).
$$

\penalty -1000
A \emph{torsion version of local duality} is the  isomorphism,
derived from~\eqref{2.2.1}:
$$
\Hom_{\D_{\!J\<}(S)}\<\<\bigl(E,\>\R\Hom_{R,J}^\bullet(S,F)\bigr)
\iso
\Hom_{\D(R)}\<\<(E,F)
\quad\ \bigl(E\in\D_{\!J\<}(S),\;F\in\D(R)\bigr).
$$

\begin{small}
Apropos of \Rref{R:Th2}(1), the functors $\varphi_{\!J}^\times$ and $\fij$
are related by
\begin{alignat*}{2}
\LL_J\R\Hom_R^\bullet(S,F)\;
  &\underset{\lower.5ex\hbox
to0pt{\hss\scriptsize\eqref{2.1.1}\hss}}{\cong}\;\fij\<F
     &&\cong\LL_J\varphi_{\!J}^\times\<F, \\
\R\GG{\<J}\>\R\Homb_R(S,F)\,
  &=\,\varphi_{\!J}^\times\<F
   && \cong\R\GG{\>J}\fij\<F.
\end{alignat*}
The first relation is the case $E=\R\GG{\<J}\<S$ of \eqref{2.2.1}, followed by
Greenlees-May duality. 
The second results, e.g.,   from the sequence of natural isomorphisms, holding for
\mbox{$G\in\D_{\!J\<}(S)$},
\mbox{$E\in\D(S)$}, and $F\in\D(R)$:
\begin{align*}
\Hom_{\D(S)}\<\<\bigl(G,\> \R\GG{\<J}\R\Hom_R^\bullet(E,F)\bigr)
&\cong
\Hom_{\D(S)}\<\<\bigl(G,\> \R\Hom_R^\bullet(E,F)\bigr)\\
&\cong
\smash{\Hom_{\D(R)}\<\<(\R\GG{\<J}S\Otimes_{\!\!S}\,G
  \Otimes_{\!\!S}\>\>E,F)}\\ 
&\cong
\Hom_{\D(S)}\<\<\bigl(G,\> \R\Hom_R^\bullet(\R\GG{\<J}E,F)\bigr)\\
&\cong
\Hom_{\D(S)}\<\<\bigl(G,  \R\GG{\<J}\R\Hom_R^\bullet(\R\GG{\<J}E,F)\bigr),
\end{align*}
which entail that the natural map is an isomorphism
$$
\R\GG{\<J}\R\Hom_R^\bullet(E,F)\iso
\R\GG{\<J}\R\Hom_R^\bullet(\R\GG{\mspace{-.5mu}J}E,F).
$$
\end{small}

\setcounter{sth}{\value{numb}}

Local Duality theorems are often formulated, as in (c) of the following,
in terms of modules and local cohomology
($\Hr_{\<\<J}^\bullet\!:=\Hr^\bullet\R\GG {\mspace{-.5mu}J}$) rather than
derived categories.

\begin{sprop}
\label{(2.2)}
Let\/ $\varphi:R\to S$ be a homomorphism of noetherian 
rings, let\/
$J$~ be an\/ $S$-ideal, and suppose that there exists a
sequence\/ ${\bf u}=(u_1,\dots,u_d)$ in\/ $J$ such that\/  $S/{\bf u}S$ is\/
$R$-finite. Then for any\/ $R$-finite module\/ $F$\textup{:}

\smallskip

\textup{(a)} $\Hr^n\fij\<F=0$ for all\/  $n<-d,$ \  
so that there is a natural\/ $\D(S)$-map 
$$
h\colon (\Hr^{-d}\fij\<F)[d\>]\to \fij\<F.
$$

\smallskip\pagebreak[2]

\textup{(b)} If\/ $\tau_F^\pd\colon \R\GG {\mspace{-.5mu}J}\fij\<F\to F$
\vspace{1pt} 
corresponds in\/ {\rm (\ref{(2.1)})} to  the identity map of\/~$\fij\<F,\,$%
\footnote
{$\,\tau_F^\pd$ may be thought of as $\text{``evaluation at 1"}\<\colon
\R\Hom_{R,J}^\bullet(S,F)\to F$.%
}
and\/ $\int=\int_{\varphi\<,J}^d(F)$ is the composed map
$$
\R\GG {\mspace{-.5mu}J}(\Hr^{-d}\fij\<F)[d\>]
 \xrightarrow{\!\R\GG {\mspace{-.5mu}J}\<(h)\>}
\R\GG {\mspace{-.5mu}J}\fij\<F
\xrightarrow{\,\tau_F^\pd\>} F,
$$
then\/ $\bigl(\Hr^{-d}\fij\<F, \int \bigr)$  represents the
functor\/ $\Hom_{\D(R)}\<\<(\R\GG {\mspace{-.5mu}J}E[d\>],F)$
of\/ $S$-modules~$E$.

\smallskip

\textup{(c)}
If\/ $J\subset \,\root\of{{\bf u}S}\>$  then there is a bifunctorial isomorphism 
\textup(with $E$, $F$ as before\textup{):}
$$
\Hom_S(E,\Hr^{-d}\fij\<F) \iso \Hom_R(\Hr_{\<\<J}^dE,F).
$$
\end{sprop}

\begin{proof}
If $\mspace{6.5mu}\hat{}\mspace{-6.5mu}\varphi$ 
is the obvious map from $R$ to the
${\mathbf u}$-adic completion $\hat S$ of~$S$, then in~$\D(S)$, 
\mbox{$\fij F =\mspace{6.5mu} \hat{}\mspace{-6.5mu}\fij F$} 
since $\R\GG {\mspace{-.5mu}J}S=\R\GG {\mspace{-.5mu}J}\hat S$.    
In proving (a), therefore, we may  assume\vadjust{\kern.5pt}
 that $S$ is ${\bf u}$-adically complete, 
so that $\varphi$ factors as 
\smash{$R\stackrel{\psi\>}{\to}R[[{\bf t}]]\stackrel{\chi\>\>}{\to} S$}
with ${\bf t}=(t_1,\dots,t_d)$ a sequence of indeterminates and
$S$ finite over~$R[[{\bf t}]]$. ($\psi$ is the natural map, and $\chi(t_i)=u_i\>$.)
In view of the easily-verified  relation 
$\fij = \ush{\chi_{\<\<J}}\smcirc\ush{\psi_{{\bf t}R[[{\bf t}]]}}$,
(\ref{2.1.2}) and~(\ref{2.1.3}) yield (a).  
Then (b) results from the natural isomorphisms
$$
\Hom_S(E, \Hr^{-d}\fij\<F)
\underset{\text{via } h }{\iso}
\Hom_{\D(S) }(E[d\>], \fij\<F)
\underset{\text{\eqref{(2.1)}}}{\iso}
\Hom_{\D(R)}\<\<(\R\Gamma_{\!\! J}^{\phantom{.}}E[d\>], F).
$$
Finally, (c) follows from (b) because
$\Hr^i_{\<\<J}E=\Hr^i_{\<\mathbf uS}E=0$ for all $i>d$
(as one sees from the usual calculation of $\Hr^i_{\<\mathbf uS}E$ via
Koszul complexes), so that the natural map is an isomorphism
$
\Hom_{\D(R)}\<\<(\R\Gamma_{\!\! J}^{\phantom{.}}E[d\>], F)
  \iso\Hom_R(\Hr_{\<\<J}^dE,F).
$
\end{proof}

\begin{sparag}\label{HuK}
Under the hypotheses of \Pref{(2.2)}(c),  the
functor~$\Hom_R(\Hr^d_{\<\<J}E,R)$ of $S$-modules~$E$ is representable.  
Under suitable extra conditions (for example,
$\hat S$~a generic local complete intersection over~$R[[{\bf t}]]$, 
H\"ubl and Kunz represent this functor by a \emph{canonical} pair
described explicitly via differential forms, residues, and certain trace maps
\cite[p.~73, Theorem~3.4]{Integration}. 
For example, with \mbox{$S=R[[{\bf t}]]$,}
$J=\mathbf tS$, and $\nu$~as in~(\ref{2.1.3}), the $S$-homomorphism
from the module~\smash{$\widehat\Omega_{S/R}^d$} of universally finite
\hbox{$d$-forms} to the relative canonical module 
$\>\omega_{R[[\mathbf t]]/R}^\pd\<=\<\Hom_R(\nu,R)$ sending the
form
$dt_1\dots dt_d$ to the $R$-homomorphism $\nu\to R$ which takes the monomial~
$t_1^{-1}\! \dots t_d^{-1}$ to~$1$ and all other monomials
$t_1^{n_1}\!\dots t_d^{n_d}$ to $0$, is clearly an isomorphism; and  the
resulting isomorphism~$\smash{\widehat\Omega_{S/R}^d}[d\>]\iso\fij R$
\emph{does not depend on the $d$-element
sequence\/~$\mathbf t$ generating} $J$---it corresponds under~(\ref{(2.1)}) to
the {\it residue map }
$$
\R\GG {\mspace{-.5mu}J}\widehat\Omega_{S/R}^d[d\>]=
 \Hr^d_{\!J}\>\widehat\Omega_{S/R}^d\to R
$$
(see, e.g., \cite[\S2.7]{Hochschild}).  Thus  $\Hom_R(\Hr_{\<\<J}^dE,R)$
\emph{is represented by\/ \smash{$\widehat\Omega_{S/R}^d$} together with the
residue map.}  The general case reduces to this one via traces of
differential forms.\index{Duality!Local|)}

\end{sparag}

\end{parag}

\smallskip
\begin{parag}
\label{sheafify}\index{Duality!Local, sheafified|(} 
(Formal sheafification of Local Duality). For 
$f\colon\X\to\Y$ as in \Tref{Th2} in
\S\ref{S:prelim}, there is  a right $\Delta$-adjoint~$\ft\<$
for  the functor\/ $\Rfs\colon\Dqct(\X)\to\Dqct(\Y)$. Furthermore, with
$\>\bj\colon\D(\Avc(\X))\to\Dvc(\X)$\index{ $\iG{\<\cJ\>}$@$\bj$} the canonical
functor, we have
$$
\Rfs\R\iGp\X\bj\>\D(\Avc(\X))\!
\!\underset{\mathstrut(\text{\ref{C:vec-c is qc})}}\subset\!
\!\Rfs\R\iGp\X\>\Dqc(\X)\!
\!\underset{\mathstrut(\text{\ref{Gamma'(qc)})}}\subset\!
\!\Rfs\Dqct(\X)\!
\!\underset{\mathstrut(\text{\ref{Rf-*(qct)})}}\subset\!
\!\Dqct(\Y)\!
\!\underset{\mathstrut(\text{\ref{C:limsub})}}\subset\!
\!\Dvc(\Y).\vspace{-1pt}
$$
It results from \eqref{adj} and  \Pref{A(vec-c)-A} that
$\Rfs\R\iGp\X\>\bj\colon\D(\Avc(\X))\to\Dvc(\Y)$ has the right
$\Delta$-adjoint~$\R Q_\X^{}\ush f\!:=\R
Q_\X^{}\R\sHomb(\R\iGp\X\cO_\X^{}\>,\ft$).

\pagebreak[4]
If, moreover, $\X$ is \emph{properly algebraic} (\Dref{D:propalg})---in particular,
if $\X$ is affine---then $\bj$ is an equivalence of categories (\Cref{corollary}),
and so the functor $\Rfs\R\iGp\X\colon\Dvc(\X)\to\Dvc(\Y)$ has a right
$\Delta$-adjoint.

For \emph{affine $f\<$}, these results are closely related to
the Local Duality isomor\-phisms \eqref{2.2.1} and ~\eqref{(2.1)}.  Recall that an
\emph{adic ring}\index{adic ring} is a pair $(R,I)$ with $R$ a ring and $I$ an
$R$-ideal such that with respect to the $I\<$-adic topology $R$ is Hausdorff and
complete. The topology on~$R$ having been specified, 
the corresponding affine formal scheme is denoted  $\Spf(\<R)$.

\begin{sprop}\label{P:affine}
Let\/ $\varphi\colon(\<R,I\>)\to(S,J\>)$ be a continuous homomorphism of
noetherian adic rings, and let\/
$\X\!:=\Spf(S)\;\smash{\stackrel{f}{\to}}\;\Spf(\<R)=:\!\Y$ be the corresponding
\textup(affine\/\textup) formal-scheme map. Let\/ $\kappa_\X^\pd\colon\X\to
X\!:=\Spec(S),$ 
$\kappa_\Y^\pd\colon\Y\to Y\!:=\Spec(R)$ be the
completion maps,  and let $^\sim={}^{\sim_S}$ denote
the standard exact functor from
$S$-modules to quasi-coherent $\cO_{\<\<X}$-modules. Then:
 
\smallskip
\textup{(a)} The restriction of\/
$\R f_{\!*}$ takes\/ $\Dqct(\X)$ to\/ $\Dqct(\Y),$\
and this restricted functor has a right adjoint 
$\ft\colon\Dqct(\Y)\to\Dqct(\X)$ given by
$$
\ft\<\F
\!:=\kappa_\X^* \bigl(\varphi_{\!J}^\times\<\R\Gamma(\Y,\>\F\>)\bigr)^\sim
\<= \kappa_\X^*\bigl(\R\Hom_{R,J}^\bullet\bigl(S,
    \>\R\Gamma(\Y,\>\F\>)\bigr)\bigr)^\sim
\qquad\bigl(\F\in\Dqct(\Y)\bigr).
$$
 
\textup{(b)} The restriction of\/
$\R f_{\!*}\>\R\iGp\X$ takes\/ $\Dvc(\X)$ to\/ $\Dvc(\Y),$\
and this restricted functor has a right adjoint 
$\ush{f_{\<\vec{\mathrm c}}}\colon\Dvc(\Y)\to\Dvc(\X)$ given by
$$
\ush{f_{\<\vec{\mathrm c}}}\F
\!:=\kappa_\X^*\bigl(\fij\R\Gamma(\Y,\>\F\>)\bigr)^\sim
\<= \kappa_\X^*\bigl(\R\Hom_R^\bullet\bigl(\R\GG {\mspace{-.5mu}J}S,
 \>\R\Gamma(\Y,\>\F\>)\bigr)\bigr)^\sim
\qquad\bigl(\F\in\Dvc(\Y)\bigr).
$$

\textup{(c)} There are natural isomorphisms
\begin{align*}
\R\Gamma(\X,\ft\<\F\>)&\iso \varphi_{\!J}^\times\R\Gamma(\Y,\>\F\>)
\qquad\bigl(\F\in\Dqct(\Y)\bigr),\\
\R\Gamma(\X,\ush{f_{\<\vec{\mathrm c}}}\F\>)&\iso \fij\>\R\Gamma(\Y,\>\F\>)
\qquad\bigl(\F\in\Dvc(\Y)\bigr).
\end{align*}

\end{sprop}

\begin{proof}
The functor~${}^\sim$ induces an equivalence of categories $\D(S)\to\Dqc(X)$,
 with quasi-inverse $\R\GG X\!:=\R\Gamma(X,-)$
(\cite[p.\,225, Thm.\,5.1]{BN}, \cite[p.\,12, Proposition~(1.3)]{AJL}); 
and \Pref{c-erator} below implies that
$\kappa_\X^*\colon\Dqc(X)\to\Dvc(\X)$ is~an 
equivalence, with quasi-inverse
$(\R\GG X\kappa_{\X*}^\pd-)^\sim=(\R\GG\X-)^\sim\<$.\kern2pt
\footnote{In checking this note that $\kappa_{\X\<*}^\pd$ 
has an exact left adjoint, hence preserves K-injectivity.%
}

It follows that \emph{the functor taking\/ \mbox{$G\in\D(S)$} to\/
$\kappa_\X^*\widetilde G$ is an
equivalence, with quasi-inverse \smash{$\R\GG\X\colon\Dvc(\X)\to\D(S)$}}, and
similarly for~$\Y$~and~$R$.  Moreover, there is an induced equivalence between
$\D_{\!J\<}(S)$ and $\Dqct(\X)$  (see \Pref{Gammas'+kappas}).
In particular, (c) follows from (a) and~(b).\looseness =-1

Corresponding to~\eqref{2.2.1} and ~\eqref{(2.1)} there are then 
functorial isomorphisms
\begin{gather*}
\Hom_{\D(\X)}\<\<(\<\E\<,\ft\<\F\>)\!\<\iso\!\< 
\Hom_{\D(\Y)}\<\<(\<\kappa_\Y^*(\R\GG\X\E)^{\sim_R}\<\<,\>\F\>)
\qquad\bigl(\<\E\<\in\<\Dqct(\X),\;\F\<\in\<\Dqct(\Y)\bigr)\<,\\
\Hom_{\D(\X)}\<\<(\<\E\<,\ush{f_{\<\vec{\mathrm c}}}\F\>)\!\<\iso\!\< 
\Hom_{\D(\Y)}\<\<(\<\kappa_\Y^*(\R\GG{\<J}\R\GG\X\E)^{\sim_R}\<\<,\>\F\>)
\qquad\bigl(\<\E\<\in\<\Dvc(\X),\medspace \>\F\<\in\<\Dvc(\Y)\bigr)\<;
\end{gather*}
and it remains to demonstrate functorial isomorphisms
\begin{alignat*}{2}
\kappa_\Y^*(\R\GG\X\E)^{\sim_R}
&\iso \R f_{\!*}\E
&&\qquad\bigl(\E\in\Dqct(\X)\bigr), \\
\kappa_\Y^*(\R\GG {\mspace{-.5mu}J}\R\GG\X\E)^{\sim_R}
&\iso \R f_{\!*}\>\R\iGp\X\E
&&\qquad\bigl(\E\in\Dvc(\X)\bigr),
\end{alignat*}
the first  a special case of the second. 

\pagebreak[3]
To prove the second,
let $E\!:=\R\GG\X\E$, let $Z\!:=\Spec(S/J\>)\subset X$, 
and let $f_0\colon X\to Y$ be
the scheme-map corresponding to~$\varphi$. The desired isomorphism
comes from the sequence of natural isomorphisms
\begin{align*}
\R f_{\!*}\>\R\iGp\X\E 
&\cong
 \R f_{\!*}\>\R\iGp\X\kappa_\X^*\widetilde E \\
&\cong
 \R f_{\!*}\kappa_\X^*\R\iG Z\widetilde E
   &&\qquad(\textup{\Pref{Gammas'+kappas}(c)})\\
&\cong
 \kappa_\Y^*\R f_{0*}\R\iG Z\widetilde E
   &&\qquad(\textup{\Cref{C:kappa-f*t}})\\
&\cong
 \kappa_\Y^*\R f_{0*}(\R\GG {\mspace{-.5mu}J}E)^\sim
   &&\qquad(\textup{\cite[p.\,9, (0.4.5)]{AJL}})\\
&\cong
 \kappa_\Y^*(\R\GG {\mspace{-.5mu}J}E)^{\sim_R}.
 \end{align*}
(The last isomorphism---well-known for bounded-below~$E$---can be checked
via~the equivalences $\R\GG X$ and $\R\GG Y\>$, which satisfy
$\R\GG Y\R f_{0*}\cong\R\GG X$ (see \cite[pp.\:142--143,
5.15(b) and~5.17]{Sp}).
\end{proof}

\medskip

\def\xf{f_{\<\<\X*}'}

\Tref{T:pf-duality} below globalizes \Pref{(2.2)}.\index{Duality!Pseudo-finite|(}
But first some preparatory remarks are needed. Recall from~\ref{maptypes}
that a map $f\colon\X\to\Y$ of noetherian formal schemes is 
\emph{pseudo\kern.6pt-finite}\index{formal-scheme
map!pseudo\kern.6pt-finite} if it is pseudo\kern.5pt-proper and has finite
fibers, or equivalently, if $f$~is pseudo\kern.5pt-proper and 
\emph{affine}. Such an~$f$ corresponds locally to a
homomorphism \mbox{$\varphi\colon(R,I\>)\to (S,J\>)$} of noetherian
adic rings such that $\varphi(I)\subset J$ and
$S/\<J$ is a finite \hbox{$R$-module.} This
$\varphi$ can be extended to a homomorphism from a power series ring 
$R[[\mathbf t]]\set R[[t_1, t_2,\dots,t_e]]$ such that the images of the
variables~$t_i$ together with~$\varphi(I\>)$ generate~$J\<$, 
and thereby $S$~becomes a finite $R[[\mathbf t]]$-module. Pseudo-finiteness is
preserved under arbitrary (noetherian) base change.

We say that a pseudo\kern.6pt-finite map $f\colon\X\to\Y$ of noetherian formal
schemes has relative dimension $\le d$ if each $y\in\Y$ has an
affine neighborhood~$\U$ such that the map~ $\varphi_\U\colon R\to  S$
of adic rings corresponding to $f^{-1}\U\to\U$ has
a continuous extension $R[[t_1,\dots,t_d]]\to S$ making $S$ into a
finite $R[[t_1,\dots,t_d]]$-module, or equivalently, there is a
topologically nilpotent sequence $\mathbf u=(u_1,\dots,u_d)$ in~$S$
(i.e., $\lim_{n\to\infty}u_i^n=0\ (1\le i\le d)$) such that $S/\mathbf uS$ is
finitely generated as an $R$-module. The \emph{relative dimension}~$\dim f$
 is defined to be the least among the integers~$d$ such that $f$ has
relative dimension $\le d$.

For any pseudo\kern.6pt-proper map $f\colon\X\to\Y$ 
of noetherian formal schemes, we have the functor~$\ush f\colon\D(\Y)\to\D(\X)$
of \Cref{C:f*gam-duality}, commuting  with open base change on~$\Y$
(\Tref{Th4}). The next Lemma is a special case of \Pref{P:coherence}.

\begin{slem}\label{L:coh}
For a  pseudo\kern.6pt-finite map \/ $f\<\colon\!\X\to\<\Y$  of noetherian formal
schemes  and for any\/
$\F\in\Dc^+(\Y),$\ it holds that
$\ush f\<\F\in\Dc^+(\X)$.
\end{slem}

\begin{proof}
Since $\ush f$ commutes with open base change, the question is local,  so we
may assume that $f$ corresponds to \mbox{$\varphi\colon(R,I\>)\to (S,J\>)$} as
above. Moreover,  the isomorphism $\ush{(gf)}\cong\ush f \<\ush g$ in
\Cref{C:f*gam-duality} allows us to assume  that \emph{either} 
\mbox{$S=R[[t_1,\dots,t_d]]$} and $\varphi$ is the natural map \emph{or} 
$S$ is a finite
$R$-module and $J=IS$. In either case $f$ is obtained by completing a proper map
$f_0\colon X \to \Spec(R)$ along a closed subscheme $Z\subset
f_0^{-1}\Spec(R/I\>)$. (In the first case, take $X$ to be the projective
space~$\mathbb P_{\!\!R}^{\mspace{.5mu}d}\supset\Spec(R[t_1,\dots,t_d])$, and
$Z\set\Spec(R[t_1,\dots,t_d]/(I,t_1,\dots,t_d))$.) The conclusion  is
given then by \Cref{C:completion-proper}.
\end{proof}

\pagebreak[3]
\begin{sth}[\textup{Pseudo-finite Duality}]\label{T:pf-duality}%
\index{Duality!Pseudo-finite|)}
 Let\/ $f\colon \X\to\Y$ be a
pseudo\kern.5pt-finite map of  noetherian formal schemes, 
and let\/ $\F$ be a coherent\/ $\cO _\Y$-module. Then:

\textup{(a)} $H^n\ush f\<\F=0$ for all\/ $n<-\dim f$.

\textup{(b)} If\/ $\dim f\le d$ and\/ $\X$ is covered by affine open
subsets with\/ \hbox{$d$-generated} defining ideals, then with\/
$\xf\set f_{\!*}\iGp\X$ and, for\/ $i\in\mathbb Z$ and\/ $\J$ a defining ideal
of~$\X,$
$$
R^i\<\<\xf\set  H^i\R\xf= H^i\Rfs\>\R\iGp\X
=\smash{\dirlm{n} H^i\Rfs\>\R\sHomb(\cO_\X/\J^n,-)},\footnotemark
$$%
\footnotetext
{The equalities hold because $\X$ being noetherian,  
any \smash{$\subdirlm{}\!\!$}\vspace{.8pt} of flasque sheaves (for example,
\smash{$\subdirlm{}\<\sHom(\cO_\X/\J^n,\E)$}\vspace{.8pt} with $\E$ an injective
$\cO_\X$-module) is
$f_{\!*}$-acyclic, and \smash{$\subdirlm{}\!\!$} commutes
with~$f_{\!*}$. (For an additive functor $\phi\colon\A(\X)\to\A(\Y)$,
an $\A(\X)$-complex $\F$ is \mbox{\emph{$\phi$-acyclic}}\index{acyclic} if the
natural map
$\phi\F\to\R\phi\F$ is a $\D(\Y)$-isomorphism. Using a standard spectral
sequence, or otherwise (cf.~\cite[(2.7.2)]{Derived categories}), one sees that any
bounded-below complex of $\phi$-acyclic $\cO_\X$-modules is $\phi$-acyclic. }%
there is, for quasi-coherent $\cO _\X$-modules\/~$\E\<,$ a functorial
isomorphism
$$
\postdisplaypenalty 10000
f_{\!*}\sHom_{\X}(\E\<, \>H^{-d}\ush f \<\F\>)\iso
\sHom_{\Y}(R^d\<\<\xf\E\<,\>\F\>).
$$
\textup(Here
 $H^{-d}\ush f\<\F$ is coherent \textup(\Lref{L:coh}\kern1pt\textup{),} and
by \textup{(a),} vanishes unless $d=\dim f$.\textup)
\end{sth}

\begin{proof} 
Since $\ush f$ commutes with open base change we may assume that $\Y$ is affine
and  that $f$ corresponds to a map $\varphi\colon(R,I\>)\to(S,J\>)$ as in
\Pref{(2.2)}. Then there is an isomorphism of functors
$$
\bj\R Q_\X^{}\ush f\cong \kappa_\X^*\bigl(\fij\R\Gamma(\Y,-)\bigr)^\sim\<,
$$
both of these functors being right-adjoint to
$\Rfs\>\R\iGp\X\colon\Dvc(\X)\to\Dvc(\Y)$ (\Pref{P:affine}(b) and 
remarks about right adjoints preceding it). Since
\mbox{$\ush f\<\F\in\Dc^+(\X)$} (\Lref{L:coh}),  therefore, by
\Cref{corollary}, the natural map is an isomorphism $\bj\R Q_\X^{}\ush
f\<\F\iso\ush f\<\F\<$; and so, since 
$\kappa_\X^*$  is exact, \Pref{(2.2)} gives~(a).

\smallskip\enlargethispage{-\baselineskip}
Next, consider the presheaf map  associating to each open
$\U\subset\Y$ the natural composition (with $\V\set f^{-1}\U$):\vspace{1.5pt}
\begin{align*}
\smash{\Hom_\V(\E\<, \>H^{-d}\ush f \<\F\>)
\underset{\textup{by (a)}}{\iso}
\Hom_{\D(\V)}\<\<(\E[d\>], \>\ush f \<\F\>)}
&\underset{\ref{C:f*gam-duality}}\iso
\Hom_{\D(\U)}\<\<(\Rfs\>\R\iGp\X\E[d\>], \F\>)
\vspace{1.5pt}\\
&\,\lra\,
\Hom_{\U}(R^d\<\<\xf\E\<,\>\F\>).
\end{align*}

\penalty -1000

\noindent To prove (b) by showing that the resulting sheaf map 
$$
f_{\!*}\sHom_{\X}(\E\<, \>H^{-d}\ush f\<\F\>)\to
\sHom_{\Y}(R^d\<\<\xf\E\<,\>\F\>)
$$
is an isomorphism, it suffices to show that  $R^i\<\<\xf\E=0$ for all~$i>d$,
a local problem for which we can (and do)  assume that $f$ corresponds to
$\varphi\colon R\to S$ as above.

 Now $\R\iGp\X\E\in\Dqct(\X)$ (\Pref{Gamma'(qc)}), so 
\Pref{Gammas'+kappas} for~\mbox{$X\set\Spec(S)$} and
$Z\set\Spec(S/J)$ gives
$\R\iGp\X\E\cong\kappa_\X^*\E_0$ with
$\E_0\set\kappa_{\X*}^\pd\R\iGp\X\E\in\DqcZ^+(X)$.\vadjust{\kern.5pt} Since $\X$
has, locally, a
$d$-generated defining ideal, we can represent~$\R\iGp\X\E$ locally  by a
\smash{$\dirlm{}\!\!$}\vadjust{\kern1pt} of   Koszul complexes on $d$ elements
\cite[p.\,18, Lemma 3.1.1]{AJL}, whence $H^i\R\iGp\X\E=0$ for all~$i>d\>$,
and so, $\kappa_{\X*}^\pd$ being exact, 
$H^i\E_0=0$. Since the map
$f_0\set\Spec(\varphi)$ is affine, it follows that  $H^i\R f_{0*}\E_0=0$,
whereupon,
$\kappa_\Y$ being flat, \Cref{C:kappa-f*t}  yields\looseness=-1
$$
\postdisplaypenalty 10000
R^i\<\<\xf\E
\cong
H^i\Rfs\kappa_\X^*\E_0
\cong 
H^i\kappa_\Y^*\R f_{0*}\E_0
\cong 
\kappa_\Y^*H^i\R f_{0*}\E_0=0\qquad(i>d),
$$
as desired. (Alternatively, use Lemmas~\ref{affine-maps} and~\ref{Gamma'+qc}.)
\end{proof}
\end{parag}\index{Duality!Local, sheafified|)}

\smallskip

\begin{parag}\label{residue thm}\index{Residue theorems|(}

Our results provide a framework for ``Residue Theorems" such
as  those appearing in \cite[pp.~87--88]{Asterisque} and \cite[pp.~750-752]{HS}
(central theorems in those papers): roughly speaking, Theorems~\ref{Th1}
and~\ref{Th2} in section~\ref{S:prelim} include both  local and
global duality, and \Cref{C:kappa-f^times-tors} expresses the
compatibility between these dualities. But the dualizing objects we deal with are
determined \vadjust{\penalty-1000}only up to isomorphism. 
The Residue Theorems run deeper
in that they include a \emph{canonical realization} of dualizing data, via
differential forms. (See the above remarks on the H\"ubl-Kunz treatment of local
duality.) This extra dimension belongs properly to a theory of the ``Fundamental
Class" of a morphism, a canonical map from relative differential forms to the
relative dualizing complex, which will be pursued in a separate paper.

\begin{sparag}

Let us be more explicit, starting with some remarks about ``Grothendieck Duality
with supports" for a map
$f\colon X\to Y$ of noetherian separated schemes  with respective closed
subschemes
$W\subset Y$ and $Z\subset f^{-1}W$. 
Via the natural  equivalence of categories $\D(\Aqc(X))\to\Dqc(X)$ (see
\S\ref{SS:Dvc-and-Dqc}), we regard the
functor~$f^\times\colon\D(Y)\to\D(\Avc(X))=\D(\Aqc(X))$ of \Tref{Th1}
as being right-adjoint to $\Rfs\colon\Dqc(X)\to\D(Y)$.%
\footnote
{For ordinary schemes, this functor $f^\times$ is well-known, and usually
denoted $f^!$ when $f$ is proper. When $f$ is an open immersion, the functors
$f^\times$ and $f^! (=f^*)$ need not agree.%
}
 The functor $\R\iGp Z$ can be regarded as being right-adjoint to the
inclusion $\D_{\!Z}(X)\hookrightarrow\D(X)$ (cf.~\Pref{Gamma'(qc)}(c));
and its restriction  to~$\Dqc(X)$ agrees naturally with that of~$\R\iG Z\>$, both
restrictions being \hbox{right-adjoint} to the inclusion
$\DqcZ(X)\hookrightarrow\Dqc(X)$. Similar statements hold for $W\subset Y\<$. 
Since $\Rfs(\DqcZ(X))\subset \D_W\<(Y)$
(cf.~ proof of \Pref{Rf-*(qct)}),
we find that the functors $\R\iG Zf^\times$ and~$\R\iG Zf^\times\R\iGp W$ 
are both right-adjoint to  $\Rfs\colon\DqcZ(X)\to \D(Y)$,  so are
isomorphic.  We
define the \emph{local integral} (a generalized residue map, 
cf.~\cite[\S4]{Integration})
$$
\rho(\G)\colon\Rfs\>\R\iG Z f^\times\G\to \R\iGp W\G
\qquad \bigl(\G\in\D(Y)\bigr)
$$
to be the natural composition
$$
\Rfs\>\R\iG Z f^\times\G
\iso
\Rfs\>\R\iG Zf^\times\R\iGp W\G
\to
\Rfs f^\times\R\iGp W\G
\to
\R\iGp W\G.
$$

Noting that for $\F\in\D_W\<(Y)$ there is
a canonical isomorphism $\R\iGp W\F\iso\F$ (proof similar to that of
\Pref{Gamma'(qc)}(a)), we have then:

\begin{sprop}[Duality with supports]\label{P:dual/supports}%
\index{Grothendieck Duality!with supports}
For\/ $\E\in\DqcZ(\<X),$ $\F\<\in\<\D_W\<(Y),$ the natural
composition
\begin{align*}
\Hom_{\DqcZ(X)}\<\<(\E\<, \R\iG Z f^\times\<\F\>)\:
&\smash{\xrightarrow[\phantom{\rho(\F\>)}]{}}
\:\Hom_{\D_W\<(Y)}\<\<(\Rfs\E\<, \Rfs\>\R\iG Zf^\times\F\>)\\
&\smash{\xrightarrow[\rho(\F\>)]{}}
\:\Hom_{\D_W\<(Y)}\<\<(\Rfs\E\<, \>\F\>)\vspace{-3pt}
\end{align*}
is an isomorphism.  
\end{sprop}
This follows from adjointness of $\Rfs$ and $f^\times\<$, via the natural
diagram 
$$
\begin{CD}
\Rfs\>\R\iG Z f^\times\G @>>>\Rfs f^\times\G\\
@V\rho(\G) VV @VVV \\
\R\iGp W\G @>>>\G
\end{CD}
\qquad\bigl(\G\in\D(Y)\bigr)\<,
$$
whose commutativity is a cheap version of the Residue Theorem 
\cite[pp.~750-752]{HS}.  

Again, however, to be worthy of the name a Residue Theorem should involve
\emph{canonical  realizations} of dualizing objects.
 For instance, when $V$ is a proper \hbox{$d$-dimensional} variety 
over a field $k$ and
$v\in V$ is a closed point, taking $X=V\<$, $Z=\{v\}$, $W=Y=\Spec(k)$,  $\G=k$, and
setting $\omega_V\!:=H^{-d}f^\times k$, we get an
 $\cO_{V\!,\>v}$-module~$\omega_{V\!,\>v}$ (commonly called
``canonical", though defined only up to isomorphism) together with the
$k$-linear map induced by $\rho(k)$:
\enlargethispage*{\baselineskip}
$$
 H_{\<v}^d(\omega_{V\!,\>v})\to k,
$$
a map whose truly-canonical realization via differentials and  residues is indicated
in  \cite[p.\,86,~(9.5)]{Asterisque}.

\pagebreak[3]

\end{sparag}

\newcommand{\Rhfs}{{\R\hat f_{\!*}}}

\begin{sparag}
\label{completion}
With preceding notation, consider the completion diagram
$$
\begin{CD}
X_{/Z}=:\,@.\X @>\kappa_\X^\pd>> X \\
@. @V \hat f VV @VV f V \\
Y_{\</W}=:\,@.\Y @>>\vbox to
0pt{\vskip-1.3ex\hbox{$\scriptstyle\kappa_\Y^\pd$}\vss}> Y 
\end{CD}
$$
\vskip-3pt
Duality with supports can be regarded\vspace{.6pt} 
more intrinsically---via $\hat f\<$ 
rather than~$f$---as a special case of the\index{Grothendieck
Duality!Torsion (global)} Torsion-Duality
\Tref{T:qct-duality} ($\>\cong\:$\Tref{Th2} of \S1) for~$\hat f$:

First of all, the local integral $\rho$
is completely determined by
$\kappa_\Y^*(\rho)$: for $\G\in\D(\Y)$, the natural map 
$\R\iGp W\G\to \kappa_{\Y*}^\pd\kappa_\Y^*\R\iGp W\G$ is an isomorphism
(\Pref{Gammas'+kappas});  and the same holds for
$\Rfs\>\R\iG Z f^\times\G\to 
\kappa_{\Y*}^\pd\kappa_\Y^*\Rfs\>\R\iG Z f^\times\G$
since as above, 
$$
\Rfs\>\R\iG Z f^\times\G\in\Rfs(\DqcZ(X))\subset\D_W\<(Y)
$$
---and so $\rho=\kappa_{\Y*}^\pd\kappa_\Y^*(\rho)$.
Furthermore, $\kappa_\Y^*(\rho)$ is determined by the ``trace" map 
$\tau_{\<\mathrm t}^{}\colon \Rhfs\hat\ft\to\mathbf 1$,
\index{ {}$\tau$ (trace map)!$\tau_{\<\mathrm t}$} as per the following
natural commutative diagram, whose rows are isomorphisms:
\begin{small}
$$
\minCDarrowwidth=18pt
\begin{CD}
\kappa_\Y^*\Rfs\>\R\iG Z f^\times\G@>\Iso>\ref{C:kappa-f*t} >
\Rhfs\kappa_\X^*\R\iG Z f^\times \kappa_{\Y*}^\pd\kappa_\Y^*\G
 @>\Iso>\ref{C:kappa-f^times-tors}>
  \Rhfs\hat\ft\kappa_\Y^*\G
   @>\Iso>\hbox to
16pt{$\scriptstyle\hss\ref{C:identities}\textup{(b\kern-.6pt)}\hss$}>
    \Rhfs\hat\ft\R\iGp\Y\kappa_\Y^*\G\\
\vspace{-21pt}\\
@V\kappa_\Y^*(\rho) VV @.  @.  @VV\tau_{\<\mathrm t}^\pd V\\
\kappa_\Y^*\R\iGp W\G
@. \hbox to 0pt{\kern75.5 pt
                         \hss$
                         \overset{\Iso}
                                      {\underset{\ref{Gammas'+kappas}}
                                                      {\hbox to 215.5 pt{\leftarrowfill}}
                                       }     
                          $\hss
                         }
 @. @.
         \R\iGp \Y\kappa_\Y^*\G  
\end{CD}
$$
\end{small}
(To see that the natural map
$\R\iG Zf^\times\G\to\R\iG Zf^\times\kappa_{\Y*}^\pd\kappa_\Y^*\G$  
is an isomorphism, replace~
$\R\iG Zf^\times$ by the isomorphic functor~$\R\iG Zf^\times\R\iGp W$ and apply
\Pref{Gammas'+kappas}.)

Finally, we have isomorphisms (for $\E\in\DqcZ(X)$, $\F\in\D_W\<(Y)$),
\begin{alignat*}{2}
\Hom_{\D(X)}\<\<(\E\<, \R\iG Z f^\times\<\F\>)
&\iso 
\Hom_{\D(\X)}\<\<(\kappa_\X^*\E\<, 
    \>\kappa_\X^*\R\iG Z f^\times\kappa_{\Y*}^\pd\kappa_\Y^*\F\>)\qquad 
&&(\ref{Gammas'+kappas})\\
&\iso
\Hom_{\D(\X)}\<\<(\kappa_\X^*\E\<,\> \smash{\hat\ft\<} \kappa_\Y^*\F\>)\qquad
&&(\ref{C:kappa-f^times-tors}) \\
&\iso\Hom_{\D(\Y)}\<\<(\R\smash{\hat f_{\!*}}\kappa_\X^*\E\<,\>
\kappa_\Y^*\F\>)\qquad &&(\ref{T:qct-duality}) \\
&\iso\Hom_{\D(\Y)}\<\<(\kappa_\Y^*\Rfs\E\<,\>
\kappa_\Y^*\F\>)\qquad &&(\ref{C:kappa-f*t}) \\
&\iso\Hom_{\D(Y)}\<\<(\Rfs\E\<,\>\F\>)\qquad &&(\ref{Gammas'+kappas}),
\end{alignat*}
whose composition can be checked, via the preceding diagram, to be the
same as the isomorphism of \Pref{P:dual/supports}.

\end{sparag}

\def\Hp#1#2{\Hr_{#1}^{\<\prime\>#2}}

\begin{sparag}\label{consequences}
\Pref{P:astrix10.2} expresses some homological consequences of
the foregoing dualities, and furnishes a general context for 
\cite[pp.~87--88, Theorem (10.2)]{Asterisque}.\looseness=-1

For any noetherian formal scheme~$\X$, $\E\in\D(\X)$, and $n\in\mathbb Z$, set
$$
\Hp\X n(\E)\!:=\Hr^n\R\Gamma(\X,\R\iGp\X\E).
$$

For instance, if 
$\X=X_{/Z}\xrightarrow{\,\kappa\,} X$ is the completion of a noetherian
scheme~$X$ along a closed~$Z\subset X$, then for
$\F\in\D(X)$, \Pref{Gammas'+kappas} yields natural isomorphisms
\begin{align*}
\R\Gamma(\X,\R\iGp\X\kappa^*\<\F\>)
  &= \R\Gamma(X,\kappa_*\R\iGp\X\kappa^*\<\F\>) \\
  &\cong \R\Gamma(X,\kappa_*\kappa^*\R\iGp Z\F\>)
    \cong \R\Gamma(X,\R\iGp Z\F\>),
\end{align*}
and so if $\F\in\Dqc(X)$, then with $\Hr_Z^\bullet$ the usual cohomology
with supports in~$Z$, 
$$
\Hp\X n(\kappa^*\<\F\>)\cong\Hr_Z^n(\F\>).
$$

Let $\J\subset\cO_\X$ be an ideal of definition.  Writing
$\Gamma_{\!\X}^{}$ for the functor $\Gamma(\X,-)$, we have a functorial
map
$$
\gamma(\E)\colon\R(\Gamma_{\!\X}^{}\<\smcirc\iGp\X)\E\to
 \R\Gamma_{\!\X}^{}\<\smcirc\R\iGp\X\>\E\qquad \bigl(\E\in\D(\X)\bigr),
$$
which is an
\emph{isomorphism} when $\E$ is bounded-below, since for any
injective $\cO_\X$-module $\cI$, \smash{$\dirlm{}\!\<_i$} of the
flasque modules $\sHom(\cO_\X/\J^i,\cI\>)$ is $\Gamma_{\!\X}^{}$-acyclic.
Whenever $\gamma(\E)$ is an isomorphism, the induced homology maps are
isomorphisms
$$
\dirlm{i}\text{Ext}^n(\cO_\X/\J^i\<,\>\E)   \iso\Hp\X n(\E).
$$

\smallskip

If $\E\in\Dqc(\X)$, then $\R\iGp\X\E\in\Dqct(\X)$ (\Pref{Gamma'(qc)}).
For any map $g\colon\X\to\Y$ satisfying the hypotheses of
\Tref{T:qct-duality}, for
$\G\in\D(\Y)$, and with $R\!:=\Hr^0(\Y,\cO_\Y)$, there~are natural maps
\begin{equation}\label{map}
\begin{aligned}
\Hom_{\D(\X)}\<\<(\R\iGp\X\E\<,\gt\<\G)
 &\iso\< \Hom_{\D(\X)}\<\<(\R\iGp\X\E\<,\gt\R\iGp\Y\G)
  &&\qquad(\ref{C:identities}(\textup b)) \\
 &\iso\< \Hom_{\D(\Y)}\<\<(\R g_*\R\iGp\X\E\<,\R\iGp\Y\G)\\
 &\,\lra\mspace{1.2mu} \Hom_R(\Hp\X n\E, \Hp\Y n\G) &&
\end{aligned}
\end{equation}
where the last map arises via the functor $\Hr^n\R\Gamma(\Y,-)\
(n\in\mathbb Z)$.

In particular, if $g=\hat f$ in the completion situation of
\S\ref{completion}, and if $\E\set\kappa_\X^*\E_0\>$, 
$\G=\kappa_\Y^*\G_0\ (\E_0\in\Dqc(X),\;\G_0\in\Dqc(Y))$, then
preceding considerations show that this composed map operates via
Duality with Supports for $f$ (\Pref{P:dual/supports}),
i.e., it can be identified with the natural composition
\begin{align*}
\Hom_{\D(X)}\<\<(\R\iG Z\E_0, \R\iG Zf^\times\G_0)
&\underset{\ref{P:dual/supports}}{\iso}\<
\Hom_{\D(Y)}\<\<(\Rfs\>\R\iG Z\E_0, \R\iG W\G_0) \\
&\,\lra\mspace{1.2mu} \Hom_{\Hr^0(Y,\cO_Y)\<}
  (\Hr_Z^n\E_0, \Hr^n_W\G_0).
\end{align*}

\penalty -1500

\begin{sparag}
Next, let $R$ be a complete noetherian local ring topologized as usual by its
maximal ideal~$I$, let
$(S,J)$ be a noetherian adic ring, let $\varphi\colon (R,I)\to (S,J)$ be a continuous
homomorphism, and let 
$$
\Y\set\Spf(S)\xrightarrow{\,f\,}\Spf(R)=:\V
$$
be the corresponding formal-scheme map. As before, $g\colon\X\to\Y$ is
a map as in \Tref{T:qct-duality}, and we set $h\set fg$. Since the
underlying space of~$\V$ is a single point, at which the stalk
of~$\cO_\V$ is just~$R$, therefore the categories of
$\cO_{\<\V}$-modules and of $R$-modules are identical, and
accordingly, for any $\E\in \D(\X)$ we can identify $\R h_*\E$ with
$\R\Gamma(\X,\E)\in\D(R)$.

Let $K$ be an injective $R$-module, and $\cK$ the corresponding injective
$\cO_{\<\V}$-module.  There exist  integers $r$, $s$ such that 
$H^i(\ush f \cK)=0$ for all $i<-r$ (resp.~$H^i(\ush h \cK)=0$ for all $i<-s$)
(\Cref{C:f*gam-duality}). Set $\omega_\Y\set H^{-r}(\ush f \cK)$ 
(resp.~$\omega_\X\set H^{-s}(\ush h \cK)$). 
\end{sparag} 

\begin{sprop}\label{P:astrix10.2}
In the preceding situation\/  $\omega_\X$
represents---via \eqref{map}---the functor\/ 
$\Hom_S(\Hp\X {\<\raisebox{.1ex}{$\scriptstyle s$}}\E\<,\Hp\Y 
{\raisebox{.1ex}{$\scriptscriptstyle 0$}}(\ush f\cK))
$ of quasi-coherent
$\cO_\X$-modules~$\E\<$. If\/~$\omega_\Y$ is the only non-zero
homology  of\/~$\ush f\cK,$\ this functor is isomorphic to\/ 
$\Hom_S(\Hp\X {\<\raisebox{.1ex}{$\scriptstyle s$}}\E\<,\Hp\Y 
{\raisebox{.1ex}{$\scriptstyle r$}}\omega_\Y)$. \looseness=2
\end{sprop}

\emph{Proof.}
There are natural maps 
$$
\postdisplaypenalty10000
\Hp\Y {\raisebox{.1ex}{$\scriptstyle r$}}(\omega_\Y)=\Hp\Y {\raisebox{.1ex}{$\scriptscriptstyle 0$}}(\omega_\Y[r])\xrightarrow{\ h\ } \Hp\Y {\raisebox{.1ex}{$\scriptscriptstyle 0$}}(\ush
f\cK)\iso\Hom_{R,J}(S,K)
$$ 
where the last isomorphism results from \Pref{P:affine}(a), in
view of the identity \mbox{$\R\iGp\Y\ush f=\ft$} (\Cref{C:identities}(a))
and the natural isomorphisms
$$
\postdisplaypenalty 10000
\R\Gamma(\Y, \kappa_\Y^*\widetilde G)\iso
\R\Gamma(Y, \kappa_{\Y*}^\pd\kappa_\Y^*\widetilde G)
\underset{\ref{Gammas'+kappas}}\iso
\R\Gamma(Y, \widetilde G)\iso G
\qquad\bigl(G\in\D_{\<J}^+(S)\bigr),
$$
for $G\set\R\Homb_{R,J}(S,\R\Gamma(\V,\cK))$. (In fact 
$\R\Gamma(\Y, \kappa_\Y^*\widetilde G)\cong G$
for any \mbox{$G\in\D(S)$,} see
\Cref{(3.2.3)} and the beginning of
\S\ref{SS:Dvc-and-Dqc}.)
In case $\omega_\Y$ is the only non-vanishing homology of $\ush f\cK$, then $h$ 
is an isomorphism too.
 
The assertions follow from the (easily checked) commutativity, for any
quasi-coherent $\cO_\X$-module~$\E\<$,  of the diagram
$$
\minCDarrowwidth=22pt
\mkern110mu
\begin{CD}
\hbox to0pt{\hss $\Hom_\X(\E\<,\omega_\X)=\!\!\<=\:$}
  \Hom_{\D(\X)}\<\<(\E[s],\>\ush g\<\ush f\cK)
    @>\<\!\!\!\textup{~\ref{C:identities}(a)}\mkern.5mu >>
      \Hom_{\D(\X)}\<\<(\R\iGp\X\E[s],\>{g_{\mathrm t}^\times}\!\ush f\cK)\\
  @V\simeq VV @VV\eqref{map}V \\
 \Hom_{\D(\V)}\<\<(\R h_*\R\iGp\X\E[s],\cK)
     @.\Hom_S\bigl(\Hp\X {\raisebox{.1ex}
        {$\scriptscriptstyle 0$}}(\E[s]),\Hp\Y {\raisebox{.1ex}{$\scriptscriptstyle
           0$}}(\ush f\cK)\bigr)\\
@| @VV\simeq V \\
\Hom_R(\Hp\X {\<\raisebox{.1ex}{$\scriptstyle s$}}\E\<,K) 
  @>\;\ \Iso\ \;>> 
\Hom_S\bigl(\Hp\X {\<\raisebox{.1ex}{$\scriptstyle s$}}\E\<,\Hom_{R,J}(S,K)\bigr)
\end{CD}
$$
\end{sparag}

\begin{sparag}
Now let us fit \cite[pp.~87--88, Theorem (10.2)]{Asterisque} into the preceding
setup. 

The cited Theorem has both local and global components. The first
deals with maps $\varphi\colon R\to S$ of local domains essentially of
finite type over a perfect field~$k$, with residue fields finite
over~$k$. To each such ring~$T$ one associates the canonical module
$\omega_T$ of ``regular" $k$-differentials of degree~$\dim T$.  Under
mild restrictions on~$\varphi$, the assertion is that the functor
$$
\postdisplaypenalty 10000
\Hom_R\bigl(\textup H_{m\<\<_{\hat S}}^{\dim\mkern-1.5mu S}G,\>
 \>\textup H_{m\<_{\lower.3ex\hbox{$\scriptscriptstyle\<R$}}}^
   {\dim\mkern-1.5mu R}\> \omega_{\<R}\bigr)
\qquad(m\set\text{maximal ideal})
$$
of $\hat S$-modules~$G$ is represented by the completion $\widehat{\omega_S}$
together with a canonical~map, the {\it relative residue}
$$
\rho_\varphi\colon \textup H_{m\<\<_{\hat S}}^{\dim\mkern-1.5mu
S}\widehat{\omega_S\>} =\textup H_{m\<_{S}}^{\dim\mkern-1.5mu S\>}\omega_S\to
 \textup H_{m\<_{R}}^{\dim\mkern-1.5mu R\>}\>\omega_{\<R}.
$$
This may be viewed as a consequence of {\it concrete\/} local duality over~$k$
(\S\ref{HuK}).

The global aspect concerns a proper map of irreducible $k$-varieties 
$g\colon V\to W$ of respective dimensions $s$ and $r$ 
with all fibers over codimension~1 points of~$W$ having
dimension $s-r$, a closed point
$w\in W\<$, the fiber $E\set g^{-1}(w)$, and the completion $\widehat
V\set V_{/E}$.  The assertion is that the functor
$$
\postdisplaypenalty 10000
\Hom_R\bigl(\textup H_{\widehat V}^{\prime s}\>\G,\mkern1.5mu
\textup H_{m\<\<_R}^{r}\omega_{\<R}\bigr)\qquad
 (R\set \cO_{W\mkern-1.5mu,\>w})
$$
of coherent $\cO_{\widehat V}$-modules~$\G$ is represented by the
completion~$\widehat{\omega_V}$ along~$E$ of the canonical
sheaf $\omega_V$ of regular differentials, together with a canonical map
$$
\theta\colon
\textup H_{\widehat V}^{\prime s}\>\widehat{\omega_V}=
\textup H_{\<E}^{s}\>\omega_V\to
\textup H_{m_{\<\<R}}^{r}\omega_{\<R}\>.
$$

Moreover, the local and global
representations are {\it compatible\/} in the sense that if
$v\in E$ is any closed point and $\varphi_v\colon R\to S\set \cO_{V\!,v}$ is
the canonical map, then the residue $\rho_v\set\rho_{\varphi_v}$
factors as the natural map
$\textup H_{m_{\<S}}^s\omega_S\to 
 \textup H_{\<E}^{s}\>\omega_V$
followed by~$\theta$. This compatibility determines $\theta$ uniquely
if the $\rho_v\ (v\in E)$ are given \cite[p.\,95, (10.6)]{Asterisque};
and of course conversely.
 
Basically, all this---\emph{without
 the explicit description of the $\omega\<$'s and the maps\/~$\rho_v$ via
differentials and residues}---is contained in \Pref{P:astrix10.2},
as follows.

In the completion situation of \S\ref{completion}, take $X$ and $Y$ to be
finite-type separated schemes over an artinian local ring~$R$, 
of respective pure dimensions $s$ and $r$, let
$W=\{w\}$ with $w$ a closed point of~$Y\<$, write $g$ in place of~$f$,
and assume that 
$Z\subset g^{-1}W$ is proper over~$R$ (which is so, e.g., if $g$ is proper and $Z$ is
closed).   Let~$K$ be an injective hull of the residue field of~$R$, and let $\cK$ be
the corresponding injective sheaf on $\Spec(R)=\Spf(R)$.  With
$f\colon Y\to
\Spec(R)$ the canonical map, and $h=fg$,  define the \emph{dualizing sheaves} 
$$
\omega_X\!:=H^{-s}h^!\cK, \qquad
\omega_Y\!:=H^{-r}f^!\cK,
$$
where $h^!$ is the Grothendieck duality functor (compatible with open
immersions, and equal to $h^\times$ when $h$ is proper), and similarly
for~$f^!\<$. It is well-known (for example via a local factorization
of $h$ as $\text{smooth}\smcirc\text{finite}$) that $h^!\cK$ has
coherent homology, vanishing in all degrees $<-s\>$; and similarly
$f^!\cK$ has coherent homology, vanishing in all degrees $<-r$.

Let 
$$
\hat f\colon\Y\set\Spf(\widehat{\cO_{W\<,\>w}})\to\Spf(R)=:\V
$$
be the completion of $f$. We may assume, after compactifying $f$ and
$g$---which\vadjust{\kern.7pt} does not affect $\hat f$ or~$\hat g$
(see~\cite{Lu}), that $f$ and $g$ are proper maps.\vadjust{\kern.4pt}
Then \Cref{C:completion-proper} shows that $\ush {\hat
h}\cK=\kappa_\X^*h^!\cK$, and so $\kappa_\X$ being flat, we see that
\begin{equation}\label{omega}
\kappa_\X^*\>\omega_X=\omega_\X
\end{equation}
where $\omega_\X$ is as in \Pref{P:astrix10.2}; and similarly
$\kappa_\Y^*\omega_Y=\omega_\Y\>$.

Once again, some form of the theory of the Fundamental Class will enable us to
represent $\omega_X$ by means of regular differential forms; and then both the
local and global components of the cited Theorem~(10.2) become special cases of
\Pref{P:astrix10.2} (modulo some technicalities 
\cite[p.\,89, Lemma (10.3)]{Asterisque} which allow a weakening of the condition
that $\omega_\Y$ be the only non-vanishing homology  of~$\ush{\hat f}\cK$). 

As for the local-global compatibility,  consider quite generally a pair of maps
$$
\X_1\xrightarrow{q\,}\X\xrightarrow{p\,}\Y
$$
of noetherian formal schemes.  In the above situation, for instance, we could
take~$p$ to be $\hat g$, $\X_1$ to be the completion of~$X$ at a closed point 
$v\in Z$, and $q$ to be the natural map. 
\Tref{Th2} gives us the adjunction\vadjust{\kern1pt}
$$
\lower.3ex\hbox{$\Dqct(\X)$}\ \,
\vbox to0pt{
\vss
\hbox{$\xrightarrow{\<\R p_*\>}$}
\vspace{-7pt}
\hbox{$\xleftarrow[\smash{\;{p_{\mathrm t}^{\<\times}}\;}]{}$}
\vss}
\ \,\lower.3ex\hbox{$\Dqct(\Y)$}.
$$
 
\medskip\vspace{1pt}

\noindent The natural isomorphism $\R(pq)_*\iso\R p_*\R q_*$ gives rise then to
an adjoint isomorphism 
$q_{\mathrm t}^{\<\times}p_{\mathrm t}^{\<\times}\iso
(pq)_{\mathrm t}^{\<\times}$; and
for $\E\in\Dqct(\Y)$ the natural map $\R(pq)_*(pq)_{\mathrm t}^{\<\times}\E\to\E$
factors as
$$
\R(pq)_*(pq)_{\mathrm t}^{\<\times}\E\iso
\R p_*\R q_*q_{\mathrm t}^{\<\times}p_{\mathrm t}^{\<\times}\E\to
\R p_*p_{\mathrm t}^{\<\times}\E\to\E.
$$
This factorization contains the compatibility between
the above maps $\theta$ and $\rho_v\>$, as one sees by interpreting them as
homological derivatives of maps of the type 
$\R p_*p_{\mathrm t}^{\<\times}\E\to\E$ (with $\E\set\R\iGp\Y\ush{\hat f}\cK)$.
Details are left to the reader.
\end{sparag}

\begin{srem}\label{R:d-vein}
In the preceding situation, suppose further that $Y=\Spec(R)$ (with $R$~artinian)
and $f=\text{identity}$, so that $h=g\colon X\to Y$ is a finite-type separated map,
$X$ being of pure dimension $s$, and $\kappa_\X\colon\X\to X$ is the completion
of~$X$ along a closed subset~$Z$ proper over~$Y\<$.   Again, $K$ is an injective
$R$-module,
$\cK$ is the corresponding $\cO_Y$-module, and $\omega_X\set H^{-s}g^!\cK$ is
a ``dualizing sheaf\kern1.5pt" on~$X$. Now~\Pref{P:astrix10.2} is just
the instance
$i=s$ of the canonical isomorphisms, for~$\E\in\Dqc(\X),\;i\in\mathbb Z$ (and with
$\Hp\X\bullet\set \Hr{}^\bullet\R\Gamma(\X,\R\iGp\X)$, see
\S\ref{consequences}, and $\hat g\set g\smcirc\kappa_\X$):
$$
\!\!\Hom_{\D(\X)}\<(\E[i],\ush {\hat g}\cK)
\underset{\textup{Thm.\,\ref{Th2}}}{\iso}
 \Hom_{\D(\Y)}\<(\R\hat g_{\<*}\R\iGp\X\E[i],\cK)\<\iso\<
   \Hom_R(\Hp\X i\>\E\<,K)
     =: \!(\Hp\X i\>\E)\mspace{.5mu}\check{}\>\>.
$$
If $X$ is Cohen-Macaulay then all the homology of $g^!\cK$ other than $\omega_X$
vanishes, so all the homology of $\ush{\hat g}\cK\cong \kappa_\X^*g^!\cK$
other than $\omega_\X=\kappa_\X^*\omega_X$ vanishes (see \eqref{omega}),
and the preceding composed isomorphism becomes
$$
\text{Ext}_\X^{s-i}(\E\<,\omega_\X)\iso (\Hp\X i\>\E)\mspace{.5mu}\check{}\>\>.
$$
 In particular, when $Z=X$ (so that $\X=X$) this is the usual duality isomorphism
$$
\text{Ext}_X^{s-i}(\E\<,\omega_X)\iso \Hr^i(X,\E)\mspace{.5mu}\check{}\>\>.
$$

If $X$ is Gorenstein and $\F$ is  a locally free $\cO_\X$-module
of finite rank, then $\omega_X$~is invertible; and taking
$\E\set\sHom_\X(\F,\omega_\X)=\check \F\otimes\omega_\X$ we get the
isomorphism
$$
\Hr^{s-i}(\X, \F\>) \iso  \bigl(\Hp\X i(\check\F\otimes\omega_\X)\bigr)
\mspace{.8mu}\check{}\>\>,
$$
which generalizes the Formal Duality theorem\index{Duality!Formal|(} 
\cite[p.\,48, Proposition~(5.2)]{De-Rham-cohomology}.
\end{srem}

\end{parag}\index{Residue theorems|)}

\smallskip
\penalty -1200

\begin{parag}
\label{bf (d)}

Both \cite[p.\,48; Proposition~(5.2)]{De-Rham-cohomology} (Formal Duality) 
and the Theorem in \cite[p.\,188]{Desingularization}  
(Local-Global Duality)\index{Duality!Local-Global} are
contained in \Pref{(2.8)}, see \cite[\S5.3]{AJL}.

 Let $R$ be a noetherian ring, discretely topologized, and set 
$$
Y\set \Spec(R)=\Spf (R)=:\<\Y.
$$
Let $g\colon X\to Y$ be a finite-type separated map, let $Z\subset X$ be
\emph{proper} over~$Y\<$, let $\kappa\colon\X=X_{/Z}\to X$ be the
completion of~$X$ along~$Z$, and set $\hat g\set g\smcirc\kappa\colon\X\to \Y$.

 Assume that $R$ has a \emph{residual complex} $\cR$ \cite[p.\,304]{H1}.
Then the corresponding quasi-coherent $\cO_Y$-complex
 \smash{${\cR}_Y\set  {\widetilde {\cR}}$} is a  \emph{dualizing
complex,} and ${\cR}_X\set g^!{\cR}_Y$ is a dualizing complex on~$X$
 \cite[p.~396, Corollary~3]{f!}.  For any $\F\in\Dc(X)$ set
$$
\F^{\>\prime}\set\R\sHomb_ X(\F,{\cR}_X)\in\Dc(X),
$$
so that  $\F\cong\F^{\>\prime}{}'=\R\sHomb_ X(\F^{\>\prime}\<, {\cR}_X)$. 

\begin{sprop}\label{(2.8)}\index{Duality!Formal|)} 
In the preceding situation, with\/  
$\Gamma_{\<\! Z}(-)\set\Gamma(X,\iG Z(-))$ 
there is a functorial isomorphism
$$
\R\Gamma(\X,\kappa^*\<\F\>)\cong
 \R\Hom_R^\bullet(\R\Gamma_{\<\! Z}\>\F^{\>\prime}\< ,\>{\cR})
\qquad \bigl(\F\in\Dc(X)\bigr).
$$
\end{sprop}

\begin{proof} Replacing $g$ by a compactification (\cite{Lu}) doesn't
affect~$\X$ or $\R\Gamma_{\<\!Z}$, so  assume that $g$  is proper. Then
\Cref{C:completion-proper} gives an isomorphism 
\mbox{$\kappa^*\cR_X\cong\ush{\hat g}\cR_Y$.} 
Now just compose the chain of functorial isomorphisms
\begin{align*}
\R\Gamma(\X, \kappa^*\<\F \>)
&\cong 
 \R\Gamma\bigl(\X,\kappa^*\R\sHomb_ X(\F^{\>\prime}\<, {\cR}_X)\bigr) &
&\textup{(see above)}\\
&\cong  
\R\Gamma\bigl(\X, \R\sHomb_{ \X}(\kappa^*\<\F^{\>\prime}\<, 
  \>\kappa^*\cR_X)\bigr)& 
&\textup{(\Lref{L:kappa*Ext})}\\
&\cong \R\Hom_\X^\bullet(\kappa^*\<\F^{\>\prime}\<, \>\ush{\hat g}{\cR}_Y))&
&\textup{(see above)}\\
&\cong \R\Hom_\Y^\bullet(\R\hat g_*\R\iGp  \X\kappa^*\<\F^{\>\prime}\<, {\cR}_Y) &
&\textup{(\Tref{Th2})}\\
&\cong \R\Hom_Y^\bullet(\R g_*\R\iG {Z_{\mathstrut}}\F^{\>\prime}\<, \cR_Y)&
&\textup{(\Pref{Gammas'+kappas})}\\ 
&\cong
 \R\Hom_Y^\bullet(\smash{\widetilde{\R\Gamma_{\<\! Z}\F}{}^{\>\prime}\<},
   \>\cR_{Y_{\mathstrut}})&
&\textup{\cite[footnote, \S5.3]{AJL}}\\
&\cong  \R\Hom_R^\bullet(\R\Gamma_{\<\! Z}\F^{\>\prime}\< ,{\cR})&
&\textup{\cite[p.\,9, (0.4.4)]{AJL}}.
\end{align*}
\vskip-3.7ex
\end{proof}

\begin{slem}\label{L:kappa*Ext}
Let\/ $X$ be a locally noetherian scheme, and let\/ $\kappa\colon\X\to X$ be its
completion along some closed subset\/~$Z$. Then for\/~$\G\in\Dqc(X)$ of finite
injective dimension and  for\/~$\F\in\Dc(X),$ the
natural map is an isomorphism
$$
\kappa^*\R\sHomb_ X(\F,\G)\iso
\R\sHomb_{\X}(\kappa^*\<\F,\kappa^*\G).
$$
\end{slem}

\begin{proof}

By \cite[p.\,134, Proposition~7.20]{H1} we may assume that
$\G$ is a bounded complex of quasi-coherent injective $\cO_{\<\<X}$-modules,
vanishing, say, in all degrees $>n$.

When $\F$ is bounded-above the (well-known) assertion is 
proved by localizing  to the affine case and applying \cite[p.\,68,
Proposition~7.1]{H1} to reduce to the trivial case 
$\F=\cO_{\<\<X}^m\ (0<m\in\mathbb Z)$. To do the same for unbounded~$\F$ we must
first  show, for fixed~$\G$, that the contravariant
functor~$\R\sHomb_{\X}(\kappa^*\<\F\<,\kappa^*\G)$ is bounded-above.

In fact we will show that if $H^i\F=0$ for all $i<i_0$ then
for all $j>n-i_0\>$, 
$$
H^j\R\sHomb_{\X}(\kappa^*\<\F\<,\kappa^*\G)=
H^j\<\kappa_*\R\sHomb_{\X}(\kappa^*\<\F\<,\kappa^*\G)=
H^j\R\sHomb_{X}(\F\<,\kappa_*\kappa^*\G)=0.
$$ 
The homology in question is the sheaf associated to the presheaf which assigns
$$
\Hom_{\D(U)}\<\bigl(\F|_{\lower.3ex\hbox{$\scriptstyle U$}}[-j], 
 (\kappa_*\kappa^*\G)|_{\lower.3ex\hbox{$\scriptstyle U$}}\bigr)=
\Hom_{\D(U)}\<\bigl(\F|_{\lower.3ex\hbox{$\scriptstyle U$}}[-j], 
 \R Q_{\lower.3ex\hbox{$\scriptstyle\<\<U$}}
  (\kappa_*\kappa^*\G)|_{\lower.3ex\hbox{$\scriptstyle U$}}\bigr)
$$ 
to each affine open subset
$U=\Spec(A)$ in~$X\<$. (Here we abuse notation by omitting~$\bj_{\<\<U}^{}$ 
in front of $\R Q_{\<\<U}^{}$, see
beginning of~\S\ref{SS:Dvc-and-Dqc}). 

\penalty-1000
Let $\U\set\kappa^{-1}U$,
and
$\hat A\set\Gamma(\U,\cO_\X)$, so that $\kappa|_\U$ factors naturally as
$$
\U=\Spf(\hat A)\xrightarrow{\kappa_1\,}U_1\set
\Spec(\hat A)\xrightarrow{k\,}\Spec(A)=U.
$$
The functors $\R Q_{\lower.3ex\hbox{$\scriptstyle\<\< U$}}k_*$ 
and\vadjust{\kern.8pt} 
$k_*\R Q_{ U_1}^{}$,  both
right-adjoint to the natural composition 
\smash{$\Dqc(U)\xrightarrow{\vbox
to0pt{\vskip-.8ex\hbox{$\scriptstyle k^*$}\vss}}
 \Dqc(U_1)\hookrightarrow\D(U_1)$,}  are isomorphic; so
there are natural isomorphisms
$$
 \R Q_{\lower.3ex\hbox{$\scriptstyle\<\<U$}}
  (\kappa_*\kappa^*\G)|_{\lower.3ex\hbox{$\scriptstyle U$}}
=
\R Q_{\lower.3ex\hbox{$\scriptstyle\<\<U$}}  k_{*}\kappa_{1*}^{}\kappa_1^*k^*
   (\G|_{\lower.3ex\hbox{$\scriptstyle U$}})
\<\iso\<
  k_{*}\R Q_{U_1}^{}
   \kappa_{1*}^{}\kappa_1^*k^*(\G|_{\lower.3ex\hbox{$\scriptstyle U$}})
\<\smash{\underset{\ref{c-erator}}{\iso}}\< 
k_{*}k^*(\G|_{\lower.3ex\hbox{$\scriptstyle U$}})
$$
and the presheaf becomes
$
U\mapsto \Hom_{\D(U)}\<\bigl(\F|_{\lower.3ex\hbox{$\scriptstyle U$}}[-j],  k_{*}k^*(\G|_{\lower.3ex\hbox{$\scriptstyle U$}})\bigr).
$

The equivalence of categories $\Dqc(U)\cong \D(\Aqc(U))=\D(A)$ indicated at the
beginning of~\S\ref{SS:Dvc-and-Dqc} yields an isomorphism
$$
\Hom_{\D(U)}\<\bigl(\F|_{\lower.3ex\hbox{$\scriptstyle U$}}[-j],  k_{*}k^*(\G|_{\lower.3ex\hbox{$\scriptstyle U$}})\bigr)\iso
\Hom_{\D(A)}\<\bigl(F[-j], G\otimes_A\hat A\bigr)
$$
where $F$ is a complex of $A$-modules\vadjust{\penalty-1000} with $\Hr^iF =0$
for
$i<i_0\>$, and  both $G$ and $G\otimes_A\hat A$ are complexes of injective
$A$-modules vanishing in all degrees $>n$\vspace{-1pt} (the latter since $\hat A$
is $A$-flat). Hence the presheaf vanishes, and the conclusion follows.
\end{proof}
\end{parag}
\penalty-1500

\begin{parag}\label{bf (e)} (Dualizing complexes.)\index{dualizing complexes|(}
Let $\X$ be  a noetherian formal scheme, and write $\D$ for $\D(\X)$, etc. The
derived functor\index{ $\iG{\raise.3ex\hbox{$\scriptscriptstyle{\ldots}$}}$
(torsion functor)!a@$\BG\set\R\iGp\X$ (cohomology colocalization)}
$\BG\set\R\iGp\X\colon\D\to\D$ (see \Sref{Gamma'1}) has a right adjoint%
\index{ $\mathbf {La}$@$\BL$ (homology localization)}
$\BL=\BL_\X\set\R\sHomb(\R\iGp\X\cO_\X^{}, -)$. This adjunction is given by
\eqref{adj}, a natural isomorphism of which we'll need the sheafified form, 
proved similarly:
\begin{equation}\label{adj0}
\R\sHomb(\M,\BL\cR)
\cong
\R\sHomb(\BG\<\M,\cR).
\end{equation}
There are natural maps
$\BG\to\mathbf1\to\BL$ inducing isomorphisms 
$\BL\BG\iso\BL\iso\BL\BL$,  $\BG\BG\iso\BG\iso\BG\BL$
(\Rref{R:Gamma-Lambda}\,(1)). \Pref{formal-GM}, a form of Greenlees-May
duality, shows that $\BL(\Dc)\subset\Dc$. (Recall that the objects of the
$\Delta$-subcategory $\Dc\subset\D$ are the complexes whose homology sheaves
are all coherent.)

Let $\Dc^*$ be the essential image of $\BG|_{\Dc}$, i.e., the full subcategory
of~$\D$ such that $\E\in\Dc^*\Leftrightarrow\E\cong\BG\F$ with
$\F\in\Dc$. \Pref{Gamma'(qc)} shows that $\Dc^*\subset\Dqct$. It follows from
the preceding paragraph that
\begin{alignat*}{2}
\E\in\Dc^*&\iff \BG\E\iso\E&&\text{ and }\,\BL\>\E\in\Dc\>,\\
\F\in\Dc&\iff \F\iso\BL\F\,&&\text{ and }\,\BG\F\in\Dc^*.
\end{alignat*}
(In particular, $\Dc^*$ is a $\Delta$-subcategory of~$\D$.) Moreover $\BG$ and
$\BL$ are quasi-inverse equivalences between the categories $\Dc$ and $\Dc^*$.
\end{parag}

\begin{sdef}\label{D:dualizing}
A complex $\cR$ is a
\emph{c-dualizing complex on} $\X\>$ if 
\begin{enumerate}
\item[(i)] $\cR\in\Dc^+(\X).$
\item[(ii)] The natural map is an isomorphism
$\,\cO_\X\iso\R\sHomb(\cR,\cR)$.
\item[(iii)] There is an integer $b$ such that for every coherent torsion sheaf~$\M$
and for every~$i>b$, $\>\Ext^i(\M,\cR)\set
H^i\>\R\sHomb(\M,\cR)=0$.
\end{enumerate}

A complex  $\cR$ is a
\emph{t-dualizing complex on} $\X\>$ if 
\begin{enumerate}
\item[(i)] $\cR\in\D_{\mathrm t}^+(\X).$
\item[(ii)] The natural map is an isomorphism
$
\,\cO_\X\iso\R\sHomb(\cR,\cR).
$
\item[(iii)] There is an integer $b$ such that for every coherent torsion sheaf~$\M$
and for  every $i>b$, $\>\Ext^i(\M,\cR)\set
H^i\>\R\sHomb(\M,\cR)=0$.
\item[(iv)] For some ideal of definition~$\J$ of~$\X$, 
$\R\sHomb(\cO_\X/\J, \cR)\in\Dc(\X).$
\item[{}](Equivalently---by simple arguments---$\R\sHomb(\M, \cR)\in\Dc(\X)$
for every coherent torsion sheaf $\M$.)
\end{enumerate}
\end{sdef}

\pagebreak[3]
\emph{Remarks.} (1) On an ordinary scheme,  (iii) signifies \emph{finite
injective dimension}  
\cite[p.\,83, Definition, and p.\,134, (iii)${}_{\textup c}$]{H1}, so both
c-dualizing and t-dualizing mean the same as what is called ``dualizing" in
\cite[p.\,258, Definition]{H1}. (For the extension to arbitrary noetherian formal
schemes, see (4) below.)

(2) By (i) and (iv), \Pref{Gamma'(qc)}(a), and \Cref{qct=plump}, any t-dualizing
complex $\cR$ is in $\Dqct^+(\X)$; and then (iii) implies that 
\emph{$\cR$~is isomorphic in $\D$ to a bounded complex of
$\Aqct$-injectives.}

To see this,  begin by imitating the proof of~\cite[p.\,80,
(iii)$\Rightarrow$(i)]{H1}, using \cite[Theorem~4.8]{Ye} and \Lref{L:Hom=RHom}
below, to reduce to showing that \emph{if\/ $\N\in\Aqct(\X)$ is such that\/
$\Ext^1(\M,\N\>)=0\>$ for every coherent torsion sheaf\/ $\M$ then\/ $\N$ is\/
$\Aqct$-injective.}

For the last assertion, suppose first that $\X$ is affine.
\Lref{Gamma'+qc} implies that  \mbox{$\sHom(\M,\N\>)\in\Avc(\X)$;} and
then $\textup{Ext}^1(\M,\N\>)=0$, by the natural exact sequence
$$
0\underset{\textup{(3.1.8)}}=\textup H^1\bigl(\X, \sHom(\M,\N\>)\bigr)
\to\textup{Ext}^1(\M,\N\>)\to \textup H^0\bigl(\X, \Ext^1(\M,\N\>)\bigr).
$$
 Since  coherent torsion sheaves
generate  $\Aqct(\X)$ (\Cref{qct=plump}, \Lref{Gamma'+qc}), a
standard argument using Zorn's Lemma shows that $\N$~is indeed
$\Aqct$-injective.

\penalty-1500

In the general case, let $\U\subset\X$ be any affine open subset. For any
coherent torsion $\cO_\U$-module~$\M_0$,  \Pref{f-*(qct)} and
\Lref{Gamma'+qc} imply there is a coherent torsion $\cO_\X$-module~$\M$
restricting on~$\U$ to~$\M_0$, whence $\>\Ext^1_\U(\M_0, \N\>|_{\U})=0$. By the
affine case, then,
$\N\>|_{\U}$ is $\Aqct(\U)$-injective, hence $\At(\U)$-injective 
\cite[Proposition~4.2]{Ye}. Finally, as in \cite[p.\,131, Lemma
7.16]{H1}, using
\cite[Lemma 4.1]{Ye},%
\footnote{where one may assume that $X$ and $\X$ have the same 
underlying space}
one concludes that $\N$ is  $\At(\X)$-injective, hence $\Aqct(\X)$-injective.

\smallskip
 (3) With (2) in mind, one finds that what is called here ``t-dualizing
complex" is what  Yekutieli\index{Yekutieli, Amnon} calls in \cite[\S5]{Ye}
``dualizing complex."

\smallskip
(4) \emph{A c-dualizing complex~$\cR$ has finite injective 
dimension}: there is an integer~$n_0$ such that for any $i>n_0$ and
any $\cO_\X$-module~$\E$, $\Hom_{\D}(\E\<, \cR[i])=0$. To see this, note 
first that
$$
\Hom_{\D}(\E\<,\cR[i])\cong \Hom_{\D}(\E\<,\BL\BG\cR[i])\cong
\Hom_{\D}(\BG\E\<,\BG\cR[i]).
$$
\Lref{L:interchange}(b) below and  (2) above show that $\BG\cR$ is isomorphic to
a bounded complex of $\Aqct$-injectives. The  complex
$\BG\E$---obtained by applying the functor $\iGp\X$ to an injective resolution
of~$\E$---consists of torsion $\cO_\X$-modules, and so as in \cite[Corollary
4.3]{Ye} (see also the proof of \Lref{L:Hom=RHom} below, with \Pref{iso-qct} in
place of \Pref{(3.2.1)}), the natural map 
$$
\textup H^i\bigl(\Homb(\BG\E\<,\BG\cR)\bigr)\to
\textup H^i\bigl(\R\Homb(\BG\E\<,\BG\cR)\bigr) =\Hom_{\D}(\BG\E\<, \BG\cR[i])
$$
is an \emph{isomorphism.} Since $\BG\E$ vanishes in degrees $<0$, the asserted
result holds for any  $n_0$  such that
$H^i(\BG\cR)=0$ for $i>n_0\>$.

\smallskip
(5) For a complex $\cR\in\Dc^+\cap\mspace{1.5mu}\Dc^-\<$, 
conditions (ii) and (iii) in
\Dref{D:dualizing} hold iff they hold stalkwise for $x\in \X$, with an
integer $b$ \emph{independent of\/~$x$.} (The idea is that such an $\cR$ is locally
resolvable by a bounded-above complex~$\F$ of finite-rank locally free
$\cO_\X$-modules, as is $\M$ in~(iii), and 
$\sHomb(\F\<,\cR)\cong\R\sHomb(\F\<,\cR)$\dots.) Proceeding as in the
proofs of~\cite{H1}, Proposition~8.2, p.\,288, and Corollary~7.2, p.\,283, one
concludes that
$\cR$ is c-dualizing iff $\>\X$ has finite Krull dimension and $\cR_x$~is a
dualizing complex for the category of
$\cO_{\X\<,\>x}$-modules for every $x\in\X$. (It is enough that the latter hold for
all \emph{closed} points~$x\in\X$.)

\begin{exams}\label{regular}
(1) If $\,\cR$ is c-(or t-)dualizing then so is $\cR\otimes\cL[n]$ for any invertible
$\cO_\X$-module and $n\in\mathbb Z$. The converse also holds, see
\Pref{P:uniqueness}.

\smallskip
(2) (Cf.~\cite[Example 5.12]{Ye}.) If $X$ is an ordinary scheme and
$\kappa\colon\X\to X$ is its completion along some closed subscheme~$Z$, then
for any dualizing $\cO_{\<\<X}$-complex~$\cR$,
$\kappa^*\cR$ is c-dualizing on ~$\X$, and 
$\BG\kappa^*\cR\cong\kappa^*\R\iG Z\cR$ (see \Pref{Gammas'+kappas}(c)) is
a t-dualizing complex lying in $\Dc^*(\X)$.

\emph{Proof.}
For $\kappa^*\cR$, conditions (i) and (ii) in the definition of c-dualizing
follow easily from the same for $\cR$ (because of
\Lref{L:kappa*Ext}). So does~(iii), after we reduce to the case $X$ affine, 
where \Pref{(3.2.1)} allows us to write $\M=\kappa^*\M_0$ with $\M_0\in\A(X)$.
(Recall from  remark (1) above that $\cR$ has finite injective dimension.)
The last assertion is given by \Lref{L:interchange}(b).

\smallskip
(3) If $\X=\Spf(A)$ where $A$ is a complete local
noetherian ring topologized by its maximal ideal~$m$---so that $\A(\X)$ is just
the category of $A$-modules---then a \mbox{c-dualizing} $\cO_\X$-complex is an
$A$-dualizing complex in the usual sense; and by~(2) (via~\cite[p.\,276, 6.1]{H1}),
or directly from \Dref{D:dualizing},
the injective hull of $A/m$ is a t-dualizing complex lying in $\Dc^*(\X)$.

\smallskip
(4) It is clear from \Dref{D:dualizing} and remark~(4)
above that $\cO_\X$ is c-dualizing iff
$\cO_\X$ has finite injective dimension over itself. 
By remark (5), $\cO_\X$ is
c-dualizing iff $\X$ is finite dimensional and
$\cO_{\X\<,\>x}$ is \emph{Gorenstein} for all $x\in\X$ \cite[p.\,295,
Definition]{H1}. 

\smallskip
(5) For instance, if the finite-dimensional  noetherian formal scheme $\>\Y$ is
\emph{regular} (i.e., the local rings
$\cO_{\Y\<,\>y}\ (y\in\Y)$ are all  regular), and 
$\cI$ is a coherent \mbox{$\cO_\Y$-ideal}, defining a closed formal
subscheme~$i\colon\X\hookrightarrow\Y$ \cite[p.\,441,(10.14.2)]{GD}, then by~
remark~(3), $\R\sHomb(i_*\cO_\X^{}, \>\cO_\Y^{})$ is  c-dualizing 
on~$\X$. So \Lref{L:interchange} gives that\looseness=-1 
$$
\R\sHomb(\cO_\Y^{}/\cI, \>\R\iGp\Y\cO_\Y^{})
\underset{\textup{\ref{R:Dtilde}(4)}}\cong 
\R\iGp\X\>\R\sHomb(i_*\cO_\X^{}, \>\cO_\Y^{})\in\Dc^*(\X)
$$
is t-dualizing on~$\X$. 
(This is also shown in \cite[Proposition 5.11,  Theorem 5.14]{Ye}.)
\end{exams}
\pagebreak[3]

\begin{slem}\label{L:interchange}
\textup{(a)} If\/ $\cR\in\Dc^*$ is t-dualizing then\/ $\BL\cR$ is c-dualizing.
\vspace{1pt}

\textup{(b)} If\/ $\cR$ is c-dualizing then\/ 
$\BG\cR$ is  t-dualizing, and lies in\/~$\Dc^*$.
\end{slem}

\begin{proof} (a) If $\;\cR\in\Dc^*$ then of course $\BL\cR\in\Dc\>$. Also,
$\BL(\D^+)\subset\D^+$ because $\R\iGp\X\cO_\X^{}$ is given
locally by a finite complex
$\cK_\infty^\bullet\>$, see proof of \Pref{Gamma'(qc)}(a). 

For condition (ii), note that if $\cR\in\Dqct^+(\X)$ then $\BG\cR\cong\cR$
(\Pref{Gamma'(qc)}), then use
the natural isomorphisms (see \eqref{adj0}:
$$
\R\sHomb(\BL\cR\>,\BL\cR)
\cong
\R\sHomb(\BG\BL\cR\>,\cR)\cong
\R\sHomb(\BG\cR\>,\cR)\cong
\R\sHomb(\cR\>,\cR). \\
$$

For (iii) note that $\BG\M\cong\M$ (\Pref{Gamma'(qc)}), then use \eqref{adj0}.

\smallskip
(b) \Pref{Gamma'(qc)} makes clear that if $\cR\in\Dc^+$
then $\BG\cR\in\Dqct^+\cap\Dc^*$. 

For (ii) use the
 isomorphisms (the second holding because $\cR\in\Dc$):
$$
\R\sHomb(\BG\cR,\BG\cR)
\underset{\textup{\ref{C:Hom-Rgamma}}}\cong
\R\sHomb(\BG\cR,\cR)
\underset{\textup{\ref{formal-GM}}}\cong
\R\sHomb(\cR,\cR).
$$

For (iii)  use the  isomorphism
$\R\sHomb(\M,\BG\cR)
\underset{\textup{\ref{C:Hom-Rgamma}}}\cong
\R\sHomb(\M,\cR).$\vspace{2pt} For (iv), note that when
$\M=\cO_\X/\J\,$ ($\J$ any ideal of definition)  this isomorphism gives
$$
\R\sHomb(\cO_\X/\J,\BG\cR)\cong\R\sHomb(\cO_\X/\J,\cR)
\underset{\textup{\ref{P:Rhom}}}\in\Dc\>,
$$ 
\vskip-4.3ex
\end{proof}
\smallskip
\pagebreak[3]
The essential \emph{uniqueness} of t-(resp.~c-)dualizing complexes is expressed
by:
\begin{sprop}\label{P:uniqueness}
\textup{(a) (Yekutieli)}\index{Yekutieli, Amnon} If\/ $\cR$ is t-dualizing then 
a complex\/ $\cR'$ is t-dualizing iff there is an invertible
sheaf\/~$\cL$ and an integer\/~$n$ such that\/
$\cR'\cong\cR\otimes\cL[n]$.

\smallskip
\textup{(b)} If\/ $\cR$ is c-dualizing then a complex\/ 
$\cR'$ is c-dualizing iff 
there is an invertible sheaf\/~$\cL$ and an integer\/~$n$ such that\/
$\cR'\cong\cR\otimes\cL[n]$.
\end{sprop}

\begin{proof}
Part (a) is proved in \cite[Theorem 5.6]{Ye}.

\smallskip
Now for a fixed invertible sheaf $\cL$ there is a natural isomorphism
of functors
\begin{equation}\label{iso}
\BL(\F\otimes\cL)\iso\BL\F\otimes\cL\qquad(\F\in\D),
\end{equation}
as one deduces, e.g., from a readily-established natural isomorphism
between the respective right adjoints
$$
\BG\E\otimes\cL^{-1}\osi\BG(\E\otimes\cL^{-1})\qquad(\E\in\D).
$$

Part (b) results, because 
 $\BG\cR'$ and~$\BG\cR$ are t-dualizing (\Lref{L:interchange}), so that by (a)
(and taking $\F\set\BG\cR[n]$ in \eqref{iso}) we have isomorphisms
\begin{flalign*}
\hskip61pt\cR'\cong\BL(\BG\cR') &\cong\BL(\BG\cR\otimes\cL[n])\qquad
&&\quad(\text{$\cL$ invertible, $n\in\mathbb Z$})\\ 
&\cong(\BL\BG\cR)\otimes\cL[n]
\cong \cR\otimes\cL[n]. \mkern-3mu 
\end{flalign*}
\vskip-3.8ex 
\end{proof}

\begin{scor}\label{C:Dc*}
If\/ $\X$ is locally embeddable
in a regular finite-dimensional formal scheme then any t-dualizing
complex on\/~$\X$ lies in\/~$\Dc^*$.
\end{scor}

\begin{proof} Whether a t-dualizing complex~$\cR$ satisfies
$\BL\cR\in\Dc$ is a local question, so we may assume that $\X$ is a closed
subscheme of a finite-dimensional regular formal scheme, and then
\Eref{regular}(5) shows that  \emph{some}---hence by
\Pref{P:uniqueness}, \emph{any}---t-dualizing complex lies in~$\Dc^*$.
\end{proof}

\begin{slem}\label{L:Hom=RHom} 
Let\/ $\X$ be a locally noetherian formal
scheme, let $\cI$
be a bounded complex of $\Aqct(\X)$-injectives,  say\/ $\cI\>^i=0$ for
all\/~$i>n,$ and
let\/ $\F\in\mathbf D^+(\X),$ say
\mbox{$H^\ell(\F)=0$} for all\/~$\ell<\!-m.$ 
Suppose there exists an open cover\/ $(\X_\alpha)$ of\/ $\X$
by completions of ordinary noetherian schemes\/~$X_\alpha$ along closed subsets,
with completion maps\/ $\kappa_\alpha\colon\X_\alpha\to X_\alpha\>,$ such that
for each\/ $\alpha$
the restriction  of\/ $\F$ to $\X_\alpha$ is\/ 
$\mathbf D$-isomorphic to\/~$\kappa_{\<\alpha}^*F^{}_{\<\alpha}$ for 
some $F_\alpha\in\mathbf D(X_\alpha).$ Then 
$$
\Ext^i(\F\<,\>\cI\>)\set H^i\R\sHomb_\X(\F\<,\>\cI\>)=0\quad\textup{for all\/~$i>m+n$.}
$$
Moreover,  if\/ $\X$ has finite Krull dimension\/~$d$ then 
$$
\textup{Ext}^i(\F\<,\>\cI\>)\set H^i\R\Homb_\X(\F\<,\>\cI\>)=0\quad\textup{for all\/~$i>m+n+d$.}
$$
\end{slem}

\emph{Remarks.} In the published version of this paper (Contemporary Math.~244) 
\Lref{L:Hom=RHom} stated:
\emph{Let\/ $\F\in\Dvc$ and let\/ $\cI$ be a bounded-below complex of\/
$\Aqct$-injectives. Then the canonical map is a\/ $\D$-isomorphism}
$$
\sHomb(\F\<,\cI\>)\iso\R\sHomb(\F\<,\cI\>).
$$
Suresh Nayak pointed out that the proof given applies
only to $\Avc$-complexes,  
not, as asserted, to arbitrary $\F\in\Dvc$. (Cf.~\cite[Corollary 4.3]{Ye}.)
\Lref{L:Hom=RHom} is used four times  in \S 2.5, so these four places need to be revisited. (There are no other references to Lemma~2.5.6 in the paper.)

\smallskip
First, in Remark (2) on p.\,24, the reference to Lemma~2.5.6 is not
necessary: the cited theorem 4.8 in \cite{Ye} (see also \Pref{1!}
below)  shows that the t-dualizing complex~
$\cR$ is $\mathbf D$-isomorphic to a bounded-below
complex~$\mathcal X'{}^\bullet$ of $\Aqct$-injectives; and then one can proceed
as indicated to show that 
for some~$n$ the (bounded) truncation~$\sigma_{\le n} \mathcal X'{}^\bullet$ is
$\Aqct$-injective and $\mathbf D$-isomorphic
to~$\mathcal X'{}^\bullet$. (To follow the details, it helps to keep in mind 
5.1.3 and~ 5.1.4 below.)

Since, 
by Remark (2), any t-dualizing complex is $\mathbf D$-isomorphic to a bounded
complex of $\Aqct$-injectives, in view of Propositions 3.3.1 and~
5.1.2 one finds that the remaining three references to Lemma 2.5.6
can be replaced by references to the present \Lref{L:Hom=RHom}.
For the reference in the proof of 2.5.7(b) this is clear. The same is true for 
Remark ~(4) on p.\,25, but $i>n_0$ at the end should be~$i>n_0 + d$, where, by Remark~(5), the Krull dimension~$d$ of $\>\X$ is finite. Finally, for the reference in the proof of 2.5.12, one can note, via 5.1.4 and~5.1.2,  that $\Dc^*\subset\Dqct\subset\Dvc\>$.\vspace{1pt} 

{\sc Proof of \ref{L:Hom=RHom}}. 
By the proof of \cite[Proposition 4.2]{Ye},  $\Aqct$-injectives are just
direct sums of sheaves of the form $\cJ(x)\ (x\in\X)$, where for any open
$\U\subset\X$,  $\Gamma(\U,\cJ(x))$~is a fixed injective hull of the residue field of $\cO_{\X\<,x}$ if $x\in\U$, and vanishes otherwise. Hence the restriction of an $\Aqct(\X)$-injective to an open $\V\subset\X$  is $\Aqct(\V)$-injective; and so the first assertion is local. Thus to prove it one may assume that $\X$ itself is a completion, with completion map 
$\kappa\colon\X\to X\<$, and that in $\D(\X)$, $\F\cong\kappa^*\<F$ for some $F\in\D(X)$.
As $\kappa^*\<$, being exact, commutes with the truncation functor~ 
$\sigma_{{\scriptscriptstyle\ge} -m}\>$, there are \mbox{$\mathbf D$-isomorphisms}
(the first as in \cite[p.\,70]{H1}): \looseness=-1
$$
\F\cong\sigma_{{\scriptscriptstyle\ge} -m}\>\F
\cong\sigma_{{\scriptscriptstyle\ge} -m}\kappa^*\<F
\cong \kappa^*\sigma_{{\scriptscriptstyle\ge} -m}\>F\>;
$$ 
so 
one can replace~$F$ by~$\sigma_{{\scriptscriptstyle\ge} -m}\>F$ and 
assume further that $F^\ell=0$ for all\/~$\>\ell<-m\>$.\vspace{1pt}

From the above description of $\Aqct$-injectives, one sees that  
$\kappa_*\cI$ is a bounded complex of
$\cO_{\<\<X}$-injectives, vanishing in degree $>n\>$. 
Since $\kappa_*$ is exact, therefore for all $i>m+n$,
\begin{align*}
\kappa_*H^i\R\sHomb_\X(\F,\>\cI\>)&\cong
H^i\kappa_*\R\sHomb_\X(\kappa^*\< F,\>\cI\>)\\
&\cong H^i\R\sHomb_X(F,\>\kappa_*\>\cI\>)
\quad\qquad\qquad\textup{\cite[p.\,147, 6.7(2)]{Sp}}\\
&\cong
H^i\sHomb_X(F,\>\kappa_*\>\cI\>)=0\>,
\end{align*}
and hence $H^i\R\sHomb_\X(\F,\>\cI\>)=0$. \vspace{1pt}

If $\X$ has Krull dimension~ $d$, and $\Gamma\set\Gamma(\X,-)$ is the global-section functor,  then by a well-known theorem of Grothendieck
the restriction of the derived \mbox{functor~$\R\Gamma$} to the category of abelian sheaves
has cohomological dimension $\le d\>$; and so since 
$\R\Homb_\X\cong \R\Gamma\>\R\sHomb_\X$ 
\cite[Exercise 2.5.10(b)]{Derived categories}, the second assertion\vspace{.6pt} 
follows from \cite[Remark 1.11.2(iv)]{Derived categories}.
\vspace{2pt}\hfill{$\square$}



\Pref{P:dualizing} below brings out the basic property of the 
\emph{dualizing functors} associated with dualizing complexes.
(For illustration, one might keep in mind the special case of \Eref{regular}(3).)

\begin{slem}\label{L:dualizing}
Let\/ $\cR$  be a c-dualizing\/  complex on\/~$\X,$ let\/ $\cRt$ be
the t-dualizing complex $\cRt\set\BG\cR,$ and for any\/
$\E\in\D$ set
$$
\cD\E\set\R\sHomb(\E\<,\>\cR),\qquad
\cDt\>\E\set\R\sHomb(\E\<,\>\cRt).
$$
\textup{(a)} There are functorial isomorphisms
$$
\BL\cDt\cong\BL\cD\cong\cD\BL\cong\cD\cong\cD\>\BG\cong\cDt\>\BG\<.
$$
\textup{(b)} For all\/ $\F\in\Dc\>,$ $\cD\F\in\Dc$ and there is a natural
isomorphism\/
$ \cDt\>\F\cong\BG\cD\F$.
\end{slem}

\begin{proof}
(a) For any $\E\in\D$, \Pref{formal-GM} gives the isomorphism
$$
\cD\E=\R\sHomb(\E\<,\cR)
\cong
\R\sHomb(\BG\E\<,\cR)
=\cD\BG\E.
$$
In particular, $\cD\BL\>\E\cong\cD\BG\<\BL\>\E\cong\cD\BG\E$.
Thus $\cD\cong\cD\>\BG\cong\cD\BL$.

\enlargethispage*{1.5\baselineskip} 
Furthermore, using that the natural map
$\smash{\BG\cO_\X\Otimes\E}\to\BG\E$ is an \emph{isomorphism}
(localize, and see \cite[p.\,20, Corollary (3.1.2)]{AJL}) we get natural isomorphisms
\vspace{-1pt}
\begin{multline*}
\R\sHomb\bigl(\BG\cO_\X\>,\>\sHomb(\E,\cR)\bigr)\iso\sHomb(\BG\E\<,\>\cR)
\underset{\textup{\ref{C:Hom-Rgamma}}}\cong
\R\sHomb(\BG\E,\BG\cR)\\
\cong
\R\sHomb\bigl(\BG\cO_\X\>,\>\sHomb(\E,\BG\cR)\bigr),
\end{multline*}
giving $\BL\cD\cong\cD\BG\cong\cDt\BG\cong\BL\cDt$.

\smallskip
(b) Given remark (2) following \Dref{D:dualizing}, \Lref{L:Hom=RHom} implies
that the functor 
$\cDt\set\R\sHomb(-, \cRt)$ is bounded on~$\Dvc\>$.
The same holds for~$\cD=\cDt\BG$ (see (a)),
 because $\BG(\Dvc)\subset\Dqct\subset\Dvc$ (\Lref{Gamma'+qc}), and
 $\BG$ is bounded. ($\BG$~is given locally by tensoring with a  bounded
flat complex $\cK_\infty^\bullet\>$, see proof of \Pref{Gamma'(qc)}(a)).

Arguing as in \Pref{P:Rhom}, we see that 
\mbox{$\cDt\>\F\set\R\sHomb(\F\<,\>\cRt)\in\Dqct$}, so that
$\BG\cDt\>\F\iso\cDt\>\F$  (\Pref{Gamma'(qc)}(a)); and similarly, 
$\cD\F\in\Dc\>$. Furthermore, the argument in \Rref{R:Dtilde}(4) gives an
isomorphism $\BG\cDt\>\F\>\cong\BG\cD\F\<.$ 
\end{proof}

\begin{sprop}\label{P:dualizing}
With notation as in \Lref{L:dualizing} we have, for\/ $\E,\F\in \D$$:$

\smallskip
\noindent\textup{(a)} $\E\in\Dc^*\!\iff\!\cDt\>\E\in\Dc$ and the natural map is an
isomorphism $\,\E\!\iso\<\cDt\cDt\>\E$.

\smallskip
\noindent\textup{(b)} $\F\in\Dc\!\iff\!\cD\F\in\Dc$ and the natural map is an
isomorphism $\,\F\!\iso\<\cD\cD\F$.

\smallskip
\noindent\textup{(c)} $\,\F\in\Dc\!\iff\!\cDt\>\F\in\Dc^*$ and the natural map is
an isomorphism $\,\F\!\iso\<\cDt\cDt\>\F$.
\end{sprop}

\begin{small}
\emph{Remark.} The isomorphism $\F\iso\cDt\cDt\>\F\>$ is a formal version  of
``Affine Duality, " see \cite[\S5.2]{AJL}.\index{Duality!Affine}
\end{small}

\smallskip
\begin{proof} For $\F\in\Dc\>$,  \Lref{L:dualizing}(b) gives $\cD\F\in\Dc$, so
\mbox{$\cDt\F\cong\BG\cD\F\in\Dc^*$.}  Moreover, from the isomorphism
$\cDt\BG\F\cong\cD\F$ of 
\Lref{L:dualizing}(a) it follows that $\cDt(\Dc^*)\subset\Dc$. The
$\Longleftarrow$ implications in (a), (b) and (c) result, as do the first parts of the
$\implies$ implications.

Establishing the  isomorphisms 
$\cD\cD\F\osi\F\iso\cDt\cDt\>\F$ is a local problem, so we may assume
$\X$ affine.  Since  the functors $\cD$ and $\cDt$ are bounded on~
$\Dvc$ (see proof of \Lref{L:dualizing}(b)), and both  take $\Dc$ into $\Dvc\>$,
therefore the functors
$\cD\cD$ and $\cDt\cDt$ are bounded on~$\Dc\>$, and so \cite[p.\,68, 7.1]{H1}
(dualized) reduces the problem to the tautological case
$\F=\cO_\X$ (cf.~\cite[p.\,258, Proposition~2.1]{H1}.)

\smallskip
For  assertion (a) one may assume that $\E=\BG\F\ (\F\in\Dc)$, so that there is a
composed isomorphism (which one checks to be the natural map):
$$
\E=\BG\F \cong  \BG\cD\cD\F 
\underset{\textup{\ref{L:dualizing}(b)}}\cong 
\cDt\cD\F
\underset{\textup{\ref{L:dualizing}(a)}}\cong \cDt\cDt\BG\F = \cDt\cDt\>\E.
\vspace{-5ex}
$$ 

\end{proof}

\smallskip
\begin{scor}\label{P:ducomp} With the preceding notation,

\smallskip
\textup{(a)} The functor\/ $\cD$ induces an 
involutive anti-equivalence of\/~$\Dc$
with itself.

\smallskip
\textup{(b)} The functor\/ $\cDt$ induces quasi-inverse anti-equivalences
between\/ $\Dc$ and\/~$\Dc^*$.
\end{scor}

\begin{slem}\label{L:building}
Let\/ $\J$ be an ideal of definition of\/~$\X$. Then a complex\/ $\cR\in\Dc$
$(\<$resp.~$\cR\in\Dqct)$ is c-dualizing 
$($resp.~t-dualizing$)$ iff for
every\/ $n>0$ the complex\/
$\R\sHomb(\cO_\X/\J^n,\>\cR)$ is dualizing on the scheme\/
$X_n\set(\X,\cO_\X/\J^n)$.
\end{slem}

\begin{proof}
Remark (1) after \Dref{D:dualizing} makes it straightforward to see that if
$\cR$ is either c- or t-dualizing on~$\X$ then $\R\sHomb(\cO_\X/\J^n,\>\cR)$ is
dualizing on $X_n$.

For the converse, to begin with, \Cref{C:Hom-Rgamma} gives
$$
\R\sHomb(\cO_\X/\J^n\<,\>\cR)=\R\sHomb(\cO_\X/\J^n\<,\>\BG\cR),
$$
and it follows from \Lref{L:interchange} that it suffices to consider the
t-dualizing case. So suppose  that $\cR\in\Dqct$ and that for
all~$n$,\vspace{.5pt}
$\R\sHomb(\cO_\X/\J^n,\>\cR)$ is dualizing on $X_n$. Taking
\smash{$\tilde\cR=\cR$} in the proof of \cite[Theorem 5.6]{Ye},\vspace{.5pt} 
one gets $\cO_\X\iso\R\sHomb(\cR, \cR)$. 

\goodbreak
It remains to check condition (iii) in \Dref{D:dualizing}.
We may assume $\cR$ to be K-injective, so that
$\cR_n\set\sHomb(\cO_\X/\J^n\<,\>\cR)$ is K-injective on~$X_n$ for
all~$n$.  Then, since $\iGp\X\cR\cong\R\iGp\X\cR\cong\cR$
(\Pref{Gamma'(qc)}(a)),  
$$
H^i\>\cR\cong
H^i\<\iGp\X\cR\cong
H^i\>\>\smash{\dirlm{n}}\!\cR_n\cong
\smash{\dirlm{n}}\!H^i\>\cR_n\qquad(i\in\mathbb Z).
$$

\smallskip\noindent
For each $n$, $\cR_n$ is quasi-isomorphic to a \emph{residual complex,} which is
an injective $\cO_{\!X_n}\<$-complex vanishing in degrees outside a certain finite
interval $I\set[a,b]$ (\cite[pp.\:304--306]{H1}). If $m\le n$, the same holds---with
the same $I$---for the complex
\mbox{$\cR_m\cong\sHom_{X_n}(\cO_\X/\J^m\<, \cR_n)$.}  It follows that
$H^i\>\cR=0$ for
$i\notin I$.  In particular, $\cR\in\Dqct^+\>$. 

So now we we may assume that
$\cR$ is a bounded-below complex of $\Aqct$-injectives \cite[Theorem 4.8]{Ye}.
Then for any coherent torsion sheaf~$\M$, the homology of
$$
\R\sHomb_\X(\M,\>\cR)
\underset{\textup{\ref{L:Hom=RHom}}}\cong
\sHomb_\X(\M,\>\cR)\cong\sHomb_\X(\M,\>\dirlm{n}\!\cR_n) 
\cong\dirlm{n}\sHomb_\X(\M,\>\cR_n)\vspace{-5pt}\\
$$
vanishes in all degrees $>b$, as required by (iii).
\end{proof}

\begin{sprop}\label{P:twisted inverse}
Let $f\colon\X\to\Y$ be a pseudo\kern.6pt-proper map of noetherian
formal schemes. \vspace{1pt}

\textup{(a)} If\/ $\cR$ is a t\kern.6pt-dualizing complex on\/ $\Y,$ then\/
$\ft\cR$ is  t\kern.6pt-dualizing  on\/~$\X$.

\smallskip
\textup{(b)} If\/ $\cR$ is a c-dualizing complex on\/ $\Y,$ 
then\/ $\ush f\cR$ is c-dualizing on\/~$\X$.
\end{sprop}

\begin{proof}
(a) Let $\J$ be a defining ideal of~$\X$, and let $\I$ be\vspace{2pt} 
a defining ideal of $\Y$ such that
$\I\cO_\X\subset\J$. Let
$X_\J\set(\X,\cO_\X/\J)\overset{{\vbox to
0pt{\vskip-5pt\hbox{$\scriptstyle \,j$}\vss}}}\hookrightarrow\X$ and 
$Y_\I\set (\Y,\cO_\Y/\I)\overset{{\vbox to
0pt{\vskip-3.5pt\hbox{$\scriptstyle i$}\vss}}}\hookrightarrow\Y$ be the resulting
closed immersions. \Eref{ft-example}(4) shows that 
$i_{\mathrm t}^{\<\times}\cR \cong \R\sHomb(\cO_\Y/\I,\>\cR),$
which is a dualizing complex on~$Y_\I$. Pseudo\kern.6pt-properness of~$f$
means  the map
$f_{\I\J}\colon X_\J\to Y_\I$ induced by
$f$ is proper, so as in \cite[p.\,396, Corollary 3]{f!}
(hypotheses about finite Krull dimension being unnecessary here for 
the existence of $f_{\<\mathrm t}^{\mkern-1.5mu\times}\!\<$, etc.),
$$
\R\sHomb(\cO_\X/\J,\>\ft\cR)\cong
j_{\mathrm t}^{\<\times}\<\ft\cR\cong
(f_{\I\J})_{\mathrm t}^{\<\times}i_{\mathrm t}^{\<\times}\cR
$$ 
is a dualizing complex on $X_\J$. The assertion is given then by \Lref{L:building}.

\smallskip
(b) By \Pref{P:coherence},
$\ush f\cR\in\Dc(\X)$.  By \Cref{C:identities}, \Lref{L:interchange}(b), and the
just-proved assertion (a),
$$
\R\iGp\X\ush f\cR\cong\ft\cR\cong\ft\<\R\iGp\Y\cR,
$$
is t-dualizing on~$\X$. So by
\Lref{L:interchange}(a), $\ush f\cR\cong\BL_\X\R\iGp\X\ush
f\cR$ is c-dualizing. 
\end{proof}

\pagebreak[3]
The following proposition generalizes \cite[p.\,291, 8.5]{H1} (see also \cite[middle
of p.\,384]{H1} and \cite[p.\,396, Corollary 3]{f!}).

\begin{sprop}\label{P:Dual!}
Let\/ $f\colon\X\to\Y$ be a pseudo\kern.6pt-proper map of noetherian formal
schemes. Suppose that\/ $\Y$ has a c-dualizing
complex\/~$\cRc,$ or equivalently 
\textup(by \Lref {L:interchange}\,\textup{),} a
t\kern.6pt-dualizing complex\/~$\cRt\in\Dc^*(\Y),$ so that $\ush f\cRc$ is
c-dualizing
$($resp.~$\ft\cRt$ is t\kern.6pt-dualizing\/$)$ on\/~$\X$ $($\Pref{P:twisted
inverse}$\mkern1.5mu)$. Define dualizing\index{dualizing functors}
functors\looseness=-1
\begin{align*}
\cDt^\Y(-)&\set\R\sHomb_\Y(-,\cRt), &
\qquad \cDc^\Y(-)&\set\R\sHomb_\Y(-,\cRc), \\
\cDt^\X(-)&\set\R\sHomb_\X(-,\ft\cRt), &
\qquad \cDc^\X(-)&\set\R\sHomb_\X(-,\ush f\cRc).\mathstrut
\end{align*}
Then there are natural isomorphisms
\begin{alignat*}{2}
\ft\<\E &\cong\cDt^\X\bL
f^*\cDt^\Y\E,\qquad&&\bigl(\E\in\Dc^*(\Y)\cap\D^+(\Y)\bigr),\\
\ush f\<\E &\cong\cDc^\X\bL f^*\cDc^\Y\E\qquad&&\bigl(\E\in\Dc^+\<(\Y)\bigr).
\end{alignat*}
\end{sprop}

\begin{proof}
When $\E\in\Dc^*(\Y)\>\cap\>\>\D^+\<(\Y)$ (resp.\ $\E\in\Dc^+(\Y))$ set
$\F\set\cDt^\Y\E$ (resp.\ $\F\set\cDc^\Y\E$). In either case, $\F\in\Dc(\Y)$
(\Pref{P:dualizing}), and also $\F\in\D^-(\Y)$---in the first case by remark (2)
following \Dref{D:dualizing} and \Lref{L:Hom=RHom}, in the second case by
remark (1) following \Dref{D:dualizing}.  
So, by \Pref{P:dualizing}, we need to find natural isomorphisms
\begin{align*}
\ft\cDt^\Y\F&\cong\cDt^\X\bL f^*\<\F, \\
\ush f\cDc^\Y\F &\cong \cDc^\X\bL f^*\<\F\<.
\end{align*}
Such isomorphisms are given by the next result---a generalization 
of~\cite[p.\,194, 8.8(7)]{H1}---for $\G\set\cRt$  (resp.~$\cRc$).%
\end{proof}

\pagebreak[3]
\begin{sprop}\label{P:Hom!}
Let\/ $f\colon\X\to\Y$ be a map of noetherian formal schemes. Then for\/ 
$\F\in\Dc^-(\Y)$ and\/ $\G\in\D^+(\Y)$ there are  natural isomorphisms
\begin{gather*}
\R\sHomb_\X(\bL f^*\<\F\<, \ft\G)\iso\ft\R\sHomb_\Y(\F\<,\G),\\
\R\sHomb_\X(\bL f^*\<\F\<, \ush f\G)\iso\ush f\R\sHomb_\Y(\F\<,\G).
\end{gather*}
\end{sprop}

\begin{proof}
The second isomorphism follows from the first, since $\ush f=\BL\ft$ and since
there are natural isomorphisms
\begin{align*}
\BL\R\sHomb_\X(\bL f^*\<\F\<, \>\ft\G)
&=\R\sHomb_\X\bigl(\R\iGp\X\cO_\X^{}\>, 
   \>  \R\sHomb_\X(\bL f^*\<\F\<,\>\ft\G)\bigl)\\
&\cong\smash{\R\sHomb_\X\bigl(\R\iGp\X\cO_\X^{}\Otimes
\bL f^*\<\F\<,\>\ft\G)\bigl)}\\ 
&\cong\R\sHomb_\X\bigl(\bL f^*\<\F\<,
  \R\sHomb_\X(\R\iGp\X\cO_\X^{},\>\ft\G)\bigl)\\
&=\R\sHomb_\X(\bL f^*\<\F\<, \>\ush f\G).
\end{align*}

For fixed $\F$ the source and target of the first isomorphism in \Pref{P:Hom!}
are functors from $\D^+(\Y)$ to $\Dqct(\X)$ (see \Pref{P:Rhom}),
right adjoint, respectively,\vspace{.7pt} to the functors
\smash{$\Rfs(\E\Otimes\bL f^*\<\F)$} and 
\smash{$\Rfs\E\Otimes\F$}
$(\E\in\Dqct(\X))$. The functorial ``projection" map  
$$
\smash{\Rfs\E\Otimes\F\to\Rfs(\E\Otimes\bL f^*\<\F),}
$$
is, by definition, adjoint to the natural composition
$$
\smash{\bL f^*(\Rfs\E\Otimes\F\>)\to\bL f^*\Rfs\E\Otimes\bL f^*\<\F
\to\E\Otimes\bL f^*\<\F;}
$$
 and it will suffice to show that this projection map is an
isomorphism. 

For this,
the standard strategy is to localize to where $\Y$ is affine,
then use boundedness of some functors, and  compatibilities with direct
sums, to reduce to the trivial case $\F=\cO_\Y$. Details appear,
e.g., in \cite[pp.\,123--125,~Proposition 3.9.4]{Derived categories}, modulo the
following substitutions: use $\Dvc$ in place of $\Dqc$, and for boundedness and
direct sums use \Lref{Gamma'+qc} and  Propositions
\ref{Rf_*bounded}(b) and~\ref{P:coprod} below.
\end{proof}\index{dualizing complexes|)}

\section{Direct limits of coherent sheaves on formal schemes.}
\label{properly}

In this section we establish, 
for a locally noetherian formal scheme~$\X$,  properties\-
of $\Avc(\X)$ needed in \S\ref{sec-th-duality} to
adapt Deligne's\index{Deligne, Pierre} proof of global Grothendieck Duality\- to the
formal context.  The basic result, \Pref{(3.2.2)}, is that
$\Avc(\X)$ is~\emph{plump} (see opening remarks in \S\ref{S:prelim}), hence
abelian, and so\ (being closed under~\smash{$\dirlm{}\!$})\vadjust{\kern.7pt}
cocomplete, i.e., it has arbitrary small colimits.  This enables us to speak
about~$\D(\Avc(\X))$, and to apply standard adjoint functor theorems to
colimit\kern.5pt-preserving functors on~$\Avc(\X)$. (See e.g.,
\Pref{A(vec-c)-A},  Grothendieck Duality for the identity map of~$\X$).

The preliminary paragraph~\ref{SS:vc-and-qc} sets up an equivalence of
categories which allows us to reduce local questions about the
(globally defined) category $\Avc(\X)$ to corresponding questions
about quasi-coherent sheaves on ordinary noetherian
schemes. Paragraph~\ref{SS:Dvc-and-Dqc} extends this equivalence to
derived categories.  As one immediate application,
\Cref{corollary} asserts that the natural functor
$\D(\Avc(\X))\to\Dvc(\X)$ is an equivalence of categories when $\X$ is
\emph{properly algebraic,}\index{properly algebraic} i.e., the $J$-adic
completion of a proper
$B$-scheme with $B$ a noetherian ring and $J$ a $B$-ideal.  This will
yield a stronger version of Grothendieck Duality on such formal
schemes---for $\Dvc(\X)$ rather than $\D(\Avc(\X))$, see
\Cref{cor-prop-duality}. We do not know whether such global
results hold over arbitrary noetherian formal schemes.

\pagebreak[3] 
Paragraph~\ref{SS:bounded} establishes boundedness 
for some derived functors, a
condition which allows us to apply them freely to  unbounded complexes,
as illustrated, e.g., in Paragraph~\ref{3.5}.

\begin{parag}\label{SS:vc-and-qc}
For $X$ a noetherian ordinary scheme,  $\Avc(X)=\Aqc\<(X)$
\cite[p.\,319,~6.9.9]{GD}. The inclusion
$j_{\lower.2ex\hbox{$\scriptstyle\<\< X$}}\colon\Aqc\<(X) \to \A(X)$ has
a right adjoint $Q_{\<\<X}\colon\A(X) \to \Aqc\<(X)$, the ``quasi-coherator,''
\index{ $\R$@ {}$Q_\X$ (quasi-coherator)\vadjust{\penalty
10000}} necessarily left exact \cite[p.\,187, Lemme 3.2]{I}. (See
\Pref{A(vec-c)-A} and~\Cref{C:Qt} for generalizations to
formal schemes.)
\begin{sprop}
\label{(3.2.1)}
Let\/ $A$ be a noetherian adic ring with ideal of definition\/~$I,$\
let\/ $f_0\colon X\to\Spec(A)$ be a proper map, 
set $Z:=f_0^{-1}\Spec(A/I),$\ and let 
$$
\kappa\colon\X= X_{/Z}\to X
$$  
be the
formal completion of\/~$X$ along\/~$Z$. Let\/ $Q\!:=Q_{\<\<X}$ be
as above.  Then\/ $\kappa^*$ induces equivalences of categories from\/
$\Aqc\<(X)$ t\/o $\Avc(\X)$ and from\/ $\Ac(X)$ to\/ $\Ac(\X),$\ both
with quasi-inverse\/~$Q\kappa_*$.
\end{sprop}

\begin{proof}
For any quasi-coherent $\cO_{\<\<X}$-module~$\G$ the canonical maps are 
\emph{isomorphisms}
\stepcounter{sth}
\renewcommand{\theequation}{\thesth} 
\begin{equation}\label{3.2.1.1}
{\rm H}^i(X\<,\G)\iso {\rm H}^i(X\<,\kappa_*\kappa^*\G) = {\rm H}^i(\X,
\kappa^*\G)\qquad(i\ge0).
\end{equation}\stepcounter{sth}%
(The equality holds because $\kappa_*$ transforms any flasque
resolution of~$\kappa^*\G$ into one of~$\kappa_*\kappa^*\G$.)

For, if $(\G_\lambda)$
is the family of coherent submodules of~$\G$, ordered by
inclusion, then $X$ and $\X$ being noetherian, one checks that
(\ref{3.2.1.1}) is the composition of the sequence of natural
isomorphisms
 \begin{alignat*}{2} 
{\rm H}^i(X\<,\G)&\iso {\rm H}^i(X,
 \>\dirlm{\lambda}\G_\lambda)&
&\hskip-60pt\mbox{\cite[p.~319,~(6.9.9)]{GD}}\\ &\iso  \dirlm{\lambda}
{\rm H}^i(X\<,\G_\lambda) &&\\
 &\iso \dirlm{\lambda} {\rm H}^i(\X,\kappa^*\G_\lambda)&
&\hskip-60pt\mbox{\cite[p.~125,~(4.1.7)]{EGA}}\\
&\iso  {\rm
H}^i(\X,\>\dirlm{\lambda}\kappa^*\G_\lambda)\\
&\iso{\rm
H}^i(\X,\kappa^*\dirlm{\lambda} \G_\lambda)
\iso   {\rm H}^i(\X,\kappa^*\G).
\end{alignat*}

Next, for any $\G$ and $\cH$ in $\Aqc\<(X)$ the natural map is an
\emph{isomorphism}
\begin{equation}
\Hom_X(\G,\cH)\iso \Hom_\X(\kappa^*\G,\>\kappa^*\cH)
\label{3.2.1.2}
\end{equation}
For, with $\G_\lambda$ as above, (\ref{3.2.1.2}) factors as the
sequence of natural isomorphisms
\begin{align*}
\Hom_{X}(\G,\cH)&\iso  
  \inlm{\lambda}\Hom_{X}(\G_\lambda\>,\cH) \\
&\iso \inlm{\lambda} {\rm H}^0\<\bigl({X\<,
\sHom_{X}(\G_\lambda\>,\cH)}\bigr) \\
&\iso \inlm{\lambda} {\rm H}^0\<\bigl({\X,\kappa^*
\sHom_{X}(\G_\lambda\>,\cH})\bigr)
 \qquad\bigl({\mbox{see }(\ref{3.2.1.1})}\bigr)\\
&\iso \inlm{\lambda} {\rm H}^0\< \bigl({\X,
\sHom_{\X}(\kappa^*\G_\lambda\>,\kappa^*\cH)}\bigr)  \\
&\iso \inlm{\lambda}\Hom_\X(\kappa^*\G_\lambda\>,\kappa^*\cH)\\ 
&\iso\Hom_\X(\dirlm{\lambda}
\kappa^*\G_\lambda\>,\kappa^*\cH) 
\iso \Hom_\X(\kappa^*\G,\>\kappa^*\cH).
\end{align*}

Finally, we show the equivalence of the following conditions, for
$\F\in\A(\X)$:
\begin{list}%
{(\arabic{t})} {\usecounter{t} \setlength{\rightmargin}{\leftmargin}}
\item \label{1} The functorial map
$\alpha(\F\>)\colon\kappa^*Q\kappa_*\F\to\F$ (adjoint to the canonical map
$Q\kappa_*\F\to \kappa_*\F\>)$ is an isomorphism.
\item \label{2} There exists an
isomorphism $\kappa^*\G\iso\F$ with  $\G\in\Aqc\<(X)$.
\item \label{3} $\F\in\Avc(\X)$.
\end{list}

Clearly $(\ref{1}) \Rightarrow (\ref{2})$; and $(\ref{2}) \Rightarrow
(\ref{3})$ because
\smash{$\dirlm{}_{\!\!{}_\lambda}\,\kappa^*\G_\lambda\iso \kappa^*\G$}
($\G_\lambda$ as before).

Since $\kappa^*$ commutes with $\smash{\dirlm{}}$\vspace{1pt} 
and induces an equivalence of
categories from
$\Ac(X)$ to $\Ac(\X)$ \cite[p.~150,~(5.1.6)]{EGA}, we see that
$(\ref{3}) \Rightarrow (\ref{2})$.\vspace{1pt}

For $\G\in\Aqc(X)$, let $\beta(\G)\colon\G\to Q\kappa_*\kappa^*\G$ be
the canonical map (the unique one whose composition with
$Q\kappa_*\kappa^*\G\to \kappa_*\kappa^*\G$ is the canonical map
$\G\to\kappa_*\kappa^*\G$). Then for any $\cH\in\Aqc\<(X)$
we have the natural commutative diagram
$$
\begin{CD}
\Hom(\cH,\G) @>\text{via }\beta>> \Hom(\cH,Q\kappa_*\kappa^*\G)\\
@V\simeq VV @VV\simeq V \\
\Hom(\kappa^*\cH,\kappa^*\G) @>\Iso>>\Hom(\cH,\kappa_*\kappa^*\G)
\end{CD}
$$

\smallskip\noindent 
where the left vertical arrow is an isomorphism by
(\ref{3.2.1.2}), the right one is an isomorphism because $Q$ is right-adjoint to
$\Aqc\<(X)\hookrightarrow\A(X)$, and the bottom arrow is an isomorphism
because $\kappa_*$ is right-adjoint to $\kappa^*$; so ``via $\beta\>$''
is an isomorphism for all~$\cH$, whence \emph{$\beta(\G)$  is an
isomorphism.}   The implication $(\ref{2}) \Rightarrow (\ref{1})$ follows now from
the easily checked fact that
$\alpha(\kappa^*\G){\<\smcirc\<}\kappa^*\<\beta(\G)$ is the identity map of
$\kappa^*\G$. 

We see also that $Q\kappa_*(\Ac(\X))\subset\Ac(X)$, since  by
\cite[p.~150,~(5.1.6)]{EGA} every \mbox{$\F\in\Ac(\X)$} is isomorphic to
$\kappa^*\G$ for some $\G\in\Ac(X)$, and $\beta(\G)$ is an
isomorphism. 

Thus we have the functors $\kappa^*\colon\Aqc\<(X)\to\Avc(\X)$ and
$Q\kappa_*\colon\Avc(\X)\to\Aqc\<(X)$, both of which preserve
coherence, and the functorial isomorphisms
$$
\alpha(\F\>)\colon\kappa^*Q\kappa_*\F \iso \F 
\ \ \bigl({\F\in\Avc(\X)}\bigr); \qquad
\beta(\G)\colon\G\iso  Q\kappa_*\kappa^*\G 
\ \ \bigl({\G\in\Aqc\<(X)}\bigr)\<\<.
$$

\Pref{(3.2.1)} results.
\end{proof}

\smallskip
Since $\kappa^*$ is right-exact, we deduce:

\begin{scor}
\label{coker}
For any affine noetherian formal scheme\/~$\X,$\
$\F\in\Avc(\X)$ iff\/ $\F$~is a cokernel of a map of free\/
$\cO_\X$-modules\/ $($i.e., direct sums of copies of\/ $\cO_\X)$.
\end{scor}

\begin{scor}\label{C:vec-c is qc}
For a locally noetherian formal scheme\/ $\X,$\
$\Avc(\X)\subset\Aqc\<(\X),$\ i.e., any\/ \smash{$\dirlm{}\!\!$} of
coherent\/ $\cO_\X$-modules is quasi-coherent.
\end{scor}

\begin{proof}
Being local, the assertion follows from \Cref{coker}.
\end{proof}

\pagebreak[3]

\begin{scor}[cf.~{\cite[3.4, 3.5]{Ye}}]
\label{C:images}
For a locally noetherian formal scheme\/~$\X$ let\/ $\F$ and\/ $\G$ be
quasi-coherent $\cO_\X$-modules. Then\/$:$
 
\textup{(a)} The kernel, cokernel, and image of any\/
$\cO_\X$-homomorphism\/ $\F\to\G$ are quasi-coherent.

\textup{(b)} $\F$ is coherent iff\/ $\F$ is locally finitely generated. 

\textup{(c)} If\/ $\F$ is coherent and\/ $\G$ is a sub- or
quotient module of\/~$\F$ then\/ $\G$ is coherent.

\textup{(d)}  If\/ $\F$ is coherent then\/
$\sHom(\F\<,\G)$ is quasi-coherent; and if also\/ $\G$ is coherent then\/ 
$\sHom(\F\<,\G)$ is coherent.
\textup(For a generalization, see \Pref{P:Rhom}.\textup)
\end{scor}

\begin{proof}
 The questions being local, we may assume
$\X=\Spf(A)$ ($A$~noetherian adic), and, by \Cref{coker}, that
$\F$ and $\G$ are in $\Avc(\X)$. Then, $\kappa^*$~being exact, 
 \Pref{(3.2.1)} with $X\set\Spec(A)$ and $f_0\set\text{identity}$
reduces the problem to noting that the corresponding statements about coherent
and quasi-coherent sheaves on~$X$ are true. 
(These statements are in \cite[p.\,217, Cor.\:(2.2.2) and~p.\,228, \S(2.7.1)]{GD}.
Observe also that if $F$ and $G$ are $\cO_{\<\<X}$-modules with $F$ coherent
then
$\sHom_\X(\kappa^*\<F,\kappa^* G)\cong \kappa^*\sHom_X(F,G)$.)
\end{proof} 

\begin{scor}\label{C:limsub}
For a locally noetherian formal scheme\/ $\X,$ any\/
\mbox{$\F\in\Avc(\X)$} is the\/ \smash{$\dirlm{}\!\<$} of its coherent\/
$\cO_\X$-submodules.
\end{scor}

\begin{proof}
Note that by \Cref{C:images}(a) and~(b)  the sum of any two coherent
submodules of~$\F$ is again coherent.  By definition,
$\F=\smash{\dirlm{}_{\!\!{}_\mu}\>\>\F_\mu}$ with $\F_\mu$ coherent, and
from~
\Cref{C:images}(a) and~(b) it follows that the canonical image of
$\F_\mu$ is a coherent submodule of~$\F\<$, whence the conclusion.
\end{proof}

\begin{scor}
\label{limit(vec-c)=qc}
For any affine noetherian formal scheme\/~$\X,$\ any\/ $\F\in\Avc(\X)$  and
any\/ $i>0,$
$$
{\rm H}^i(\X,\F\>)=0.
$$
\end{scor}

\begin{proof}
Taking  $f_0$ in \Pref{(3.2.1)} to be the identity map,
we have \mbox{$\F\cong\kappa^*\G$} with $\G$
quasi-coherent; and so by \eqref{3.2.1.1},
$\mathrm H^i(\X,\F\>)\cong \mathrm H^i(\Spec(A),\G)=0$.
\end{proof}
\end{parag}

\smallskip
\begin{parag}
\Pref{(3.2.1)} will now be used to show, for locally noetherian formal
schemes~$\X$, that $\Avc(\X)\subset\A(\X)$ is plump, and that this inclusion has a
right adjoint, extending to derived categories.

\begin{slem}
\label{L:Ext+lim}
Let\/ $\X$ be a noetherian formal scheme, let\/ 
$\F\in\Ac(\X),$\ and let\/ 
$(\G_\alpha\>,\gamma_{\alpha\beta}\colon\G_\beta\to
  \G_\alpha)_{\alpha,\>\beta\in \Omega}$
be a directed system in\/~$\Ac(\X)$. Then
for every\/~$q\ge0$ the natural map is an isomorphism
$$
\dirlm{\alpha}\mathrm {Ext}^q(\F\<,\>\G_\alpha)\iso \mathrm
{Ext}^q(\F\<,\>\dirlm{\alpha}\G_\alpha).
$$
\end{slem}

\begin{proof}
For an $\cO_\X$-module $\M$, let $\mathrm E(\M)$ denote the usual
spectral sequence\looseness=-1
$$
\mathrm  E\mspace{.5mu}_2^{pq}(\M)\!:=\mathrm
H^p\bigl(\X,\>\Ext^{\>q}(\F\<,\M)\bigr)
\Rightarrow
\mathrm {Ext}^{p+q}(\F\<,\M).
$$
It suffices that the natural map of spectral sequences be an
isomorphism 
$$
\smash{\dirlm{} \mathrm E(\G_\alpha)\iso 
  \mathrm E(\>\dirlm{} \G_\alpha)\qquad(\dirlm{}\!\<:=\dirlm{\alpha}\!),}
$$ 

\smallskip\noindent
and for that we need only check out the $\mathrm E\mspace{.5mu}_2^{pq}$ terms,
i.e., show that the natural maps
$$
\postdisplaypenalty 10000
\smash{\dirlm{}\mathrm H^p\bigl(\X,\>\Ext^{\>q}(\F\<, \G_\alpha)\bigr)\to
\mathrm H^p\bigl(\X,\>\dirlm{}\Ext^{\>q}(\F\<, \G_\alpha)\bigr)\to
\mathrm H^p\bigl(\X,\>\Ext^{\>q}(\F\<, \dirlm{}\G_\alpha)\bigr)}
$$

\smallskip\noindent
are isomorphisms. 
The first one is, because $\X$ is noetherian.  So we
need only show that the natural map is an isomorphism
$$
\dirlm{}\Ext^{\>q}(\F\<,\>\G_\alpha)\iso \Ext^{\>q}(\F\<,\>\dirlm{}\G_\alpha).
$$

For this localized question we may assume that $\X=\Spf(A)$ with $A$ a
noetherian adic ring. By \Pref{(3.2.1)} (with $f_0$ the
identity map of $X\!:=\Spec(A)$) there is a coherent
$\cO_{\<\<X}$-module~$F$ and a directed system 
$(G_\alpha\>,g_{\alpha\beta}\colon
 G_\beta\to G_\alpha)_{\alpha,\>\beta\in\Omega}$  
of coherent $\cO_{\<\<X}$-modules such that
$\F=\kappa^*F\<$,
$\G_\alpha=\kappa^*G_\alpha\>$, and 
$\gamma_{\alpha,\>\beta}=\kappa^*\<g_{\alpha,\>\beta}$. 
Then the well-known natural isomorphisms (see
\cite[(Chapter 0), p.\,61, Prop.\,(12.3.5)]{EGA}---or the proof of
\Cref{(3.2.3)} below)
\begin{multline*}
\dirlm{}\Ext_\X^q(\F\<,\>\G_\alpha)
\iso\dirlm{}\kappa^*\Ext_{\<\<X}^{\>q}(F\<,G_\alpha)\iso
\kappa^*\dirlm{}\Ext_{\<\<X}^{\>q}(F\<,G_\alpha) \\
\iso\kappa^*\Ext_{\<\<X}^{\>q}(F\<,\>\dirlm{}G_\alpha)\iso \Ext_\X^q(\kappa^*F\<,
\>\kappa^*\dirlm{} \<G_\alpha)\iso \Ext_{\<\<X}^{\>q}(\F\<, \>\dirlm{} \G_\alpha)
\end{multline*}
\penalty10000
give the desired conclusion.
\end{proof}

\begin{sprop}
\label{(3.2.2)}
Let\/ $\X$ be a locally noetherian formal scheme.  If\/
$$
\F_1\to\F_2\to\F\to\F_3\to\F_4
$$ 
is~an exact sequence of\/
$\cO_\X$-modules and if\/ $\F_1\>,$ $\F_2\>,$ $\F_3$ and\/ $\F_4$ are all
in\/ $\Aqc\<(\X)$ $($resp.~$\Avc(\X))$ then\/ $\F \in \Aqc\<(\X)$
$($resp.~$\Avc(\X))$. Thus\/ $\Aqc\<(\X)$ and $\Avc(\X)$ are
plump---hence abelian---subcategories of\/~$\A(\X),$\ and both
$\Dqc\<(\X)$ and its subcategory\/~$\Dvc(\X)$ are triangulated
subcategories of\/ $\D(\X)$. Furthermore,  $\Avc(\X)$ is closed under
arbitrary small\/ $\A(\X)$-colimits.
\end{sprop}

\begin{proof}
Part of the $\Aqc$ case is covered by \Cref{C:images}(a), and all of it
by \cite[Proposition 3.5]{Ye}. At any rate, since every quasi-coherent
$\cO_\X$-module is locally in~$\Avc\subset\Aqc$ (see
Corollaries~\ref{coker}  and~\ref{C:vec-c is qc}), it suffices to treat 
the $\Avc$ case. 

Let us first show that the kernel~$\mathcal K$ of an $\Avc$ map
$$
\mspace{160mu}\psi\colon\smash{\dirlm{}_{\<\<\!\beta}}\>\>\H_\beta=\H\to\G
=\smash{\dirlm{}_{\!\!\alpha}}\>\>\G_\alpha\qquad
(\G_\alpha\>,\H_\beta\in\Ac(\X)) 
$$
is itself in $\Avc(\X)$.  It will
suffice\vadjust{\kern1.5pt} to do so for the kernel~$\mathcal K_\beta$
of the composition 
$$
\smash{\psi_\beta\colon\H_\beta\xrightarrow{\text
{natural}\,}\H\xrightarrow{\psi\>}\G,}
$$ 
since $\mathcal K=\smash{\dirlm{}_{\<\<\!\beta}\>\>\mathcal K_\beta\>}$. 

\pagebreak[2]
By
the case\vadjust{\kern.75pt} $q=0$ of \Cref{L:Ext+lim}, there is
an $\alpha$ such that $\psi_\beta$ factors as\vadjust{\kern-3pt}
$$
\H_\beta\xrightarrow{\psi_{\beta\alpha}\>}\G_\alpha \xrightarrow{\text
{natural}\,}\G\>;
$$
and then with $\mathcal K_{\beta\alpha'}\ (\alpha'>\alpha)$ the
(coherent) kernel of the composed map
$$
\H_\beta\xrightarrow{\psi_{\beta\alpha}\>}\G_\alpha \xrightarrow{\text
{natural}\,}\G_{\alpha'}
$$
we have $\mathcal K_\beta =\dirlm{}_{\<\<\!\alpha'}\>\>\mathcal
K_{\beta\alpha'}\in\Avc(\X)$.

Similarly, we find that $\text{coker}(\psi)\in\Avc(\X)$.  Being
closed under  small direct sums, then, $\Avc(\X)$ is closed
under arbitrary small $\A(\X)$-colimits
\cite[Corollary\,2, p.\,109]{currante}. 

Consideration of the exact sequence
$$
0\longrightarrow\text{coker}(\F_1\to\F_2)\longrightarrow\F\longrightarrow
\ker(\F_3\to\F_4)\longrightarrow 0
$$
now reduces the original question to where $\F_1=\F_4=0$. Since $\F_3$
is the $\smash{\dirlm{}}$ of its coherent submodules
(\Cref{C:limsub}) and $\F$ is the $\smash{\dirlm{}}$ of the
inverse images of those submodules, we need only show that each such
inverse image is in $\Avc(\X)$.  Thus we may assume
$\F_3$ coherent (and $\F_2=\smash{\dirlm{}_{\!\!\alpha}}\>\G_\alpha$
with $\G_\alpha$ coherent).\vspace{1pt}

\penalty-1000

The exact sequence $0\to\F_2 \to \F\to \F_3 \to 0$ represents an
element 
$$
\eta\in\mathrm{Ext}^1(\F_3,\>\F_2)=
\mathrm{Ext}^1(\F_3,\>\smash{\dirlm{}_{\!\!\alpha}}\>\G_\alpha);
$$ 
and by
\Cref{L:Ext+lim}, there is an~$\alpha$ such that $\eta$~is
the natural image of an element
$\eta_\alpha\in\mathrm{Ext}^1(\F_3,\>\G_\alpha)$, represented by an
exact sequence~$0\to\G_\alpha \to \F_\alpha\to \F_3 \to
0$.  Then
$\F_\alpha$~is coherent, and by \cite[p.\,66, Lemma 1.4]{sM75}, we
have an isomorphism
$$
\postdisplaypenalty 10000
\F\iso \F_2\oplus_{\G_\alpha} \<\F_\alpha\>.
$$
Thus $\F$ is the cokernel of a map in~$\Avc(\X)$, and so as
above, $\F\in \Avc(\X)$.
\end{proof}

\begin{sprop}
\label{A(vec-c)-A}\index{ $\R$@ {}$Q_\X$ (quasi-coherator)\vadjust{\penalty
10000}} On a locally noetherian formal scheme\/~$\X,$\ the
inclusion functor\/ 
$j_{\lower.2ex\hbox{$\scriptstyle\X$}}\colon \Avc(\X) \to \A(\X)$ has a right
adjoint\/ $Q_{\X}\colon\A(\X) \to \Avc(\X);$ and\/ $\R Q_\X^{} $ is right-adjoint to
the natural functor $\D(\Avc(\X))\to \D(\X)$.
In particular, if\/ $\kappa\colon\X\to X$
is as in \Pref{(3.2.1)} then\/
$Q_{\X}\cong\kappa^*Q_{\<\<X}\kappa_*$ and\/
$\R Q_{\X}\cong\kappa^*\R Q_{\<\<X}\kappa_*\>$.
\end{sprop}

\begin{proof}
Since $\Avc(\X)$ has a small family of (coherent) generators, and
is closed under arbitrary small
$\A(\X)$-colimits,  the existence of~$Q_\X^{}$ follows  from the Special
Adjoint Functor Theorem\index{Special Adjoint Functor Theorem} 
(\cite[p.\,90]{pF1964} or \cite[p.\,126, Corollary]{currante}).%
\footnote{It follows that $\Avc(\X)$ is closed under \emph{all}
$\A(\X)$-colimits (not necessarily small): if
$F$ is any functor into $\Avc(X)$ and
$\F\in\A(X)$ is a colimit of 
$j_{\<\lower.2ex\hbox{$\scriptscriptstyle\X$}}\smcirc F\<$, then $Q_{\X}\F$
is a colimit of $F\<$, and the natural map is an isomorphism $\F\iso
j_{\<\lower.2ex\hbox{$\scriptscriptstyle\X$}}Q_{\X}\F$.  (Proof:
exercise, given in dual form in \cite[p.\,80]{pF1964}.) \looseness=-1}

In an abelian category~$\A$,
a complex~$J$ is, by definition, K-injective
if for each exact $\A$-complex $G$, 
the complex $\Homb_{\A}(G, J)$ is exact too. Since $j_\X^{}$
is exact, it follows that its right adjoint $Q_\X^{}$ 
transforms K-injective
$\A(\X)$-complexes into \mbox{K-injective} $\Avc(\X)$-complexes,
whence the derived functor $\R Q_\X^{}$ is right-adjoint to the natural
functor $\D(\Avc\<(\X)) \to \D(\X)$ (see \cite[p.\,129,
Proposition~1.5(b)]{Sp}). 

The next assertion is a corollary of~\Pref{(3.2.1)}: any
$\M\in\Avc(\X)$ is isomorphic to $\kappa^*\G$ for some $\G\in\Aqc\<(X)$,
and then for any $\N\in\A(\X)$ there are natural isomorphisms\looseness=1
\begin{align*}
\Hom_\X(j_{\lower.2ex\hbox{$\scriptstyle\X$}}\M,\>\N\>)&\cong
 \Hom_\X(j_{\lower.2ex\hbox{$\scriptstyle\X$}}\kappa^*\G,\>\N\>)\\
&\cong
  \Hom_X(j_{\lower.2ex\hbox{$\scriptstyle\<\<X$}}\G,\> \kappa_*\>\N\>)\cong
   \Hom_{\Aqc\<(X)}(\G,\>Q_{\<\<X}\kappa_*\>\N\>)\\
&\cong
    \Hom_{\Avc(\X)}(\kappa^*\G,\>\kappa^*Q_{\<\<X}\kappa_*\>\N\>)\cong
     \Hom_{\Avc(\X)}(\M,\>\kappa^*Q_{\<\<X}\kappa_*\>\N\>).
\end{align*}
Moreover, since $\kappa_*$ has an exact left adjoint (viz.~$\kappa^*$), therefore,
as above, $\kappa_*$~transforms K-injective $\A(\X)$-complexes into K-injective
$\A(X)$-complexes, and it follows at once that 
$\R Q_\X^{}\cong\kappa^*\R Q_{\<\<X}\kappa_*$. 
\end{proof}

Let\/ $\X$ be a locally noetherian formal scheme. A property~$\mathbf P$ 
of sheaves of modules is \emph{local} if it is defined on~$\A(\U)$
for arbitrary open subsets\/ $\U$ of\/~$\X$,  and is such that\/ for any
$\E\in\A(\U)$ and any open covering\/
$(\U_\alpha)$ of\/~$\U,$\ $\mathbf P(\E)$~holds 
iff\/ $\mathbf P(\E|_{\U_\alpha}\<)$ holds for all~$\alpha$.

For example, coherence  and quasi-coherence are both local properties---to which 
by \Pref{(3.2.2)}, the following Proposition applies.

\begin{sprop}\label{P:Rhom}
Let\/ $\X$ be a locally noetherian formal scheme, and let\/~$\mathbf P$ be a
local property of\/ sheaves of modules.
Suppose further that for all open\/~$\U\subset\X$ the full subcategory\/
$\A_{\mathbf P}\<(\U)$ of\/~$\A(\U)$ whose objects are all the $\E\in\A(\U)$ for
which\/ $\mathbf P(\E)$ holds is a
\emph{plump} subcategory of~$\>\A(\U)$. Then for all\/ $\F\in\Dc^-\<(\X)$ and\/
$\G\in\D_{\mathbf P}^+\<(\X),$\
it holds that ~$\R\sHomb(\F,\G)\in\D_{\mathbf P}^+\<(\X)$.
\end{sprop}

\begin{proof}
Plumpness implies that $\D_{\mathbf P}\<(\X)$ is a triangulated
subcategory of~$\D(\X)$, as is $\Dc(\X)$, so \cite[p.\,68, Prop.\,7.1]{H1} gives
a ``way-out'' reduction to where $\F$ and~$\G$ are
$\cO_\X$-modules. The question being local on $\X$,   we may assume
$\X$ affine and replace
$\F$ by a quasi-isomorphic bounded-above complex $\F\>^\bullet$ of finite-rank
free $\cO_\X$-modules, see \cite[p.\,427, (10.10.2)]{GD}. Then 
$\R\sHomb(\F\>^\bullet\<,\>\G)=\sHomb(\F\>^\bullet\<,\>\G)$, and the conclusion
follows easily.
\end{proof}

\end{parag}
\smallbreak

\begin{parag}
\label{SS:Dvc-and-Dqc}
\Pref{A(vec-c)-A} applies 
in particular to any noetherian scheme~$X\<$.
When $X$ is separated, $j_{\lower.2ex\hbox{$\scriptstyle\<\<X$}}$ induces an
\emph{equivalence of categories} $\bj_{\!X}\colon\D (\Aqc\<(X)) \cong \Dqc\<(X)$,%
\index{ $\iG{\<\cJ\>}$@$\bj$} 
with quasi-inverse $\R Q_{\<\<X}|_{\Dqc\<(X)}$.  
(See \cite[p.\,133, Corollary 7.19]{H1} for
bounded-below complexes, and \cite[p.\,230, Corollary~5.5]{BN} or
\cite[p.\,12, Proposition~(1.3)]{AJL} for the general case.) We do not know if
such an equivalence, with ``$\vec{\mathrm c}\,$" in place of~``qc," always holds for
separated noetherian formal schemes. The next result will at least take care of the
``properly algebraic" case, see \Cref{corollary}.

\begin{sprop}
\label{c-erator}
In \Pref{(3.2.1)}\textup{,} the functor\/ 
$\kappa^*\colon\D(X)\to \D(\X)$ induces equivalences from\/ 
$\Dqc\<(X)$ to\/ $\Dvc(\X)$ and from\/
$\Dc(X)$ to\/ $\Dc(\X),$\ both with quasi-inverse\/~$\R Q\kappa_*$
$($where\/~$\R Q$ stands 
for\/~$\bj_{\!X}\mspace{-1.5mu}\smcirc\<\R Q_{\<\<X}^{})$.
\end{sprop}

\begin{proof}
Since $\kappa^*$ is exact, \Pref{(3.2.1)} implies that
$\kappa^*(\Dqc\<(X))\subset\Dvc(\X)$ and $\kappa^*(\Dc(X))\subset
\Dc(\X)$.  So it will be enough to show that:

\smallskip
(1) If $\F\in\Dvc(\X)$ then the functorial $\D(\X)$-map $\kappa^*\R
Q\kappa_*\F\to \F$ adjoint to the natural map $\R Q\kappa_*\F\to
\kappa_*\F$ is an isomorphism.

(2) If $\G\in\Dqc\<(X)$ then the natural map 
$\G\iso  \R Q\kappa_*\kappa^*\G$ is an
isomorphism.

(3) If $\F\in\Dc(\X)$ then $\R Q\kappa_*\F\in\Dc(X)$.

\smallskip

Since $\Dvc(\X)$ is triangulated (\Pref{(3.2.2)}), we can use
way-out reasoning \cite[p.~68, Proposition~7.1 and p.~73,
Proposition~7.3]{H1} to reduce to where $\F$ or $\G$ is a single
sheaf.  (For bounded-below complexes we just need the obvious facts
that $\kappa^*$ and the restriction of $\R Q\kappa_*$ to~$\Dvc(\X)$
are both bounded-below ($=$~way-out right) functors.  For unbounded
complexes, we need those functors to be bounded-above as well, which
is clear for the exact functor $\kappa^*\<$, and will be shown for $\R
Q\kappa_*|_{\Dvc(\X)}$ in \Pref{(3.2.7.1)} below.)

Any $\F\in\Avc(\X)$ is isomorphic to $\kappa^*\G$ for some $\G\in\Aqc\<(X)$; and
one checks that the natural composed map
$\kappa^*\G\to  \kappa^*\R Q\kappa_*\kappa^*\G\to\kappa^*\G$
is the identity, whence $(2)\Rightarrow(1)$. 
Moreover,  if $\F\in\Ac(\X\>)$ then $\G\cong
Q\kappa_*\F\in\Ac(X)$, whence $(2)\Rightarrow(3)$.

Now a map $\varphi:\G_1\to\G_2$ in $\Dqc^+(X)$ is an isomorphism iff

\begin{quote}
\hskip-2.25em$(*)\colon\!$ the induced map 
$\Hom_{\D(X)}(\E[-n], \>\G_1)\to\Hom_{\D(X)}(\E[-n], \>\G_2)$ 
is an\newline isomorphism for every $\E\in\Ac(X)$ and every $n\in\mathbb Z$.
\end{quote}
(For, if  $\cV$ is the vertex of a triangle with base~$\varphi$, then $(*)$ says that 
for all~\mbox{$\E$, $n$,} $\Hom_{\D(X)}(\E[-n], \V)=0$; but if $\varphi$ is not an
isomorphism, i.e., $\V$ has non-vanishing homology, say
$H^n(\cV) \neq 0$ and $H^i(\cV)= 0$ for all $i<n$, then the inclusion~into~
$H^n(\cV)$ of any coherent non-zero submodule~
$\E$ gives a non-zero map $\E[-n]\to\V$.) So for~(2) it's enough to check
that the natural composition\looseness=-2
\begin{align*}
\Hom_{\D(X)}(\E[-n], \>\G) 
 &\longrightarrow \Hom_{\D(X)}(\E[-n], \R Q\kappa_*\kappa^*\G)\\ 
&\!\iso  \Hom_{\D(X)}(\E[-n],\kappa_*\kappa^*\G)
 \iso 
  \Hom_{\D(\X)}(\kappa^*\E[-n], \kappa^*\G)
\end{align*}
is the \emph{isomorphism} 
${\rm Ext}_X^n(\E\<,\G)\<\iso \<  {\rm Ext}_{\X}^n(\kappa^*\E\<,\kappa^*\G)$ 
in the following consequence of ~\eqref{3.2.1.1}:

\begin{scor}
\label{(3.2.3)}
With $\kappa \colon\X \to X$ as in \Pref{(3.2.1)} and\/ 
\mbox{$\cL\in \Dqc(X),$} the natural map\/ $\R \Gamma(X,\cL) \to
\R\Gamma(\X,\kappa^* \cL)$ is an isomorphism. In particular, for\/ $\E
\in \Dc^- (X)$ and\/ $\G \in \Dqc^+(X)$ the natural map\/ ${\rm
Ext}_X^n(\E\<,\G) \to {\rm Ext}_{\X}^n(\kappa^*\E\<,\kappa^*\G)$ is an
isomorphism.
\end{scor}
{\it Proof.}  After  ``way-out'' reduction to the case where
$\cL \in \Aqc\<(X)$ 
(the $\R\Gamma$'s are bounded, 
by \Cref{Rf_*bounded}(a) below), the
first assertion is given by~\eqref{3.2.1.1}.
To get the second assertion, take $\cL \!:=\R\sHomb_{\<\<X}(\E\<,\G)$
(which is in $\Dqc^{\raise.2ex\hbox{$\scriptscriptstyle+$}}(X)$,
\cite[p.~92,~Proposition~3.3]{H1}), so that $\kappa^*\cL \cong
\R\sHomb_{\X}(\kappa^*\E\<,\kappa^*\G)$ (as one sees easily after way-out
reduction to where $\E$ and $\G$ are $\cO_{\<\<X}$-modules, and further
reduction to where $X$ is affine, so that $\E$ has a resolution by finite-rank
free modules\dots\!).
\end{proof}

\begin{sdef}
\label{D:propalg}\index{properly algebraic}
A formal scheme
$\X$ is said to be \emph{properly algebraic} if there exist 
a noetherian ring~$B$, a
$B$-ideal~$J\<$, a proper $B$-scheme $X\<$, and an isomorphism
from~$\X$ to the $J$-adic completion of~$X\<$.
\end{sdef}

\begin{scor}
\label{corollary}
On a properly algebraic formal scheme\/~$\X$ the natural functor\/
$\bj_{\!\X}\colon\D(\Avc(\X))\to\Dvc(\X)$\index{ $\iG{\<\cJ\>}$@$\bj$}
is an equivalence of categories 
with quasi-inverse\/~$\R Q_{\X}^{}\>;$\ and therefore\/
$\bj_{\!\X}\mspace{-1.5mu}\smcirc\<\R Q_{\X}^{}$ 
is right-adjoint to the inclusion\/ $\Dvc(\X)\hookrightarrow\D(\X)$.
\end{scor}

\begin{proof}
\smallskip
If $\X$ is properly algebraic, then with $A\!:=J$-adic
completion of~$B$ and $I\!:=JA$, it holds that $\X$ is the
$I\<$-adic completion of~$X\otimes_{B}A$,  and so
we may assume the hypotheses and conclusions of \Pref{(3.2.1)}.
We have also, as above, the equivalence of categories
$\bj_{\!X}\colon\D(\Aqc(X))\to\Dqc(X)$; and so the assertion follows from
 Propositions~\ref{c-erator} and~\ref{A(vec-c)-A}.
\end{proof}

\begin{sprop}\label{P:Lf*-vc}
For a map\/ $g\colon\Z\to\X$  of locally noetherian formal schemes,
$$
\bL g^*\<(\Dvc(\X))\subset\Dqc(\Z).
$$
If\/ $\X$ is properly algebraic, then
$$
\bL g^*\<(\Dvc(\X))\subset\Dvc(\Z).
$$

\end{sprop}

\begin{proof}
The first assertion, being local on~$\X$, follows from the second. Assuming $\X$
properly algebraic we may, as in the proof of \Cref{corollary}, place
ourselves in the situation of  \Pref{(3.2.1)}, so that any $\G\in\Dvc(\X)$
is, by \Cref{corollary} and \Pref{(3.2.1)}, isomorphic to 
$\kappa^*\E$ for some $\E\in\Dqc(X)$. By \cite[p.\,10, Proposition~(1.1)]{AJL}),
$\E$ is isomorphic to a
\smash{$\dirlm{}\!\<$}\vspace{1pt} of bounded-above quasi-coherent flat
complexes (see the very end of the proof of~\emph{ibid.}); and therefore
$\G\cong\kappa^*\E$ is isomorphic to a K-flat complex of $\Avc(\X)$-objects.
Since $\bL g^*$~agrees with~$g^*$ on K-flat complexes, and
$g^*(\Avc(\X))\subset\Avc(\Z)$, we are done.
\end{proof}

\begin{srems}
\label{(3.2.4.1)}

(1) Let $\X$ be a properly algebraic formal scheme (necessarily
noetherian) with ideal of definition~$\I$, and set
 $I\!:={\rm H}^0(\X,\I)\subset A\!:= {\rm H}^0(\X,\cO_\X)$.
Then {\it $A$ is a noetherian $I\<$-adic ring, and $\X$ is 
$\Spf(A)$-isomorphic to the $I\<$-adic completion of a proper
$A$-scheme}. Hence $\X$ is proper over~$\Spf(A)$, via the  canonical map 
given by \cite[p.~407,~(10.4.6)]{GD}.

\penalty -1000

Indeed, with $B$, $J$ and $X$ as in
\Dref{D:propalg},
\cite[p.\,125, Theorem~(4.1.7)]{EGA} implies that the topological ring
$$
A=\inlm{\,n>0} {\rm H}^0(\X, \cO_\X/\I^n\cO_\X) = \inlm{\,n>0} {\rm H}^0(X,
\cO_{\<\<X}/I^n\cO_{\<\<X}) 
$$
is the $J$-adic completion of the noetherian $B$-algebra $A_0:=
H^0(X,\cO_{X})$, and that the $J$-adic and $I\<$-adic topologies on~$A$ are
the same; and then $\X$ is the $I\<$-adic completion of~$X\otimes_{A_0}A$.

(2) It follows that a quasi-compact formal scheme $\X$ 
is properly algebraic iff so~is each of its
connected components.

(3) While (1) provides a less relaxed characterization of 
properly algebraic formal schemes  than \Dref{D:propalg},
\Cref{(3.5.2)} below provides a more relaxed one.
\end{srems}

\begin{slem}
\label{(3.5.1)}
Let\/ $X$ be a locally noetherian scheme, $\cI_1\subset\cI_2$ 
quasi-coherent\/ $\cO_{\<\<X}$-ideals, $Z_i$ the support of\/ $\cO_{\<\<X}/\cI_i\>,$\
and\/ $\X_i$ the completion\/ $X_{\!/Z_i}\ (i=1,2)$. Suppose that\/
$\cI_1\cO_{\X_2}$ is an
ideal of definition of\/ $\X_2$. Then\/ $\X_2$ is a union of connected
components of\/ $\X_1$ $($with the induced formal-subscheme structure$)$.
\end{slem}
\begin{proof}
We need only show that $Z_2$ is open in~$Z_1$.
Locally we have a noetherian ring $A$ and $A$-ideals
$I\subset J$ equal to their
own radicals such that with ${\hat A}$  the $J$-adic completion, 
$J^n{\hat A}\subset I{\hat A}$ for some $n>0$; and we want the
natural map $A/I\twoheadrightarrow A/J$ to be {\it flat}.
(For then with $L\set J/I$, 
$L/L^2={\rm Tor}_1^{A/I}(A/J,A/J)=0$, whence
\mbox{$(1-\ell)L = (0)$} for some $\ell\in L$, whence $\ell=\ell^2$ and
$L=\ell( A/I)$, so that \mbox{$A/I\cong L\times (A/J)$} and
$\Spec(A/J) \hookrightarrow\Spec(A/I)$ is open.)

So it suffices that the localization
$(A/I\>)_{1+J}\to (A/J\>)_{1+J} = A/J$ by the
multiplicatively closed set~$1+J\>$ be an isomorphism, i.e., that
its kernel $J (A/I\>)_{1+J}$ be  nilpotent (hence (0), since $A/I$ is reduced.)  But
this is so because the natural map $A_{1+J}\to {\hat A}$ is faithfully
flat, and therefore $J^n\< A_{1+J}\subset IA_{1+J}$.
\end{proof}

\vspace{-5pt}
\begin{scor}
\label{(3.5.2)}
Let\/ $A$ be a noetherian ring, let\/ $I$ be an\/ $A$-ideal, and let\/ ${\hat A}$ be
the\/
$I\<$-adic completion of\/ $A$. Let\/ $f_0\colon X \to\Spec(A)$ be a separated
finite-type  scheme-map, let\/
$Z$ be a closed subscheme of\/ $f_0^{-1}(\Spec(A/I)),$\ let\/ $\X=X_{\</Z}$ be the
completion of\/ $X$ along\/~$Z,$\ and let\/ $f\colon \X \to\Spf({\hat A})$ be the
formal-scheme map induced by\/ $f_0\>$\textup{:}
$$
\begin{CD}
\X\!:=X_{\</Z} @>>> X \\
\vspace{-22pt}\\
@VfVV @VVf_0 V \\
\vspace{-22pt}\\
\Spf(\hat A) @>>> \Spec(A)
\end{CD}
$$

\smallskip\noindent
If\/ $f$ is proper \textup(see\/ \textup{\S\ref{maptypes})} 
then\/ $\X$ is properly algebraic.
\end{scor}

\begin{proof}
Consider a compactification of $f_0$ (see \cite[Theorem~3.2]{Lu}): 
$$
X 
{\begin{array}[t]{c} 
  \hookrightarrow \\[-2.5 mm]
  \mbox{\tiny open} 
\end{array}}
{\>\>\overline {\<\<X}} 
{\begin{array}[c]{c} 
  \mbox{\scriptsize $\bar f_{\scriptscriptstyle\<0}^{}$}\\[-2 mm]
  \longrightarrow \\[-2.5 mm]
  \mbox{\tiny proper} 
\end{array}} \Spec(A).
$$
Since $f$ is proper, therefore $Z$ is proper over $\Spec(A)$, hence closed in
$\>\>{\overline {\<\<X}}\<$.  Thus we may replace\- $f_0$~by~${\bar
f_0}\>$, i.e., we may assume $f_0$ proper. Since $f\<$, being proper, is adic,
\Lref{(3.5.1)}, with \hbox{$Z_2\set Z$} and
$Z_1\set f_0^{-1}(\Spec(A/I))$, shows that $\X$ is a union of connected
components\- of the properly algebraic formal scheme~$X_{\<\</Z_1}$.
Conclude by \Rref{(3.2.4.1)}(2).
\end{proof}

\end{parag}

\pagebreak[3]

\begin{parag}
\label{SS:bounded} To deal with unbounded complexes we need the
following boundedness results on certain derived functors. (See, e.g., 
Propositions~\ref{P:proper f*} and~\ref{P:kappa-f*} below.)

\begin{sparag}\label{note1} Refer to \S\ref{maptypes} for the definitions of
separated, resp.~affine, maps.

 A formal scheme~$\X$ is \emph{separated}\index{formal scheme!separated} if the
natural map \hbox{$f_\X\colon \X\to\text{Spec}({\mathbb Z})$} is separated, i.e.,
for some---hence any---ideal of definition~$\J$, the scheme
$(\X,\cO_\X/\J)$ is separated. For example, any locally noetherian
affine formal scheme is separated.

A locally noetherian formal scheme $\X$ is affine if and only if the
map~$f_\X$ is affine, i.e., for some---hence any---ideal of
definition~$\J$, the scheme $(\X,\cO_\X/\J)$ is affine. Hence the
intersection $\V\cap\V'$ of any two affine open subsets of a separated
locally noetherian formal scheme~$\Y$ is again affine.  In other
words, the inclusion $\V\hookrightarrow\Y$ is an affine map.
More generally, if $f\colon\X\to \Y$ is a map of locally noetherian
formal schemes, if $\Y$~is separated, and if $\V$ and $\V'$ are affine
open subsets of~$\Y$ and~$\X$ respectively, then $f^{-1}\V\>\cap\V'$ is
affine \cite[p.\,282, (5.8.10)]{GD}.
\end{sparag}

\begin{slem} \label{affine-maps}
If\/ $g\colon \X \to \Y$ is an affine map\/ of locally noetherian
formal schemes, then every\/ $\M\in\Avc(\X)$ is $g_*\<$-acyclic, i.e.,
$R^ig_*\M=0$ for all\/ $i>0$. More generally, if\/ $\G\in\Dvc(\X)$
and\/ $e\in\mathbb Z$ are such that\/ $H^i(\G)=0$ for all\/ $i\ge e$,
then\/ $H^i(\R g_*\G)=0$ for all\/ $i\ge e$.
\end{slem}
\begin{proof}
$R^ig_*\M$ is the sheaf associated to the presheaf $\U \mapsto
{\rm H}^i({g^{-1}(\U),\M})$, ($\U$~open in~$\Y$)
\cite[Chap.\,0, (12.2.1)]{EGA}.  If $\>\U $ is affine then so is
$g^{-1}(\U)\subset \X$, and \Cref{limit(vec-c)=qc} gives ${\rm
H}^i({g^{-1}(\U),\M})=0$ for all $i>0$.

Now consider in $\K(\X)$ a quasi-isomorphism $\G\to I$ where $I$ is a
``special" inverse limit of injective resolutions~$I_{-e}$ of the
truncations $\G^{{\scriptscriptstyle\ge}e}$ (see \eqref{trunc}), so that 
$H^i(\R g_*\G)$ is the sheaf associated to the presheaf $\>\U\mapsto \mathrm
H^i(\Gamma(g^{-1}\U,I))$, see
\cite[p.\,134, 3.13]{Sp}. If $C_{\<-e}$ is the kernel of the split surjection
$I_{-e}\to I_{\>1-e}$ then
$C_{\<-e}[e]$ is an injective resolution of $H^e(\G)\in\Avc(\X)$, and
so for any affine open $\U\subset\Y$ and~$i>e$, 
\mbox{$\mathrm H^i(\Gamma(g^{-1}\U,C_{\<-e}))=0$.} Applying \cite[p.\,126,
Lemma]{Sp}, one finds then that for $i\ge e$ the natural map $\mathrm
H^i(\Gamma(g^{-1}\U,I))\to \mathrm H^i(\Gamma(g^{-1}\U,I_{-e}))$ is an
isomorphism.  Consequently if $H^i(\G)=0$ for all~$i\ge e$ (whence
$I_{-e}\cong \G^{{\scriptscriptstyle\ge}e}=0$ in~$\D(\X)$) then
$\mathrm H^i(\Gamma(g^{-1}\U,I))=0$.
\end{proof}

\begin{sprop} 
\label{Rf_*bounded} Let\/ $\X$ be a noetherian formal scheme.
Then:

\textup{(a)} The functor\/ $\R\Gamma(\X,-)$
is bounded-above on\/~$\Dvc(\X)$. In other words,
there is an integer\/ $e\ge 0$ such that if\/ $\G\in\Dvc(\X)$ and\/
$H^i(\G\>)=0$ for all $i\ge i_0$ then\/ $\mathrm H^i(\R\Gamma(\X,-))=0$ 
for all\/ $i\ge i_0+e$.

\textup{(b)} For any formal-scheme map\/ $f\colon\X\to\Y$
with\/ $\Y$ quasi-compact,
the functor\/~$\Rfs$ is bounded-above on\/~$\Dvc(\X),$\ i.e.,
there is an integer\/ $e\ge 0$ such that if\/~$\G\in\Dvc(\X)$ and\/
$H^i(\G\>)=0$ for all $i\ge i_0$ then\/ $H^i(\Rfs\G\>)=0$ for all\/
$i\ge i_0+e$.
\end{sprop}

\begin{proof} Let us prove (b). (The proof of (a) is the same, \emph{mutatis
mutandis.}) Suppose first that $\X$ is separated, see~\S\ref{note1}.  Since $\Y$
has a finite affine open cover and $\Rfs$ commutes with open base change, we
may assume that $\Y$ itself is affine. Let $n(\X)$ be the least positive
integer~$n$ such that there exists a finite affine open cover
$\X=\cup_{i=1}^n \X_i\>$, and let us show by induction on~$n(\X)$ that
$e\!:=n(\X)-1$ \emph{will do.}

\enlargethispage*{.6\baselineskip}
The case $n(\X)=1$ is covered by \Lref{affine-maps}. So assume
that $n\!:=n(\X)\ge 2$, let $\X=\cup_{i=1}^n \X_i$ be an affine
open cover, and let $u_{1}\colon\<\X_1\hookrightarrow \X$,\,\
$u_2\colon\!\cup_{i=2}^n \X_i \hookrightarrow \X$, \,\
$u_3\colon\!\cup_{i=2}^n (\X_1 \cap \X_i) \hookrightarrow \X$ be the
respective inclusion maps. Note that $\X_1 \cap\> \X_i$ is affine
because $\X$ is separated. So by the inductive\vadjust{\penalty-500}
hypothesis, the assertion holds for the~maps
$f_i\!:=f\smcirc u_i\ (i=1,2,3)$.

\pagebreak[3]

Now apply the $\Delta$-functor~$\Rfs$ to the ``Mayer\kern.5pt-Vietoris" triangle%
\index{Mayer-Vietoris triangle}
$$
\G \longrightarrow \R u_{1*}^{}u_1^*\G \oplus \R
u_{2*}^{}u_{2}^*\G\longrightarrow \R u_{3*}^{}u_{3}^*\G
\stackrel{+1\>}{\longrightarrow}
$$
(derived from the standard exact sequence 
$$
0\to\E\to u_{1*}^{}u_1^*\E
\oplus u_{2*}^{}u_2^*\E\to u_{3*}^{}u_3^*\E\to 0
$$ 
where \hbox{$\G\to\E$} is a K-injective resolution) to get the $\D(\Y)$-triangle
$$
\Rfs\G \longrightarrow \R f_{1*}^{}u_1^*{\G} \oplus \R
f_{2*}^{}u_{2}^*\G\longrightarrow \R f_{3*}^{}u_{3}^*\G
\stackrel{+1\>}{\longrightarrow}
$$
whose associated long exact homology sequence yields the assertion
for~$f$.

The general case can now be disposed of with a similar Mayer\kern.5pt-Vietoris
induction on the least number of \emph{separated} open subsets needed
to cover~$\X$.
\end{proof}

\begin{sprop}
\label{(3.2.7.1)}
Let\/ $X$ be a separated noetherian scheme, let\/ $Z\subset X$ be a closed
subscheme, and let\/ $\kappa_\X^\pd\colon\X=X_{/Z} \to X$ 
be the completion map.
Then the~functor\/~$\R Q_{\<\<X} \kappa_*$ is bounded-above on\/ $\Dvc(\X)$.
\looseness =1
\end{sprop}
\begin{proof}
Set $\kappa\set\kappa_\X^{}$. Let $n(X)$ be the least number of affine
open subschemes needed to cover~$X\<$. When $X$ is affine, $Q_{\<\<X}$
is the sheafification of the global section functor, and since
$\kappa_*$ is exact and, being right adjoint to the \emph{exact}
functor~$\kappa^*\<$, preserves K-injectivity, we find that for any
$\F\in \D(\X)$, $\R Q_{\<\<X}\kappa_*\F$ is the sheafification of the
complex $\R \Gamma(X, \kappa_*\F\>)=\R \Gamma(\X\<, \F\>)$. Thus
\Pref{Rf_*bounded}(a) yields the desired result for
$n(X)=1$.

Proceed by induction when $n(X)>1$, using a ``Mayer\kern.5pt-Vietoris''
argument as in the  proof of \Pref{Rf_*bounded}.  The enabling points
are that if
$v\colon V\hookrightarrow X$ is an open immersion with $n(V)<n(X)$, giving
rise to the natural commutative diagram
$$
\CD
V_{\<\</Z\cap V}=:\:@.\V @>\kappa_\V^\pd>> V \\
@.@V\hat v VV @VVvV \\
@.\X@>>\vbox to
0pt{\vskip-1.3ex\hbox{$\scriptstyle\kappa_\X^\pd$}\vss}> X
\endCD
$$
then there are natural isomorphisms, for $\F\in\Dvc(\X)$ and 
$v_*^{\rm qc}\colon\Aqc\<(V)\to \Aqc\<(X)$ the restriction of $v_*\>$:
$$
\R Q_{\<\<X} \kappa_{\X*}^\pd\R {\hat v}_*{\hat v}^*\F \cong \R Q_{\<\<X} \R
v_*\kappa_{\V*}^\pd{\hat v}^*\F
\cong \R v_*^{\rm qc} \R Q_V \kappa_{\V*}^\pd{\hat v}^*\F,
$$ 
and the functor $\R Q_V\kappa_{\V*}^\pd{\hat v}^*$ is bounded-above, by the
inductive hypothesis on $n(V)<n(X)$, as is $\R v_*^{\rm qc}$, by the proof of
\cite[p.\,12, Proposition (1.3)]{AJL}.
\end{proof}
\end{parag}

\medskip
\pagebreak[3]
\begin{parag}\label{3.5}
Here are some examples of how boundedness is used.
\begin{sprop}\label{P:proper f*}
Let\/ $f\colon\X\to\Y$ be a proper map of noetherian formal schemes. Then
$$
\R f_{\!*}\Dc(\X)\subset\Dc(\Y)\quad\textup{and}\quad
\R f_{\!*}\Dvc(\X)\subset\Dvc(\Y).
$$
\end{sprop}

\begin{proof}
For a coherent\vspace{.4pt} $\cO_\X$-module $\M$, 
$\R f_{\!*}\M\in\Dc(\Y)$ \cite[p.\,119, (3.4.2)]{EGA}. Since $\X$ is noetherian,
the homology functors $H^i\R f_{\!*}$ commute with\vspace{.9pt}
\smash{$\dirlm{}\!\!$} on $\cO_\X$-modules, 
whence\vspace{.6pt} $\R f_{\!*}\N\in\Dvc(\Y)$ for all $\N\in\Avc(\X)$. $\R
f_{\!*}$ being bounded on~$\Dvc(\X)$ (\Pref{Rf_*bounded}(b)),  
way-out reasoning \cite[p.\,74, (iii)]{H1} completes the proof.\looseness=-1
\end{proof}

\begin{sprop}\label{P:coprod}
Let\/ $f\colon\X\to\Y$ be a map of quasi-compact formal schemes, with\/
$\X$ noetherian. Then the functor\/ $\Rfs|_{\Dvc(\X)}$ commutes with small 
direct sums, i.e., for any small family\/ $(\E_\alpha)$ in\/~$\Dvc(\X)$
the natural map 
$$
\oplus_{\<\alpha}( \Rfs\E_\alpha) \to\Rfs(\oplus_{\<\alpha}\>\E_\alpha)
$$
is a\/ $\D(\Y)$-isomorphism.
\end{sprop}

\begin{proof}
It suffices to look at the induced homology maps in each degree, i.e., setting
$R^{\>i}\<\<f_{\!*}\set H^i\Rfs\ (i\in\mathbb Z)$, we need to show that \emph{the
natural map}
$$
\oplus_{\<\alpha}( R^{\>i}\<\<f_{\!*}\E_\alpha)\iso
R^{\>i}\<\<f_{\!*}(\oplus_{\<\alpha}\>\E_\alpha).
$$
\emph{is an isomorphism.} 

For any $\F\in\Dvc(\X)$ and any integer $e\ge0$, the
vertex~$\G$ of a triangle based on the natural map~$t_{i-e}$ from $\F\>$ to the
truncation $\F^{\>{\scriptscriptstyle\ge}i-e}$ (see~\eqref{trunc}) satisfies 
$H^j(\G)=0$ for all $j\ge i-e-1$; so if $e$ is the integer
in~\Pref{Rf_*bounded}(b), then
$R^{\>i-1}\<\<f_{\!*}\G=R^{\>i}\<f_{\!*}\G=0$, and the map induced by
$t_{i-e}$  is an \emph{isomorphism}\looseness=-1
$$
 R^{\>i}\<\<f_{\!*}\F\iso
R^{\>i}\<\<f_{\!*}\F^{\>{\scriptscriptstyle\ge}i-e}\<.
$$
We can therefore replace each $\E_\alpha$ by
$\E_\alpha^{{\scriptscriptstyle\ge}i-e}$, i.e., we may assume that the
$\E_\alpha$ are uniformly bounded below.

We may  assume further that each complex~$\E_\alpha$ is injective, hence 
$f_{\!*}$-acyclic (i.e., the canonical map is an \emph{isomorphism} 
$f_{\!*}\E_\alpha\iso \Rfs\E_\alpha$). Since $\X$ is noetherian,
$R^{\>i}\<\<f_{\!*}$ commutes with direct
sums; and so each component of
$\oplus_{\<\alpha}\>\E_\alpha$ is an $f_{\!*}$-acyclic $\cO_\X$-module. This
implies  that the bounded-below complex $\oplus_{\<\alpha}\>\E_\alpha$ is itself 
$f_{\!*}$-acyclic. Thus in the natural commutative diagram
$$
\begin{CD}
 \oplus_{\<\alpha}(
f_{\!*}\E_\alpha)@>\Iso>>f_{\!*}(\oplus_{\<\alpha}\>\E_\alpha)\\ 
@V\simeq VV @VV\simeq V\\
\oplus_{\<\alpha}( \Rfs\E_\alpha)@>>>\Rfs(\oplus_{\<\alpha}\>\E_\alpha)
\end{CD}
$$
the top and both sides are isomorphisms, whence so is the bottom.
\end{proof}

The following Proposition generalizes \cite[p.\,92, Theorem (4.1.5)]{EGA}.
 
\begin{sprop}\label{P:kappa-f*}
Let\/ $f_0\colon\! X\to Y$ be a proper map of locally noetherian schemes, 
let\/ $W\subset Y$ be a closed subset, let\/ $Z\!:=f_0^{-1}W,$\ let\/
$\kappa_\Y^\pd\colon\Y=Y_{/W}\to Y$ and\/~ 
$\kappa_\X^\pd\colon\X=X_{/Z}\to X$  be the respective \(flat\) completion maps,
and let\/ $f\colon\X\to\Y$ be the map induced by~$f_0\>$. Then for
$\E\in\Dqc(X)$ the map\/ $\theta_{\<\E}$ adjoint to the natural composition
$$
\R f_{\<0*}^{}\E\longrightarrow\R f_{\<0*}^{}\kappa_{\X*}^\pd\kappa_\X^*\E\iso
\kappa_{\Y*}^\pd\R f_{\!*}\kappa_\X^*\E
$$
is an \emph{isomorphism}
$$
\theta_{\<\E}\colon \kappa_\Y^*\R f_{\<0*}^{}\E\iso \R f_{\!*}\kappa_\X^*\E\<.
$$
\end{sprop}

\begin{proof}
We may assume $Y$ affine, say
$Y=\Spec(A)$, and then $W=\Spec(A/I)$ for some $A$-ideal~$I\<$. Let $\hat A$ be
the $I\<$-adic completion of~$A$, so that there is a natural cartesian diagram
$$
\begin{CD}
X\otimes_A\hat A=:\>@. X_1 @>k_{\<X}^{}>> X \\
@.@V f_1 VV @VV f_0 V  \\
\Spec(\hat A)=:\>@. Y_1 @>>\vbox to
0pt{\vskip-1.3ex\hbox{$\scriptstyle k_Y^{}$}\vss}> Y
\end{CD}
$$
Here $k_Y^{}$ is flat, and the natural map is an isomorphism
$k_Y^*\R f_{\<0*}^{}\E\iso \R {f_{\<1*}^{}}k_{\<X}^*\E\colon$
since $\R  f_{\<0*}^{}$ (resp.~$\R {f_{\<1*}^{}})$ is bounded-above on~$\Dqc(X)$
(resp.~$\Dqc(X_1)$), see \Pref{Rf_*bounded}(b), way-out reasoning
reduces this assertion to the well-known case where $\E$ is a single
quasi-coherent
$\cO_X$-module. Simple considerations show then that we can replace
$f_0$ by $f_1$ and $\E$ by $k_{\<X}^*\E$; in other words, we can
assume \mbox{$A=\hat A$}.

 From \Pref{P:proper f*} it follows that 
$\R f_{\<0*}^{}\E\in\Dqc(Y)$ and $\R f_{\!*}\kappa_\X^*\E\in \Dvc(\Y)$. 
Recalling the equivalences in~\Pref{c-erator}, we see that any
$\F\in\Dvc(\Y)$ is isomorphic to~$\kappa_\Y^*\F_0$ for some
$\F_0\in\Dqc(Y)$ (so that $\bL f_{\<0}^*\F_0\in\Dqc(X)$),
and that there is a sequence of natural
isomorphisms
\begin{align*}
 \Hom_\Y(\F\<, \>\kappa_\Y^*\R f_{\<0*}^{}\E)
&\iso 
 \Hom_Y(\F_0\>, \>\R f_{\<0*}^{}\E) \\
&\iso
 \Hom_X(\bL f_{\<0}^*\F_0\>, \>\E) \\
&\iso
 \Hom_\X(\kappa_\X^*\bL f_{\<0}^*\F_0\>, \>\kappa_\X^*\E) \\
&\iso
 \Hom_\X(\bL f^*\<\kappa_\Y^*\F_0\>, \>\kappa_\X^*\E) 
\iso
\Hom_\Y(\F\<, \>\R f_{\!*}\kappa_\X^*\E).
\end{align*}
The conclusion follows. 
\end{proof}

\end{parag}

\section{Global Grothendieck Duality.}
\label{sec-th-duality}
\index{Grothendieck Duality!global}

\begin{thm}
\label{prop-duality}
Let\/ $f\colon\<\X \to \Y$ be a map of quasi-compact formal schemes, with 
$\X$~noetherian, and
let\/~$\bj\colon\<\D(\Avc(\X))\index{ $\iG{\<\cJ\>}$@$\bj$}
\to\D(\X)$ be the natural functor.  Then the\/
\hbox{$\Delta$-functor\/}\vadjust{\kern.4pt} 
$\Rfs \<\smcirc\>\bj\>$  has a right\/ $\Delta$-adjoint.  
In fact there is a bounded-below\/ \hbox{$\Delta$-functor\/}
$f^{\times}\<\colon\D(\Y)
\to\D\left({\Avc(\X)}\right)\index{ $\iG$@$f^{{}^{\>\ldots}}$ (right adjoint of
$\R f_{\<\<*}\cdots$)!$f^\times\<\<$}$\vadjust{\kern.3pt}
and a map of\/ $\Delta$-functors 
$\tau\colon\Rfs \>\bj f^{\times}\to {\bf 1}$\index{ {}$\tau$ (trace map)}
such that for all\/
$\G\in\D(\Avc(\X))$ and\/ $\F\in\D(\Y),$\ the composed map 
$($in\/ the derived category of abelian groups\/$)$
\begin{align*}
\R\Homb_{\Avc(\X)\<}(\G,\>f^\times\<\<\F\>)
&\xrightarrow{\mathrm{natural}\,}
   \R\Homb_{\A(\Y)\<}(\Rfs\>\bj\>\G,
     \>\Rfs \>\bj f^{\times}\<\F\>)\\ 
&\xrightarrow{\;\>\mathrm{via}\ \tau\ }
  \R\Homb_{\A(\Y)\<}(\Rfs \>\bj\>\G,\>\F\>)
\end{align*}
is an isomorphism.
\end{thm}

With \Cref{corollary} this gives: 

\begin{scor}\label{cor-prop-duality}
If\/ $\X$ is properly algebraic, the restriction of\/~$\Rfs$ 
to\/~$\Dvc(\X)$ has a right\/ $\Delta$-adjoint \textup(also to be denoted
$f^\times$ when no confusion results\/\textup).
\end{scor}

\noindent\emph{Remarks.} 1.~Recall that over any abelian category~$\A$ in which
each complex~$\F$ has a K-injective resolution~$\rho(\F\>)$, we can
set
$$
\R\Homb_\A(\G,\F\>)\!:=\Homb_\A\bigl(\G,\rho(\F\>)\bigr) \qquad
\bigl(\G,\F\in\D(\A)\bigr);\index{ $\R$ (right-derived functor)!$\R\Homb$}
$$
and there are natural isomorphisms
$$
\mathrm H^i\R\Homb_\A(\G,\F\>)\cong \mathrm{Hom}_{\D(\A)}\bigl(\G,\F[i]\bigr)
\qquad(i\in\mathbb Z).
$$

2.~Application of homology to the second assertion in the Theorem
reveals that it~is equivalent to the first one.

3.~We do not know in general (when $\X$ is not properly algebraic)
that the functor~$\bj$ is fully faithful---$\bj$\index{ $\iG{\<\cJ\>}$@$\bj$} 
has a right adjoint
$(\text{identity})^{\<\times}\cong\R Q_\X^{}$ (see
\Pref{A(vec-c)-A}), but it may be that for some 
$\E\in\Avc(\X)$ the natural map 
$\>\E\to\R Q_\X^{}\>\bj\>\E$ is  not an isomorphism. 

4. For a  \emph{proper} map $f_0\colon X\to Y$ of \emph{ordinary} schemes it is
customary to write~$f_0^!$\index{ $\iG$@$f^{{}^{\>\ldots}}$ (right adjoint of
$\R f_{\<\<*}\cdots$)!$\mathstrut \fs\<$} instead
of~$f_0^\times\<$. (Our extension of this notation to maps of formal
schemes---introduced immediately after
\Dref{D:basechange}---is not what would be expected here.) 

5.~\Tref{prop-duality} includes the case when $\X$ and~$\Y$ are ordinary
noetherian schemes.
(In fact the proof below applies with minor changes to arbitrary maps of
quasi-compact, quasi-separated schemes, 
cf.~\cite[Chapter~4]{Derived categories}.)
The next Corollary relates the formal situation to the
ordinary one.

\begin{scor}\label{C:kappa-f^times}
Let\/ $A$ be a noetherian adic ring with ideal of definition\/~$I,$\ 
set\/ $Y\!:=\Spec(A)$ and $W\!:=\Spec(A/I)\subset Y\<$. 
Let\/ $f_0\colon X\to Y$ be a proper map and set~
 $Z\!:=f_0^{-1}W,$ so that there is a commutative diagram
$$
\begin{CD}
\X\!:=@.X_{/Z} @>\kappa_\X^\pd>> X  \\
@.@V f VV @VV f_0 V \\
\Y\!:=\:@.\Spf(A)@>>\vbox to
0pt{\vskip-1.3ex\hbox{$\scriptstyle\kappa_\Y^\pd$}\vss}> Y
\end{CD}
$$
with\/ $\kappa_\X^\pd$ and\/ $\kappa_\Y^\pd$  the respective
\textup(flat\textup) completion maps, and
$f$ the \textup(proper\textup) map induced by\/~$f_0\>$. 

\pagebreak[3]
Then the map
adjoint to the natural composition
$$
\Rfs\kappa_\X^*f_{\<0}^!\kappa_{\Y*}^\pd
\xrightarrow{\ref{P:kappa-f*}\>}
\kappa_\Y^*\R f_{0*} f_{\<0}^!\kappa_{\Y*}^\pd
\longrightarrow 
\kappa_\Y^*\kappa_{\Y*}^\pd\longrightarrow\mathbf 1
$$
is an isomorphism of
functors---from\/ $\D(\Y)$ to\/ $\Dvc(\X),$\ see
\Cref{cor-prop-duality}---
$$
\kappa_\X^*f_{\<0}^!\kappa_{\Y*}^\pd\iso f^\times\<.
$$

\end{scor}

\begin{proof}
For any $\E\in\Dvc(\X)$ set $\E_0\!:=\bj_{\!X}\R
Q_{\<\<X}^{}\kappa_{\X*}^\pd\E\in\Dqc(X)$ (see \Sref{SS:Dvc-and-Dqc}). Using
\Pref{c-erator} we have then for any $\F\in\D(\Y)$ the 
natural isomorphisms
\begin{align*}
 \Hom_{\D(\X)}(\E,\>\kappa_\X^*f_{\<0}^!\kappa_{\Y*}^\pd\F\>)
&\iso
 \Hom_{\D(X)}(\E_0\>,\>f_{\<0}^!\kappa_{\Y*}^\pd\F\>) \\
&\iso
 \Hom_{\D(Y)}(\R f_{\<0*}^{}\E_0\>,\>\kappa_{\Y*}^\pd\F\>) \\
&\iso
 \Hom_{\D(\Y)}(\kappa_\Y^*\R f_{\<0*}^{}\E_0\>,\>\F\>) \\
&\underset{\ref{P:kappa-f*}}{\iso}
 \Hom_{\D(\Y)}(\R f_{\!*}\kappa_\X^*\E_0\>,\>\F\>) 
\iso
 \Hom_{\D(\Y)}(\R f_{\!*}\E,\>\F\>). 
\end{align*}
Thus $\kappa_\X^*f_{\<0}^!\kappa_{\Y*}^\pd$ is right-adjoint 
to~$\Rfs|_{\Dvc(\X)}\>$, whence the conclusion.
\end{proof}

\smallskip
\begin{proof}[Proof of \Tref{prop-duality}.]
1. Following Deligne\index{Deligne, Pierre} \cite[p.\,417, top]{H1}, we begin by
considering  for
$\M\in\A(\X)$ the functorial flasque \emph{Godement resolution}
$$
0\to\M\to G^0(\M)\to G^1(\M)\to\cdots\,.
$$
Here, with $G^{-2}(\M)\!:=0$, $G^{-1}(\M)\!:=\M$, and for 
$i\ge 0$,  
$K^i(\M)$ the cokernel of
$G^{i-2}(\M)\to G^{i-1}(\M)$, the sheaf
$G^i(\M)$ is specified inductively by
$$
G^i(\M)\bigl(\U\bigr)\!:=\prod_{x\in\U}\,K^i(\M)_x
 \qquad(\U\text{ open in }\X).
$$
One shows by induction on~$i$ that all the functors $G^i$
and $K^i$ (from $\A(\X)$ to itself) are \emph{exact.}
Moreover, for $i\ge 0$, $G^i(\M)$, being flasque, is
$f_{\!*}$-\emph{acyclic,} i.e.,
$$
R^j\<\<f_{\!*}G^i(\M)=0\quad\text{for all }j>0.
$$

 The category $\Avc(\X)$ has small colimits
(\Pref{(3.2.2)}), and is generated~by its coherent
members, of which there exists a small set containing
representatives of every isomorphism class.  The Special
Adjoint Functor Theorem\index{Special Adjoint Functor Theorem}
\mbox{(\cite[p.\,90]{pF1964}} or~\cite[p.\,126, Corollary]{currante}) guarantees
then that a right-exact functor $F$ from~$\Avc$ into~an  abelian category~$\A'$
has a right adjoint iff $F$ is  
\emph{continuous}\index{continuous functor} in the sense that it commutes with
filtered direct limits, i.e., for any small directed system 
\mbox{$(\M_\alpha\>,\,\varphi_{\alpha\beta}\colon \M_\beta\to 
 \M_\alpha)$} in~$\Avc\>$, with 
$\dirlm{\alpha}\M_\alpha=(\M,\,\varphi_\alpha\colon \M_\alpha\to \M)$
it holds that
$$
\bigl(F(\M), F(\varphi_\alpha)\bigr)=\dirlm{\alpha}\bigl(F(\M_\alpha),
F(\varphi_{\alpha\beta})\bigr).
$$
Accordingly, for constructing right adjoints we need to replace the 
restrictions of~$G^i$ and~$K^i$ to $\Avc(\X)$ by continuous functors.

\begin{slem}\label{L:vc-functor}
Let\/ $\X$ be a locally noetherian formal scheme and let\/
$G$ be a functor from $\Ac(\X)$ to a category\/~$\A'$ in which direct
limits exist for all small directed systems.  
Let\/ $j\colon\Ac(\X)\hookrightarrow\Avc(\X)$ be the inclusion functor. Then:

\smallskip
\textup{(a)} There exists a continuous functor\/~$G_\vc\>\colon\Avc(\X)\to\A'$ 
and an isomorphism of functors\/ 
$\varepsilon\colon G\iso G_\vc\smcirc\< j$ such that 
for any map of functors\/ $\psi\colon G\to F\<\smcirc\< j$ with\/ $F$~
continuous, there is a unique map of functors\/ $\psi_\vc\>\colon G_\vc\to F$
such that\/ $\psi$ factors as
$$
G \xrightarrow{\,\varepsilon\,} 
G_\vc\smcirc\< j \xrightarrow{\textup{via }\psi_\vc\,} F\<\smcirc\< j\>.
$$

\smallskip
\textup{(b)} Assume that $\A'$ is abelian, and has exact filtered
direct limits $($i.e., satisfies Grothendieck's axiom\/
\textup{AB5)}. Then if\/ $G$ is exact, so is $G_\vc\>$.
\end{slem}

\begin{proof}
(a) For $\M\in\Avc(\X)$, let $(\M_\alpha)$ be the directed system of coherent
$\cO_\X$-submodules of\/~$\M,$\ and set 
$$
G_\vc(\M)\set\smash{\dirlm{\alpha}}G(\M_\alpha).
$$

\smallskip\noindent
For any $\Avc(\X)$-map $\nu\colon\M\to\N$ and any~$\alpha$,
there exists a coherent
submodule~$\N_\beta\subset\N$ such that $\nu|_{\M_\alpha}$ factors as
$\M_\alpha\to\N_\beta\hookrightarrow\N$
(\Cref{C:limsub} and \Lref{L:Ext+lim}, with
$q=0$); and the resulting composition 
$$
\nu_\alpha'\colon G(\M_\alpha)\to G(\N_\beta)\to G_\vc(\N\>)
$$ 
does not depend on the choice of~$\N_\beta\>$. We define the map
$$
G_\vc(\nu)\colon G_\vc(\M)=\smash{\dirlm{\alpha}G(\M_\alpha)}\to G_\vc(\N\>)
$$

\smallskip\noindent
to be the unique one whose composition with $G(\M_\alpha)\to
G_\vc(\M)$ is $\nu_\alpha'$ for all~$\alpha$. Verification of
the rest of assertion (a) is straightforward.

(b) Let $0\to\M\to\N\xrightarrow{\,\pi\,}\mathcal{Q}\to 0$ be an exact
sequence in $\Avc(\X)$. Let $(\N_\beta)$ be the filtered system of
coherent submodules of~$\N\<$, so that $\N=\smash{\dirlm{}\<\N_\beta}$
(\Cref{C:limsub}). Then
$(\M\cap\N_\beta)$ is a filtered system of coherent $\cO_\X$-modules whose
$\smash{\dirlm{}}$ is~$\M$, and
$(\pi\N_\beta)$ is a filtered system of coherent $\cO_\X$-modules whose
$\smash{\dirlm{}}\vspace{.7pt}$ is~$\mathcal Q$ (see \Cref{C:images}). The
exactness of~$G_\vc$ is then made apparent by application
of~\smash{$\dirlm{}_{\<\<\!\beta}$} to the system of exact sequences
$$
0\to G(\M\cap\N_\beta)\to G(\N_\beta)\to G(\pi\N_\beta)\to 0.
$$
\vskip-3.8ex
\end{proof}

\smallskip
Now for  $\M \in\Avc(\X)$,  
the $\smash{\dirlm{}}$\vspace{1pt} of the system of Godement resolutions of
all the coherent submodules~$\M_\alpha\subset\M$ is a functorial resolution 
$$
0\to\M\to G_\vc^0(\M)\to G_\vc^1(\M)\to\cdots;
$$
and the cokernel of $G_\vc^{i-2}(\M)\to G_\vc^{i-1}(\M)$ is\vspace{1pt}
$K_\vc^i(\M)\set\smash{\dirlm{}\!}K^i(\M_\alpha)$.
By~(b) above (applied to the exact functors $G^i$ and $K^i$),\vspace{.8pt} the
continuous functors $G_\vc^i$ and $K_\vc^i$ are exact;  and
$G_\vc^i(\M)=\smash{\dirlm{}}\>G^i(\M_\alpha)$\vspace{1.5pt} is
$f_{\!*}$-acyclic since
$G^i(\M_\alpha)$ is, and---$\X$ being noetherian---the functors
$R^j\<f_{\!*}$\vspace{1pt} commute with~$\smash{\dirlm{}\!}$.\vspace{1.2pt}  
\Pref{Rf_*bounded}(b) implies\- then that there is an
integer~$e\ge 0$ such that for all $\M\in\Avc(\X)$,
$K_\vc^e(\M\>)$ is $f_{\!*}$-acyclic. \vspace{1pt}

So if we define the exact
functors~$\cD^i\colon \Avc(\X)\to\A(\X)$ by
$$
\cD^i(\M\>)\!=
\begin{cases}
G_\vc^i(\M\>)\qquad &(0\le i< e) \\ 
K_\vc^e(\M\>) \qquad&(i=e) \\ 
0\qquad&(i>e)
\end{cases}
$$
then for $\M\in\Avc(\X)$, each $\cD^i(\M)$ is $f_{\!*}$-acyclic 
and the natural sequence
$$
0\longrightarrow\M\xrightarrow{\delta(\M)\>}\cD^0(\M\>)
\xrightarrow{\delta^0(\M)\>}\cD^1(\M\>)
\xrightarrow{\delta^1(\M)\>}\cD^2(\M\>)
\longrightarrow\cdots\longrightarrow\cD^e(\M\>) 
\longrightarrow 0
$$ 
is exact. In short, the sequence
$\cD^0\to\cD^1\to\cD^2\to\cdots\to\cD^e\to 0$ 
is an \emph{exact, continuous,
$f_{\!*}\<$-acyclic, finite resolution of the inclusion functor
$\Avc(\X)\hookrightarrow\A(\X)$.}
\smallskip

\pagebreak[3]

2. We have then a $\Delta$-functor
$({\cD}^{\bullet}\<,\text{Id})\colon\K(\Avc(\X) )\to \K(\X)$ which assigns
an \mbox{$f_{\!*}$-acyclic} resolution to each $\Avc(\X)$-complex~$\G
=(\G^{\>p})_{p\in\lower.1ex\hbox{$\scriptstyle\mathbb Z$}}\,$:
$$
(\cD^\bullet\<\G\>)^m\!:=\bigoplus_{p+q=m}\cD^q(\G^{\>p})\qquad
(m\in\mathbb Z,\ 0\le q\le e),
$$
the differential $(\cD^\bullet \G\>)^m\to(\cD^\bullet \G\>)^{m+1}$
being defined on $\cD^q(\G^{\>p})$ $(p+q=m)$ to be $d'+(-1)^pd''$ where
$d'\colon\cD^q(\G^{\>p})\to\cD^q(\G^{\>p+1})$ comes from the
differential in $\G$ and
$d''=\delta^q(\G^{\>p})\colon\cD^q(\G^{\>p})\to\cD^{q+1}(\G^{\>p})$.

It is elementary to check that the natural map
$\delta(\G\>)\colon\G\to {\cD}^{\bullet}\G$ is a
\emph{quasi-isomorphism}. The canonical maps are 
\emph{$\D(\Y)$-isomorphisms}
\begin{equation}\label{f*D}
f_{\!*}{\cD}^{\bullet}(\G) \iso \Rfs{\cD}^{\bullet}(\G) \;\underset{\Rfs\delta(\G)}{\osi}\;\Rfs\G,
\end{equation}
i.e., the natural map $\alpha^i\colon
H^i\bigl(f_{\!*}\cD^\bullet(\G)\bigr)\to H^i\bigl(\Rfs\cD^\bullet(\G)\bigr)$ is an isomorphism for all~$i\in\mathbb Z\>\>$:
this holds for bounded-below $\G$ because ${\cD}^{\bullet}(\G)$ is a
complex of $f_{\!*}$-acyclic objects; and for arbitrary~$\G$ since for
any $n\in\mathbb Z$, with $\G^{{\scriptscriptstyle\ge}n}$ denoting the
truncation\looseness=-1
\stepcounter{numb}
\begin{equation}\label{trunc}
\cdots \to 0\to 0 \to
\textup{coker}(\G^{n-1}\to\G^n)\to\G^{n+1}\to\G^{n+2}\to\cdots
\end{equation}
there is a natural commutative diagram
$$
\CD 
H^i\bigl(f_{\!*}\cD^\bullet(\G)\bigr)@>\alpha^i>>H^i\bigl(\Rfs\cD^\bullet(\G)\bigr) \\ @V\beta_n^i VV @VV\gamma_n^i V \\
H^i\bigl(f_{\!*}\cD^\bullet(\G^{{\scriptscriptstyle\ge}n})\bigr)@>>\alpha_n^i>
H^i\bigl(\Rfs\cD^\bullet(\G^{{\scriptscriptstyle\ge}n})\bigr)
\endCD
$$
in which, when $n\ll i$, $\beta_n^i$ is an isomorphism (since $\G$ and
$\G^{{\scriptscriptstyle\ge}n}$ are identical in all degrees
$>n$), $\gamma_n^i$~is an isomorphism (by \Pref{Rf_*bounded}(b)
applied to the mapping cone of the natural composition
$\cD^\bullet(\G)\iso\G\longrightarrow\G^{{\scriptscriptstyle\ge}n}
\iso\cD^\bullet(\G^{{\scriptscriptstyle\ge}n})$), 
and $\alpha_n^i$ is an isomorphism (since
$\G^{{\scriptscriptstyle\ge}n}$ is bounded below).

Thus we have realized $\Rfs\smcirc\>\>\bj$ at the homotopy
level, via the functor $\C^\bullet\!:=f_{\!*}\cD^\bullet\>$; and our task is
now to find a right adjoint at this level.
\smallskip

3. Each functor ${\C}^p=f_{\!*}{\cD}^p\colon\Avc(\X)\to\A(\Y)$ is
exact,\vspace{.3pt} since $R^1\<\<f_{\!*}(\cD^p(\M\>))=0$ for all
$\M\in\Avc(\X)$. ${\C}^p$ is continuous, since $\cD^p$ is\vspace{.6pt}
and, $\X$ being noetherian,
$f_{\!*}$ commutes with $\smash{\dirlm{}}\!$.  As before, the Special
Adjoint Functor Theorem\index{Special Adjoint Functor Theorem}\vspace{1.2 pt}
yields that \emph{$\C^p$~has a right adjoint
$\C_p\colon\A(\Y) \to\Avc(\X)$.}\vadjust{\kern.4 pt}

\penalty-1000
For each $\A(\Y)$-complex
$\F=(\F^p)_{p\in\lower.1ex\hbox{$\scriptstyle\mathbb Z$}}\>$ let
$\C_\bullet \>\F$ be the $\Avc(\X)$-complex with
$$
(\C_\bullet \>\F\>)^m\!:=\prod_{p-q\,=\,m}\C_q\F^p\qquad (m\in\mathbb Z,
0\le q\le e),
$$
and with differential $(\C_\bullet \>\F\>)^m\to(\C_\bullet \>\F\>)^{m+1}$
the unique map making the following diagram commute for all $r, s$
with $r\mspace{-1.5mu}-\mspace{-1.5mu}s\>=\>m\mspace{-1.5mu}+\!1\>$:
$$
\CD 
\underset{p-q\,=\,m}{\prod} \C_q\F\>^p 
 @>\phantom{d_\prime\>+\>(-1)^rd_{\prime\prime}}>> 
   \underset{p-q\,=\,m+1}{\prod} \C_q\F\>^p\\ 
@VVV @VVV\\ 
\C_s\F\>^{r-1}\oplus\C_{s+1}\F\>^r 
 @>>d_\prime\>+\>(-1)^rd_{\prime\prime}>
   \C_s\F\>^r 
\endCD
$$
where: 

(i) the vertical arrows come from projections, 

(ii) $d_\prime\colon\C_s\F\>^{r-1}\to\C_s\F\>^r$ corresponds to the
differential in $\F$, and 

 (iii) with $\delta_s\colon\C_{s+1}\to \C_s$
corresponding by adjunction to $f_{\!*}(\delta^s)\colon \C^s\to\C^{s+1}$,
$$
d_{\prime\prime}\!:=(-1)^s\delta_s(\F\>^r)\colon\C_{s+1}\F\>^r\to\C_s\F\>^r.
$$
This construction leads naturally to a $\Delta$-functor
$({\C}_{\bullet}\>,\text{Id})\colon\K(\Y)\to \K(\Avc(\X))$.
The adjunction isomorphism
$$
\Hom_{\Avc(\X)}(\M\<,{\C}_p\>\N\>) \iso \Hom_{\A(\Y)}(\C^p\<\M\<,\>\N\>)
\qquad \bigl(\M\in \Avc(\X),\ \N\in \A(\Y)\bigr)
$$ 
applied componentwise produces an isomorphism of complexes of
abelian groups
\stepcounter{sth}
\begin{equation}\label{Deligne}
\Homb_{\Avc(\X)}(\G,\>\C_{\bullet}\F\>) \iso
\Homb_{\A(\Y)}({\C}^{\bullet}\G,\>\F\>)
\end{equation}
for all $\Avc(\X)$-complexes $\G$ and $\A(\Y)$-complexes~$\F\<$.

\smallskip
4. The isomorphism \eqref{Deligne} suggests that we use $\>\C_\bullet$
to construct $f^\times\<$, as follows.  Recall~that a complex
$\J\in\K(\Avc(\X))$ is K-injective iff for each exact complex
\hbox{$\G\in\K(\Avc(\X))$}, the complex $\Homb_{\Avc(\X)}(\G,\J)$ is exact
too.  By~\eqref{f*D},
$\C^\bullet\G$ is exact if $\G$ is; so it follows from~\eqref{Deligne}
that \emph{if $\F$ is K-injective in~$\K(\Y)\mspace{-.6mu}$ then
$\C_\bullet\F$~is
\hbox{K-injective} in~$\K(\Avc(\X))$.} Thus if 
$\K_{\text{\textbf I}}(-)\subset\K(-)$%
\index{ $\K$ (homotopy category)!a@$\K_{\text{\textbf I}}$}
 is the full subcategory of all
\mbox{K-injective} complexes, then we have a $\Delta$-functor
\hbox{$(\C_\bullet\>,\text{Id})\colon
\K_{\text{\textbf I}}(\Y)\to\K_{\text{\textbf I}}(\Avc(\X))$.} 
Associating a K-injective resolution to each complex in~$\A(\Y)$ leads to 
a $\Delta$-functor
\hbox{$(\rho, \Theta)\colon \D(\Y)\to\K_{\text{\textbf I}}(\Y)$}.%
\footnote
{In fact $(\rho, \Theta)$ is an equivalence of
$\Delta$-categories, see \cite[\S1.7]{Derived categories}.
But note that $\Theta$ need not be the identity morphism, i.e., one may
not be able to find a complete family of K-injective resolutions
commuting with translation. For example, we do not know that
every periodic complex has a periodic K-injective resolution.%
} 
This $\rho$ is bounded below: an $\A(\Y)$-complex~$\E$ such that $H^i(\E)=0$
for all $i<n$ is quasi-isomorphic to its truncation
$\E^{{\scriptscriptstyle\ge}n}\<$ (see~\eqref{trunc}),
which is quasi-isomorphic to an
injective complex~$\F$ which vanishes in all degrees below~$n$. 
(Such an~$\F$ is K-injective.)

Finally, one can define~$f^\times$ to be the composition of the functors
$$
\D(\Y)\xrightarrow{\,\rho\,}
\K_{\text{\textbf I}}(\Y)\xrightarrow{\C_\bullet\>} 
\K_{\text{\textbf I}}(\Avc(\X))
\xrightarrow{\text{natural}\,} \D(\Avc(\X)),
$$
and check, via \eqref{f*D} and \eqref{Deligne} that
\Tref{prop-duality} is satisfied. (This involves some
tedium with respect to $\Delta$-details.)
\end{proof}

\section{Torsion sheaves.}
Refer to \S\ref{Gamma'} for notation and first sorites regarding 
torsion sheaves.

Paragraphs~\ref{tors-sheaves} and~\ref{tors-D} develop 
properties of quasi-coherent torsion sheaves and their derived categories on
locally noetherian formal schemes---see e.g.,
Propositions~\ref{Gamma'(qc)}, \ref{Gammas'+kappas}, \ref{Rf-*(qct)}, 
and \Cref{C:f* and Gamma}. (There is some overlap
here with \S4 in \cite{Ye}.) Such properties will be needed throughout the rest of
the paper.  For instance,
Paragraph~\ref{tors-eqvce}  establishes for a noetherian formal scheme~$\X$,
 either separated or finite-dimensional,  an
\emph{equivalence of categories} $\D(\Aqct(\X)){{\mkern8mu\longrightarrow
\mkern-25.5mu{}^\approx\mkern17mu}}\Dqct(\X)$, thereby enabling the use of~
$\Dqct(\X)$---rather than $\D(\Aqct(\X))$---in~\Tref{T:qct-duality}
($\:\cong\:$\Tref{Th2} of \Sref{S:prelim}). Also, \Lref{Gam as holim},
identifying the derived functor $\R\iG\cJ(-)$ (for any $\cO_{\<\<X}$-ideal~$\cJ$,
where $X$ is a ringed space) with the homotopy colimit\- of the functors
$\R\sHomb(\cO_{\<\<X}/\<\cJ^n\<,-)$, plays a key role in the proof of the Base Change
 \Tref{T:basechange} ($\:\cong\:$\Tref{Th3}).

\pagebreak[3]
\begin{parag}\label{tors-sheaves} 
This paragraph deals with categories of quasi-coherent torsion sheaves on locally
noetherian formal schemes.
\end{parag}

\begin{sprop}
\label{f-*(qct)}
Let $f\colon \X\to \Y$ be a map of noetherian formal schemes, and let\/ 
$\M\in\Aqct(\X)$.  Then\/ $f_{\!*}\>\M \in \Aqct(\Y)$. Moreover, if\/ $f$ is
pseudo\kern.6pt-proper\/ \(see~\textup{\S\ref{maptypes}}\) and\/ $\M$ is
coherent then $f_{\!*}\>\M$ is coherent.
\end{sprop}

\begin{proof} 
Let $\J\subset\cO_\X$ and $\I\subset\cO_\Y$ be ideals of definition
such that $\I\cO_\X\subset\J$, and let\vspace{-1pt} 
$$
X_{n}\!:=(\X,\cO_\X/\J^n)
 \xrightarrow{f_{\<n}^{}\>}(\Y,\cO_\Y/\I^n)=:Y_{n}\qquad(n>0)
$$ 
be the scheme-maps induced by $f\<$, so that if  $j_n$ and $i_n$ are the canonical
closed immersions then $fj_n=i_nf_{\<n}^{}$. 
Let $\M_{n}\set\sHom(\cO_{\X}/\J^n\<,\>\M)$, so that 
$$
\M=\iGp{\X}\M= \dirlm{n}\<\M_n=\dirlm{n}j_{n*}j_n^*\>\M_n\>.
$$
Since $\J^{n}$ is  a coherent $\cO_\X$-ideal \cite[p.\,427]{GD},  therefore
$\M_n$ is quasi-coherent (\Cref{C:images}(d)), and it is
straightforward to check that
$i_{n*}f_{\<n*}^{}j_n^*\>\M_n\in\Aqct(\Y)$. Thus, $\X$~being noetherian, and by
\Cref{qct=plump} below, 
$$
f_{\!*}\>\M  =  f_{\!*}\>\> \dirlm{n}\<\M_n\cong
 \dirlm{n} f_{\!*}j_{n*}j_n^*\>\M_n   =  
\dirlm{n} i_{n*}f_{\<n*}^{}j_n^*\>\M_n  \in\Aqct(\Y).
$$

When $f$ is pseudo\kern.6pt-proper every
$f_{\<n}^{}$ is proper; and if $\M\in\Aqct(\X)$ is coherent then so is $f_{\!*}\>\M$, 
because for some $n$,
$f_{\!*}\>\M=f_{\!*}j_{n*}j_{n}^{*}\>\M_n=i_{n*}f_{\<n*}^{}j_n^*\>\M_n$.
\end{proof}

\begin{sprop}
\label{iso-qct}
Let\/ $Z$ be a closed subset of a locally noetherian scheme\/~$X\<,$ and let\/
$\kappa\colon\X \to X$ be the completion of\/ $X$ along $Z$. Then the functors\/
$\kappa^*$ and~$\kappa_*$ restrict to inverse isomorphisms between the
categories\/ $\A_Z(X)$ and\/ $\At(\X),$\ and  between the categories\/ $\AqcZ(X)$
and\/ $\Aqct(\X);\>$\ and  if\/  $\M\in\Aqct(\X)$ is coherent, then so is\/
$\kappa_*\>\M$.
\end{sprop}

\begin{proof}
Let $\cJ$ be a quasi-coherent $\cO_X$-ideal such that the support of~$\cO_X/\cJ$
is~$Z$. 
Applying $\dirlm{n}\!\!$ to the natural isomorphisms\vadjust{\kern-2pt} 
$$
\postdisplaypenalty 1000
\quad\kappa^*\sHom_X(\cO_X/\cJ^n\<\<,\>\>\N\>)\iso
\sHom_\X(\cO_\X/\cJ^n\cO_\X\>,\>\kappa^*\N\>)
 \qquad (\N\in\A(X),\ n>0)
$$
we get a functorial isomorphism $\kappa^*\<\iGp{Z}\iso\iGp{\X}\kappa^*\<$,
and hence $\kappa^*(\A_Z(X))\subset\A_t(\X)$. 
 Applying $\dirlm{n}\!\!$ to the natural isomorphisms\vadjust{\kern-1pt}  
$$
\sHom_X(\cO_X/\cJ^n\<\<,\>\kappa_*\M)\iso
\kappa_*\sHom_\X(\cO_\X/\cJ^n\cO_\X\>,\>\M)\qquad(\M\in\A(\X),\ n>0)
$$
we get a functorial isomorphism
$\iGp{Z}\kappa_*\iso\kappa_*\iGp{\X}\>$, and hence
$\kappa_*(\A_t(\X))\subset\A_Z(X)$.

\goodbreak
As $\kappa$ is a
pseudo\kern.6pt-proper map of locally noetherian formal schemes 
((0) being an ideal of definition of~$X$), 
we see as in the proof of \Pref{f-*(qct)} that
for~\mbox{$\M\in\Aqct(\X)$,} $\kappa_*\>\M$ is a \smash{$\dirlm{}\!\!$}
of quasi-coherent $\cO_X$-modules,  so is itself quasi-coherent, and
$\kappa_*\>\M$ is coherent whenever $\M$ is.%
\footnote{The noetherian assumption in \Lref{f-*(qct)} is needed
only for commutativity of $f_{\!*}$ with $\smash{\subdirlm{}}\!$, a
condition clearly satisfied by $f=\kappa$ in the present situation.\par}

Finally, examining  stalks (see \S\ref{Gamma'}) we find that 
the natural transformations ${\rm 1} \to \kappa_* \kappa^*$ and 
$\kappa^* \kappa_* \to 1$ induce isomorphisms
\begin{align*}
\iGp{Z}\N \iso 
 \kappa_*\kappa^*\iGp{Z}\N \qquad & \bigl(\N \in \A(X)\bigr)\<, \\
\kappa^* \kappa_*\iGp{\X}\M\iso
 \iGp{\X}\M \qquad  & \bigl(\M \in\A(\X)\bigr)\<.
\end{align*}

\vspace{-6.2ex}
\phantom{xxx}
\end{proof}

\begin{scor}\label{qct=plump}
If\/ $\X$ is a locally noetherian formal scheme then\/~$\Aqct(\X)$ is
plump in\/~$\A(\X)$ and closed under small\/  $\A(\X)$-colimits.%

\end{scor}

\begin{proof}
The assertions are local, and so,  since $\At(\X)$ is plump (\S\ref{Gamma'1}),
\Pref{iso-qct} (where
$\kappa^*$ commutes with~\smash{$\dirlm{}\!$}) enables reduction to
well-known facts about
$\AqcZ(X)\subset\A(X)$ with $X$ an affine noetherian (ordinary) scheme.
\end{proof}

\begin{slem}
\label{Gamma'+qc}
Let\/ $\X$ be a locally noetherian formal scheme. 
If\/ $\M$ is a quasi-coherent\/ $\cO_\X$-module then\/
$\iGp{\X}\M\in\Aqct(\X)$ is the\/ \smash{$\dirlm{}\!\!$} of its coherent
submodules.  In particular, $\Aqct(\X)\subset\Avc(\X)$.
\end{slem}

\begin{proof}
 Let $\J$ be an ideal of definition of $\X\>$. For any positive integer~$n$,
let  $X_{n}$ be the scheme~$(\X,
\cO_{\X}/\J^n)$,  let $j_n\colon X_{n}\to\X$  be the canonical closed
immersion, and  let
$\M_n\set\sHom(\cO_{\X}/\J^n\<,\>\M)\subset\iGp \X(\M)$, so that
$\M_n\in\Aqct(\X)$ (\Cref{C:images}(d)). Then the quasi-coherent
$\cO_{\!X_{n}}\<$-module $j_n^*\M_n$ is the 
$\smash{\dirlm{}}\!\!$\vspace{1pt} of its
coherent submodules \cite[p.\,319, (6.9.9)]{GD}, hence so is
$\>\M_n=j_{n*}j_n^*\M_n$\vspace{.6pt}  
(since $j_n^*$ and $j_{n*}$ preserve both $\smash{\dirlm{}}\!\!$ and
coherence\vadjust{\kern1pt}
\cite[p.\,115, (5.3.13) and~(5.3.15)]{GD}), and therefore so is
\mbox{$\iGp{\X}\M=\smash{\dirlm{n}\!\M_n\>}$.} That
$\smash{\dirlm{n}\<\<\M_n\>}\in\Aqct(\X)$\vspace{1.5pt}  results
from \Cref{qct=plump}.
\end{proof}

\begin{scor}\label{C:Qt}
For a locally noetherian formal scheme\/~$\X,$  the inclusion functor\/
$j^{\mathrm t}_{\<\X}\colon\Aqct(\X)\hookrightarrow \A(\X)$ has a right
adjoint\/~$\Qt_{\<\X}$.\index{ $\R$@ {}$Q_\X$
(quasi-coherator)\vadjust{\penalty 10000}!$\Qt_{\<\X}$|(} If\/ moreover\/ $\X$
is noetherian then\/ $\Qt_{\<\X}$ commutes with\/
\smash{$\dirlm{}\!\<.$} 

\end{scor}

\begin{proof}
To show that $j^{\mathrm t}_{\<\X}$ has a right adjoint one can, in
view of \Cref{qct=plump} and \Lref{Gamma'+qc}, simply
apply the Special Adjoint Functor theorem.  

More specifically, since
$\iGp{\X}$ is right-adjoint to the inclusion
$\At(\X)\hookrightarrow\A(\X)$, and $\Avc(\X)\subset\Aqc(\X)$
(\Cref{C:vec-c is qc}), it follows from \Lref{Gamma'+qc}
that the restriction of~$\iGp{\X}$ to $\Avc(\X)$ is right-adjoint to
$\Aqct(\X)\hookrightarrow\Avc(\X)$; and by
\Pref{A(vec-c)-A}, $\Avc(\X)\hookrightarrow\A(\X)$ has a
right adjoint~$Q_\X^{}\>$; so $\Qt_{\<\X}\!:=\iGp{\X}\smcirc Q_\X^{}$
is right-adjoint to~$j^{\mathrm t}_{\<\X}\>$.  (Similarly,
$Q_\X^{}\smcirc \iGp \X$ is right-adjoint to~$j^{\mathrm
t}_{\<\X}\>$.)

Commutativity with \smash{$\dirlm{}\!\!$}\vspace{1pt} means that for any small
directed system~$(\G_\alpha)$ in~$\A(\X)$ and any $\M\in\Aqct(\X)$,
the natural map
$$
\phi\colon\Hom(\M,\>\dirlm\alpha\<\Qt_{\<\X}\>\G_\alpha)\to
\Hom(\M,\Qt_{\<\X}\>\dirlm\alpha\<\G_\alpha)
$$
is an \emph{isomorphism}. This follows from \Lref{Gamma'+qc}, which allows us
to assume that $\M$ is coherent, in which case $\phi$ is isomorphic to
the natural composed isomorphism
$$
\dirlm\alpha\Hom(\M,\>\Qt_{\<\X}\>\G_\alpha)\iso
\dirlm\alpha\Hom(\M,\>\G_\alpha)\iso
\Hom(\M,\>\dirlm\alpha\<\G_\alpha).
$$
\vskip-4ex
\end{proof}
\smallskip

\emph{Remark.} For an ordinary noetherian scheme~$X$ we have 
$\Qt_{\<\<X}=Q_{\<\<X}^{}$\index{ $\R$@ {}$Q_\X$
(quasi-coherator)\vadjust{\penalty 10000}!$\Qt_{\<\X}$|)} (see
\S\ref{SS:vc-and-qc}). More generally,  if $\kappa\colon\X\to X$ is as in
\Pref{iso-qct}, then 
$\Qt_{\X}=\kappa^*\<\iG {Z\>} Q_{\<\<X}^{}\kappa_*$. Hence 
\Pref{f-*(qct)} (applied to open immersions $\X\hookrightarrow\Y$
with $\X$ affine) lets us  construct  the functor $\Qt_{\Y}$ for
any noetherian formal scheme~$\Y$ by mimicking the construction for
ordinary schemes  (cf.~ \cite[p.\, 187, Lemme 3.2]{I}.)

\begin{parag}\label{tors-D}

The preceding results carry over to  derived categories. 

From \Cref{qct=plump} it follows that on a locally noetherian formal
scheme~$\X$,  $\Dqct(\X)$ is a triangulated subcategory of $\D(\X)$, closed under
direct sums.

\penalty-1000
\begin{sprop}\label{Gamma'(qc)}
\hskip-1pt For a locally noetherian formal scheme\/ $\X,$ set\/ $\At\!:=\At(\X),$\
the category of torsion\/ $\cO_\X$-modules,  and let\/ 
$\bi\colon\D(\At)\to \D(\X)$ be the natural \mbox{functor.}
Then$\>:$

\textup{(a)} An\/ $\cO_\X$-complex\/~$\E$ is in\/~$\Dt(\X)$ iff the natural
map\/ $\bi\R\iGp\X\E\to\E$ is a\/ $\D(\X)$-isomorphism.

\textup{(b)} If\/ $\E\in\Dqc(\X)$ then\/ $\R\iGp{\X}\E\in\Dqc(\At)$.

\textup{(c)} The functor $\bi$ and its right adjoint $\R\iGp{\X}$ induce
quasi-inverse equivalences between\/ $\D(\At)$ and\/~$\Dt(\X)$ and between
$\Dqc(\At)$  and\/~$\Dqct(\X)$.%
\footnote{We may therefore sometimes abuse notation and write 
$\R\iGp\X$ instead of $\>\bi\R\iGp\X$; but the meaning should be clear
from the context.}
\end{sprop}

\begin{proof}
(a) For $\F\in\D(\At)$ (e.g., $\F\set\R\iGp\X\E$), any complex isomorphic to 
$\bi\F$ is clearly in~$\Dt(\X)$.

Suppose conversely that $\E\in\Dt(\X)$. The assertion that
$\bi\R\iGp{\X}\E\cong\E$ is local, so we may assume that 
$\X=\Spf(A)$ where $A=\Gamma(\X,\>\cO_\X)$ is a noetherian adic ring, so that
any defining ideal~$\J$ of~$\X$ is  generated by a finite sequence in~$A$. Then
$\bi\R\iGp{\X}\E\cong
\cK_\infty^\bullet\otimes\,\E$, where $\cK_\infty^\bullet$%
\index{ $\K$ (homotopy category)@$\cK_\infty^\bullet$ (limit of Koszul
complexes)} is a bounded flat complex---a
\smash{$\>\dirlm{}\!\!$}\vadjust{\kern1pt} of Koszul complexes on powers of the
generators of~$\J$---see \cite[p.\,18, Lemma 3.1.1]{AJL}. 

So $\bi\R\iGp{\X}$ is a bounded functor, and the usual way-out argument reduces
the question to where
$\E$ is a single torsion sheaf.  But then it is immediate from the construction of
$\cK_\infty^\bullet$ that $\cK_\infty^\bullet\<\otimes\E=\E$.
\smallskip

(b) Again, we can assume that $\X=\Spf(A)$ and  $\R\iGp{\X}$ is bounded,
and since $\Aqc(\X)$ is plump in~$\A(\X)$ (\Pref{(3.2.2)}) we can reduce to where
$\E$ is a single quasi-coherent $\cO_\X$-module, though it is better to assume
only that
$\E\in\Dqc^+(\X)$, for then we may also assume $\E$ injective, so that
$$
\R\iGp{\X}\E\cong\iGp{\X}\E=\dirlm{n>0\,\,\>}\sHom(\cO/\J^n\<,\>\E).
$$
From \Cref{C:images}(d) it follows that
$\sHom(\cO/\J^n\<,\>\E)\in\Dqct(\X)$---for this assertion another way-out
argument reduces us again to where $\E$ is a single quasi-coherent
$\cO_\X$-module---and since homology commutes with \smash{$\dirlm{}\!\!$} and
$\Aqct$ is closed under \smash{$\dirlm{}\!\!$} (\Cref{qct=plump}),
therefore $\R\iGp{\X}\E$ has quasi-coherent homology.
\smallskip

Assertion (c) results now from the following simple lemma. 
\end{proof}

\enlargethispage*{3pt}

\begin{slem}\label{L:j-gamma-eqvce}
Let\/ $\A$ be an abelian category,  let\/
$j\colon\A_\flat\to\A$ be the inclusion of a plump subcategory such that $j$ has a right
adjoint\/~$\varGamma\<\<,$\  and let\/
$\bj\colon\D(\A_\flat)\to\D(\A)$\index{ $\iG{\<\cJ\>}$@$\bj$} be the
derived-category extension of\/ $j$. Suppose that every\/ $\A$-complex has a
K-injective resolution, so that the derived functor\/
$\R\varGamma\colon\D(\A)\to\D(\A_\flat)$ exists.  Then $\R \varGamma$ is
right-adjoint to~$\bj$.  Furthermore, the following conditions are equivalent.
\begin{enumerate}
\item[(1)] $\bj$ induces an equivalence of categories from\/ $\D(\A_\flat)$
to\/ $\D_\flat(\A),$ with quasi-inverse\/~
$\R_\flat\varGamma\!:=\R\varGamma|_{\D_\flat(\A)}$.\vadjust{\kern1pt}
\item[(2)] For every\/ $\E\in \D_\flat(\A)$ the natural map\/ $\bj\R
\varGamma\E\to\E$  is an isomorphism.\vadjust{\kern1pt}
\item[$(3)$] The functor\/~$\R_\flat\varGamma$ is
bounded,  and for\/ $\E_0\in
\A_\flat$ the natural map\/ $\bj\R \varGamma\E_0\to\E_0$ is a\/
$\D(\A)$-isomorphism.
\end{enumerate}
When these conditions hold,  every \/ $\A_\flat$-complex
has a K-injective resolution.
\end{slem}

\begin{proof}
 Since $\varGamma$ has an exact left adjoint, it takes K-injective
$\A$-complexes to K-injective $\A_\flat$-complexes, whence there is a
bifunctorial isomorphism in the derived category of abelian groups
$$
\R\Homb_\A(\bj\G,\>\E)\iso\R\Homb_{\A_\flat}(\G,\>\R\varGamma\E)
\qquad\bigr(\G\in\D(\A_\flat),\ \E\in\D(\A)\bigl).
$$
(To see this, one can assume $\E$ to be K-injective, and then drop the
$\R$'s\dots\!\!) \:Apply homology $\mathrm H^0$ to this isomorphism to get
adjointness of $\bj$ and~$\R\varGamma\<$. 

The implications $(1)\!\Rightarrow\!(3)\!\Rightarrow\!(2)$ are straightforward.
For $(2)\!\Rightarrow\!(1)$, one
needs that for $\G\in\D(\A_\flat)$ the natural map
$\G\to\R\varGamma\bj\G$ is an isomorphism, or equivalently (look at
homology), that the corresponding map
$\bj\G\to\bj\R\varGamma\bj\G$ is an isomorphism. But the composition of
this last map with the isomorphism $\bj\R\varGamma\bj\G\iso\bj\G$ 
(given by~(2)) is the identity, whence the conclusion.

Finally, if $\G$ is an $\A_\flat$-complex  and $j\G\to \cJ$ is a K-injective
$\A$-resolution, then as before $\varGamma\<\cJ$ is a K-injective
$\A_\flat$-complex; and (1) implies that the  natural composition\looseness=-1
$$
\G\to\varGamma\< j\G\to\varGamma\< \cJ\ (\:\cong \R\varGamma\bj \G)
$$ 
is a $\D(\A_\flat)$-isomorphism, hence  an $\A_\flat$-K-injective resolution.
\end{proof}

\begin{scor}\label{C:Hom-Rgamma}
For any complexes\/ $\E\in\Dt(\X)$ and\/ $\F\in\D(\X)$ 
the natural map\/ $\R\iGp\X\F\to\F$ induces an isomorphism
$$
\R\sHomb(\E, \R\iGp\X\F\>)\iso\R\sHomb(\E, \F\>).
$$
\end{scor}

\begin{proof}
Consideration of homology presheaves shows it sufficient that for each affine
open $\U\subset\X$, the natural map 
$$
\Hom_{\D(\U)}\bigl(\E|_\U\>,\>(\R\iGp\X\F\>)|_\U\bigr) 
\to\Hom_{\D(\U)}\bigl(\E|_\U\>,\>\F|_\U\bigr)
$$
be an isomorphism. But since $\R\iGp{\X}$ commutes with restriction to~$\U$,
that is a direct consequence of \Pref{Gamma'(qc)}(c) (with $\X$
replaced by~$\U$).
\end{proof}

Parts (b) and (c) of the following Proposition will be generalized in
parts (d) and~(b), respectively, of \Pref{P:f* and Gamma}.

\begin{sprop}\label{Gammas'+kappas}
Let\/ $Z$ be a closed subset of a locally noetherian scheme\/~$X\<,$ and let\/
$\kappa\colon\X \to X$ be the completion of\/ $X$ along $Z$. Then$\>:$
\vadjust{\kern1pt}

\textup{(a)} The exact functors\/
$\kappa^*$ and~$\kappa_*$ restrict to inverse isomorphisms between the
categories\/ $\D_{\<Z}(X)$ and\/ $\Dt(\X),$\ and  between the categories\/
$\DqcZ(X)$ and\/ $\Dqct(\X);\>$\ and  if\/  $\M\in\Dqct(\X)$ has coherent
homology, then so does\/ $\kappa_*\>\M$.
\vadjust{\kern1pt}

\textup{(b)} There is a unique derived-category isomorphism
$$
\R\iGp Z\kappa_*\E\iso \kappa_*\R\iGp \X\E\qquad\:\bigl(\E\in\D(\X)\bigr)
$$
whose composition with the natural map 
$\kappa_*\R\iGp\X\E\to\kappa_*\E$ is just the natural map 
\mbox{$\R\iGp Z\kappa_*\E\to\kappa_*\E$.}
\vadjust{\kern1pt} 

\textup{(c)} There is a unique derived-category isomorphism
$$
\kappa^*\R\iGp Z\F \iso\R\iGp \X\kappa^*\<\F\qquad\bigl(\F\in\D(X)\bigr)
$$
whose composition with the natural map $\R\iGp\X\kappa^*\<\F\to\kappa^*\<\F$ 
is just the natural map \mbox{$\kappa^*\R\iGp Z\F\to\kappa^*\<\F$.}
\end{sprop}

\pagebreak[3]

\begin{proof}
The assertions in (a) follow at once from \Pref{iso-qct}. 
\vadjust{\kern1pt}

 (b) Since $\kappa_*$ has an exact left adjoint (namely~$\kappa^*$),
therefore $\kappa_*$ transforms \mbox{K-injective} $\A(\X)$-complexes into
K-injective $\A(X)$-complexes, and consequently the isomorphism in~(b) results
from the isomorphism $\iGp Z\kappa_*\iso \kappa_*\iGp\X$ in the proof
of \Pref{iso-qct}.  That the composition in~(b) is as
asserted comes down then to the elementary fact that the natural
composition
$$
\postdisplaypenalty5000
\sHom_X(\cO_X/\cJ^n\<\<,\>\kappa_*\M)\iso
\kappa_*\sHom_{\>\X}(\cO_\X/\cJ^n\cO_\X\>,\>\M)\lra \kappa_*\M
$$
(see proof of \Pref{iso-qct}) is just the obvious map.
Since $\kappa_*\R\iGp\X\E\in\D_{\<Z}(X)$ (by~(a) and 
\Pref{Gamma'(qc)}(a)), the uniqueness
assertion (for the inverse isomorphism) 
results from adjointness of~$\R\iGp Z$ and the inclusion
$\D_{\<Z}(X)\hookrightarrow \D(X)$. (The proof is similar to that of 
\Pref{Gamma'(qc)}(c)).
\vadjust{\kern1pt}

(c) Using (b), we have the natural composed map
$$
\kappa^*\R\iGp Z\F \to \kappa^*\R\iGp Z\kappa_*\kappa^*\<\F\iso
\kappa^*\kappa_*\R\iGp\X\kappa^*\<\F\to\R\iGp\X\kappa^*\<\F.
$$
Showing this  to be an isomorphism is a local problem, so 
assume  $X=\Spec(A)$  with $A$  a noetherian adic ring. Let
$K_\infty^\bullet$ be the usual
\smash{$\dirlm{}\!\!$}\vspace{.8pt} of Koszul
complexes on powers of a finite system of
generators\vadjust{\kern.4pt} of an ideal of definition of~$A$ 
(\cite[\S3.1]{AJL}); and let \smash{$\widetilde K_\infty^\bullet$} be the 
corresponding quasi-coherent complex on $\Spec(A)$, so that the complex
$\cK_\infty^\bullet$%
\index{ $\K$ (homotopy category)@$\cK_\infty^\bullet$ (limit of Koszul
complexes)} in the proof of \Pref{Gamma'(qc)}(a) is
just~\smash{$\kappa^*\widetilde K_\infty^\bullet$.} Then one checks via~
\cite[p.\,18, Lemma~(3.1.1)]{AJL} that the map in question is
isomorphic to the natural isomorphism of complexes
$$
\kappa^*(\widetilde K_\infty^\bullet\otimes_{\cO_{\<\<X}}\<\F\>)\iso
\kappa^*\widetilde K_\infty^\bullet\otimes_{\cO_{\X}}\<\kappa^*\<\F.
$$

That the composition in (c) is as asserted results from the following
natural commutative diagram, whose bottom row composes to the identity:
$$
\minCDarrowwidth=22pt
\begin{CD}
\kappa^*\R\iGp Z\F
 @>>> \kappa^*\R\iGp Z\kappa_*\kappa^*\<\F
  @>\Iso>>\kappa^*\kappa_*\R\iGp\X\kappa^*\<\F
   @>>>\R\iGp\X\kappa^*\<\F \\
@VVV @VV\hskip3.5em\textup{(b)}V @VVV @VVV \\
\kappa^*\<\F
 @>>>\kappa^*\kappa_*\kappa^*\<\F 
  @= \kappa^*\kappa_*\kappa^*\<\F
   @>>>\kappa^*\<\F
\end{CD}
$$
Uniqueness is shown as in (b).
\end{proof}
\smallskip

\begin{scor}
\label{C:Gammas'+kappas}
The natural maps are isomorphisms
\begin{alignat*}{2}
\Hom_{\<X\<}(\E,\>\F\>)
  &\cong\Hom_{\<X\<}(\E,\kappa_*\kappa^*\<\F\>)\cong
   \Hom_\X(\kappa^*\E,\kappa^*\<\F\>)
\quad&&\bigl(\E\in\D_{\<Z}(X),\, \F\in\D(X)\bigr), \\
\Hom_{\<X\<}(\E,\>\F\>)
   &\cong\Hom_{\<X\<}(\E,\kappa_*\kappa^*\<\F\>)\cong
    \Hom_\X(\kappa^*\E,\kappa^*\<\F\>)
\quad&&\bigl(\E\in\D(X),\, \F\in\D_{\<Z}(X)\bigr), \\
\Hom_\X(\G,\>\H\>)
  &\cong\Hom_\X(\kappa^*\kappa_*\G,\>\H\>)\cong
    \Hom_{\<X\<}(\kappa_*\G,\kappa_*\H\>)
\quad&&\bigl(\G\in\Dt(\X),\, \H\in\D(\X)\bigr).
\end{alignat*}
\end{scor}

\begin{proof}
For the first line, use
\Pref{Gamma'(qc)} and its analogue for $\D_{\<Z}(X)$,
\Lref{L:j-gamma-eqvce}, and \Pref{Gammas'+kappas} to get the
equivalent sequence of natural isomorphisms
\begin{align*}
\Hom_{\<X\<}(\E,\F\>)&\cong \Hom_{\<X\<}(\E,\R\iGp Z\F\>)\\
&\cong \Hom_\X(\kappa^*\E,\kappa^*\R\iGp Z\F\>)\\
&\cong \Hom_\X(\kappa^*\E,\R\iGp \X\kappa^*\<\F\>)\\
&\cong\Hom_\X(\kappa^*\E, \kappa^*\<\F\>)\\
&\cong\Hom_{\<X\<}(\E,\kappa_*\kappa^*\<\F\>).
\end{align*}

The rest is immediate from \Pref{Gammas'+kappas}(a).
\end{proof}
\medskip
\pagebreak[3]
The next series of results concerns the behavior of $\Dqct$ with respect to maps
of formal schemes.

\begin{sprop}
\label{Rf-*(qct)}
Let $f\colon \X \to \Y$ be a map of noetherian formal
schemes. Then\/ $\Rfs|_{\Dqct(\X)}$ is bounded, and
$$
\R f_{\!*} \bigl(\Dqct(\X)\bigr) \subset \Dqct(\Y).
$$
Moreover, if\/ $f$ is
pseudo\kern.6pt-proper and\/ $\F\in\Dt(\X)$ has coherent homology, then\/ 
so does $\R f_{\!*}\>\F\in\Dt(\Y)$.
\end{sprop}

\begin{proof}
Since $\Dqct(\X)\subset\Dvc(\X)$ (\Lref{Gamma'+qc}), the boundedness
assertion is given by \Pref{Rf_*bounded}(b). (Clearly, $\Rfs$ is
bounded-below.) It suffices then for the next assertion (by the usual way-out
arguments \cite[p.\,73, Proposition 7.3]{H1}) to show for any $\M\in\Aqct(\X)$ that
$\R f_{\!*}\M\in\Dqct(\Y)$. 

Let $\E$ be an injective resolution of~$\M$, let $\J$ be an ideal of
definition of~$\X$, and let $\E_n$ be the flasque complex $\E_n\set
\sHom(\cO/\J^n\<,\>\E)$. Then by \Pref{Gamma'(qc)}(a),
\mbox{$\M\cong\R\iGp \X\M\cong \smash{\dirlm{}\!\<_n\>\>\E_n}\>$.}
 Since $\X$ is noetherian, 
\smash{$\dirlm{}\!\!$}'s of flasque sheaves are \mbox{$f_{\!*}$-acyclic} and
\smash{$\dirlm{}\!\!$} commutes with~$f_{\!*}\>$; 
so with  notation as in the proof of \Pref{f-*(qct)},

$$
\Rfs \M\cong\R f_{\!*}\>\R\iGp \X\M
\cong
f_{\!*}\>\dirlm{n}\<\E_n
\cong
\dirlm{n}\<f_{\!*}j_{n*}j_n^*\E_n
\cong 
\dirlm{n} i_{n*}f_{\<n*}^{}j_n^*\E_n\>.
$$
Since $\E\in\Dqc^+(\X)$, therefore
$$
j_{n*}j_n^*\E_n=\sHom(\cO/\J^n\<,\>\E)\in\Dqc(\X),
$$
as we see by way-out reduction to where $\E$ is a single quasi-coherent
sheaf and then by \Cref{C:images}(d);  and hence 
$j_n^*\E_n\in\Dqc(X_{n})$ (see \cite[p.\,115, (5.3.15)]{GD}).
Now $j_n^*\E_n$~is a flasque bounded-below $\cO_{\!X_{\<n}}$-complex,  so
by way-out reduction to (for example) \cite[p.\,643, corollary~11]{Ke},
$$
f_{\<n*}^{}j_n^*\E_n\cong
\R f_{\<n*}^{}j_n^*\E_n\in\Dqc(Y_{n});
$$
and finally, in view of \Cref{qct=plump}, 
$$
\R f_{\!*}\M\cong i_{n*}\>\dirlm{n}\< f_{\<n*}^{}j_n^*\E_n
\in\Dqct(\Y).
$$

For the last assertion, we reduce as before to showing for each coherent
torsion $\cO_\X$-module~$\M$ and each $p\ge0$ that 
${R}^p\!f_{\!*}\M\set H^p\Rfs\M$ is a
coherent $\cO_\Y$-module. With notation remaining as in \Pref{f-*(qct)},
the maps $i_n$ and $j_n$ 
are exact, and for some $n$, \mbox{$\M=j_{n*}j_n^*\M_n$.} So 
$$
R^p\!f_{\!*}\>\M=R^p\!f_{\!*}\>j_{n*}j_n^*\M_n=
i_{n*}R^p\!f_{n*}\>j_n^*\M_n,
$$
which is coherent since $j_n^*\M_n$ is a coherent $\cO_{\!X_n}$-module
and $f_n\colon X_n\to Y_n$ is a proper scheme-map.
\end{proof}

\begin{scor}[cf.~\Cref{P:kappa-f*}]\label{C:kappa-f*t}
Let\/ $f_0\colon X\to Y$ be a  map of locally noetherian schemes, let\/
$W\subset Y$ and\/ $Z\subset f_0^{-1}W$ be closed subsets, with associated
\(\kern.5pt flat\) completion
maps\/ $\kappa_\Y^{\pd}\colon\Y=Y_{\</W}\to Y\<, \,$ 
$\kappa_\X^{\pd}\colon\X=X_{\</Z}\to X\<,$   and let\/
$f\colon\X\to\Y$ be the map induced by~$f_0\>$. For\/
$\E\in\D(X)$ let 
$$
\theta_{\<\E}\colon\kappa_\Y^*\R f_{\<0*}^{}\E\to \R f_{\!*}\kappa_\X^*\E
$$
be the map  adjoint to the natural composition
$$
\R f_{\!0*}\E\longrightarrow\R f_{\<0*}^{}\kappa_{\X*}^{\phantom*}\kappa_\X^*\E\iso
\kappa_{\Y*}^{\phantom*}\R f_{\!*}\kappa_\X^*\E.
$$
Then\/ $\theta_{\<\E}$ is an isomorphism for all\/ $\E\in\DqcZ(X)$ .
\end{scor}

\begin{proof}
 $\theta_{\<\E}$ is the composition of the natural maps
$$
\kappa_\Y^*\R f_{\!0*}\E\to
\kappa_\Y^*\R f_{\<0*}^{}\kappa_{\X*}^{\phantom*}\kappa_\X^*\E\iso
\kappa_\Y^*\kappa_{\Y*}^{\phantom*}\R f_{\!*}\kappa_\X^*\E
\to \R f_{\!*}\kappa_\X^*\E.
$$
By \Pref{Gammas'+kappas}, the first map and (in view of
\Pref{Rf-*(qct)})  the third~map are both isomorphisms.
\end{proof}

\begin{sprop}\label{P:f* and Gamma}
Let\/ $f\colon\X\to\Y$ be a map of locally noetherian formal schemes.
Let\/ $\I$ be a coherent\/ $\cO_\Y$-ideal, and let $\D_\I(\Y)$ be the
triangulated subcategory of\/~$\D(\Y)$ whose objects are the
complexes\/~$\F$ with\/ $\I$-torsion homology\/ \textup(i.e.,
$\iG\I\<\<H^i\F=H^i\F$ for all\/ $i\in\mathbb Z$---see
\S\S\textup{\ref{S:prelim}} and\/~\textup{\ref{Gamma'1}).}  Then$\>:$

\smallskip
\textup{(a)} $\bL f^*\<(\D_\I(\Y))\subset\D_{\I\cO_\X}\<(\X)$.

\smallskip
\textup{(b)}  There is a unique
functorial isomorphism 
$$
\hskip100pt \xi(\E)\colon\bL f^*\R\iG\I\E\iso \R\iG{\I\cO_\X}\bL
f^*\E
\qquad \ \bigl(\E\in\D(\Y)\bigr)
$$
whose composition with the natural map\/ 
$\R\iG{\I\cO_\X}\bL f^*\E\to\bL f^*\E$ is
the natural map\/ $\bL f^*\R\iG\I\E\to\bL f^*\E\<$. 

\smallskip
\textup{(c)}  The natural map is an isomorphism
$$
\R\iGp\X\> \bL f^*\R\iGp\Y\E\iso \R\iGp\X\> \bL f^*\E
\qquad \ \bigl(\E\in\D(\Y)\bigr).
$$

\smallskip
\textup{(d)}  If\/ $\X$ is noetherian, there is a unique
functorial isomorphism
$$
\hskip100pt\R\iG\I\<\Rfs\>\G \iso \Rfs\R\iG{\I\cO_\X}\G\qquad
\ \bigl(\G\in\D^+(\X)\bigr)
$$
whose composition with the natural map\/ 
$\Rfs\R\iG{\I\cO_\X}\G\to\Rfs\>\G$ is
the natural map\/ $\R\iG\I\<\Rfs\>\G\to\Rfs\>\G$.
\end{sprop}

\begin{proof}
(a) Let $\F\in\D_\I(\Y)$.  To show that $\bL
f^*\<\F\in\D_{\I\cO_\X}\<(\X)$ we may assume that $\F$~is K-injective.
Let $x\in\X$, set $y\set f(x)$, and let $P_{\<x}^\bullet$ be a flat
resolution of the $\cO_{\Y\<,y}$-module~$\cO_{\X\<,x}\>$.  Then, as in
the proof of \Pref{Gamma'(qc)}(a), there is a canonical
$\D(\Y)$-isomorphism
$$
\dirlm{n} \sHomb(\cO_\Y/\I^n\<,\>\F\>)=\iG\I\<\F=\R\iG\I\<\F\iso\F,
$$
and it follows that for any~$i$ the stalk at~$x$ of the homology~$H^i\bL
f^*\<\F$ is 
$$
\textup H^i\bigl(P_{\<x}^\bullet\otimes_{\cO_{\Y\<,y}}
  \F^{\phantom{.}}_{\!y} \bigr)
=\dirlm{n}\textup H^i\bigl(P_{\<x}^\bullet\otimes_{\cO_{\Y\<,y}}
 \Homb_{\cO_{\Y\<,y}}(\cO^{\phantom{.}}_{\Y\<,y}/\I^n_{\!y}\>,
 \>\F^{\phantom{.}}_{\!y})\bigr).
$$
Hence each element of the stalk is annihilated by a power
of~$\I\cO_{\X\<,x}\>$, and (a) results.

\smallskip
(b) The existence and uniqueness of a functorial map $\xi(\E)$ satisfying
everything except the isomorphism property result from (a) and the fact that
$\R\iG{\I\cO_\X\<\<}$ is right-adjoint to the inclusion
$\D_{\I\cO_\X}\<\<(\X)\hookrightarrow\D(\X)$.

To show that $\xi(\E)$ is an isomorphism we may assume that $\Y$ is affine and
that $\E$ is K-flat, and then proceed as in the proof of (the special case)
\Pref{Gammas'+kappas}(c), via the bounded flat complex $K_\infty^\bullet\>$.

\smallskip
(c) Let $\I$, $\J$  be  defining ideals of~$\Y$ and~$\X$ respectively, so that
$\cK\set\I\cO_\X\subset\J$.  The natural map
$
\R\iGp\X\>\R\iG\cK:=\R\iG\J\R\iG\cK\to\R\iG\J=:\R\iGp\X
$
 is an \emph{isomorphism,} as one checks locally via
~\cite[p.\,20, Corollary~(3.1.3)]{AJL}. So for any $\E\in\D(\Y)$, 
(b) gives
$$
\R\iGp\X\>\bL f^*\<\E
\cong \R\iGp\X\>\R\iG\cK\bL f^*\<\E
\cong \R\iGp\X\>\bL f^*\>\R\iGp\Y\>\E.
$$

\smallskip
(d) $\G$ may be assumed bounded-below and injective, so that
$$
\G_n\set\sHomb(\cO_\X/\I^n\cO_\X\>,\>\G)
$$ 
is flasque.  

Then, since $\X$ is noetherian, 
$\iG{\I\cO_\X}\G=\smash{\dirlm{}\!\mkern-1.5mu _n\>\>\G_n}\>$
is flasque too, and
$$
\Rfs\>\R\iG{\I\cO_\X}\G\cong
\Rfs\iG{\I\cO_\X}\G\cong f_{\!*}\mkern1.5mu{\dirlm{n}\G_n}
\cong {\dirlm{n}\<f_{\!*}\>\G_n}\in\D_\I(\Y).
$$  
By \Lref{L:j-gamma-eqvce}, $\R\iG\I$ (resp.~$\R\iG{\I\cO_\X\<\<}$) is
right-adjoint to the inclusion
\mbox{$\D_\I(\Y)\hookrightarrow\D(\Y)$}
(resp.~$\D_{\I\cO_\X}(\X)\hookrightarrow\D(\X)$), 
whence, in particular, the uniqueness in~(d).
Moreover, in view of~(a), for any
$\E\in\D_\I(\Y)$ the natural maps are isomorphisms
\begin{multline*}
 \Hom_\Y(\E,\>\R\iG\I\<\Rfs\>\G)
\iso
 \Hom_\Y(\E\<,\>\Rfs\>\G)
\iso
 \Hom_\X(\bL f^*\<\E\<,\>\G) \\
\iso
 \Hom_\X(\bL f^*\<\E\<,\>\R\iG{\I\cO_\X}\G) 
\iso
  \Hom_\Y(\E\<,\>\Rfs\>\R\iG{\I\cO_\X}\G).
\end{multline*}
It follows formally that the image under this composed isomorphism of the
identity map of~$\R\iG\I\<\Rfs\>\G$
is an isomorphism as asserted. (In fact this
isomorphism is adjoint to the composition
$\bL f^*\R\iG\I\<\<\Rfs\G
\xrightarrow[\xi(\Rfs\G)]{}
\R\iG{\I\cO_\X}\bL f^*\Rfs\G
\xrightarrow[\textup{nat'l}]{}
\R\iG{\I\cO_\X}\G.\>)$
\end{proof}

\begin{sdef}\label{D:Dtilde}\index{ $\D$ (derived category)!z@${ \widetilde
{\vbox to5pt{\vss\hbox{$\mathbf D$}}}_{\mkern-1.5mu\mathrm {qc}} }$}
For a locally noetherian formal scheme~$\X$, 
$$
\wDqc(\X)\set\R\iGp\X{}^{-1}(\Dqc(\X))
$$
is the $\Delta$-subcategory of $\D(\X)$ whose objects are those complexes~$\F$
such that $\R\iGp\X\F\in\Dqc(\X)$---or equivalently, $\R\iGp\X\F\in\Dqct(\X)$.
\end{sdef}

\begin{srems}\label{R:Dtilde}
(1) By \Pref{Gamma'(qc)}(b), $\Dqc(\X)\subset\wDqc(\X)$. Hence
$$
\R\iGp\X\bigl(\>\wDqc(\X)\bigr)\subset\wDqc(\X).
$$

(2) Since $\R\iGp\X$ is idempotent (see \Pref{Gamma'(qc)}), the vertex  of any
triangle based on the canonical map $\R\iGp\X\E\to\E\ (\E\in\D(\X))$ is
annihilated by $\R\iGp\X$. It follows that $\wDqc(\X)$ is the smallest
$\Delta$-subcategory of~$\D(\X)$ containing $\Dqct(\X)$ and all complexes~$\F$
such that $\R\iGp\X\F=0$.\vspace{1.5pt}

\smallskip
(3) The functor $\R\iGp\X\colon\D(\X)\to \D(\X)$ has a
right adjoint\index{ $\mathbf {La}$@$\BL$ (homology localization)}
$$
\BL_\X(-)\set\R\sHomb(\R\iGp\X\cO_\X^{}\>,-).
$$
Indeed, there are  natural functorial isomorphisms for  $\E,\F\in\D(\X)$, 
\begin{equation}\label{adj}
\begin{aligned}
\Hom_{\D(\X)\<}(\R\iGp\X\E\<,\>\F\>)
&\iso\Hom_{\D(\X)\<}(\E\Otimes\R\iGp\X\cO_\X^{}\>,\>\F\>)\\
&\iso\Hom_{\D(\X)\<}\bigl(\E,\>\R\sHomb(\R\iGp\X\cO_\X^{}\>,\>\F\>)\bigr).
\end{aligned}
\end{equation}
(Whether the natural map \smash{$\E\Otimes\R\iGp\X\cO_\X^{}\iso\R\iGp\X\E$}
is an isomorphism\vspace{1.3pt} is a local question, dealt with e.g., in
\cite[p.\,20, Corollary~(3.1.2)]{AJL}. \vspace{.8pt} The second isomorphism
is given, e.g., by \cite[p.\,147, Proposition 6.6\,(1)]{Sp}.)\vspace{.5pt}

\pagebreak[3]
There is a natural isomorphism $\R\iGp\X\iso\R\iGp\X\BL_\X^{}$ (see (d) in 
\Rref{R:Gamma-Lambda} below), and consequently 
$$
\BL_\X^{}\bigl(\>\wDqc(\X)\bigr)\subset\wDqc(\X).
$$

\smallskip
(4) \emph{If\/ $\E\in\Dc^-(\X)$ and\/ $\F\in\wDqc(\X)$ then\/
$\R\sHomb(\E\<,\>\F\>)\in\wDqc(\X)$, and hence\/
$\R\sHomb(\R\iGp\X\E\<,\>\F\>)\in\wDqc(\X)$.} Indeed, the natural map
$$\postdisplaypenalty10000
\R\iGp\X\>\R\sHomb(\E\<,\R\iGp\X\F\>)\to\R\iGp\X\>\R\sHomb(\E\<,\>\F\>)
$$
is an \emph{isomorphism,} since for any
$\G$ in $\Dt(\X)$, $\>\smash{\G\Otimes\E}\in\Dt(X)$ (an assertion which can
be checked locally, using
\Pref{Gamma'(qc)}(a) and the complex~$\cK_\infty^\bullet$ in its
proof), so that there is a sequence of natural isomorphisms (see
\Pref{Gamma'(qc)}(c)):
\begin{align*}
\Hom\bigl(\G, \R\iGp\X\>\R\sHomb(\E\<,\R\iGp\X\F\>)\bigr)
&\iso
\Hom\bigl(\G, \R\sHomb(\E\<,\R\iGp\X\F\>)\bigr) \\
&\iso
\Hom\bigl(\smash{\smash{\G\Otimes\E}}\<, \R\iGp\X\F\bigr) \\
&\iso
\Hom\bigl(\smash{\G\Otimes\E}\<, \F\bigr) \\
&\iso
\Hom\bigl(\G, \R\sHomb(\E\<,\>\F\>)\bigr)\\
&\iso
\Hom\bigl(\G, \R\iGp\X\>\R\sHomb(\E\<,\>\F\>)\bigr).
\end{align*}
Since $\Aqct(\X)$ is plump in~$\A(\X)$ (\Cref{qct=plump}),
\Pref{P:Rhom} shows that
\mbox{$\R\iGp\X\>\R\sHomb(\E\<,\R\iGp\X\F\>)\in\Dqct(\X)$,} whence
$\R\sHomb(\E\<,\>\F\>)\in\wDqc(\X)$.

From (3) and the natural isomorphisms
$$
\R\sHomb(\R\iGp\X\E\<,\>\F\>)\cong
\R\sHomb(\R\iGp\X\cO_\X^{}\Otimes\E\<,\>\F\>)\cong
\BL_\X^{}\R\sHomb(\E\<,\>\F\>)
$$
we see then that 
$$
\R\sHomb(\R\iGp\X\E\<,\>\F\>)\in\wDqc(\X).
$$

\smallskip

\pagebreak[3]
(5) For $\F\in\D(\X)$ it holds that
$$
\F\in\wDqc(\X)\iff\R\sHomb(\cO_\X/\J,\>\F\>)\in\Dqct(\X) 
\text{ for all defining
ideals~$\J$ of~$\X$.}
$$
The implication $\implies$ is given, in view of \Cref{C:Hom-Rgamma},
by \Pref{P:Rhom};
and the converse is given by \Lref{Gam as holim},
since \Cref{qct=plump} implies that $\Dqct(\X)$ is a $\Delta$-subcategory of
$\D(\X)$ closed under direct sums.

\smallskip
(6)\vspace{.7pt} Let $f\colon\X\to\Y$ be a map of locally noetherian formal
schemes. For any $\F\in\wDqc(\Y)$, 
\Lref{Gamma'+qc} and \Pref{P:Lf*-vc} give
$$
\bL f^*\>\R\iGp\Y\>\F\in \bL f^*(\Dqct(\Y))\subset\bL
f^*(\Dvc(\Y))\subset\Dqc(\X)\subset\wDqc(\X),
$$
and so 
$\R\iGp\X\>\bL f^*\<\F
\underset{\textup{\ref{P:f* and Gamma}(c)}}\cong
\R\iGp\X\>\bL f^*\R\iGp\Y\F\in\Dqct(\X).$
Thus\vspace{-5pt} 
$$
\bL f^*\bigl(\>\wDqc(\Y)\bigr)\subset\wDqc(\X).
$$
\end{srems}

\smallskip
\begin{scor}\label{C:f* and Gamma}
Let\/ $f\colon\X\to\Y$ be an adic map of locally noetherian formal schemes.
Then$\>:$

\smallskip
\textup{(a)} $\bL f^*\<(\Dt(\Y))\subset\Dt(\X)$.

\smallskip
\textup{(b)} $\bL f^*\<(\Dqct(\Y))\subset\Dqct(\X)$.

\smallskip
\pagebreak[3]

\textup{(c)} There is a unique
functorial isomorphism
$$
\bL f^*\R\iGp\Y\E\iso \R\iGp\X\>\bL f^*\E
\qquad \ \bigl(\E\in\D(\Y)\bigr)
$$
whose composition with the natural map\/ 
$\R\iGp\X\>\bL f^*\E\to\bL f^*\E$ is
the natural map\/ $\bL f^*\R\iGp\Y\E\to\bL f^*\E$. There results a conjugate 
isomorphism of  right-adjoint functors
$$
\Rfs\BL_\X\>\G\iso \BL_\Y\Rfs\G\qquad\bigl(\G\in\D(\X)\bigr).
$$
whose composition with the natural map\/ $\Rfs\G\to\Rfs\BL_\X\>\G$
is the natural map\/ $\Rfs\G\to\BL_\Y\Rfs\G$.

\smallskip
\goodbreak
\textup{(d)} If\/ $\X$ is noetherian then there is a unique
functorial isomorphism
$$
\hskip100pt\R\iGp\Y\>\Rfs\>\G \iso \Rfs\R\iGp\X\>\G\qquad
\ \bigl(\G\in\D^+(\X)\textup{ or }\>\G\in\wDqc(\X)\bigr)
$$
whose composition with the natural map\/ $\Rfs\R\iGp\X\>\G\to\Rfs\>\G$ is
the natural map\/ $\R\iGp\Y\>\Rfs\>\G\to\Rfs\>\G$.

\smallskip
\textup{(e)} If\/ $\X$ is noetherian then\/ 
$\Rfs\<\bigl(\>\wDqc(\X)\bigr)\subset\wDqc(\Y).$
\end{scor}

\begin{proof}
To get (a) and (c) take $\I$ in \Pref{P:f* and
Gamma} to be an ideal of definition of~$\Y$.  
(The second assertion in (c) is left to the reader.) As
$\Dqct(\Y)\<=\>\Dvc(\Y)\cap\>\Dt(\Y)$ (\Cref{C:vec-c is qc} and
\Lref{Gamma'+qc}), (b) follows from (a) and \Pref{P:Lf*-vc}.
The same choice of~$\>\I$ gives (d) for $\G\in\D^+(\X)$---and the
argument also works for  \smash{$\G\in\wDqc(\X)$} once one notes that
$$
\Rfs\>\R\iGp\X\<\bigl(\>\wDqc(\X)\bigr)
\subset\Rfs\<\bigl(\Dqct(\X)\bigr)
\underset{\mathstrut\text{\ref{Rf-*(qct)}}}\subset \Dqct(\Y)\subset\Dt(\Y).
$$
The isomorphism in (d) gives (e) via \Pref{Rf-*(qct)}.
\end{proof}

\begin{scor} \label{C:kappa-f*t'}
In \Cref{C:kappa-f*t}\textup{,} if\/ $\X$ is noetherian and\/
 $Z= f_0^{-1}W$ then for all\/ $\F\in\Dqc(X)$ the map\/
$\theta_{\<\F}'\set\R\iGp\Y(\theta_{\<\F}\<)$ is an isomorphism
$$
\theta_{\<\F}'\colon\R\iGp\Y\kappa_\Y^*\R f_{\<0*}^{}\>\F\iso
\R\iGp\Y\>\Rfs\kappa_\X^*\>\F.
$$
\end{scor}

\begin{proof} 
Arguing as in \Pref{Gamma'(qc)}, we find that 
$\R\iG Z\>\F\in\DqcZ(X)$, so that we have the isomorphism 
$\theta_{\R\iG Z\F}$ of \Cref{C:kappa-f*t}.\vspace{1pt} 

Imitating the proof of~ \Cref{C:f* and Gamma}, we get an
isomorphism 
$$
\alpha_{\<\F}^{}\colon\R f_{\<0*}^{}\R\iG Z\>\F\iso\R\iG W\R f_{\<0*}^{}\>\F
$$
whose composition with the natural map 
$\R\iG W\R f_{\<0*}^{}\>\F\to\R f_{\<0*}^{}\>\F$ is the natural map
$\R f_{\<0*}^{}\R\iG Z\>\F\to\R f_{\<0*}^{}\>\F$.

Consider then the diagram
$$
\begin{CD}
\kappa_\Y^*\R f_{\<0*}^{}\R\iG Z\>\F @>\Iso>
\smash{\kappa_\Y^*(\<\alpha_{\<\<\F}^{}\<)}>
\kappa_\Y^*\R\iG W\R f_{\<0*}^{}\>\F
@>\Iso>\smash{\textup{\ref{Gammas'+kappas}(c)}}>
\R\iGp\Y\kappa_\Y^*\R f_{\<0*}^{}\>\F@>\textup{nat'l}>>
\kappa_\Y^*\R f_{\<0*}^{}\>\F\\
@V\theta_{\R\iG Z\F} V \simeq V  @. @V(1)\hskip8.8em V\theta_{\!\F}'V 
    @VV\theta_{\!\F}^{} V\\
\Rfs\kappa_\X^*\R\iG Z\>\F @>\Iso>\textup{\ref{Gammas'+kappas}(c)}>
\Rfs\R\iGp\X\kappa_\X^*\>\F @>\Iso>\textup{\ref{C:f* and Gamma}(d)}>
\R\iGp\Y\>\Rfs\kappa_\X^*\>\F @>>\textup{nat'l}>
\Rfs\kappa_\X^*\>\F\\
\end{CD}
$$
It suffices to show that subdiagram (1) commutes; and 
since $\R\iGp\Y$ is right-adjoint to the inclusion $\Dt(\Y)\hookrightarrow\D(\Y)$
it follows that it's enough to show that the outer border of the diagram commutes.
But it is straightforward to check that the top and bottom rows compose to the
maps induced by the natural map $\R\iG Z\to\bf 1$, whence the conclusion.
\end{proof}
\end{parag}

\begin{parag}\label{tors-eqvce}

From the following key \Pref{1!}---generalizing the noetherian
case of
\cite[p.\,12, Proposition~(1.3)]{AJL}---there will result, for complexes
with quasi-coherent torsion homology, a stronger version of the
Duality \Tref{prop-duality}, see \Sref{S:t-duality}.

Recall what it means for a noetherian formal scheme~$\X$ to
be \emph{separated} (\S\ref{note1}).
Recall also from \Cref{C:Qt} that the inclusion functor\/
$j^{\mathrm t}_{\<\X}\colon\Aqct(\X)\hookrightarrow \A(\X)$ has a right
adjoint\/~$\Qt_{\<\X}$.

\pagebreak[3] 
\begin{sprop}
\label{1!}Let\/ $\X$ be a noetherian formal scheme.

\textup{(a)} The  extension
of\/~$j^{\mathrm t}_{\<\X}$ induces an \emph{equivalence of
categories}\index{ $\iG{\<\cJ\>}$@$\bj$!$\bj^{\mathrm t}$} 
$$
\bj^{\mathrm t}_{\<\<\X}\colon\D^+(\Aqct(\X)){{\mkern8mu\longrightarrow
\mkern-25.5mu{}^\approx\mkern17mu}}\Dqct^+(\X),
$$
with bounded quasi-inverse\/  $\R\Qt_{\X}|_{\Dqct^+(\X)}$.\vspace{1pt}

\textup{(b)} If\/ $\X$ is separated, or of finite Krull dimension, then
the  extension of\/~$j^{\mathrm t}_{\<\X}$ induces an \emph{equivalence of
categories}
$$
\bj^{\mathrm t}_{\<\<\X}\colon\D(\Aqct(\X)){{\mkern8mu\longrightarrow
\mkern-25.5mu{}^\approx\mkern17mu}}\Dqct(\X),
$$
with bounded quasi-inverse\/  $\R\Qt_{\X}|_{\Dqct(\X)}$.
\end{sprop}

\begin{proof}
(a) The asserted equivalence is 
given by \cite[Theorem 4.8]{Ye}.  The idea is that  $\Aqct(\X)$ contains enough
$\At(\X)$-injectives \cite[Proposition 4.2]{Ye}, so by
\cite[p.\,47, Proposition 4.8]{H1}, $\D^+(\Aqct(\X))$ is equivalent to
$\Dqc^+(\At(\X))$, which is  equivalent to~$\Dqct^+(\X)$ (\Pref{Gamma'(qc)}(c)). 

Since $\R\Qt_{\X}$\vspace{.4pt} is right-adjoint to
$\bj^{\mathrm t}_{\<\<\X}$ (\Lref{L:j-gamma-eqvce}), 
its restriction to $\Dqct^+(\X)$ is quasi-inverse to 
$\bj^{\mathrm t}_{\<\<\X}|_{\D^+(\Aqct(\X))}$.
From the resulting isomorphism 
$$
\iota^{}_\E\colon\bj^{\mathrm t}_{\<\<\X}\R\Qt_{\X}\E\iso\E\qquad
(\E\in\Dqct^+(\X))
$$
we see that if $H^i\E=0$ then $H^i\> \R\Qt_{\X}\E=0$, so that
$\R\Qt_{\X}|_{\Dqct^+(\X)}$ is bounded.

(b) By \Lref{L:j-gamma-eqvce}, and having the isomorphism~$\iota^{}_\E$,
we need only show that 
$\R\Qt_\X$~is bounded on~$\Dqct(\X)$. 

Suppose that $\X$ is the completion of a 
separated ordinary noetherian scheme~$X$ along some closed subscheme, and let
$\kappa\colon\X\to X$ be the completion map, so that 
$\Qt_{\X}=\kappa^*\<\iG Z Q_{\<\<X}^{}\kappa_*^{}$ (see remark
following \Cref{C:Qt}). The exact functor~$\kappa_*^{}$ preserves K-injectivity,
since it has an exact left adjoint, namely~$\kappa^*\<$. Similarly
$Q_{\<\<X}$~transforms K-injective
$\A(X)$-complexes into K-injective
$\Aqc(X)$-complexes. Hence 
$\R\Qt_\X\cong \kappa^*\R\iG Z^{\!\!\textup{qc}}\>\R Q_{\<\<X}^{}\kappa_*^{}$,
where $\iG Z^{\!\!\textup{qc}}\colon\Aqc(X)\to\AqcZ(X)$ is the restriction 
of~$\iG Z\>$.
Now by the proof of 
\cite[p.\,12, Proposition~(1.3)]{AJL}, $\R Q_{\<\<X}^{}$ is bounded on 
$\Dqc(X)\supset\kappa_*{}\Dqct(\X)$ (\Pref{Gammas'+kappas}). Also, by
\cite[p.\,24, Lemma (3.2.3)]{AJL}, $\R\iG Z$ is
bounded; and hence by \cite[p.\,26, Proposition (3.2.6)]{AJL}, so is 
$\iG Z^{\!\!\textup{qc}}\<$. Thus $\R\Qt_\X$ is bounded on~$\Dqct(\X)$.

In the general separated case, one proceeds by induction on the least number\- of
affine open subsets covering $\X$, as in 
the proof of \cite[p.\,12, Proposition~(1.3)]{AJL} (which is \Pref{1!} for
$\X$  an ordinary scheme), \emph{mutatis mutandis}---namely, substitute ``$\X\>$"
for ``$\<X\<$,\!" ``qct" for ``qc," ``$\Qt$" for ``$Q$,\!" and recall for 
a map~\mbox{$v\colon\V\to\X$} of noetherian formal schemes that
\mbox{$v_*(\Aqct(\V))\subset\Aqct(\X)$} (\Pref{f-*(qct)}), and
furthermore that if $v$ is affine then $v_*|_{\Aqct(\V)}$ is \emph{exact}
(Lemmas~\ref{Gamma'+qc} and~\ref{affine-maps}).

A similar procedure works when the Krull dimension $\dim\X$ is
finite, but now the induction is
on $n(\X)\set{}$least $n$ such that $\X$ has an open covering
$\X=\cup_{i=1}^n\U_i$ where for each~$i$ there is a  
separated ordinary noetherian scheme $U_i$ such that $\U_i$~is isomorphic
to the completion of $U_i$ along
one of its closed subschemes. (This property of $\U_i$ is inherited by any of its
open subsets). 

\newcommand{\vqct}{v_*^{\mathrm{qct}}}

The case $n(\X)=1$ has just been done. Consider, for any open
immersion $v\colon\V\hookrightarrow \X$, the functor 
$\vqct\set v_*|_{\Aqct(\V)}$.
To complete the induction as in the proof of 
\cite[p.\,12, Proposition~(1.3)]{AJL}, 
one needs to show that \emph{the derived functor
$\R \vqct\colon\D(\Aqct(\V))\to\D(\X)$
is bounded above.}

For $\>\N\in\Aqct(\V)$, let $\>\N\to\cJ^\bullet\>$ be an 
$\Aqct$-injective---hence flasque---resolution
\cite[Proposition 4.2]{Ye}. Now
$H^i\R\vqct(\N)$ is the sheafification of the presheaf sending~
an
open
$\W\subset\X$ to
$\textup{H}^i\Gamma(\W\cap\V,\> \cJ^\bullet)=\textup H^i(\W\cap\V,\>\N)$,
which vanishes when $i>\dim\X$, whence the conclusion (\cite[Proposition
(2.7.5)]{Derived categories}). 
\end{proof}

\end{parag}

\begin{parag}

Let $(X, \cO_{\<\<X})$ be a ringed space, and let $\cJ$ be an
$\cO_{\<\<X}$-ideal. The next Lemma, 
expressing $\R\iG{\<\cJ\>}$ as a ``homotopy colimit,"
\index{ $\iG{\raise.3ex\hbox{$\scriptscriptstyle{\ldots}$}}$ (torsion
functor)!$\R\iG{\raise.3ex\hbox{$\scriptscriptstyle{\ldots}$}}$ as homotopy
colimit|(}  lifts back to $\D(X)$ the well-known relation
$$
H^i\R\iG{\<\cJ\>}\G=
\dirlm{n}\text{\emph{$\E\<$xt}}^i_{\cO_X}(\cO_X/\<\cJ^n\<,\>\G)
\qquad\bigl(\G\in \D(X)\bigr).
$$

Define $\bh_n\colon\D(X)\to\D(X)$~by
$$
\bh_n(\G)\set\R\sHomb(\cO_{\<\<X}/\<\cJ^n\<,\>\G)\qquad\bigl(n\ge1,\
\G\in\D(X)\bigr).
$$
There are  natural functorial maps $s_n\colon \bh_n\to \bh_{n+1}$ and
$\varepsilon_n\colon \bh_n\to\R \iG{\<\cJ\>}$, satisfying
\mbox{$\varepsilon_{n+1}s_n=\varepsilon_n$.}
The family
$$
(1,-s_m)\colon \bh_m\to \bh_m\oplus \bh_{m+1} \subset \oplus_{n\ge1}\>\bh_n
\qquad(m\ge 1)
$$
defines a natural map $s\colon {\oplus_{n\ge1}\>\bh_n}\to\oplus_{n\ge1}\>\bh_n$.
There results,  for each $\G\in\D(X)$, a map of triangles
$$
\begin{CD}
\oplus_{n\ge1}\>\bh_n\>\G @>s>> \oplus_{n\ge1}\>\bh_n\>\G @>>> \textup{??} @>+>> \\
@VVV @VV \<\sum\!\varepsilon_n V @VV\overline\varepsilon V \\
0 @>>> \R \iG{\<\cJ\>}\G @= \R \iG{\<\cJ\>}\G @>+>> 
\end{CD}
$$

\begin{slem}
\label{Gam as holim}
The map\/ $\overline{\varepsilon}$ is a\/ $\D(X)$-isomorphism, and so 
we have a triangle
$$
\begin{CD}
\oplus_{n\ge1}\>\bh_n\>\G @>s>> \oplus_{n\ge1}\>\bh_n\>\G 
@>\sum\!\varepsilon_n>> \R\iG{\<\cJ\>}\G @>+>> 
\end{CD}
$$
\end{slem}

\begin{proof}
In the exact homology sequence
$$
\cdots\to H^i\bigl(\!\oplus_{n\ge1}\>\bh_n\>\G\bigr)
\overset{\sigma^i}{\lra}  H^i\bigl(\!\oplus_{n\ge1}\>\bh_n\>\G\bigr)
\lra H^i\bigl(\textup{??}\bigr)
\lra H^{i+1}\bigl(\!\oplus_{n\ge1}\>\bh_n\>\G\bigr)
\to\cdots
$$
the map $\sigma^i$ is injective, as can be verified stalkwise at each
$x\in X$.  Assuming, as one may, that $\G$ is K-injective, one deduces
that
$$
H^i(\textup{??}) = \dirlm{n}\! H^i(\bh_n\>\G)=H^i\>\>\dirlm{n}\! (\bh_n\>\G) 
= H^i\>\>\dirlm{n}\!\sHomb(\cO_{\<\<X}/\<\cJ^n\<,\>\G)=H^i(\R\iG{\<\cJ\>}\G),
$$
whence the assertion.
\end{proof}
\index{ $\iG{\raise.3ex\hbox{$\scriptscriptstyle{\ldots}$}}$ (torsion
functor)!$\R\iG{\raise.3ex\hbox{$\scriptscriptstyle{\ldots}$}}$ as homotopy
colimit|)} 
\end{parag}

\section{Duality for torsion sheaves.}
\label{S:t-duality}

Paragraph~\ref{T:qct-duality} contains the proof of \Tref{Th2} 
(section~\ref{S:prelim}), that is, of two
essentially equivalent forms of Torsion Duality%
\index{Grothendieck Duality!Torsion (global)} on formal
schemes---\Tref{T:qct-duality}  and~
\Cref{C:f*gam-duality}. The rest of the paragraph deals with numerous relations
among the functors which have been introduced, and with compatibilities among
dualizing functors occurring before and after completion of maps of ordinary
schemes.

More can be said for complexes with coherent homology, thanks to Greenlees-May
duality.\index{Greenlees-May Duality} This is done in
paragraph~\ref{coherent}.

\pagebreak[3]
Paragraph~\ref{SS:Gam-Lam} discusses additional relations involving 
$\R\iGp\X\colon\D(\X)\to\D(\X)$ and its right adjoint
$\R\sHomb(\R\iGp\X\cO_\X^{}\>,-)$ on a locally noetherian formal
scheme~$\X$.

\begin{thm}\label{T:qct-duality}
\textup{(a)} Let\/ $f\colon\<\X \to \Y$ be a  map of noetherian formal schemes.
Assume  that\/ $f$ is separated, or\/ $\X$ has finite Krull
dimension, or else restrict~to bounded-below complexes.
 Then the\/  \hbox{$\Delta$-functor\/}
\mbox{$\>\Rfs\colon\<\Dqct(\X)\<\xrightarrow{\!\ref{Rf-*(qct)}\>}\<\Dqct(\Y)
  \hookrightarrow\D(\Y)$}  has a right\/ $\Delta$-adjoint.   

In~fact there is a\/ bounded-below $\Delta$-functor
$\ft\colon\D(\Y)\to\Dqct(\X)$\vadjust{\kern.3pt}%
\index{ $\iG$@$f^{{}^{\>\ldots}}$ (right adjoint of
$\R f_{\<\<*}\cdots$)!$\ft\<\<$} and a map of\/ $\Delta$-functors 
$\tau_{\<\mathrm t}^\pd\colon\R f_{\!*} \ft\to {\bf 1}$%
\index{ {}$\tau$ (trace map)!$\tau_{\<\mathrm t}$} 
such that for all\/ $\G\in\Dqct(\X)$ and\/ $\F\in\D(\Y),$\ the composed map 
$($in\/ the derived category of abelian groups\/$)$
\begin{align*}
\R\Homb_\X(\G,\>\ft\<\F\>)
 &\xrightarrow{\mathrm{natural}\,}
  \R\Homb_\Y(\R f_{\!*}\G, \>\R f_{\!*}\ft\<\F\>) \\ 
&\xrightarrow{\:\mathrm{via}\ 
\tau_{\<\mathrm t}^{\phantom{.}}\;\mkern1.5mu}\R\Homb_\Y(\R f_{\!*}\G,\F\>)
\end{align*}
is an isomorphism. 

\smallskip
\textup{(b)} If\/ $g\colon \Y\to\Z$ is another such map  then there is a natural
isomorphism\/ 
$$
(gf)_{\mathrm t}^{\<\times}
\iso\ft\< g_{\mathrm t}^\times\<.
$$
\end{thm}

\noindent\emph{Proof.} Assertion (b) follows from (a),
which easily implies that $(gf)_{\mathrm t}^\times$ and 
$\ft\< g_{\mathrm t}^\times$ are both right-adjoint to the restriction
of~$\R(gf)_*=\R g_*\Rfs$ to~$\Dqct(\X)$.
 
As for (a), assuming first that $\X$ is separated or finite-dimensional, or that
only bounded-below complexes are considered, we can replace~$\Dqct(\X)$ by the
\emph{equivalent} category
$\D(\Aqct(\X))$ (\Pref{1!}).  The inclusion\mbox{
$k\colon\Aqct(\X)\hookrightarrow\Avc(\X)$} has the right adjoint $\iGp
\X$. (\mbox{$\iGp\X(\Avc(\X))\subset\Aqct(\X)$,} by
\Lref{Gamma'+qc} and~\Cref{C:vec-c is qc}.) 
So for all
$\Aqct(\X)$-complexes~$\G'$ and $\Avc(\X)$-complexes~$\F^{\>\prime}$ there~is a natural
isomorphism of abelian-group complexes
$$
\Homb_{\Aqct}(\G'\<,\>\iGp\X \F^{\>\prime}\>)
\iso
\Homb_{\Avc}(k\G'\<,\>\F^{\>\prime}\>).
$$
Note that if $\F^{\>\prime}$ is K-injective over $\Avc(\X)$ then $\iGp\X\F^{\>\prime}$ is K-injective
over~$\Aqct(\X)$, because $\iGp\X$ has an exact left adjoint. Combining this
isomorphism with the isomorphism (\ref{Deligne}) in the proof of
\Tref{prop-duality}, we can conclude just as in part~4 at the end of
that proof, with the functor~$\ft$ defined to be the composition
$$
\D(\Y)\xrightarrow{\,\rho\,}
\K_{\text{\textbf I}}(\Y)\xrightarrow{\C_\bullet\>} 
\K_{\text{\textbf I}}(\Avc(\X))\xrightarrow{\iGp\X\>}
\K_{\text{\textbf I}}(\Aqct(\X))
\xrightarrow{\text{natural}\,} \D(\Aqct(\X)).
$$

(We have in mind here simply that the natural
functor $\D(\Aqct(\X))\to\D(\Avc(\X))$ has a right adjoint.  That is easily seen
to be true once one knows the existence of K-injective resolutions
in~$\D(\Avc(\X))$; but we don't know how to prove the latter other than
by quoting the generalization to arbitrary Grothendieck categories
\cite[Theorem~2]{BR}, \cite[Theorem 5.4]{AJS}.  The preceding argument avoids this
issue. One could also apply Brown Representability\index{Brown Representability}
directly, as in the proof of
\Tref{Th1} described in the Introduction.)

\smallskip

Now suppose only that the map $f$ is separated.  If $\Y$ is separated 
then so is~$\X$, and the preceding argument holds. For arbitrary
noetherian~$\Y$ the existence of a bounded-below right adjoint for $\R
f_{\!*}\colon\Dqct(\X)\to\D(\Y)$ results then from the
following Mayer\kern.5pt-Vietoris pasting argument, by induction on the least
number of separated open subsets needed to cover~$\Y$. 
Finally,  to dispose of the assertion about the $\R\Homb\>$'s apply homology to
reduce it to $\ft$ being a right adjoint.

To reduce clutter, we will abuse  notation---but only in the rest of the proof
of \Tref{T:qct-duality}---by writing  ``$f^\times\>$" in place of~``$\ft\<$."

\begin{slem}\label{L:pasting}
Let\/ $f\colon\X\to\Y=\Y_1\cup\Y_2$ $(\Y_i$ open in\/~$\Y)$ be a  map
of formal schemes, with\/ $\X$ noetherian.  Consider the commutative~diagrams
$$
\begin{CD}
\X_{12}\!:=@.\;\X_1\cap \X_2 @>q_i>> \X_i @>x_i>> \X \\
@.@Vf_{12}VV @Vf_iVV @VV f\hbox to 0pt{\quad \qquad$(i=1,2)$\hss}V  \\
\Y_{12}\!:=@.\Y_1\cap \Y_2 @>>p_i> \Y_i @>>y_i> \Y
\end{CD}
$$
where\/ $\X_i\!:=f^{-1}\Y_i$\ and all the horizontal arrows represent inclusions.
Suppose that for\/ $i=1,2,12,$ the functor\/ 
$\R f_{i*}\colon\Dqct(\X_i)\to\D(\Y_i)$ has a right adjoint\/~$f_i^\times\<$.
Then  $\R f_{\!*}\colon\Dqct(\X)\to\D(\Y)$ has a right adjoint $f^\times;$\ and
with the inclusions\/ $y_{12}\!:=y_i\smcirc p_i\>,$\ $x_{12}\!:=x_i\smcirc
q_i\>,$\ there is for each\/ $\F\in\D(\Y)$ a natural\/ $\D(\X\>)$-triangle
$$
f^\times\<\F\to \R\>x_{1*}^{}f_1^\times y_1^*\F    \oplus
  \R\>x_{2*}^{}f_2^\times y_2^*\F\xrightarrow{\lambda_\F}
\R\>x_{12*}^{}f_{12}^\times y_{12}^*\F \to (f^\times\<\F\>)[1]\,.
$$
\end{slem}

\emph{Remark.}
If we expect $f^\times$ to exist, and the natural maps 
$x_i^*f^\times \to f_i^\times y_i^*$ to be isomorphisms, then there should
be such a triangle---the Mayer\kern.5pt-Vietoris triangle
\index{Mayer-Vietoris triangle} of~$f^\times\<\F$.
This suggests  we first define $\lambda_\F\>$, then let $f^\times\<\F$ be
the vertex of a triangle based on $\lambda_\F\>$, and verify~\dots

\smallskip

\begin{proof}
There are  natural maps
$$
\tau_1\colon\R f_{1*}^{}f_1^\times\to \mathbf1,\qquad
\tau_2\colon\R f_{2*}^{}f_2^\times\to \mathbf1,\qquad
\tau_{12}\colon\R f_{12*}^{}f_{12}^\times\to \mathbf1.
$$
For $i=1,2$, define the ``base-change" map 
$\beta_i\colon q_i^*\<f_i^\times\to f_{12}^\times\> p_i^*$ to be 
adjoint under \Tref{T:qct-duality} to the  map of functors
$$
\R f_{12*}^{}q_i^*f_i^{\times}\underset{\>\text{natural}\>}{\xrightarrow{\ \Iso\ }}
 p_i^*\R f_{i*}f_i^{\times} \xrightarrow{\;\tau_i\;} p_i^*.
$$
This $\beta_i$ corresponds  to a functorial map
$\beta_i'\colon f_i^\times\to \R q_{i*}f_{12}^\times\> p_i^*$,
from which we obtain a functorial map
$$
\nopagebreak
\R\>x_{i*}f_i^\times y_i^*
 \longrightarrow\R\>x_{i*} \R q_{i*}f_{12}^\times\> p_i^*y_i^*
   \iso\R\>x_{12*}^{}f_{12}^\times\> y_{12}^*\>\>,
$$
and hence a natural map, for any $\F\in\D(\Y)$:
$$
\check D^0(\F\>)\!:=
 \R\>x_{1*}^{}f_1^\times y_1^*\F    \oplus
  \R\>x_{2*}^{}f_2^\times y_2^*\F    \xrightarrow{\lambda_\F\>}
   \R\>x_{12*}^{}f_{12}^\times\> y_{12}^*\F=: \check D^1(\F\>)\>.
$$
Embed this map  in a triangle $\check D(\F\>)$, and denote the
 third vertex by $f^\times(\F\>)$:
$$
\check D(\F\>)\colon\ f^\times\<\F\to\check D^0(\F\>)\xrightarrow{\lambda_\F\>}
\check D^1(\F\>)\to (f^\times\<\F\>)[1]\>. 
$$
Since $\check D^0(\F )$ and $\check D^1(\F )$ are in $\Dqct(\X)$ (see
\Pref{Rf-*(qct)}), therefore so is $f^{\times} \F$
(\Cref{qct=plump}). 

This is the triangle  in \Lref{L:pasting}. 
Of course we must still show that this~$f^\times$ is functorial, and right-adjoint
to~
$\R f_{\!*}$. (Then by uniqueness of adjoints such a triangle will
exist no matter which right adjoint~$f^\times$ is used.)

\goodbreak

Let us next construct a map
$
\tau_{\<\F}^{}\colon \R
f_{\!*}f^\times\<\F\to\F\ (\F\in\D(\Y)).
$
Set
$$
\check C^0(\F\>)\set\R y_{1*}^{} y_1^*\F\oplus \R y_{2*}^{} y_2^*\F ,\qquad
\check C^1(\F\>)\set \R y_{12*}^{} y_{12}^*\F.
$$
We have then the Mayer\kern.5pt-Vietoris $\D(\Y)$-triangle
$$
\check C(\F\>)\colon\  \F \to \check C^0(\F\>)  \xrightarrow{\mu_\F\>}
\check C^1(\F\>)
\to \F\>[1],
$$
arising from the usual exact sequence (\v Cech resolution)
$$
0\to \F \to y_{1*}^{}y_1^*\F\oplus y_{2*}^{}y_2^*\F \to y_{12*}^{}y_{12}^*\F  \to 0,
$$
where $\F$ may be taken to be K-injective.
Checking commutativity of the following natural diagram is a purely
category-theoretic exercise (cf.~\cite[Lemma (4.8.1.2)]{Derived categories} :
$$
\CD
\quad\R f_{\!*} \check D^0(\F\>)\quad
 \rlap{$\overset{\R f_{\!*}\lambda_\F}{\hbox
   to 94pt{\rightarrowfill}}$}
   @.@.
  \R f_{\!*} \check D^1(\F\>)\\
@| @. @|  \\
\R f_{\!*}(\R\>x_{1*}^{}f_1^\times y_1^*\F  \oplus
            \R\>x_{2*}^{}f_2^\times y_2^*\F\>)
  @.@.
   \R f_{\!*}\>\R\>x_{12*}^{}f_{12}^\times\> y_{12}^*\F\\
@V\simeq VV @. @VV\simeq V  \\
\R y_{1*}^{}\R f_{1*}^{}f_1^\times y_1^*\F \oplus
  \R y_{2*}^{}\R f_{2*}^{}f_2^\times y_2^*\F
   @.\hbox to36pt{}@.
    \R y_{12*}^{}\R f_{12*}^{}f_{12}^\times\> y_{12}^*\F\\
@V \tau_{\<1}^{}\oplus\tau_2^{} VV @.@VV \tau_{\<12}^{} V \\
\R y_{1*}^{} y_1^*\F \oplus
  \R y_{2*}^{} y_2^*\F  @.@.
   \R y_{12*}^{} y_{12}^*\F \\
@| @. @|  \\
\quad\ \ \,\check C^0(\F\>)\quad\ \ \,
 \rlap{$\underset{\mu_\F}{\hbox
   to 94pt{\rightarrowfill}}$}
   @.@.
   \check C^1(\F\>)
\endCD
$$
This commutative diagram extends to a map $\check\tau_{\<\F}^{}$ of triangles:
$$
\begin{CD}
\R f_{\!*}f^\times\<\F@>>>\R f_{\!*}\check D^0(\F\>)
   @>>>\R f_{\!*} \check D^1(\F\>)  @>>>\R f_{\!*}f^\times\<\F[1] \\
@V\tau_{\<\<\F}^{} VV @VVV @VVV @VV\tau_{\<\<\F}^{}[1] V \\
\F@>>> \check C^0(\F\>) @>>> \check C^1(\F\>) @>>>\F\>[1]
\end{CD}
$$
 
The  map $\tau_{\<\F}^{}$ is not necessarily unique.  But the next
Lemma will show, for fixed~$\F\<$, that \emph{the pair 
$(f^{\times}\<\F, \>\tau_{\<\F}^{})$ represents the functor} 
$$
\Hom_{\D(\Y)} (\R f_{\!*} \E, \F\>)\qquad\bigl(\E\in\Dqct(\X)\bigr).
$$
 It follows formally that one can make $f^{\times}$ 
into a functor and $\tau\colon \R f_{\!*}f^{\times}\to \mathbf 1$ into a morphism
of functors in such a way that the pair $(f^\times\<, \tau)$ is a right adjoint for
$\R f_{\!*}\colon\Dqct(\X) \to \D(\Y)$ (cf.~\cite[p. 83, Corollary~2]{currante});
and that there is a unique  isomorphism of functors $\Theta\colon  f^\times
T_2\iso T_1f^\times$ (where
$T_1$ and~$T_2$ are the respective translations on $\Dqct(\X)$ and~$\D(\Y)$)
such that 
$(f^\times\<,\Theta)$  is a $\Delta$-functor $\Delta$-adjoint to~$\R f_{\!*}$
(cf.~\cite[Proposition (3.3.8)]{Derived categories}). That will complete the proof of
\Lref{L:pasting}.
\end{proof}

\begin{slem}\label{L:f^times}
For\/ $\E\in \Dqct(\X),$ and with\/ $f^\times\<\F,$  $\tau_{\<\F}^{}$ as above,  the
composition
$$
\Hom_{\Dqct(\X)}\<(\E,\,f^\times\<\F\>) \xrightarrow{\!\textup{natural}\,}
 \Hom_{\D(\Y)}\<(\R f_{\!*}\E,\>\R f_{\!*}f^\times\< \F\>)
\xrightarrow{\textup{via}\;\tau_{\<\<\F}^{}}
\Hom_{\D(\Y)}\<(\R f_{\!*}\E,\>\F\>)
$$
is an isomorphism.
\end{slem}

\begin{proof}
In the following diagram,  to save space we 
write  $H_\X$ for $\Hom_{\Dqct(\X)}$, 
$H_\Y$ for $\Hom_{\D(\Y)}$, and $f_{\!*}$ for $\R f_{\!*}\>$:
$$
\nopagebreak
\def\H#1,{H_\X\bigl(\E,\>#1\bigr)}
\def\h#1,{H_\Y\bigl(f_{\!*}\E,\>#1\bigr)}
\CD
\H  (\check D^0\F\>)[-1], 
 @>>> \h f_{\!*}\bigl((\check D^0\F\>)[-1]\bigr), 
  @>>>  \h (\check C^0\F\>)[-1], \\ 
@VVV @VVV @VVV\\
\H  (\check D^1\F\>)[-1], 
 @>>> \h f_{\!*}\bigl((\check D^1\F\>)[-1]\bigr), 
  @>>>  \h (\check C^1\F\>)[-1], \\
@VVV @VVV @VVV\\
\H f^\times \F, @>>> \h f_{\!*}f^\times \F, @>>>  \h \F,\\
@VVV @VVV @VVV\\
\H \check D^0\F, 
 @>>> \h f_{\!*}\check D^0\F, 
  @>>>  \h \check C^0\F, \\
@VVV @VVV @VVV\\
\H \check D^1\F, 
 @>>> \h f_{\!*}\check D^1\F, 
  @>>>  \h \check C^1\F,
\endCD
$$

The first column maps to the second via  functoriality of
$f_{\!*}\>$, and the second to the third via the above triangle map 
$\check \tau_{\<\F}^{}\>$; so the diagram commutes. The 
columns are exact \cite[p.\,23, Prop.\,1.1\,b)]{H1}, and thus  
if  each of the first two  and last two rows is shown to compose to an
isomorphism, then the same holds for the middle row, proving
\Lref{L:f^times}.

Let's look at the fourth  row. With notation as in \Lref{L:pasting} (and again,
with all the appropriate $\R$'s omitted),
we want the left column  in the following natural diagram to compose to an
isomorphism:
$$
\def\H#1,#2,#3,{H_{#1}(#2,\>#3\>)}
\begin{CD}
\H \X,\E,x_{i*}f_i^\times y_i^*\F, 
 @>\Iso>>
 \H \X_i,x_i^*\E,f_i^\times y_i^*\F, \\ 
@VVV  @VVV   \\ 
\H \Y, f_{\!*}\E,f_{\!*}x_{i*}f_i^\times y_i^*\F, @.  
 \H \Y_i,f_{i*}x_i^*\E, f_{i*}f_i^\times y_i^*\F, \\
@V\simeq VV @VV \simeq V  \\
\H \Y, f_{\!*}\E,y_{i*}f_{i*}f_i^\times y_i^*\F, @>\Iso>> 
 \H \Y_i,y_i^*f_{\!*}\E, f_{i*}f_i^\times y_i^*\F,  \\
@V\text{via }  \tau_i VV  @VV\text{via } \tau_i V  \\
\H \Y, f_{\!*}\E,y_{i*}y_i^*\F, @>\Iso>>
 \H \Y_i,y_i^*f_{\!*}\E, y_i^*\F, 
\end{CD}
$$
Here the horizontal arrows represent adjunction isomorphisms.  Checking that
the diagram commutes is a straightforward category-theoretic exercise. By
hypothesis, the right column composes to an isomorphism. Hence so does 
the left one.

\enlargethispage{-.2\baselineskip}
The argument for the fifth row is similar. Using the (easily checked) fact that 
the morphism $f_{\!*}\check D^0\to \check C^0$ is 
$\Delta$-functorial, we find that  the first row is, up to
isomorphism, the same as the fourth row with $\F[-1]$ in place of $\F$, so it too
composes to an isomorphism. Similarly, isomorphism for the second row follows
from that for the fifth.
\end{proof}

\penalty-2000

\begin{exams}\label{ft-example}
(1) Let $f\colon\X\to\Y$ be a
map of quasi-compact formal schemes with  $\X$
\emph{properly algebraic,} and  let $f^\times$ be the right
adjoint given by \Cref{cor-prop-duality}.
Using \Pref{Gamma'(qc)} we find then
that $\ft:=\R\iGp\X\smcirc \<f^\times$ is a right
adjoint for the restriction of $\>\R f_{\!*}$ to~$\Dqct(\X)$.\vadjust{\kern1.5pt}

(2) For a noetherian formal scheme~$\X$, \Tref{T:qct-duality} gives
a right adjoint~\mbox{$\mathbf1^{\<!}\set\mathbf1^{\!\times}_{\mathrm t}$} to the
inclusion
$\Dqct(\X)\hookrightarrow\D(\X)$. If $\G\in\wDqc(\X)$ (i.e.,
$\R\iGp\X\>\G\in\Dqct(\X)$, see \Dref{D:Dtilde}),
then the natural $\Dqct(\X)$-map 
$\R\iGp\X\>\G\to \mathbf1^{\<!}\G$ 
(corresponding to the natural
$\D(\X)$-map $\R\iGp\X\>\G\to \G$) is an
\emph{isomorphism,} see \Pref{Gamma'(qc)}.\vadjust{\kern1.5pt} 

(3) If $\X$ is \emph{separated} or if $\X$ is \emph{finite-dimensional,} then we
have the equivalence
$\bj^{\mathrm t}_{\<\<\X}\colon\D(\Aqct(\X)){{\mkern8mu\longrightarrow
\mkern-25.5mu{}^\approx\mkern17mu}}\Dqct(\X)$ of \Pref{1!}, and
we can take 
$\mathbf1^{\<!}\set\bj^{\mathrm t}_{\<\<\X}
\smcirc\R\Qt_{\<\X}$,  see \Cref{C:Qt} and
\Lref{L:j-gamma-eqvce}.\vadjust{\kern1.5pt}

(4) Let $f\colon\X\to\Y$ be a closed immersion of noetherian formal schemes
(see \cite[p.\,442]{GD}). The functor~$f_{\!*}\colon\A(\X)\to\A(\Y)$ is exact,
so
$\Rfs=f_{\!*}$. Let $\I$ be the kernel of the surjective map\vadjust{\kern.5pt}
$\cO_\Y\twoheadrightarrow f_{\!*}\cO_\X$\vspace{1pt} and let
$\overline{\Y\<}\>$  be the ringed space~$(\Y, \cO_\Y/\I)$, so
that
$f$~factors naturally as 
$\X\stackrel
{\vbox to 0pt{\vskip-4.5pt\hbox{$\scriptscriptstyle\bar{\! f}\>$}\vss}}
\to\overline{\Y\<}\>\stackrel
{\vbox to
0pt{\vskip-3pt\hbox{$\scriptscriptstyle i\,$}\vss}}\to\Y$,\vadjust{\kern1.2pt} the
map~$\,\>\bar{\<\!f}\>$ being flat. The\vspace{2.5pt} inverse isomorphisms
$
\A(\X)\,
\begin{minipage}[b][10pt][c]{22pt}
$$
\begin{CD}
_{_{\bar{\!f}_{\!\<*}^{}}}\\
\vspace{-31pt}\\
\xrightarrow{\ \ \ }\\
\vspace{-35.5pt}\\
\xleftarrow{\ \ \,}\\
\vspace{-30pt}\\
^{^{\,\bar {\!f}^{\mkern-.5mu*}}}\\
\vspace{-23.5pt}
\end{CD}
$$
\end{minipage}
\<\A(\Y)
$
extend to inverse isomorphisms 
$
\D(\X)\,
\begin{minipage}[b][10pt][c]{22pt}
$$
\begin{CD}
_{_{\bar{\!f}_{\!\<*}^{}}}\\
\vspace{-31pt}\\
\xrightarrow{\ \ \ }\\
\vspace{-35.5pt}\\
\xleftarrow{\ \ \,}\\
\vspace{-30pt}\\
^{^{\,\bar {\!f}^{\mkern-.5mu*}}}\\
\vspace{-23.5pt}
\end{CD}
$$
\end{minipage}
\<\D(\Y)
$
\medskip

The functor $\cH_\I\colon \A(\Y)\to\A(\>\overline{\<\Y\<\<}\>\>)$ 
defined by\vspace{.6pt}
$
\cH_\I(F\>)\set\sHom(\cO_\Y/\I,\>F\>)
$
has an exact left adjoint, namely
$i_*\colon\A(\>\overline{\<\Y\<\<}\>\>)\to\A(\Y)$, so $\cH_\I$ preserves
K-injectivity and
$\R\cH_\I$ is right-adjoint to 
$i_*\colon\D(\>\overline{\<\Y\<\<}\>\>)\to\D(\Y)$ (see proof
of~\Lref{L:j-gamma-eqvce}). Hence the functor
$f^\natural\colon\D(\Y)\to\D(\>\overline{\<\Y\<\<}\>\>)$ defined by
\begin{equation}\label{f^natl}
f^\natural(\F\>)\set\bar{f}^*\R\cH_\I(\F\>)= \bar{f}^*
\R\sHomb(\cO_\Y/\I,\>\F\>)\qquad\bigl(\F\in\D(\Y)\bigr)
\end{equation}\index{ $\iG$@$f^{{}^{\>\ldots}}$ (right adjoint of
$\R f_{\<\<*}\cdots$)!{}$f^\natural\<$} 
is right-adjoint to $f_{\!*}=i_*\bar f_{\!*}$, and 
$f_{\!*}\colon\Dqct(\X)\to\D(\Y)$ has the right adjoint%
\index{ $\iG$@$f^{{}^{\>\ldots}}$ (right adjoint of
$\R f_{\<\<*}\cdots$)!$\mathstrut \fs\<$}
$$
\ft\set f^!\set \mathbf1^{\<!}\smcirc f^\natural.
$$

We recall that  $\G\in\A(\X)$ is quasi-coherent iff 
$\bar f_{\!*}\G\in\Aqc(\>\overline{\<\Y\<\<}\>)$  iff $f_{\!*}\G\in\Aqc(\Y)$, see
\cite[p.\,115, (5.3.15), (5.3.13)]{GD}. Also, by looking at stalks (see
\S\ref{Gamma'1}) we find that \mbox{$f_{\!*}\G\in\At(\Y)\Rightarrow
\G\in\At(\X)$.} Hence  \Rref{R:Dtilde}(4) together with the isomorphism
$\Rfs\R\iGp\X\cong\R\iGp\Y\Rfs$ of \Cref{C:f* and Gamma}(d) yields
that $f^\natural\wDqcp(\Y)\subset\wDqcp(\X)$;
and given \Cref{qct=plump}, \Pref{P:Rhom} yields
$f^\natural\Dqct^+(\Y)\subset\Dqct^+(\X)$. Thus if
\mbox{$\F\in\wDqcp(\Y)$} then by~(2) above, 
$f^!\<\F\cong \R\iGp \X f^\natural\F\>$; and if $\F\in\Dqct^+(\Y)$
then\vspace{1.5pt} $f^!\<\F\cong  f^\natural\F$. 

(5) Let $f\colon\X\to\Y$ be  any map satisfying the hypotheses of~
\Tref{T:qct-duality}.  Let $\J\subset\cO_\X$ and $\I\subset\cO_\Y$ be ideals of
definition such that $\I\cO_\X\subset\J$, and let 
$$
X_n\set(\X,\cO_\X/\J^n)
 \xrightarrow{f_n^{}\>}(\Y,\cO_\Y/\I^n)=:Y_n\qquad(n>0)
$$ 
be the scheme-maps induced by $f\<$, so that each $f_n$ also satisfies the
hypotheses of~\kern.5pt\Tref{T:qct-duality}. As the target of the
functor $(f_n)_{\mathrm t}^{\<\times}$ is $\Dqct(X_n)=\Dqc(X_n)$,
we write $f_n^\times$ for $(f_n)_{\mathrm t}^{\<\times}$
(see~(1) above).
If  $\>j_n\colon X_n\hookrightarrow\X$ and
$i_n\colon Y_n\hookrightarrow\Y$ are the canonical closed immersions then
$fj_n=i_nf_n$, and so $\>j_n^!\ft\<=f_n^\times i_n^!$.

\vspace{.8pt}

The functor
$j_n^\natural\colon\D(\X)\to\D(X_{n})$ being  as
in~\eqref{f^natl}, we have, using~(4),
$$
\bh_n\G\set\R\sHomb(\cO_{\<\<X}/\<\J^n\<,\>\G)
=j_{n*}j_n^\natural\G\cong j_{n*}j_n^!\G
\qquad\bigl(\G\in\wDqcp(\X)\bigr). 
$$
Hence for $\G\set\ft\<\F\ (\F\in\D^+(\Y))$, \Lref{Gam as holim}  gives a%
\index{ $\iG$@$f^{{}^{\>\ldots}}$ (right adjoint of
$\R f_{\<\<*}\cdots$)!$\ft\<\<$!as homotopy colimit} ``homotopy colimit"
triangle\looseness=-1
$$
\oplus_{n\ge1}\,
  j_{n*}f_n^\times i_n^!\>\F 
\lra 
\oplus_{n\ge1}\,j_{n*}f_n^\times i_n^!\>\F  
\lra
\ft\<\F 
\overset{+}\lra
$$
\end{exams}

\goodbreak

Once again, $\wDqc(\X)\set(\R\iGp\X)^{-1}\Dqct(\X)$  (\Dref{D:Dtilde}).

\begin{scor}\label{C:f*gam-duality}
 \textup{(a)} Let\/ $f\colon\<\X \to \Y$ be a  map of noetherian formal schemes.
Suppose that\/ $f$ is separated or that\/ $\X$ has finite Krull
dimension, or else restrict~to bounded-below complexes. 
Let\/ $\BL_\X\colon\D(\X)\to\D(\X)$  be the bounded-below\/ $\Delta$-functor%
\index{ $\mathbf {La}$@$\BL$ (homology localization)}
$$
\BL_\X(-)\set\R\sHomb(\R\iGp\X\cO_\X^{}\>,-),
$$
and let\/ $\ush f\colon\D(\Y)\to\wDqc(\X)$  be the\/ $\Delta$-functor
$\ush f\!:=\BL_\X\ft$\index{ $\iG$@$f^{{}^{\>\ldots}}$ (right adjoint of
$\R f_{\<\<*}\cdots$)!$\ush f$}
\textup(see Example~\textup{\ref{R:Dtilde}(3)).}

 The functor\/ $\ush f$ is bounded-below, 
and is right-adjoint to 
$$
\Rfs\R\iGp\X\colon\wDqc(\X)\xrightarrow{\ref{Rf-*(qct)}\,}\Dqct(\Y)
\hookrightarrow\D(\Y).
$$
\textup(In particular with\/
$\bj\colon\D(\Avc(\X))\to\D(\X)$\index{ $\iG{\<\cJ\>}$@$\bj$} the natural
functor, the~functor\/
$$
\R f_{\!*}\>\R\iGp\X\>\bj\colon\D(\Avc(\X))\to\D(\Y)
$$
has the bounded-below right adjoint\/~$\R Q_\X^{} \ush f\<\<$---see
Proposition~\textup{\ref{A(vec-c)-A}.)}

In fact there is a map of\/ $\Delta$-functors%
\index{ {}$\tau$ (trace map)!$\ush\tau$} 
$$
\ush\tau\colon\R f_{\!*}\R\iGp\X \ush f\to {\bf 1}
$$ 
such that for all\/ $\G\in\wDqc(\X)$ and\/ $\F\in\D(\Y),$\ the~composed map 
\begin{align*}
\R\Homb_\X(\G,\>\ush f\F\>)
 &\xrightarrow{\mathrm{natural}\,}
  \R\Homb_\Y(\R f_{\!*}\R\iGp\X \G, \>\R f_{\!*}\R\iGp\X \ush f\F\>) \\ 
&\xrightarrow{\:\mathrm{via}\ 
\ush \tau}\R\Homb_\Y(\R f_{\!*}\R\iGp\X \G,\F\>)
\end{align*}
is an \emph{isomorphism.} 

\smallskip
\textup{(b)} If\/ $g\colon \Y\to\Z$ is another such map then there is a natural
isomorphism\/ 
$$
\ush{(gf)}\iso\ush f\< \ush g.
$$
\end{scor}

\begin{proof}
 (a) The functor 
$\BL_\X$
is bounded below because
$\R\iGp\X\cO_\X$ is locally
isomorphic to the bounded complex~$\cK_\infty^\bullet$  in the proof of
\Pref{Gamma'(qc)}(a), hence homologically bounded-above. Since
$\BL_\X$ is right-adjoint to $\R\iGp\X$ (see
\eqref{adj}), (a) follows directly from \Tref{T:qct-duality}.

\smallskip
(b) Propositions~\ref{Rf-*(qct)} and~\ref{Gamma'(qc)}(a) show that
for any $\G\in\Dqct(\X)$ we have \mbox{$\R\iGp\Y\>\Rfs\G\cong \Rfs\G$,} 
and hence the functors
$\ft\<\<\BL_\Y$ and~$\ft$ are both right-adjoint to~$\Rfs|_{\Dqct(\X)}$, so they
are isomorphic. Then \Tref{T:qct-duality}(b) yields
functorial isomorphisms
$$
\ush{(gf)}=\BL_\X(gf)_{\mathrm t}^{\<\times}\iso
\BL_\X\ft\<g_{\mathrm t}^\times\iso
\BL_\X\ft\<\BL_\Y g_{\mathrm t}^\times=
\ush f\<\ush g.
$$
\end{proof}

\medskip
Here are some ``identities" involving the dualizing functors $f^\times$
(\Tref{prop-duality}), $\ft$~(\Tref{T:qct-duality}), and 
$\ush f\!:=\BL_\X\ft$ (\Cref{C:f*gam-duality}). 

Note that $\BL_\X$ is
right-adjoint to~$\R\iGp\X$, see~\eqref{adj}.
Simple arguments show that the natural maps are isomorphisms 
$\BL_\X\iso\BL_\X\BL_\X\>$, $\R\iGp\X\iso\R\iGp\X\BL_\X\>$,  see~(b) and~(d)
in \Rref{R:Gamma-Lambda}(1).

\pagebreak[3]
\begin{scor}\label{C:identities}
With the notation of \Cref{C:f*gam-duality}\kern.5pt\textup{,} 
 
\smallskip
\textup{(a)} There are natural  isomorphisms
\begin{alignat*}{2}
\R\iGp \X\ush f&\iso\ft\<,\qquad & \ush f&\iso \BL_\X\ft\<,\\
\R\iGp\X\ft&\iso\ft\<,\qquad & \ush f&\iso\BL_\X\ush f\<.
\end{alignat*}

\smallskip
\textup{(b)} The natural
functorial maps\/ $\R\iGp\Y\to\mathbf 1\to\BL_\Y$ induce isomorphisms
\begin{gather*}
\ft\R\iGp\Y\iso\ft\iso\ft\<\<\BL_\Y, \\
\ush f\R\iGp\Y \iso\ush f \iso \ush f\<\<\BL_\Y.
\end{gather*}

\smallskip
\textup{(c)} There are natural pairs of maps
\begin{gather*}
\ft\xrightarrow{\alpha_1\>} \R\iGp\X\>\bj f^\times
 \xrightarrow{\alpha_2\>}\ft\<, \\
\ush f\xrightarrow{\beta_1\>} \BL_\X\>\bj f^\times
 \xrightarrow{\beta_2\>}\ush f\<, 
\end{gather*}
each of which composes to an identity map.  If\/ $\X$ is properly algebraic
then all of these maps are isomorphisms.
 
\smallskip
\textup{(d)} If\/ $f$ is\/ \emph{adic} then the  isomorphism\/\vspace{.5pt}
$\Rfs\>\R\iGp\X\>\bj\osi\R\iGp\Y\>\Rfs\>\bj$ in\/~\textup{\ref{C:f* and
Gamma}(d)} indu\-ces an isomorphism of the  right adjoints 
$(\<$see \Tref{prop-duality}\textup{,} \Pref{A(vec-c)-A}$)$
$$
f^\times\<\<\BL_\Y\iso\R Q_\X^{}\ush f\<.
$$
\end{scor}

\begin{proof}
(a) The second isomorphism (first row) is the
identity map. \Pref{Gamma'(qc)} yields the third. The first is the
composition 
$$
\R\iGp \X\ush f=\R\iGp \X\BL_\X\ft\iso\R\iGp \X\ft\iso\ft\<.
$$
The fourth is the composition
$$
\hskip21pt\ush f=\BL_\X\ft\iso\BL_\X\BL_\X\ft\iso\BL_\X\ush f\<.
$$

(b) The first isomorphism results from 
$\R\iGp\Y$ being  right adjoint to the
inclusion $\Dt(\Y)\hookrightarrow\D(\Y)$ (see \Pref{Gamma'(qc)}(c)). 
For the second, check that $\ft$ and~$\ft\<\<\BL_\Y$ are
both right-adjoint to  $\R f_{\!*}|_{\Dqct(\X)}\>\dots$
(Or, consider the composition
$\ft\iso\ft\R\iGp\Y\iso\ft\R\iGp\Y\BL_\Y\iso\ft\<\BL_\Y$.) 
Then apply $\BL_\X$ to the first row to get the second row.

\smallskip
(c) With $\boldsymbol k\colon\D(\Aqct(\X))\to\D(\Avc(\X))$ the natural functor,
let $$\alpha\colon\boldsymbol k\R\Qt_\X\ft\to f^\times$$ be adjoint to
$\R f_{\!*}\bj\boldsymbol k\R\Qt_\X\ft
 \stackrel{\ref{1!}}{=}\R f_{\!*}\ft
\xrightarrow{\tau_{\textup t}^{\phantom{.}}} \mathbf 1.$
By \Cref{C:Hom-Rgamma}, $\bj(\alpha)\colon\ft\to\bj \< f^\times$
factors naturally as 
$$
\ft\xrightarrow{\alpha_{\<1}^{}\>}\R\iGp\X\bj \< f^\times\to\bj \< f^\times.
$$
Let $\alpha_2$ be the map adjoint to the natural composition
$\R f_{\!*}\>\R\iGp\X\bj \< f^\times\to\R f_{\!*}\bj \< f^\times\to\mathbf 1$.
One checks that 
$\tau_{\textup t}\smcirc\R f_{\!*}(\alpha_2\alpha_1)=\tau_{\textup t}$
($\tau_{\textup t}$ as in \Tref{T:qct-duality}),
whence  $\alpha_2\alpha_1=\text{identity}$.

The pair $\beta_1\>,\>\beta_2$ is obtained from $\alpha_1\>,\>\alpha_2$
by application of the functor~$\BL_\X$---see \Cref{C:Hom-Rgamma}.
(Symmetrically, the pair $\alpha_1\>,\>\alpha_2$ is obtained from
$\beta_1\>,\>\beta_2$ by application of the functor~$\R\iGp\X$.)

When $\X$ is properly algebraic, the functor~$\bj$ is fully faithful
(\Cref{corollary}); and it follows that $\R\iGp\X\>\bj \< f^\times$
and~$\ft$ are both right-adjoint to $\R f_{\!*}|_{\Dqct(\X)}$.

\smallskip
(d) Straightforward.
\end{proof}

\pagebreak[3]
The next three corollaries deal with compatibilities between formal (local) and
ordinary (global) Grothendieck duality.

\begin{scor}\label{C:kappa-f^times-tors}
Let\/ $f_0\colon X\to Y$ be a  map of  noetherian ordinary schemes. Suppose either
that\/ $f_0$ is separated or that\/ $X$ is finite-dimensional, or else restrict to
bounded-below complexes. Let\/
$W\subset Y$ and\/~ $Z\subset f_0^{-1}W$  be closed subsets,
$\kappa_\Y^\pd\colon\Y=Y_{/W}\to Y$ and\/ 
$\kappa_\X^\pd\colon\X=X_{/Z}\to X$  the respective  completion
maps, and\/ $f\colon\X\to\Y$ the map induced by~$f_0.$ 
\vadjust{\penalty-750}
$$
\begin{CD}
\X@.:=X_{\</Z} @>\kappa_\X^\pd>> X  \\
@V f VV @. @VV f_0^{} V \\
\Y@.:=Y_{/W}@>>\vbox to
0pt{\vskip-1.3ex\hbox{$\scriptstyle\kappa_\Y^\pd$}\vss}> Y
\end{CD}
$$
With $f_{\<0}^{\<\times}\set(f_0^{})_{\textup t}^{\<\<\times}$ right-adjoint to
$\Rfs\colon\Dqc(X)\to\D(Y),$
let $\tau_{\<\mathrm t}'$ be the composition
$$
\R f_{\!*}\kappa_\X^*\R\iG{Z}f_{\<0}^{\<\times}\kappa_{\Y*}^\pd
\underset{\ref{C:kappa-f*t}\>}{\iso}
\kappa_\Y^*\R f_{\<0*}^{}\R\iG{Z}f_{\<0}^{\<\times}\kappa_{\Y*}^\pd
\longrightarrow 
\kappa_\Y^*\R f_{\<0*}^{}f_{\<0}^{\<\times}\kappa_{\Y*}^\pd
\longrightarrow 
\kappa_\Y^*\kappa_{\Y*}^\pd\longrightarrow\mathbf 1.
$$
Then for all $\E\in\Dqct(\X)$ and $\F\in\D(\Y),$ the composed map
\begin{align*}
\alpha(\E\<,\F\>)\colon\Hom_{\D(\X)}\<(\E\<,
\kappa_\X^*\R\iG{Z}f_{\<0}^{\<\times}\kappa_{\Y*}^\pd
\F\>) &\lra \Hom_{\D(\Y)}\<(\Rfs\E\<,
\>\Rfs\kappa_\X^*\R\iG{Z}f_{\<0}^{\<\times}\kappa_{\Y*}^\pd\F\>) \\
&\,\<\underset{\textup{via }\tau_{\<\mathrm t}'}{\lra}
\Hom_{\D(\Y)}\<(\Rfs\E\<,\>\F\>)
\end{align*}
\noindent is an isomorphism. Hence the map adjoint 
to~$\tau_{\<\mathrm t}'$ 
is an isomorphism of functors
$$
\kappa_\X^*\R\iG{Z}f_{\<0}^{\<\times}\kappa_{\Y*}^\pd 
\iso \ft\<.
$$
\end{scor}

\begin{proof}
For any $\E\in\Dqct(\X)$,  
set $\E_0\!:=\kappa_{\X*}^\pd\E\in\DqcZ(X)$ (\Pref{Gammas'+kappas}). 
\Pref{Gammas'+kappas} and 
\cite[p.\,7, Lemma (0.4.2)]{AJL} give natural isomorphisms\vspace{-3pt}
\begin{multline*}
 \Hom_{\D(\X)}\<(\E,\>\kappa_\X^*\R\iG Z\G)
\iso \Hom_{\D(X)}\<(\E_0\>,\>\R\iG Z\G)
\iso \Hom_{\D(X)}\<(\E_0\>,\>\G)\\
\bigr(\G\in\Dqc(X)\bigl).
\end{multline*}
\vspace{-15pt}

\noindent
(In other words, $\kappa_\X^*\R\iG Z\G=(\kappa_\X^\pd)_{\mathrm
t}^{\<\times}\G$.) One checks then that the map
$\alpha(\E\<,\F\>)$ factors as the sequence of natural isomorphisms
\begin{align*}
\Hom_{\D(\X)}\<(\E,\>\kappa_\X^*\R\iG{Z}f_{\<0}^{\<\times}
   \kappa_{\Y*}^\pd\F\>)
&\iso
 \Hom_{\D(X)}\<(\E_0\>,\>f_{\<0}^{\<\times}\kappa_{\Y*}^\pd\F\>) \\
&\iso
 \Hom_{\D(Y)}\<(\R f_{\<0*}^{}\E_0\>,\>\kappa_{\Y*}^\pd\F\>) \\
&\iso
 \Hom_{\D(\Y)}\<(\kappa_\Y^*\R f_{\<0*}^{}\E_0\>,\>\F\>) \\
&\iso
 \Hom_{\D(\Y)}\<(\R f_{\!*}\kappa_\X^*\E_0\>,\>\F\>)
              \qquad\textup{(\Cref{C:kappa-f*t})}\\ 
&\iso
 \Hom_{\D(\Y)}\<(\R f_{\!*}\E,\>\F\>). 
\end{align*}
\vskip-3.8ex
\end{proof}
\vskip1pt

\begin{scor}\label{C:kappa+duality}
With  hypotheses as in \Cref{C:kappa-f^times-tors}\kern.5pt\textup{:}

\textup{(a)} There are natural isomorphisms
\begin{align*}
\R\iGp\X\kappa_\X^*f_{\<0}^{\<\times}\kappa_{\Y*}^\pd
&= (\kappa_\X^\pd)_{\mathrm t}^{\<\times}\<f_{\<0}^{\<\times}
       \kappa_{\Y*}^\pd
 \iso \ft\<, \\
\BL_\X\kappa_\X^*f_{\<0}^{\<\times}\kappa_{\Y*}^\pd
&=\kappa_\X^{\textup{\texttt\#}}f_{\<0}^{\<\times}
\kappa_{\Y*}^\pd
\iso \ush f;
\end{align*}
and if $f_0$ is proper, $Y=\Spec(A)$  \($A$ adic\)$,$
$Z=f_0^{-1}W,$
then with $f^{\<\times}\!$  as in
\Cref{cor-prop-duality}\kern.5pt\textup{:}\vspace{-1.4ex}
$$
\kappa_\X^*f_{\<0}^{\<\times}\kappa_{\Y*}^\pd\iso f^\times\<.\vspace{4pt}
$$

\enlargethispage*{\baselineskip}
\textup{(b)} The functor 
$
f_{\<0,Z}^\times\!:=\R\iG Z f_{\<0}^{\<\times}\colon\D(Y)\to\DqcZ(X)
$ 
is right-adjoint to the functor\/ $\Rfs|_{\DqcZ(X)}\>;$\ and
there is an isomorphism\vspace{-1.5pt}
$$
\kappa_\X^*f_{\<0,Z}^\times\kappa_{\Y*}^\pd\iso \ft\<.
$$

\pagebreak[3]

\textup{(c)} If\/ $X$ is separated then, with notation as in
\Sref{SS:Dvc-and-Dqc}\textup{,} the functor 
$$
\ush{f_{\<0,Z}}\!:=
\bj_{\!X}^{}\R Q_{\<\<X}^{}\R\sHomb_X(\R\iG Z\cO_{\<\<X}\>,\>f_{\<0}^{\<\times}-)
 \colon\D(Y)\to\Dqc(X)
$$  
is right-adjoint to\/~$\R f_{\<0*}^{}\>\R\iG Z|_{\Dqc(X)}\>;$\ \vspace{.6pt}and
if\/
$\X$ is properly algebraic, so that we have the equivalence\/~
$\bj_{\!\X}\colon\D(\Avc(\X))\to\Dvc(\X)$ 
\(\kern-1pt \Cref{corollary}\kern.7pt\), there\vspace{.6pt} is an
isomorphism\vspace{-2pt}
$$
\kappa_\X^*\ush{f_{\<0,Z}}\kappa_{\Y*}^\pd\iso \bj_{\!\X}\R Q_\X^{}\ush f\<.
$$
\end{scor}

\pagebreak[3]
\begin{proof}
(a) The first isomorphism combines \Cref{C:kappa-f^times-tors} (in
proving which we noted that $\kappa_\X^*\R\iG Z\G
=(\kappa_\X^\pd)_{\mathrm t}^{\<\times}\G$ for $\G\in\Dqc(X)$) and
\Pref{Gammas'+kappas}.
 The second follows from 
$\ush f=\BL_\X\ft\<$. The third is
\Cref{C:kappa-f^times}.

(b) The first assertion is easily checked; and the isomorphism is given by
\Cref{C:kappa-f^times-tors}.

(c) When $X$ is separated, $\bj_{\!X}^{}$ is an equivalence \cite[p.\,12,
Proposition (1.3)]{AJL}, and then the first assertion is easily checked.  

From
\Cref{C:kappa-f^times-tors} and \Pref{Gammas'+kappas} we
get an isomorphism
$$
\R\iG{Z}f_{\<0}^{\<\times}\kappa_{\Y*}^\pd\iso 
\kappa_{\X*}^\pd\ft\<.
$$
As in \Cref{C:Hom-Rgamma}, the natural map is an isomorphism
$$
\R\sHomb_X(\R\iG Z\cO_{\<\<X}\>,\G)
 \iso \R\sHomb_X(\R\iG Z\cO_{\<\<X}\>, \R\iG Z\G)\qquad\bigl(\G\in\Dqc(X)\bigr).
$$
When $\X$ is properly algebraic, 
$\bj_{\!\X}^{}\R Q_\X^{}\cong \kappa_\X^*\>\bj_{\!X}^{}\R
Q_{\<\<X}^{}\kappa_{\X*}^{}$ 
(\Pref{A(vec-c)-A}). So then we have a sequence of
natural isomorphisms
\begin{align*}
\kappa_\X^*\ush{f_{\<0,Z}}\kappa_{\Y*}^\pd
& \ \>\raise.2ex\hbox{\EQAL{17}}\;\>
 \kappa_\X^*\>\bj_{\!X}^{}\R Q_{\<\<X}^{}\R\sHomb_X(\R\iG Z\cO_{\<\<X},
    \>f_{\<0}^{\<\times}\kappa_{\Y*}^\pd -) \\
&\iso
\kappa_\X^*\>\bj_{\!X}^{}\R Q_{\<\<X}^{}\R\sHomb_X(\R\iG Z\cO_{\<\<X},
    \>\R\iG Z f_{\<0}^{\<\times}\kappa_{\Y*}^\pd- ) \\
&\iso
\kappa_\X^*\>\bj_{\!X}^{}\R Q_{\<\<X}^{}\R\sHomb_X(\R\iG Z\cO_{\<\<X},
    \>\kappa_{\X*}^\pd\ft -) \\
&\iso
\kappa_\X^*\>\bj_{\!X}^{}\R Q_{\<\<X}^{}\kappa_{\X*}^\pd
  \R\sHomb_\X(\kappa_\X^*\R\iG Z\cO_{\<\<X},\>\ft -) \\
&\iso \bj_{\!\X}^{}\R Q_\X^{}\R\sHomb_\X(\R\iGp \X\cO_{\X}\>,\>\ft-)\\
& \ \>\raise.2ex\hbox{\EQAL{17}}\;\>\bj_{\!\X}^{}\R Q_\X^{}\ush f.
\end{align*}
\vskip-3.8ex
\end{proof}

The following instance of ``flat base change" will be needed in the proof of
the general base-change \Tref{Th3}.
\begin{scor}\label{C:compln+basechange}
In \Cref{C:kappa-f^times-tors}\textup{,} assume further that $Z=f_{\<0}^{-1}W\<$.
Then the natural map is an isomorphism 
$$
\R\iG Zf_{\<0}^{\<\times}\<\F\iso\R\iG
Zf_{\<0}^{\<\times}\<\kappa_{\Y*}\kappa_\Y^*\F
\qquad\bigl(\F\in\D(Y)\bigr),
$$
and so  there is a composed isomorphism
$$
\zeta\colon\R\iGp\X\kappa_\X^*f_{\<0}^{\<\times}\<\F
\underset{\textup{\ref{Gammas'+kappas}(c)}}\iso
\kappa_\X^*\R\iG Z f_{\<0}^{\<\times}\<\F
\iso
\kappa_\X^*\R\iG Z f_{\<0}^{\<\times}\<\kappa_{\Y*}\kappa_\Y^*\F
\underset{\textup{\ref{C:kappa+duality}(b)}}\iso
\ft\<\kappa_\Y^*\F.
$$
\end{scor}

\begin{proof}
First, $\R f_{\<0*}^{}(\DqcZ(X))\subset\Dqc{}_W(Y)$. For, by
\cite[Proposition~(3.9.2)]{Derived categories}, 
$\R f_{\<0*}^{}(\Dqc(X))\subset\Dqc(Y)$; and then the assertion follows from the
natural isomorphism of functors (from
$\Dqc(X)$ to $\Dqc(Y)$)
$\R\iG W\R  f_{\<0*}^{}\cong\R  f_{\<0*}^{}\R\iG {f^{-1}W}\>,$ because $\G\in\DqcZ(X)$
(resp.~$\H\in\Dqc{}_W(Y)$) iff $\R\iG Z\G\cong\G$ (resp.~$\R\iG W\H\cong\H$),
cf.~\Pref{Gamma'(qc)}(a) and its proof. (The said functorial
isomorphism arises from the corresponding one without the $\R$'s, since
 $\R  f_{\<0*}^{}$ preserves K-flabbiness, see \cite[5.12, 5.15(b), 6.4, 6.7]{Sp}.

Now \Cref{C:Gammas'+kappas} gives that the natural map
is an isomorphism
$$
\Hom_{\D(Y)}\<(\R f_{\<0*}^{}\E, \F\>)\iso
\Hom_{\D(Y)}\<(\R f_{\<0*}^{}\E,\kappa_{\Y*}\kappa_\Y^*\F\>)
\qquad\bigl(\E\in\DqcZ(X)\bigr),
$$
and the conclusion follows from the adjunction in
\Cref{C:kappa+duality}(b).
\end{proof}

\medskip
\begin{parag}\label{coherent}
The next Proposition is a special case of Greenlees-May Duality for formal
schemes\index{Greenlees-May Duality} (see \cite[Proposition 0.3.1]{AJL$'$}). It is
the key to many statements in this paper concerning complexes with coherent
homology.

\begin{sprop}
\label{formal-GM}
Let\/ $\X$ be a locally noetherian formal scheme, $\E\in\D(\X)$. 
Then for all\/ $\F\in\Dc(\X)$ the natural map\/ $\R\iGp{\X}\E\to \E$ induces an
isomorphism
$$
\R\sHomb(\E, \>\F\>) \iso \R\sHomb(\R\iGp{\X}\E,\>\F\>).
$$
\end{sprop}

\begin{proof} The canonical isomorphism (cf.~\eqref{adj})
$$
\R\sHomb(\R\iGp\X\E\<,\>\F\>)
 \iso\R\sHomb\bigl(\E,\>\R\sHomb(\R\iGp\X\cO_\X^{}\>,\>\F\>)\bigr)
$$
reduces the question to where $\E=\cO_{\X}\>$.
It suffices then---as in the proof of
\Cref{C:Hom-Rgamma}---that for affine~$\X=\Spf(A)$, the natural map
be an isomorphism
$$
\Hom_{\D(\X)}(\cO_{\X}\>, \>\F\>)
\iso \Hom_{\D(\X)}(\R\iGp{\X}\cO_{\X}\>, \>\F\>)
\qquad \bigl(\F \in \Dc(\X)\bigr).
$$

Let $I$ be an ideal of definition of the adic ring $A$, set
$Z\set\Supp(A/I)$, and let $\kappa\colon\X\to X\!:=\Spec(A)$ be the
completion map. The categorical equivalences in
\Pref{c-erator} and the isomorphism $\kappa^*\R\iGp
Z\cO_X\!\iso\!\R\iGp\X\cO_\X^{}$ in \Pref{Gammas'+kappas} make the
problem whether for all $F\in\Dc(X)$ (e.g.,
$F=\R Q\>\kappa_*^{}\F\set\bj_{\!X}^{}\R Q_{\<\<X}^{}\kappa_*^{}\F\>$) 
the~natural
map is an isomorphism
$$
\Hom_{\D(X)}(\cO_{\<\<X}, \>F)
\iso \Hom_{\D(X)}(\R\iG Z\cO_{\<\<X}, \>F).
$$

Now, the canonical functor
$\bj_{\!X}^{}\colon\D(\Aqc(X))\to \D(X)$ induces an equivalence of
categories
\mbox{$\D(\Aqc(X)){{\mkern6mu\longrightarrow
\mkern-25.5mu{}^\approx\mkern15mu}}\Dqc(X)$} (see beginning of
\S\ref{SS:Dvc-and-Dqc}), and so we may assume that $F$~is a K-flat 
quasi-coherent complex.
\Lref{L:j-gamma-eqvce} shows that
$\bj_{\!X}^{}\R Q_{\<\<X}^{}$ is right-adjoint to the inclusion
$\Dqc(X)\hookrightarrow\D(X)$. The natural
map 
$$
\R\sHomb(\cO_X, F)\to\R\sHomb(\R\iG Z\cO_X, F)\vspace{-5pt}
$$
factors then as\vspace{3pt}
\begin{equation}\label{nat}
\begin{aligned} 
\R\sHomb(\cO_X, F)=F
&\underset{\ref{c-erator}}\iso
\bj_{\!X}^{}\R Q_{\<\<X}^{}\kappa_*^{}\kappa^*F \\
&\,\lra\,\mathstrut\kappa_*^{}\kappa^*F\\
&\underset{\vbox to 0pt{\vss\hbox{$\scriptstyle\lambda$}\vskip1pt}}\iso
\,\inlm{n}F/(I\cO_X)^nF
\underset{\vbox to 0pt{\vss\hbox{$\scriptstyle\Phi$}\vskip1pt}}{\iso}
\R\sHomb(\R\iG Z\cO_X, F),
\end{aligned}
\end{equation}
where the map $\lambda$, obtained by applying $\kappa_*^{}$ to the natural map
from $\kappa^*F$ to the completion~$F_{/Z}$, 
is a $\D(X)$-isomorphism by \cite[p.\,6, Proposition~(0.4.1)]{AJL};
and $\Phi$~is the isomorphism $\Phi(F,\cO_X)$ of
\cite[\S2]{AJL}. (The fact that $\Phi$ is an isomorphism is essentially the main
result in \cite{AJL}.)  
Also, by adjointness, the natural map is an isomorphism
$$
\Hom_{\D(X)}(\cO_{\<\<X}, \>\bj_{\!X}^{}\R Q_{\<\<X}\kappa_*^{}\kappa^*F)
\iso
\Hom_{\D(X)}(\cO_{\<\<X},\>\kappa_*^{}\kappa^*F).
$$
Conclude now by applying the functor
$\textup H^0\R\Gamma(X,-)$ to~\eqref{nat}.
\end{proof}

\begin{scor}\label{C:coh-dual}
Let\/ $f\colon\X\to\Y$ be as in  \Cref{C:f*gam-duality}\textup{,}
and assume further that\/ $f$~is adic.
 Then for all\/ $\F\in\Dc(\Y)$ the map 
corresponding to the natural composition\/
$\R f_{\!*} \R\iGp\X\>\bj \< f^\times\<\F
\to \R f_{\!*} \bj \< f^\times\<\F
\to \F
$
\(see \Tref{prop-duality}\/\)
is an isomorphism
$$
f^\times\<\F\iso\R Q_\X^{}\ush f\<\F.
$$
\end{scor}

\begin{proof}
By ~\Pref{formal-GM},
$\F\cong\BL_\Y\>\F\set\R\sHomb(\R\iGp\Y\cO_\Y\>,\>\F\>)$; so this Corollary
is a special case of \Cref{C:identities}(d).
\end{proof}

\begin{scor}\label{C:completion-proper}
In \Cref{C:kappa-f^times-tors}\textup{,} suppose\/ $Y=\Spec(A)$ \($A$ adic\/\)
and\vadjust{\kern.5pt} that the the map\/~$f_0$ is
\emph{proper.} Then with the customary\vadjust{\kern.5pt}
notation\/~$f_{\<0}^!\>$ for\/
$f_{\<0}^{\<\times}$ we have, for any\/ $\F\in\Dc^+\<(\Y),$
a natural isomorphism
$$
 \kappa_\X^*f_{\<0}^!\kappa_{\Y*}^\pd\F\iso
\ush f\F \in\Dc^+\<(\X).
$$
\end{scor}

\begin{proof}
The natural  map $f_{\<0}^!\mkern1.5mu\bj_{\<Y}\<\R Q_{\<Y}\kappa_{\Y*}^\pd\to 
f_{\<0}^!\kappa_{\Y*}^\pd$
is an isomorphism of functors from $\D(\Y)$ to~$\Dqc(X)$, both 
being right-adjoint to~$\kappa_\Y^*\R f_{\<0*}^{}$.
\Pref{c-erator} gives 
$\bj_{\<Y}\R Q_{\<Y}\kappa_{\Y*}^\pd\F\in\Dc^+\<(Y)$; so by
\cite[p.\,396, Lemma~1]{f!},
$f_{\<0}^!\kappa_{\Y*}^\pd\F\in\Dc^+\<(X)$.%
\footnote{%
For $\G\in\Dc^{\scriptscriptstyle+}(Y)$ one has
$f_0^!\G\in\Dc^{\scriptscriptstyle+}(X)$: The question being local on $X$ one
reduces to where\vadjust{\kern.6pt} \emph{either} $X$ is a projective space
${\mathbf P}^n_Y$  and $f_0$ is  projection, so that~ 
{$f_0^!\G=f_0^*\G\otimes\Omega^n_{X/Y}[n]\in
\Dc^{\scriptscriptstyle+}(X)$}, \emph{or} $f_0$ is a closed immersion and
$f_{0*}f_0^!\G=\R\sHomb_Y(f_{0*}\cO_X,\F\>)\in\Dc^{\scriptscriptstyle+}(Y)$
\cite[p.\,92, Proposition~3.3]{H1} whence, again,
$f_0^!\G\in\Dc^{\scriptscriptstyle+}(X)$ \cite[p.~115, (5.3.13)]{GD}.%
\vadjust{\kern 3pt}%
}
Hence \Pref{formal-GM} and 
Corollary ~\ref{C:kappa+duality}(a) yield isomorphisms\looseness=-1
$$
\kappa_\X^*f_{\<0}^!\kappa_{\Y*}^\pd\F
\iso\R\sHomb(\R\iGp\X\cO_\X^{}\>,\>\kappa_\X^*f_{\<0}^!\kappa_{\Y*}^\pd\F\>)
=:\<\<\BL_\X\kappa_\X^*f_{\<0}^!\kappa_{\Y*}^\pd\F
\iso\ush f\F\<.
$$
\vskip-3.8ex
\end{proof}

\smallskip
\end{parag}
\begin{parag} \label{SS:Gam-Lam}
More relations, involving the functors $\R\iGp\X$ and~
$\BL_\X\set\R\sHomb(\R\iGp\X\cO_\X^{},-)$\index{ $\mathbf {La}$@$\BL$
(homology localization)} on a locally noetherian formal scheme~$\X$, will now be
summarized.

\begin{small}

\begin{srems}\label{R:Gamma-Lambda}
Let $\X$ be a locally noetherian formal scheme.

(1) The functor\index{ $\iG{\raise.3ex\hbox{$\scriptscriptstyle{\ldots}$}}$
(torsion functor)!a@$\BG\set\R\iGp\X$ (cohomology colocalization)}
$\BG\!:=\R\iGp\X\colon\D(\X)\to\D(\X)$ admits a natural map
\mbox{$\BG\xrightarrow{\gamma\>}\mathbf 1$}, which induces a functorial
isomorphism
\begin{equation} 
\Hom(\BG\E, \BG\F\>)\iso  \Hom(\BG\E,\F\>)\qquad\bigl(\E,\F\in\D(\X)),
\tag{A}
\end{equation}
see \Pref{Gamma'(qc)}(c). Moreover $\BG$ has a right adjoint, viz.~
$\BL\!:=\BL_\X$ (see~\eqref{adj}). 

The rest of (1)
consists of (well-known) formal consequences of these properties.

\smallskip
Since $\gamma$ is functorial, it holds  that
$\gamma(\F\>)\smcirc\gamma(\BG\F\>)
=\gamma(\F\>)\smcirc\BG(\gamma(\F\>))\colon\BG\BG\F\to\F\<$, so 
injectivity of the map in (A) (with $\E=\BG\F\>$) yields 
$
\gamma(\BG\F\>)=\BG(\gamma(\F\>))\colon\BG\BG\F\to\BG\F;
$
and one finds after
setting $\F=\BG\G$  in (A) that this functorial map is an \emph{isomorphism}
\begin{equation} 
\gamma(\BG)=\BG(\gamma)\colon\BG\BG\iso  \BG.
\tag{a}
\end{equation}

Conversely, given (a) one  can deduce that the
map in (A) is an isomorphism, whose inverse takes 
$\alpha\colon\BG\E\to\F\>$ to the composition
$\BG\E\iso\BG\BG\E\xrightarrow{\BG\alpha\,}\BG\F\>$.%
\footnote
{The \emph{idempotence} of~$\BG\<$, expressed by
(a) or (A),\vadjust{\kern1.3pt}  can be interpreted as follows.  

Set $\D\set\D(\X)$,
$\mathbf S\set\{\,\E\in\D\mid
\BG(\E)=0\,\}\<$,\vspace{.6pt} so that
$\BG$ factors uniquely as \smash{$\D
\overset{\vbox to0pt{\vss 
                                 \hbox{$\scriptstyle q\>$}
                                 \vskip-.35ex}
              }\to
\D/\mathbf S
\overset{\vbox to0pt{\vss 
                                 \hbox{$\scriptstyle\,\overline{\!\BG\<}\>\>$}
                                 \vskip-.35ex}
              }\to\D$}
where $q$ is the ``Verdier quotient" functor. Then
\emph{$\,\overline{\!\BG}$ is left-adjoint to~$q$}, so that $\mathbf
S\subset\D$ admits a ``Bousfield colocalization."\index{Bousfield
colocalization}  It follows from (c) and (d) below that  
$\mathbf S=\{\,\E\in\D\mid \BL(\E)=0\,\}$, and (b) below means that
\emph{the functor $\bar\BL\colon\D/\mathbf S\to\D$
 defined by $\BL=\bar\BL\smcirc q$
is right-adjoint to~$q\>$}; thus $\mathbf S\subset\D$ also admits a ``Bousfield
localization." 
And $\D/\mathbf S$ is equivalent, via $\,\overline{\!\BG}$ and $\bar\BL$
respectively, to the categories $\Dt\subset\D$
and~$\D\>\>\hat{}\>\subset\D$ introduced below---categories denoted by
$\mathbf S^\perp$ and ${}^\perp\mathbf S$ in~\mbox{\cite[Chapter
8]{TC}.}%
}
The  composed functorial map
$\lambda\colon\mathbf1\to\BL\BG\xrightarrow{\BL(\gamma)}\BL$
 induces an isomorphism
\begin{equation} 
 \Hom(\BL\E,\BL\F\>)\iso \Hom(\E, \BL\F\>) \qquad\bigl(\E,\F\in\D(\X)),
\tag{B}
\end{equation}
or equivalently (as above), $\lambda$ induces an isomorphism
\begin{equation} 
\lambda(\BL)=\BL(\lambda)\colon\BL\iso\BL\BL.
\tag{b}
\end{equation}

Moreover, the isomorphism (A) transforms via adjointness to an isomorphism
$$
\Hom(\E, \BL\BG\F\>)\iso  \Hom(\E,\BL\F\>)\qquad\bigl(\E,\F\in\D(\X)),
$$
whose meaning is that $\gamma$ induces an isomorphism
\begin{equation} 
\BL\BG\iso  \BL.
\tag{c}
\end{equation}
Similarly,  (B) means that $\lambda$ induces the conjugate isomorphism
\begin{equation} 
\BG\BL\osi  \BG.
\tag{d}
\end{equation}
Similarly,  that $\BL(\lambda(\F\>))$---or $\gamma(\BG(\E))$---is an
isomorphism (respectively that 
$\lambda(\BL(\F\>))$---or~$\BG(\gamma(\E))$---is an isomorphism) is equivalent
to the first (respectively the second) of the following maps (induced by
$\lambda$ and $\gamma$ respectively) being an isomorphism:
\begin{equation}
\Hom(\BG\E\<,\F\>)
\iso\Hom(\BG\E\<,\BL\>\F\>) 
\osi\Hom(\E\<,\BL\>\F\>).
\tag{AB}
\end{equation}

That (c) is an isomorphism also means that the functor~$\BL$ factors, via
$\BG\<$, through the essential image~$\Dt(\X)$ of~$\BG\<$ (i.e., the  full
subcategory~$\Dt(\X)$ whose objects are isomorphic to $\BG\E$ 
for some~$\E$); and similarly (d) being an isomorphism means that $\BG$ factors,
via~$\BL$, through the essential image $\D\>\>\hat{}\>(\X)$ of~$\BL$; and the
isomorphisms
\mbox{$\BG\BL\BG\cong\BG$} and $\BL\BG\BL\cong\BL$ deduced from~(a)--(d)
signify that $\BL$ and $\BG$ induce quasi-inverse equivalences between the
categories
$\Dt(\X)$ and~$\D\>\>\hat{}\>(\X)$.

\smallskip

(2) If $\X$ is properly algebraic, the natural functor 
$\bj\colon\D(\Avc(\X))\to\Dvc(\X)$ is an
\emph{equivalence,} and the inclusion
$\Dvc(\X)\hookrightarrow\D(\X)$ has a
right adjoint~$\mathbf Q:=\bj\R Q_\X^{}$ (\Cref{corollary}.) Then (easily
checked, given \Cref{C:vec-c is qc} and \Pref{Gamma'(qc)}) all
of~(1) holds with $\D$,~$\Dt\>$, and~$\BL$ replaced by $\Dvc\>$, $\Dqct\>$,
and~$\BLc\!:=\mathbf Q\BL$, respectively. 

\smallskip

(3) As in (1), $\BL$ induces an
equivalence from $\Dqct(\X)$ to $\Dqc\mspace{-11mu}\hat{}\mspace{10mu}(\X)$,
the essential image of $\BL|_{\Dqct(\X)}$---or, since $\BL\cong\BL\BG\<$, of
$\BL|_{\Dqc(\X)}$ (\Pref{Gamma'(qc)}).  So for\vadjust{\kern.7pt}
any  $f\colon\X\to\Y$ as in \Cref{C:identities}, the functor 
$$
\BL_\Y\R f_{\<\<*}\>\R\iGp\X\colon 
\Dqc\mspace{-11mu}\hat{}\mspace{10mu}(\X) \to
\Dqc\mspace{-11mu}\hat{}\mspace{10mu}(\Y)
$$
has the right adjoint $\BL_\X\ft\R\iGp\Y=\BL_\X\ft=\ush f\<$. There result two
``parallel" adjoint pseudofunctors
\cite[(3.6.7)(d)]{Derived categories} (where ``3.6.6" should be ``3.6.2"): 
$$
(\R f_{\<\<*}\>,\>\ft)\text{ (on $\Dqct$)\quad and\quad }
(\BL_\Y\R f_{\<\<*}\>\R\iGp{\X},\>\ush f\>)
\text{ (on $\Dqc\mspace{-11mu}\hat{}\mspace{10mu}\>$)}.
$$
Both of these correspond to the same adjoint pseudofunctor on the
quotient~$\Dqc/(\mathbf S\cap\Dqc)$, see footnote
under (1).

If $f$ is \emph{adic} then $\Rfs\BL_\X\cong\BL_\Y\Rfs$ (\Cref{C:f* and
Gamma}(c)), and so \Pref{Rf-*(qct)} gives that
$\Rfs(\Dqc\mspace{-11mu}\hat{}\mspace{10mu}\>(\X))\subset
\Dqc\mspace{-11mu}\hat{}\mspace{10mu}\>(\Y)$. Moreover, there are functorial
isomorphisms
$$
\BL_\Y\R f_{\<\<*}\>\R\iGp\X\BL_\X\cong
\R f_{\<\<*}\>\BL_\X\R\iGp\X\BL_\X\cong
\R f_{\<\<*}\>\BL_\X\>.
$$
Thus for adic $f\<$, $\BL_\Y\R f_{\<\<*}\>\R\iGp\X$ can be replaced above
by~$\Rfs$.\vspace{1pt}

When $f$ is \emph{proper} more can be said, see \Tref{T:properdual}.

\end{srems}
\end{small}
\end{parag}

\section{Flat base change.}
\label{sec-basechange}
\renewcommand{\theequation}{\thesth}

A \emph{fiber square}\index{fiber square} of adic 
formal schemes is a commutative diagram
$$
\begin{CD}
\V@>v>>\X \\
@VgVV @VVfV \\
\U@>>\vbox to 0pt{\vskip-1ex\hbox{$\scriptstyle u$}\vss}>\Y
\end{CD}
$$
such that the natural map is an \emph{isomorphism} 
$\V\iso\X\times_\Y\U$. If $\I$, $\J$, $\mathscr K$ are 
ideals of definition of $\Y$, $\X$, $\U$
respectively, then $\cL\set\J\cO_\V+\mathscr K\cO_\V$ is an ideal of definition
of~$\>\V$, and the scheme $V\set(\V,
\cO_\V/\cL)$ is the fiber product of the
$(\Y,\cO_\Y/\I)$-schemes $(\X,\cO_\X/\J)$ and $(\U,\cO_\U/\mathscr K)$, see
\cite[p.\,417, Proposition~(10.7.3)]{GD}. By
\cite[p.\,414, Corollaire~(10.6.4)]{GD}, if $V$ is locally
noetherian and the $\cO_V$-module~$\cL/\cL^2$ is of finite type
then $\V$ is locally noetherian. That happens whenever
$\X$, $\Y$ and $\U$ are locally noetherian and either $u$ or $f$ is of
pseudo\kern.6pt-finite type.\looseness=-1

Our goal is to prove \Tref{T:basechange} (=\:\Tref{Th3} of the
Introduction).
That is, given a fiber square  as above, with $\X$, $\Y$, $\U$ and
$\V$ noetherian, $f$ \emph{pseudo\kern.6pt-proper,}
and $u$ \emph{flat,} we want to establish a functorial isomorphism
$$
\beta_\F\colon\R\iGp\V\>v^*\<\ft\<\F \iso 
\gt\>\R\iGp\U u^*\F\ (\cong \gt\<u^*\<\F\>)
\qquad \bigl(\F\in\wDqcp(\Y)\bigr).
$$ 

Some consequences of this theorem will be given in \Sref{Consequences}.

In order to define $\beta_\F$ (\Dref{D:basechange})
we first need to set up a canonical isomorphism 
$\R\iGp\U u^*\Rfs\iso\R\iGp\U \R g_*v^*\<$. This is done
in \Pref{uf=gv}. (When $u$~is \emph{adic} as well as flat, $\R\iGp\U$
can be omitted.)

Our proof of \Tref{T:basechange} has the weakness that it \emph{assumes} the case
when $f$ is a proper map of noetherian ordinary schemes. As far as
we know, the published proofs of this latter result  make use of
finite-dimensionality hypotheses on the schemes involved  (see \cite[p.\, 392,
Thm.\,2]{f!}, \cite[p.\,383, Cor.\,3.4]{H1}), or projectivity hypotheses on~$f$
\cite[p.\,191, 5]{H1}). There is however an outline of a proof for the general
case, even without noetherian hypotheses, in
\cite{Non noetherian}---see Corollary 4.3 there.%
\footnote{
Details may eventually appear in \cite{Derived categories}.
It is quite possible
that the argument can be adapted to give a direct proof for formal
schemes too.}

\medskip
To begin with, here are several properties of formal-scheme maps
(see \S\ref{maptypes})
which propagate across fiber squares.

\begin{prop}
\label{P:basechange}
\textup{(a)} Let\/ $f\colon \X\to\Y$ and\/ $u:\U\to\Y$ be maps of
locally noetherian formal schemes, such that the fiber product\/
$\X\times_{\Y} \U$ is locally noetherian \textup(a condition which
holds, e.g., if either \/ $f$ or\/ $u$ is of pseudo-finite type, see
\cite[p.\,414, Corollaire (10.6.4)]{GD}\textup).  If\/
$f$ is\/ \emph{separated} \textup(resp.~\emph{affine,}
resp.~\emph{pseudo\kern.6pt-proper,} resp.~\emph{pseudo\kern.6pt-finite,}
resp.~\emph{of pseudo\kern.6pt-finite type,} 
resp.~\emph{adic}\textup)  then so is the projection\/ $\X\times_{\Y}\U\to\U$. 
\vspace{1pt}

\textup{(b)} With\/ $f\colon \X\to\Y$ and  $u:\U\to\Y$ as in\/ \textup{(a),}
assume either that\/  $u$ is adic or that\/ $f$ is of
pseudo\kern.6pt-finite type. 
If\/ $u$ is flat then so is the projection\/ 
$\X\times_{\Y}\U\to\X$.\vadjust{\kern1pt}

\textup{(c)} Let\/ $f\colon \X\to\Y$,  $u:\U\to\Y$ be maps of
locally noetherian formal schemes, with\/ $u$ flat and locally
over\/~$\Y$ the completion of a finite-type map of ordinary
schemes. Then\/ $\X\times_{\Y}\U$ is locally noetherian, and  
the projection\/ $\X\times_{\Y}\U\to\X$ is flat. 
\end{prop}

\begin{proof}
(a) The adicity assertion is obvious, and the rest follows from
corresponding assertions for the ordinary schemes obtained by
factoring out defining ideals.\vadjust{\kern1pt}

(b) It's enough to treat the case when $\Y$, $\X$, and $\U$
are the formal spectra, respectively, of noetherian adic rings
$(A,I)$, $(B,J)$ and $(C,K)$ such that $B$ and~$C$ are $A$-algebras
with \mbox{$J\supset IB$} and $K\supset IC$, and such that 
$B\, \widehat{\otimes}_{\<A} \>\>C$ is
noetherian  (since $\X\times_{\Y} \U$ is
locally noetherian, see \cite[p.\,414, Corollaire (10.6.5)]{GD}). 
By the following \Lref{(4.1.2)}, the problem is to show that if
$C$ is $A$-flat and \emph{either} $K=IC$ ($u$ adic), \emph{or} $B/\<J$ is a
finitely-generated $A$-algebra ($f$ of pseudo-finite type), 
then $B\, {\widehat \otimes}_{\<A}\>\>C$ is $B$-flat. 

The local criterion of
flatness \cite[p.\,98, \S5.2, Thm.\,1 and p.\,101, \S5.4, Prop.\,2]{Bou}
reduces the problem further to showing that for all $n>0$,
$
(B\, {\widehat\otimes}_{\<A}
\>\>C)/J^n(B\, {\widehat\otimes}_{\<A} \>\>C)$ is $(B/\<J^n)$-flat, i.e., that
$(B/\<J^n)\, {\widehat\otimes}_{\<A} \>\>C$
is $(B/\<J^n)$-flat.
But, $C$ being $A$-flat, if $K=IC$ then
$(B/\<J^n)\, {\widehat\otimes}_{\<A} \>\>C
=(B/\<J^n)\otimes_{A/I^n} (C/I^nC)
$ 
is clearly $B/\<J^n$-flat; while if $B/J$ is a finitely-generated
$A$-algebra, then
$(B/\<J^n)\>\otimes_A \>C$ is noetherian and $(B/\<J^n)$-flat,
whence so is its $K$-adic completion $(B/\<J^n)\, {\widehat\otimes}_{\<A} \>\>C$. 
\vadjust{\kern1pt}

(c) Proceeding as in the proof of (b), we may assume~$C$ to be the
$K'$-adic completion of a finite\kern.5pt-type $A$-algebra~$C'$ ($K'$
a $C'$-ideal). If $C$ is $A$-flat then by \cite[\S 5.4, Proposition
4]{Bou}, the localization $C'{}'\set C'[(1+K')^{-1}]$ is $A$-flat, so
the noetherian $B$-algebra $ B\otimes_{\<A} C'{}' $ is $B$-flat, as is
its (noetherian) completion $B\,\widehat {\otimes}_{\<A}\, C$.
\end{proof}

\begin{slem}
\label{(4.1.2)}
Let $\varphi:A \to C$ be a continuous homomorphism of noetherian adic
rings. Then $C$ is $A$-flat iff the corresponding map
$\Spf(\varphi)\colon \Spf(\<C)\to\Spf(A)$ is~flat, i.e., iff for each
open prime~$q\subset C,$ $C_{\{q\}}$ is
$A_{\{\varphi^{-1}q\}}\<$-\kern.5pt flat.
\end{slem}
\begin{proof} Recall that if $K$ is an ideal of definition of~$C$ and
$q\supset K$ is an open prime ideal in~$C$, then with $C\setminus q$ 
ordered by divisibility,
$$
C_{\{q\}}\set\cO_{\Spf(C)\<,\>q}=
\dirlm{{\displaystyle\mathstrut}\hbox to 0pt
{\hss$\scriptstyle f\in \>C\setminus q\;$\hss}}C_{\{f\}}\vspace{3pt}
$$
\vspace{1pt}%
where $C_{\{f\}}$ is the $K$-adic completion of the localization~$C_{\<f}\>$.

Now for each
$f\notin q$ and $n>0$ the canonical map\vspace{1pt} $C_{\<f}/\<K^n\<C_{\<f}
\to C_{\{f\}}/\<K^n\<C_{\{f\}}$~is bijective,\vadjust{\kern.7pt}
so the \smash{$\dirlm{}\!\!$} of these maps is an
isomorphism\vadjust{\kern.8pt} $C_q/\<K^n\<C_q\iso C_{\{q\}}/\<K^n\<C_{\{q\}}$,
whence so is
 the $K$-adic completion \smash{$\widehat{C_q}
\iso \widehat{C_{\{q\}}}$} of the canonical map $C_q\to C_{\{q\}}$. 
We can therefore apply \cite[\S5.4, Proposition 4]{Bou} twice to get
that $C_q$ is $A_{\varphi^{-1}q}$-flat iff $C_{\{q\}}$ is
$A_{\{\varphi^{-1}q\}}$-flat. So if $C$ is $A$-flat then
$\Spf(\varphi)$ is flat; and the converse holds because $C$ is
$A$-flat iff $C_m$ is $A_{\varphi^{-1}m}$-flat for every maximal
ideal~$m$ in~$C$, and every such $m$ is open since $C$ is complete.
\end{proof}

\begin{prop}
\label{uf=gv}

\textup{(a)} Consider a fiber square of noetherian formal schemes
$$
\begin{CD}
\V@>v>>\X \\
@VgVV @VVfV \\
\U@>>\vbox to 0pt{\vskip-1ex\hbox{$\scriptstyle u$}\vss}>\Y
\end{CD}
$$
with $u$ and $v$ flat.  Let 
$$
\psi_{\<\G}^{}\colon \R g_*\R\iGp\V\> v^*\<\G\to \R\iGp\U\> \R g_*v^*\<\G
\qquad\bigl(\G\in\Dqc(\X)\bigr)
$$
be the unique map
whose composition with the natural map 
$\R\iGp\U\> \R g_*v^*\<\G\to\R g_*v^*\<\G$
is the natural map
$\R g_*\R\iGp\V\> v^*\<\G\to\R g_* v^*\<\G$. \textup(The existence
of\/~$\psi_{\<\G}^{}$ is given by Propositions~\ref{Gamma'(qc)} and
~\ref{Rf-*(qct)}.\textup) Then for all\/ $\E\in\Dqct(\X),$
$\psi_\E$ is an \emph{isomorphism}.

In particular, if\/ $u$ \textup(hence $v)$ is \emph{adic} then\/
$\psi_{\E}$ can be identified with the identity map of\/ $ \R g_*v^*\<\E$.

\vspace{2pt}

\textup{(b)}  Let\/ $\X,$\ $\Y,$\ $\U$ be  noetherian formal schemes, let\/
$f\colon\X\to\Y$ and\/ 
$u\colon\U\to\Y$ be maps, with\/ $u$ flat, and assume further that one of the
following holds:\vadjust{\kern1pt}

\item[\hspace{2.87em}(i)] $u$ is adic, and\/ $\V\set\X\times_\Y\U$ is
noetherian,

\item[\hspace{2.6em}(ii)] $f$ is of pseudo\kern.6pt-finite type,

\item[\hspace{2.3em}(iii)] $u$ is locally the completion of 
a finite-type map of ordinary schemes;\vadjust{\kern1.5pt}

\noindent 
so that by \Pref{P:basechange} we have a fiber square as
in\/ \textup{(a)}. Let
$$
\theta_\G\colon u^*\>\Rfs\G\to \R g_*v^*\G\qquad\bigl (\G\in\D(\X)\bigr) 
$$
be  adjoint\vadjust{\kern.5pt} to the canonical map\/
$\Rfs \G\to\Rfs  \R v_* v^*\G=\R u_*\R g_*v^*\G$.\vadjust{\kern1pt}

Then for all\/ $\E\in\Dqct(\X),$\ the map
$\theta_{\!\E}'\set\R\iGp\U(\theta_{\<\E})$\index{ {}$\theta'\<$} is an
\emph{isomorphism}
$$\postdisplaypenalty10000
\theta_{\!\E}'\colon\R\iGp\U u^*\>\Rfs\E\iso\R\iGp\U\> \R g_*v^*\<\E.
$$

In particular, if\/ $u$ \textup(hence $v)$ is \emph{adic} then\/
$\theta_{\<\E}$ itself is an isomorphism.

\smallskip
\textup{(c)} Under the hypotheses of\/ \textup{(a)}\vspace{.6pt} resp.~\textup{(b),}
if\/
$f$ \textup(hence $g)$ is adic then\/ $\psi_\E$ resp.~$\theta_{\!\E}'$ is an
isomorphism for all\/ $\E\in\wDqc(\X)$ \textup(see
\Dref{D:Dtilde}\kern.5pt\textup{).}

\end{prop}

\begin{proof}
(a) Let $\J$ be an ideal of definition of~$\X$, and $\cK$ of~$\U$, so that
$\J\cO_\V+\cK\cO_\V$ is an ideal of definition of~$\V$. The obvious equality
 $\iG{\J\cO_\V+\cK\cO_\V}=\iG{\cK\cO_\V}\iG{\J\cO_\V}$, 
applied to K-injective $\cO_\V$-complexes, leads to a natural
functorial map 
$$
\R\iGp\V \underset{\ref{Gamma'1}}{\overset{\textup{def}}=}
\R\iG{\J\cO_\V+\cK\cO_\V}
\lra\R\iG{\cK \cO_\V}\R\iG{\J \cO_\V}
$$
which is an \emph{isomorphism,} as one checks locally via
~\cite[p.\,20, Corollary~(3.1.3)]{AJL}.  Also, there are natural
isomorphisms
$$
\R\iG{\J\cO_\V}v^*\<\E
\underset{\textup{\ref{P:f* and Gamma}(b)}}{\iso}
v^*\>\R\iGp\X\>\E=v^*\>\R\iG\J\E
\underset{\textup{\ref{Gamma'(qc)}(a)}}{\iso} v^*\<\E
\qquad\bigl(\E\in\Dqct(\X)\bigr).
$$
Thus the natural map $\R\iGp\V\to\R\iG{\cK \cO_\V}$ induces an
\emph{isomorphism}---the composition
$$
\R g_*\R\iGp\V\mkern.6mu v^*\<\E
\iso\R g_*\R\iG{\cK \cO_\V}\R\iG {\J \cO_\V}v^*\<\E
\iso\R g_*\R\iG{\cK \cO_\V} v^*\<\E.
$$

Since $(*)\<\colon$\!
$\<\R g_*\R\iG{\cK \cO_\V} v^*\<\E\cong\R g_*\R\iGp\V\mkern.6mu v^*\<\E
\<\in\<\Dt(\U)$ (Propositions~\ref{Gamma'(qc)} and ~\ref{Rf-*(qct)}) 
therefore we can imitate the proof of \Pref{P:f* and Gamma}(d)---\emph{without}
the boundedness imposed there on~$\G$, since that would be needed only to
get $(*)$---to see that
the map $\R g_*\R\iG{\cK \cO_\V} v^*\<\E\to\R g_*v^*\<\E$ 
induced by  $\R\iG{\cK \cO_\V}\to \mathbf1$
factors uniquely as
$$
\R g_*\R\iG{\cK \cO_\V} v^*\<\E\iso\R\iGp\U\> \R g_*v^*\<\E\lra\R g_*v^*\<\E,
$$
with the first map an isomorphism. 
It follows  that $\psi_\E$ is the composed isomorphism
$$
\R g_*\R\iGp\V\mkern.6mu v^*\<\E
\iso\R g_*\R\iG{\cK \cO_\V}\R\iG {\J \cO_\V}v^*\<\E
\iso\R g_*\R\iG{\cK \cO_\V} v^*\<\E
\iso\R\iGp\U\> \R g_*v^*\<\E.
$$

The last statement in (a) (for adic~$u$) results then from
\Cref{C:f* and Gamma}(b) and
Propositions~\ref{Rf-*(qct)} and \ref{Gamma'(qc)}(a).

\smallskip
(b) Once $\theta_{\!\E}'$ is shown to be an isomorphism,
the last statement in (b) (for adic~$u$)
follows from \Cref{C:f* and Gamma}(b),
and  Propositions~\ref{Rf-*(qct)} and \ref{Gamma'(qc)}(a). 

To show that $\theta_{\!\E}'$ is an isomorphism, it suffices to show
that the composition
$$
\psi_\E^{-1}\theta_{\!\E}'\colon\R\iGp\U u^*\>\Rfs\E
\to \R g_*\R\iGp\V\> v^*\<\E
\qquad\bigl(\E\in\Dqct(\X)\bigr).
$$
is an isomorphism. We use \Lref{Gam as holim} to
reduce the problem, as follows.
 
First, the functors $u^*\<$, $v^*\<$, $\R\iGp\U$ and $\R\iGp\V$ are bounded,
and commute with direct\- sums: for $u^*\<$ and $v^*\<$ that is clear, and
for $\R\iGp\U$ and $\R\iGp\V$ it holds because they can be realized
locally by tensoring with a bounded flat complex (see proof of
\Pref{Gamma'(qc)}).  Furthermore, \Lref{Gamma'+qc},
\Pref{Gamma'(qc)}, and~\Pref{P:Lf*-vc} show that
$\R\iGp\V\>v^*\>\Dqct(\X)\subset\Dqct(\V)$; and the functor $\R g_*$
(resp.~$\Rfs$) is bounded on, and commutes with direct sums in,
$\Dqct(\V)$ (resp.~$\Dqct(\X)$), see Propositions~\ref{Gamma'+qc},
~\ref{Rf_*bounded} and~\ref{P:coprod}.  Hence, standard
way-out reasoning allows us to assume that $\E\in\Dqct^+(\X)$.

Next, let $\J$ be
an ideal of definition of~$\X$, $X_{n}\ (n>0)$ the scheme
$(\X,\cO_\X/\<\J^n)$, and $j_n\colon X_{n}\hookrightarrow \X$ the
associated closed immersion. The functor
$j_{n*}\colon\A(X_{n})\to\A(\X)$ is exact, so it extends to a functor
$\D(X_{n})\to\D(\X)$. The functor
$j_n^\natural\colon\D(\X)\to\D(X_{n})$ being defined as
in~\eqref{f^natl}, we have
$$
\bh_n(\G)\set\R\sHomb(\cO_\X/\<\J^n\<,\>\G)
=j_{n*}j_n^\natural\G\qquad\bigl(\G\in\D(\X)\bigr). 
$$

If $\E\in\Dqct^+(\X)$ then $\E=\R\iG{\<\J}\<\E$
(\Pref{Gamma'(qc)}(a)), and, as noted just
after~\eqref{f^natl}, $j_n^\natural\E\in\Dqc(X_{n})$. 
Hence, from the triangle in
\Lref{Gam as holim} (with $\G$~replaced by an $\E\in\Dqct^+(\X)$)
we derive a diagram of triangles
$$
\minCDarrowwidth=18pt
\begin{CD}
\R\iGp\U u^*\>\Rfs(\oplus_{n\ge1}\>\bh_n\>\E) 
  @>>> \R\iGp\U u^*\>\Rfs(\oplus_{n\ge1}\>\bh_n\>\E) 
    @>>> \R\iGp\U u^*\>\Rfs(\R\iG{\<\J}\<\E) @>+>> \\
@V\simeq VV @V\simeq VV @VV\simeq V \\
\oplus_{n\ge1}\>\R\iGp\U u^*\>\Rfs \bh_n\>\E 
  @>>> \oplus_{n\ge1}\>\R\iGp\U u^*\>\Rfs \bh_n\>\E 
    @>>> \R\iGp\U u^*\>\Rfs\E @>+>> \\
@V\oplus V\psi_{\bh_{\<n}\E}^{-1}\>\theta_{\< \bh_{\<n}\E}'V 
 @V\oplus V\psi_{\bh_{\<n}\E}^{-1}\>\theta_{\<\bh_{\<n}\E}'V 
  @VV\psi_\E^{-1}\theta_{\!\E}'V \\
\oplus_{n\ge1}\R g_*\R\iGp\V\> v^*\<\bh_n\>\E 
  @>>> \oplus_{n\ge1}\R g_*\R\iGp\V\> v^*\<\bh_n\>\E 
    @>>> \R g_*\R\iGp\V\> v^*\<\E @>+>> \\
@V\simeq VV @V\simeq VV @VV\simeq V \\
\R g_*\R\iGp\V\> v^*(\oplus_{n\ge1}\>\bh_n\>\E)
  @>>> \R g_*\R\iGp\V v^*(\oplus_{n\ge1}\>\bh_n\>\E) 
    @>>> \R g_*\R\iGp\V\> v^*(\R\iG{\<\J}\<\E) @>+>>
\end{CD}
$$
From this diagram we see that if each 
$\psi_{\bh_n\E}^{-1}\>\theta_{\<\bh_n\E}'$ 
is an isomorphism, then so is~$\psi_\E^{-1}\>\theta_\E'$.
So we need only prove (b) when $\E=j_{n*}\>\F$ with
$\F\set j_n^\natural\E\in\Dqc(X_{n})$.\vspace{1pt}
Let us show  that in fact \emph{for any $n>0$ and any $\F\in\Dqc(X_{n})\<,$
$\theta_{\!j_{\<n\<*}\F}'$ is an isomorphism.}

\smallskip

The assertion~(b)  is local both on~$\Y$ and on $\U$. Indeed, for (b)
to hold 
it  suffices, for every diagram of fiber squares
$$
\CD
\V@<j'<<\V'@>v'>>\X' @>j>>\X \\
@VgVV@Vg'VV @VVf'V @VVfV\\
\U@<<i'<\U'@>>u'>\Y'@>>i>\Y
\endCD
$$
where $\Y'$ ranges over a base of open subsets of~$\Y$,
$\U'$ ranges over a base of open subsets of~$u^{-1}\Y'\<$, 
$u'$ is induced by~$u$, and  $i$, $i'$ are the inclusions, that
$i'{}^{\<*}\<\theta_{\!\E}'$ \mbox{$(=\theta_{\!\E}'|_{\U'})$} be an isomorphism. 
Now when $u$ is an open immersion,
$\theta_\G$ is an isomorphism for all $\G\in\D(\X\>)$. (One may
assume $\G$ to be K-injective and note that $v^*\<$, having the exact left
adjoint ``extension by zero," preserves K-injectivity, so that $\theta_\G$
becomes the usual isomorphism $u^*\<\<f_{\!*}\>\G\iso g_*v^*\G$). Thus
there are functorial isomorphisms
$i'{}^*\>\R g_*\iso\R g_*'\>j'{}^*$ and $i^*\>\Rfs\iso\R f_{\<\<*}'\>j^*$;
and similarly there is an isomorphism  
$i'{}^*\>\R\iGp\U\iso\R\iGp{\U'} i'{}^*\<$.
So it  suffices  that the composition
$$
i'{}^*\>\R\iGp\U u^*\>\Rfs \E\xrightarrow{i'{}^{\<\<*}\<\<\theta_{\!\E}'\>}
 i'{}^*\>\R\iGp\U\>\R g_*v^*\<\E\iso
\R\iGp{\U'} i'{}^*\>\R g_*v^*\<\E\iso
 \R\iGp{\U'} \R g_*'\>j'{}^*v^*\<\E
$$
be an isomorphism; and with a bit of patience one identifies this 
composition with
$$
 \R\iGp{\U'}u'{}^*i^*\>\Rfs \E\iso
 \R\iGp{\U'} u'{}^*\>\R f_{\<\<*}'\>j^*\<\E
  \xrightarrow{\theta_{\<j^{\<\<*}\!\E}'\>\>}
   \R\iGp{\U'}\R g_*'v'{}^*\<j^*\<\E,  
$$
thereby reducing to showing that $\theta_{\<j^{\<\<*}\<\<\E}'$ is an
isomorphism. Thus one may
assume that both $\Y$ and $\U$ are affine, say $\Y=\Spf(A)$ and
$\U=\Spf(C)$ with $C$ a flat $A$-algebra (\Lref{(4.1.2)}).

Suppose next that $\X$ and $\Y$ are ordinary schemes, so that 
$\Y=\Spec(A)$. In cases~(i) and (ii) of (b), set $C'{}'\set C$,
and in case (iii) let $C'{}'$ be as in the proof of part~(c) of
\Pref{P:basechange}. In any case, $C'{}'$ is $A$-flat, $C$
is the $K$-adic completion of~$\>C'{}'$ for some $C'{}'$-ideal~$K\<$,
$\>\X\times_\Y \Spf(C)$ is the $K$-adic completion of~$\X\times_\Y
\Spec(C'{}')$,
and we have a natural commutative diagram
$$
\begin{CD}
\X\times_\Y \Spf(C)@>v_2>>\X\times_\Y \Spec(C'{}')@>v_1>>\X  \\
@VgVV @Vg_1VV @VVfV \\
\Spf(C) @>>u_2> \Spec(C'{}') @>>u_1> \Y
\end{CD}
$$
With $\iGp{}$ denoting $\iGp{\Spf(C)}$,
$\theta_{\!\E}'=:\theta'(\E,f,u)$ factors naturally as the composition
$$
\R\iGp{}u_2^*u_1^*\>\Rfs\E 
\xrightarrow{\R\iGp{}u_2^*(\theta(\E\<\<,\>f\<,\>u_1))}
\R\iGp{}u_2^*\>\R g_{1\<*}^{}v_1^*\E
\xrightarrow{\theta'(v_1^*\<\E\<\<,\>\>g_1\<,\>\>u_2)}
\R\iGp{}v_2^*\>\R g_*v_2^*v_1^*\E.
$$
Here $\theta(\E,f,u_1)$ is an isomorphism because all the schemes
involved are ordinary schemes. (One argues as in \cite[p.\,111,
Prop.\,5.12]{H1}, using \cite[p.\,35, (6.7)]{AHK}; for a fussier
treatment see
\cite[Prop.\,(3.9.5)]{Derived categories}.) 
Also, $\theta'(v_1^*\E,g_1,u_2)$ is an isomorphism, in case (i) of (b)
since then $u_2$ and $v_2$ are identity maps, and in cases (ii) and
(iii) by \Cref{C:kappa-f*t'} since then
$\X\times_\Y\Spec(C'{}')$ is noetherian.  Thus:

\begin{slem}\label{L:ordinary}
\Pref{uf=gv} holds when\/ $\X$ and\/ $\Y$ are both ordinary
schemes.
\end{slem}

We will also need the following special case of \Pref{uf=gv}:

\begin{slem}\label{L:closed}
Let\/ $\I$ be an ideal of definition of\/ $\Y,$\  $Y_n$ the scheme
$(\Y,\cO_\Y/\I^n),$\ and 
$i_n\colon Y_n\hookrightarrow\Y$  the canonical closed immersion. Let
$u_n\colon Y_n\times_\Y\U\to Y_n$ and 
$p_n\colon Y_n\times_\Y\U\to\U$ be the projections \textup(so that\/ $u_n$
is flat and\/ $p_n$ is a closed immersion, see\/ \textup{\cite[p.\,442,
(10.14.5)(ii)]{GD})}. Then the natural map is an isomorphism
$$
u^*i_{n*}\G\iso p_{n*}u_n^*\G\qquad\bigl(\G\in\Dqc(Y_n)\bigr).
$$
\end{slem}

\begin{proof}
Since the functors $u^*\<$, $i_{n*\>}$, $ p_{n*\>}$, and $u_n^*$ are
all exact,\vspace{.5pt} we may assume that $\G$ is a quasi-coherent
$\cO_{Y_n}$-module; and since those functors commute with\vspace{1pt}
\smash{$\dirlm{}\!\!$} we may further assume $\G$ coherent, and then
refer to \cite[p.\,443, (10.14.6)]{GD}.
\end{proof}

Finally, for general noetherian formal schemes $\X$ and $\Y$,
and  $\I$ and $Y_n$ as above, let
$\J\supset\I\cO_\X$ be an ideal of definition of\/ $\X$, let $X_n$ be the scheme
$(\X,\cO_\X/\<\J^n),$\ and let $f_n\colon X_n\to Y_n$ be the map induced
by $f$.  Then for any $\F\in\Dqc(X_n)$, it holds that $\R f_{n*}\F\in\Dqc(Y_n)$. 
(See \Pref{Rf-*(qct)}---though the simpler case
$\F\in\Dqc^+(X_n)$ would do for proving \Pref{uf=gv}.)
Associated to the natural diagram\looseness=-1

\bigskip
\xdef\RestoreCatCode{\catcode`\noexpand\@=\the\catcode`\@}
\catcode`@=11
\expandafter\ifx\csname graph\endcsname\relax
\alloc@4\box\chardef\insc@unt\graph\fi
\expandafter\ifx\csname graphtemp\endcsname\relax
\alloc@1\dimen\dimendef\insc@unt\graphtemp\fi
\RestoreCatCode
\setbox\graph=\vtop{%
  \vbox to0pt{\hbox{%
    \special{pn 8}%
    \special{pn 2}%
    \special{ia 275 1889 140 140 0 6.28319}%
    \graphtemp=.6ex \advance\graphtemp by 1.889in
    \rlap{\kern 0.275in\lower\graphtemp\hbox to 0pt{\hss $\U$\hss}}%
    \special{ia 1740 1889 140 140 0 6.28319}%
    \graphtemp=.6ex \advance\graphtemp by 1.889in
    \rlap{\kern 1.74in\lower\graphtemp\hbox to 0pt{\hss $\Y$\hss}}%
    \special{ia 975 1389 140 140 0 6.28319}%
    \graphtemp=.6ex \advance\graphtemp by 1.389in
    \rlap{\kern 0.975in\lower\graphtemp \hbox to 0pt%
{\hss $Y_n\!\<\times_{\<\Y}\!\U$\hss}}%
    \special{ia 2440 1389 140 140 0 6.28319}%
    \graphtemp=.6ex \advance\graphtemp by 1.389in
    \rlap{\kern 2.44in\lower\graphtemp\hbox to 0pt{\hss $Y_n$\hss}}%
    \special{pa 0 780}%
    \special{pa 0 500}%
    \special{pa 550 500}%
    \special{pa 550 780}%
    \special{pa 0 780}%
    \special{ip}%
    \graphtemp=.6ex \advance\graphtemp by 0.64in
    \rlap{\kern 0.275in\lower\graphtemp\hbox to 0pt{\hss
$\V$\hss}}%
    \special{ia 1740 640 140 140 0 6.28319}%
    \graphtemp=.6ex \advance\graphtemp by 0.64in
    \rlap{\kern 1.74in\lower\graphtemp\hbox to 0pt{\hss $\X$\hss}}%
    \special{pa 700 260}%
    \special{pa 700 0}%
    \special{pa 1250 0}%
    \special{pa 1250 280}%
    \special{pa 700 280}%
    \special{ip}%
    \graphtemp=.6ex \advance\graphtemp by 0.14in
    \rlap{\kern 0.975in\lower\graphtemp\hbox to 0pt
{\hss$X_n\!\<\times_{\<\Y}\! \U\ $\hss}}%
    \special{ia 2440 140 140 140 0 6.28319}%
    \graphtemp=.6ex \advance\graphtemp by 0.14in
    \rlap{\kern 2.44in\lower\graphtemp\hbox to 0pt{\hss $X_n$\hss}}%
    \graphtemp=0.3\baselineskip \advance\graphtemp by .6ex
    \advance\graphtemp by 1.889in
    \rlap{\kern 1.007in\lower\graphtemp\hbox to 0pt{\hss $\scriptscriptstyle u$\hss}}%
    \special{pa 415 1889}%
    \special{pa 1600 1889}%
    \special{fp}%
    \special{pa 1600 1889}%
    \special{pa 1575 1914}%
    \special{fp}%
    \special{pa 1575 1864}%
    \special{pa 1600 1889}%
    \special{fp}%
    \graphtemp=-0.5\baselineskip \advance\graphtemp by .6ex
    \advance\graphtemp by 1.55in
    \rlap{\kern 1.405in\lower\graphtemp\hbox to 0pt{\hss $\scriptscriptstyle\ \quad
u_{\<n}$\hss}}%
    \special{pa 1250 1389}%
    \special{pa 1695 1389}%
    \special{fp}%
    \special{pa 1785 1389}%
    \special{pa 2300 1389}%
    \special{fp}%
    \special{pa 2300 1389}%
    \special{pa 2275 1414}%
    \special{fp}%
    \special{pa 2275 1364}%
    \special{pa 2300 1389}%
    \special{fp}%
    \graphtemp=-0.5\baselineskip \advance\graphtemp by .6ex
    \advance\graphtemp by 1.639in
    \rlap{\kern 0.6in\lower\graphtemp\hbox to 0pt
{\hss {$\scriptscriptstyle p_{n}\ $}\hss}}%
    \special{pa 389 1808}%
    \special{pa 837 1514}%
    \special{fp}%
    \special{pa 389 1808}%
    \special{pa 396 1773}%
    \special{fp}%
    \special{pa 389 1808}%
    \special{pa 424 1815}%
    \special{fp}%
    \graphtemp=0.5\baselineskip \advance\graphtemp by .6ex
    \advance\graphtemp by 1.639in
    \rlap{\kern 2.09in\lower\graphtemp\hbox to 0pt{\hss
\raise1ex\hbox{$\scriptscriptstyle\,\ \quad i_{\<n}$}\hss}}%
    \special{pa 1854 1808}%
    \special{pa 2326 1471}%
    \special{fp}%
    \special{pa 1854 1808}%
    \special{pa 1861 1773}%
    \special{fp}%
    \special{pa 1854 1808}%
    \special{pa 1889 1815}%
    \special{fp}%
    \graphtemp=-0.5\baselineskip \advance\graphtemp by .6ex
    \advance\graphtemp by 0.64in
    \rlap{\kern 1.075in\lower\graphtemp\hbox to 0pt{\hss $\scriptscriptstyle\qquad v$\hss}}%
    \special{pa 415 640}%
    \special{pa 1600 640}%
    \special{fp}%
    \special{pa 1600 640}%
    \special{pa 1575 665}%
    \special{fp}%
    \special{pa 1575 615}%
    \special{pa 1600 640}%
    \special{fp}%
    \graphtemp=-0.5\baselineskip \advance\graphtemp by .6ex
    \advance\graphtemp by 0.401in
    \rlap{\kern 0.564in\lower\graphtemp\hbox to 0pt{\hss {$\scriptscriptstyle\quad\! q_{n}\quad
$}\hss}}%
    \special{pa 389 560}%
    \special{pa 861 254}%
    \special{fp}%
    \special{pa 389 560}%
    \special{pa 396 525}%
    \special{fp}%
    \special{pa 389 560}%
    \special{pa 424 567}%
    \special{fp}%
    \graphtemp=-0.5\baselineskip \advance\graphtemp by .6ex
    \advance\graphtemp by 0.14in
    \rlap{\kern 1.775in\lower\graphtemp\hbox to 0pt{\hss $\scriptscriptstyle v_{\<n}$\hss}}%
    \special{pa 1250 140}%
    \special{pa 2300 140}%
    \special{fp}%
    \special{pa 2300 140}%
    \special{pa 2275 165}%
    \special{fp}%
    \special{pa 2275 115}%
    \special{pa 2300 140}%
    \special{fp}%
    \graphtemp=0.5\baselineskip \advance\graphtemp by .6ex
    \advance\graphtemp by 0.39in
    \rlap{\kern 2.09in\lower\graphtemp\hbox to 0pt{\hss
\raise1ex\hbox{$\scriptscriptstyle\,\ \quad j_{\<n}$}\hss}}%
    \special{pa 1854 558}%
    \special{pa 2326 221}%
    \special{fp}%
    \special{pa 1854 558}%
    \special{pa 1861 523}%
    \special{fp}%
    \special{pa 1854 558}%
    \special{pa 1889 565}%
    \special{fp}%
    \graphtemp=.6ex \advance\graphtemp by 1.264in
    \rlap{\kern 1.74in\lower\graphtemp\hbox to 0pt{\raise1.5ex\hbox{$\scriptscriptstyle \ f$}\hss}}%
    \special{pa 1740 780}%
    \special{pa 1740 1749}%
    \special{fp}%
    \special{pa 1740 1749}%
    \special{pa 1715 1724}%
    \special{fp}%
    \special{pa 1765 1724}%
    \special{pa 1740 1749}%
    \special{fp}%
    \graphtemp=.6ex \advance\graphtemp by 0.665in
    \rlap{\kern 2.44in\lower\graphtemp%
\hbox to 0pt{$\scriptscriptstyle \ f_{\<n}$\hss}}%
    \special{pa 2440 260}%
    \special{pa 2440 1249}%
    \special{fp}%
    \special{pa 2440 1249}%
    \special{pa 2415 1224}%
    \special{fp}%
    \special{pa 2465 1224}%
    \special{pa 2440 1249}%
    \special{fp}%
    \graphtemp=.6ex \advance\graphtemp by 1.264in
    \rlap{\kern 0.275in\lower\graphtemp%
\hbox to 0pt{\hss $\scriptscriptstyle g\ $}}%
    \special{pa 275 780}%
    \special{pa 275 1749}%
    \special{fp}%
    \special{pa 275 1749}%
    \special{pa 250 1724}%
    \special{fp}%
    \special{pa 300 1724}%
    \special{pa 275 1749}%
    \special{fp}%
    \special{pa 975 280}%
    \special{pa 975 590}%
    \special{fp}%
    \graphtemp=.6ex \advance\graphtemp by 0.97in
    \rlap{\kern 0.975in\lower\graphtemp\hbox to 0pt{\hss
\raise2.5ex\hbox{$\scriptscriptstyle g_{\<n}^\pd\; $}}}%
    \special{pa 975 690}%
    \special{pa 975 1249}%
    \special{fp}%
    \special{pa 975 1249}%
    \special{pa 950 1224}%
    \special{fp}%
    \special{pa 1000 1224}%
    \special{pa 975 1249}%
    \special{fp}%
    \kern 2.58in
  }\vss}%
  \kern 2.088in
}

\stepcounter{sth}
\renewcommand{\theequation}{\thesparag}

\vspace{45pt}
\begin{equation}\label{cube}
\end{equation}

\vspace{-95pt}

\centerline{\quad\ \box\graph}

\bigskip
\noindent
there is a composed isomorphism
\begin{align*}
\R\iGp\U u^*\>\Rfs j_{n*}\>\F
&\iso
\R\iGp\U u^*i_{n*}\R f_{\<\<n*}\>\F  
 &&\bigl (\F\in\Dqc(X_n)\bigr) \\
&\iso
\R\iGp\U \>\>p_{n*}u_n^*\>\R f_{\<\<n*}\>\F  
    &&(\textup{\Lref{L:closed})} \\
&\iso
\R\iGp\U\>\> p_{n*}\R g_{n*}v_n^*\F 
  && (\textup{\Lref{L:ordinary})}\\
&\iso
\R\iGp\U\>\R g_{*}q_{n*}v_n^*\F \\
&\iso
\R\iGp\U\>\R g_{*}v^*\<j_{n*}\>\F
  &&(\textup{\Lref{L:closed}),}
\end{align*}
which---the conscientious reader will verify---is just
$\theta_{\!j_{\<n\<*}\F}'\>$. 

Thus $\theta_{\!j_{\<n\<*}\F}'\>$ is indeed an isomorphism.

\smallskip
(c)  By definition $\R\iGp\X(\>\wDqc(\X))\subset \Dqct(\X)$, and so by (a) and~(b)
it's enough to see, as follows, that the natural map $\R\iGp\X\E\to\E$ induces
isomorphisms of the source and target of both $\psi_\E$ and $\theta_{\!\E}'\>$.

\Pref{P:f* and Gamma}(c) gives the isomorphism 
$\R g_*\R\iGp\V\> v^*\R\iGp\X\E \iso \R g_*\R\iGp\V\> v^*\<\E$,
as well as the second of the following isomorphisms, the first and
third of which follow from \Cref{C:f* and Gamma}(d):
$$
\R\iGp\U\> \R g_*v^*\R\iGp\X\>\E\cong
\R g_*\R\iGp\V\> v^*\R\iGp\X\>\E\cong
\R g_*\R\iGp\V\> v^*\E\cong
\R\iGp\U\> \R g_*v^*\E.
$$
Likewise, there are natural isomorphisms
$$
\R\iGp\U u^*\>\Rfs\R\iGp\X\E\cong
\R\iGp\U u^*\>\R\iGp\Y\>\Rfs\E\cong
\R\iGp\U u^*\>\Rfs\E.\vspace{-3.7ex}
$$
\end{proof}

\goodbreak
\smallskip

Notation and assumptions stay as in \Pref{uf=gv}(a). 
Assume that $f$ and~$g$ satisfy the hypotheses of \Tref{T:qct-duality},
so that the functor \mbox{$\Rfs\colon\Dqct(\X)\to\D(\Y)$}
has a right adjoint~$\ft\<$, and similarly for~$g$.
Recall from \Cref{C:identities}(b) that there is a natural isomorphism
$\gt \R\iGp\U\iso\gt$.

\begin{defi}\label{D:basechange}
\hskip-1pt With conditions as in \Pref{uf=gv}(b), the
\kern-.5pt\emph{base-change~map}
$$\index{base-change map}
\beta_\F\colon\R\iGp\V\> v^*\<\ft\<\<\F
\to\gt\>\R\iGp\U u^*\<\F
\qquad\bigl(\F\in\D(\Y)\bigr)
$$
is defined to be the map adjoint to the natural composition
$$
\R g_*\R\iGp\V\>\>v^*\<\<\ft\<\<\F
\underset{\psi}{\iso}
\R\iGp\U\>\R g_* v^*\<\<\ft\<\<\F
\underset{\theta'{}^{-1}}{\iso}
\R\iGp\U u^*\>\Rfs\ft\<\<\F
\to
\R\iGp\U u^*\<\F
$$
where $\psi\set\psi_{\<\ft\!\F}$ and $\theta'\set \theta'_{\!\ft\!\F}\>\>$.
In particular,  if $u$ (hence $v$) is \emph{adic} then $$
\beta_\F\colon v^*\<\<\ft\<\<\F\to \gt u^*\<\F
$$
is the map adjoint to the natural composition
$$
\R g_* v^*\<\<\ft\<\<\F
\underset{\theta^{-1}}{\iso}
u^*\>\Rfs\ft\<\<\F\to u^*\<\F
$$
where $\theta\set \theta_{\!\ft\!\F}\>\>$.
\end{defi}

\emph{Notation.} For a pseudo\kern.6pt-proper (hence separated) map~$f$ (see
\S\ref{maptypes}), we write~$f^!$\index{ $\iG$@$f^{{}^{\>\ldots}}$ (right adjoint of
$\R f_{\<\<*}\cdots$)!$\mathstrut \fs\<$} instead of
$\ft\<$.

\pagebreak[3]

\smallskip
\begin{thm}\label{T:basechange}\index{base-change isomorphism}
Let\/ $\X,$ $\Y$ and\/ $\U$ be noetherian formal schemes, let\/ $f\colon\X\to\Y$
be a pseudo\kern.6pt-proper map, and let\/ $u\colon \U\to\Y$ be flat, so that in any fiber
square
$$
\begin{CD}
\V@>v>>\X \\
@VgVV @VVfV \\
\U@>>\vbox to 0pt{\vskip-1ex\hbox{$\scriptstyle u$}\vss}>\Y
\end{CD}
$$
the formal scheme\/ $\V$ is noetherian, $g$ is pseudo\kern.6pt-proper, and $v$ is
flat \textup(\Pref{P:basechange}\kern.5pt\textup). Then for all\/
$\F\in\wDqcp(\Y)\set\wDqc(\Y)\cap\D^+(\Y)$ the base-change map\/ $\beta_\F$ is
an
\emph{isomorphism}
$$
\beta_\F\colon\R\iGp\V\>v^*\fs\F \iso 
\gs\>\>\R\iGp\U u^*\F\ (\cong \gs\<u^*\<\F\>).
$$
\end{thm}

\begin{small}
\emph{Remark.} In \cite[p.\,233, Example 6.5]{N1} Neeman\index{Neeman, Amnon}
gives an example where $f$ is a finite\- map of ordinary
schemes,  $u$ is an open immersion,  $\F\in\Dc^-(\Y)$, and 
$\beta_\F$  is \emph{not} an isomorphism.
\end{small}

\begin{proof}
Recall diagram \eqref{cube}, in which, 
$\I$ and $\J\supset\I\cO_\X$ being
defining ideals of $\Y$ and $\X$ respectively,
$Y_n$ is the scheme $(\Y,\cO_\Y/\I^n)$ and $X_n$ is the
scheme $(\X,\cO_\X/\J^n)$. Let
$\cK\supset\I\cO_\U$ be a defining ideal of~$\U$, let
$\cL\set\J\cO_\V+\cK\cO_\V\>$, a defining ideal of~$\V$, let
$V_n\ (n>0)$ be the scheme $(\V,\cO_\V/\cL^n)$, and let
$
l_n\colon V_n\hookrightarrow\V
$
be the canonical closed immersion. Then by \Eref{ft-example}(4), 
$$
l_{n*}l_n^!\G=l_{n*}l_n^\natural\G=\R\sHom(\cO_\V/\cL^n\<,\G\>)=:\,\bh_n(\G\>)
\qquad \bigl(\G\in\Dqct^+(\V)\bigr).
$$
So in view of the natural  isomorphism
$\R\iGp\V\gs u^*\F\<\iso\<\gs u^*\F$ 
(\Pref{Gamma'(qc)}(a)), \Lref{Gam as holim} shows
it sufficient to prove that the maps
$$
\bh_n(\beta_\F\>)\colon l_{n*}l_n^!\R\iGp\V v^*\<\<f^!\F\to
l_{n*}l_n^!g^!u^*\F
\qquad(n>0)
$$
are all isomorphisms.

\pagebreak[3]
Moreover, the closed immersion $l_n$ factors uniquely as
$$
V_n\xrightarrow{r_n\>\>} X_n\times_\Y\U\xrightarrow{q_n\>\>}\V,
$$so we can replace $l_n^!$ by~$r_n^!q_n^!$ (\Tref{T:qct-duality}(b)).
Thus \emph{it will suffice to prove that the maps
$$
q_n^!(\beta_\F\>)\colon q_n^!\R\iGp\V v^*\<\<f^!\F\to
q_n^!g^!u^*\F
\qquad(n>0)
$$
are all isomorphisms.}

\smallskip

In the cube \eqref{cube}, the front, top,
rear, and bottom faces are fiber squares,  denoted, respectively, 
by $\square$, $\square_{\textup t}\>$, $\square_{\textup r}\>$ and 
$\square_{\textup b}\>$; and we have the ``composed" fiber square
$\square_{\textup c}\>$:
$$
\begin{CD}
X_n\times_\Y\U @>v_n>> X_n \\
@Vp_ng_nV \!=\,gq_n V @Vi_nf_nV \!=\,fj_n V \\
\U @>>\vbox to 0pt{\vskip-1ex\hbox{$\scriptstyle u$}\vss}> \Y
\end{CD}
$$
The proper map  $f_n$
and the closed immersions~$i_n$ and~$j_n$ are all of pseudo-finite type. Also, it
follows from \Pref{P:basechange}(b) that in addition to~$u$, the maps
$u$, $u_n$, $v$ and~$v_n$ are all flat.  So corresponding to the  fibre
squares~$\square_\bullet$ we have base-change maps~$\beta_\bullet\>$.  

\goodbreak

Consider the following diagram of functorial  maps where, to save space, we set
$\blacktriangle\set X_n\!\times_\Y\<\U$ and 
$\blacktriangledown\set Y_n\!\times_\Y\<\U$.

$$
\minCDarrowwidth=15pt
\begin{CD}
q_n^!\R\iGp\V v^*\<\<f^!\<   @<\;\beta_{\textup t}\; <<
\<\R\iGp\blacktriangle v_n^*j_n^!f^!  @<\Is<<
\R\iGp\blacktriangle v_n^*(fj_n)^!  @=
\R\iGp\blacktriangle v_n^*(i_nf_n)^!\< @>\Is>>
\<\R\iGp\blacktriangle v_n^*f_n^!i_n^! \\
@V q_n^!(\beta) VV @. @V\beta_{\textup c} VV @. @VV\beta_{\textup r} V \\
 q_n^!g^!u^*  @>\Is>>
(gq_n)^!u^*  @=
 (p_ng_n)^!u^*  @<\Is<<
g_n^!p_n^!u^*  @<<\hbox to0pt{\hss$\scriptstyle g_n^!(\beta_{\textup b})$\hss}<
g_n^!\R\iGp{\blacktriangledown} u_n^*i_n^!
\end{CD}
$$
As above, we want to see that $q_n^!(\beta)$ is an isomorphism (in the category of
functors from $\wDqcp(\Y)$ to $\D(X_n\!\<\times_\Y\<\<\U)$). For that the
following assertions clearly suffice:\vspace{2pt}

(a) The preceding diagram commutes.

(b)  The base-change maps $\beta_{\textup t}$ and $\beta_{\textup b}$ are
isomorphisms.

(c) The base-change map $\beta_{\textup r}$ is an isomorphism.

\smallskip\noindent
Assertion (a) results from  part (b) of the transitivity lemma~\ref{L:trans}
below. 
Since $i_n$ and $j_n$ are closed immersions, 
assertion~(b) results from \Lref{L:immbc}, which is just
\Tref{T:basechange} for the case when $f$ is a closed immersion.
Since $f$~is pseudo\kern.6pt-proper therefore $f_n$ is proper, and assertion~(c)
is essentially the case of \Tref{T:basechange}---established in
\Lref{L:ordbc}---when $\X$ and $\Y$ are ordinary schemes.  

Thus these three Lemmas will complete the proof of
\Tref{T:basechange}.
\end{proof}

\begin{parag}\label{S:trans}

We will need some ``transitivity" properties\index{transitivity} of the maps
$\theta_{\!\E}'$ and
$\beta_\F$ relative to horizontal and vertical composition of fiber squares
of noetherian formal schemes,
i.e., diagrams of the form
\begin{subequations}\label{E:trans}
 \begin{equation}\label{E:transh}
  \begin{CD}
\V@>v_2>> \V_1 @>v_1>>\X \\
@VgVV @Vg_1VV @VVfV \\
\U @>>\vbox to 0pt{\vskip-1.3ex\hbox{$\scriptstyle u_2$}\vss}> \U_1 @>>\vbox to
0pt{\vskip-1.3ex\hbox{$\scriptstyle u_1$}\vss}> \Y
  \end{CD}
 \end{equation}

\vspace{10pt}

\begin{equation}\label{E:transv}
  \begin{CD}
\V @>v>> \X \\
@Vg_2 VV @VVf_2V \\
\W @>w>>\Z \\
@Vg_1 VV @VVf_1V \\
\U @>>\vbox to 0pt{\vskip-1.3ex\hbox{$\scriptstyle u$}\vss}>\Y
  \end{CD}
 \end{equation}
\end{subequations}
where all squares are fiber squares, and the maps $u$, $u_i\>$, $v$, $v_i\>$,
and~$w$ are all flat.

As we will be dealing with several fiber squares simultaneously we will
indicate the square with which, for instance, the map $\theta_\G$ in
\Pref{uf=gv} is associated, by writing
$\theta_{\<\<f\<,u}(\G)$ instead. 

The transitivity properties\index{transitivity} begin with:

\begin{slem}\label{L:transtheta}
Coming out of the fiber square diagrams\/~\eqref{E:transh} and\/
~\eqref{E:transv}\textup, the following natural diagrams commute for all
$\G\in\D(\X)\colon$
$$
\begin{CD}
(u_1u_2)^*\>\Rfs\G
 @.\!\!\overset{\theta_{\<\<f\<,\>u_{\mkern-1.5mu 1}^{}\!\<u_{\mkern-1.5mu
2}^{}}(\G)}{\Rarrow{158pt}}@.
   \R g_*(v_1v_2)^*\G\\
@V\simeq VV @. @VV\simeq V \\
u_2^*u_1^*\>\Rfs\G
 @>>u_2^*(\theta_{\<\<f\<,\>u_{\<1}^{}}(\G))>
  u_{2}^*\>\R g_{1\<*}^{}v_1^*\G
   @>>\theta_{\<g_{\<1}^{}\<,\>u_{\<2}^{}}(v_1^*\G)>
    \R g_*v_2^*v_1^*\G
\end{CD}
$$

\vspace{10pt}

$$
\begin{CD}
u^*\>\R (f_1f_2)_*\G
 @.\quad\,\overset{\theta_{\<\<f_1\<f_2,\>u}(\G)}{\Rarrow{189pt}}\ @.
   \R (g_1g_2)_{\<*}\>v^*\G\\
@V\simeq VV @. @VV\simeq V \\
u^*\>\R f_{1\<*}\>\R f_{\<2*}^{}\G
 @>>\theta_{\<\<f_{1}\<,\>u}(\R f_{\<2*}^{}\G)>
 \R g_{1\<*}^{}\>w^*\>\R f_{\<2*}^{}\G
   @>>\R g_{1\<*}^{}(\theta_{\<\<f_2,\>w}(\G))>
    \R g_{1\<*}^{}\>\R g_{2*}^{}v^*\G
\end{CD}
$$
\end{slem}

\begin{proof}
This is a formal exercise, based on adjointness of $u^*$ and $\R
u_*\>$, etc.  Details are left to the reader.
\end{proof}

\begin{slem}\label{L:trans}\index{transitivity}
\textup{(a)} In the fiber square diagram\/~\eqref{E:transh}\vspace{.6pt} \(with\/
$u_1\>,$
$v_1\>,$ $u_2$~and $v_2$ flat\), let\/ $\F\in\D(\Y)$ be such that the maps\/
$\theta_{\<\<1}'\set\theta'_{\<\<f\<,\>u_1^{}}(\ft\<\F\>),$
\mbox{$\theta_{\<2}'\set
  \theta'_{\<g_1,u_2^{}}((g_1^{}\<)_{\textup t}^{\!\times}\<u_1^*\F\>)$} 
and\/ $\theta_{\<2}'{}\!'\set
\theta'_{\<g_1,u_2^{}}(\R\iGp{\V_1}\<v_1^*\ft\<\F\>)$
of \Pref{uf=gv} are isomorphisms.
Then the map\/ $\theta'\set\theta'_{\<\<f\<,\>u_1^{}\<u_2^{}}(\ft\<\F\>)$
is an isomorphism, so the base-change maps\/
$\beta_1\set\beta_{\<\<f\<,\>u_1^{}}(\F\>),$
$\beta_2\set\beta_{g_1^{},\>u_2^{}}(u_1^*\F\>)$ and\/
$\beta\set\beta_{\<\<f\<,\>u_1^{}\<u_2^{}}(\F\>)$ can all be defined as in
\Dref{D:basechange};\vspace{1pt} and the following natural diagram, all of whose
uparrows are isomorphisms, commutes$\>:$
$$
\begin{CD}
\R \iGp\V(v_1v_2)^*\<\ft\<\F
   @.\hskip-15pt\overset{\beta}{\Rarrow{180pt}}   @.
    \gt\<\R\iGp\U(u_1u_2)^*\F\\
@A\simeq AA @.  @AA\simeq A \\
\R \iGp\V v_2^*v_1^*\<\ft\<\F 
 @.  
  \R\iGp\V \>v_2^*(g_1^{}\<)_{\textup t}^{\!\times} u_1^*\F 
   @>\beta_2>>
     \gt\<\R\iGp\U u_2^*u_1^*\F\\
@A\simeq A\textup{\ref{P:f* and Gamma}(c)} A 
 @A\simeq A\textup{\ref{C:identities}(b)}A 
   @A\textup{\ref{P:f* and Gamma}(c)}A\simeq A \\
\R\iGp\V \>v_2^*\>\R\iGp{\V_1}\<v_1^*\ft\<\F
 @>> \R\iGp\V \>v_2^*(\beta_1)>
 \R\iGp\V \>v_2^*(g_1^{}\<)_{\textup t}^{\!\times} \R\iGp{\U_1}u_1^*\F
   @>>\beta_2>
     \gt\<\R\iGp\U u_2^*\R\iGp{\U_1}u_1^*\F
\end{CD}
$$

\textup{(b)} In the fiber square diagram\/~\eqref{E:transv}---where\/ $u,$
$v$ and~$w$ are assumed flat---set\/ \mbox{$f\set f_{\<1}f_2$} and\/ $g\set
g_1g_2$. Let\/ $\F\in\D(\Y)$ be such that the maps\/
$\theta_{\<\<1}'\set\theta'_{\<\<f_{\<1}\<,\>u}
     ((f_{\<1}^{}\<)_{\textup t}^{\!\times}\F\>),$
\mbox{$\theta_{\<2}'\set
  \theta'_{\<\<f_2,w}(\ft\<\<\F\>)$} 
and\/
$\theta'\set\theta'_{\<\<f\<,\>u}(\ft\<\<\F\>)$
of \Pref{uf=gv} are isomorphisms, 
so that the base-change maps\/
$\beta_1\set\beta_{\<\<f_{\<1},\>u}(\F\>),$
$\beta_2\set\beta_{\<\<f_2,w}((f_{\<1}^{}\<)_{\textup t}^{\!\times}\F\>)$ and\/
$\beta\set\beta_{\<\<f,\>u}(\F\>)$ are all defined. Then 
the following diagram, whose two uparrows are
isomorphisms, commutes\/\textup{:}
$$
\begin{CD}
\R\iGp\V v^*\< \ft\<\F
 @.
   \overset{\beta}{\hskip9.8pt\Rarrow{167pt}}
    @.
        \gt\R\iGp\U u^*\F
      \\
@A\simeq AA @. @AA\simeq A
      \\
\R\iGp\V v^* (f_2^{})\<_{\textup t}^{\!\times}\<(f_{\<1}^{}\<)\<_{\textup
t}^{\<\times}\F
 @>>\beta_2>
   (g_2^{})\<_{\textup t}^{\<\times}\R\iGp\W \>w^*\<(f_{\<1}^{}\<)\<_{\textup
t}^{\<\times}\F
     @>> (g_2^{})\<_{\textup t}^{\<\times}\<(\beta_1) >
       (g_2^{})\<_{\textup t}^{\<\times}\<(g_1^{})\<_{\textup t}^{\<\times} \R\iGp\U
u^*\<\F 
\end{CD}
$$
\end{slem}

\pagebreak[3]
\def\ftt{f^{\<\times}}

\begin{proof} 
(a)  The map 
$$
\gamma\set \R\iGp\U u_2^*(\theta_{\<f,u_{\<1}^{}}\!(\ft\<\F\>))\colon
 \R\iGp\U u_2^*u_1^*\Rfs \ft\<\F\lra  
\R\iGp\U u_2^*\R g_{1\<*}^{}v_1^*\<\ft\<\F
$$
 is isomorphic, by \Pref{P:f* and Gamma}(c), to 
$$
\R\iGp\U u_2^*(\theta_{\<\<1}')\colon
 \R\iGp\U u_2^*\R\iGp{\U_1}\<u_1^*\Rfs \ft\<\F\lra  
\R\iGp\U u_2^*\R\iGp{\U_1}\<\R g_{1\<*}^{}v_1^*\<\ft\<\F\<,
$$ 
and so is an isomorphism (since $\theta_{\<\<1}'$ is).  

The map 
$$
\theta_{\<g_1^{}\<\<,u_2^{}}'(v_1^*\ft\<\F\>)\colon
 \R \iGp \U u_2^*\>\R g_{1\<*}^{}v_1^*\<\ft\<\F \to
   \R \iGp\U \>\R g_* v_2^*v_1^*\<\ft\<\F
$$
is also an isomorphism, as it is isomorphic to
$$
\theta_{\<2}'{}\!'\colon
 \R \iGp \U u_2^*\>\R g_{1\<*}^{}\R\iGp{\V_1}\<v_1^*\<\ft\<\F \to
   \R \iGp\U \>\R g_* v_2^*\R\iGp{\V_1}\<v_1^*\<\ft\<\F\<,
$$
because the natural map 
$\R \iGp \U u_2^*\>\R g_{1\<*}^{}\R\iGp{\V_1}\<v_1^*\<\ft\<\F \to
\R \iGp \U u_2^*\>\R g_{1\<*}^{}v_1^*\<\ft\<\F$ is
the composed isomorphism
\begin{multline*}
\smash{\R \iGp \U\< u_2^*\>\R g_{1\<*}^{}\R\iGp{\V_1}\<v_1^*\<\ft\<\F
\xrightarrow[\!\!\R \iGp \U u_2^*
\psi_{\! f_{\<\<\mathrm t}^{\!\times}\!\F}\!\!]{\Iso}
\R \iGp \U u_2^*\>\R\iGp{\U_1}\<\R
g_{1\<*}^{}v_1^*\<\ft\<\F}{}_{\displaystyle\mathstrut}
\\
\xrightarrow[\!\textup{\ref{P:f* and Gamma}(c)}\!]{\Iso}
\R \iGp \U u_2^*\>\R g_{1\<*}^{}v_1^*\<\ft\<\F{}
\end{multline*}

\smallskip
\noindent(see \Pref{uf=gv}(a)); and because
$\R \iGp\U \>\R g_*v_2^*\R\iGp{\V_1}\<v_1^*\<\ft\<\F \to
 \R \iGp\U \>\R g_* v_2^*v_1^*\<\ft\<\F$ is one of the maps in the
commutative diagram (B) below, all of whose other maps are isomorphisms.

Thus in the next diagram, 
whose commutativity results easily from that of the first diagram in
\Lref{L:transtheta},  all the maps other than
$\theta'$ are isomorphisms, whence so is $\theta'$.
$$
\begin{CD}
 \R\iGp\U \>\R g_*  (v_1v_2)^{\<*}\<\ft\<\F 
  @<\mkern48mu\theta'\mkern48mu<<
    \R \iGp \U u_2^*u_1^*\Rfs \ft\<\F
                   \\
@V\simeq V\hbox{\hskip67.5pt{\footnotesize(A)}}V 
    @V\simeq V\gamma V 
                    \\
\R \iGp\U \>\R g_* v_2^*v_1^*\<\ft\<\F
    @<\Iso<\theta_{\<\<g_1^{}\<\<,u_2^{}}'\<\<(\<v_1^*\ft\!\F \>)<
    \R \iGp \U u_2^*\>\R g_{1\<*}^{}v_1^*\<\ft\<\F 
\end{CD}
$$

Now it suffices to show that the diagram which is
\emph{adjoint} to the diagram in~(a) without its southeast (bottom right)
corner, commutes. That adjoint diagram is the outer border of the following one,
where, to reduce clutter, we omit all occurrences of the symbols~$\R$ and $\F$, 
write
$\ftt$ for
$\ft\<$, etc., and leave some obvious maps unlabeled:
$$
\minCDarrowwidth=21.5pt
\begin{CD}
g_*\iGp\V\> (v_1v_2)^{\<*}\<\<\ftt
 @>\psi>>
  \iGp\U \>g_*  (v_1v_2)^*\<\<\ftt
    @>\mkern33mu\theta'{}^{-1}\mkern33mu>>
     \iGp \U u_2^*u_1^*f_{\<\<*} \ftt
       @>>>
         { \iGp \U u_2^*u_1^* }
             \\
@A AA
  @A A\hbox{\hskip54.5pt{\footnotesize(A)}}A 
    @A \hbox to 0pt{$\scriptstyle\hskip7pt
          \gamma^{-1}
         \hss$} 
      A \hskip76pt\raise4.5ex
           \UnderElement{\hbox{\footnotesize(C)\hskip22pt}}{\|}{50pt}{} 
      A      \\
g_*\iGp\V\>  v_2^*v_1^*\<\ftt\hbox to0pt{\hskip17pt {\footnotesize(B)}\hss}
 @.
  \iGp\U \>g_* v_2^*v_1^*\<\ftt
    @<<\theta_{\<g_1^{}\<\<,u_2^{}}'(v_1^*\ftt\<)<
     \iGp \U u_2^*\>g_{1\<*}^{\phantom{.}}v_1^*\<\ftt
              \\
 \vspace{-22pt}
               \\
@A\simeq A\textup{\ref{P:f* and Gamma}(c)}A  @AAA  @AAA
              \\     
g_*\iGp\V\>  v_2^*\iGp{\V_1}\<\<v_1^*\<\ftt
 @>>\psi>
  \iGp\U \>g_* v_2^*\iGp{\V_1}\<\<v_1^*\<\ftt
    @>>\mkern31mu\theta_{\<2}'\!{\mkern-1.5mu}_{\phantom A}'\!\!^{-1}
       \mkern31mu >
     \iGp \U u_2^*\>g_{1\<*}^{\phantom{.}}\iGp{\V_1}\<\<v_1^*\<\ftt
               \\
\vspace{-22pt}
                \\
@V\beta_1VV  @V\beta_1VV  @VV\beta_1V
                \\
g_*\iGp\V\>  v_2^*g_1^{\<\times}\< u_1^*
 @>> \psi >
  \iGp\U \>g_* v_2^*g_1^{\<\times}\< u_1^*
    @>>\mkern32mu \theta'_{\<2}{}^{-1} \mkern32mu >
     \iGp \U u_2^*\>g_{1\<*}^{\phantom{.}}\>g_1^{\<\times}\< u_1^*
       @>>>
         \iGp \U u_2^*u_1^*
\end{CD}
$$
It suffices then that each one of the subrectangles
commute. 

\pagebreak[2]
For the three unlabeled subrectangles commutativity is clear.

As before, commutativity of subrectangle~(A) follows from that of the
first diagram in \Lref{L:transtheta}.

Commutativity of~(B) is easily checked after composition with 
the natural map
$\iGp\U \>g_*  (v_1v_2)^*\<\<\ftt\to g_*  (v_1v_2)^*\<\<\ftt\<\<$. (See the
characterization of~$\psi$ in~\Pref{uf=gv}(a).) 

Commutativity of~(C) results from that of the following diagram:
$$
\begin{CD}
g_{1\<*}^{}\<v_1^*\<\ftt 
 @=
   g_{1\<*}^{}\<v_1^*\<\ftt
     @< \theta_{\<\<f\<,\>u_1^{}} <<
       u_1^*f_{\<\<*}\ftt
         @>>>  u_1^*
             \\
@A  A \hbox{\hskip33pt{\footnotesize(D)}} A  
 @AAA   @AAA  @AAA
             \\
g_{1\<*}^{}\iGp{\V_1}\<\<v_1^*\<\ftt
  @> > \psi >
    \iGp{\U_1} \>g_{1\<*}^{}v_1^*\<\ftt
      @>>   \theta_{\<1}'{}^{-1} >
        \iGp{\U_1}\<\<u_1^*f_{\<\<*} \ftt      
          @>>>  \iGp{\U_1}\<\<u_1^*
              \\
@V \beta_1 VV  @. @. @V \hbox{\footnotesize (E)\hskip109pt} VV
              \\
g_{1\<*}g_1^{\<\times} u_1^*  
 \hbox to 0pt{\hskip113.5pt\Rarrow{205pt}\hss}  
    @. @.  @. u_1^*
\end{CD}
$$
Here subrectangle~(D) commutes by the characterization of~$\psi$
in~\Pref{uf=gv}(a); and (E) commutes by the very definition
of the base-change map~$\beta_1$. 

\smallskip

(b) As in (a), we consider the \emph{adjoint} diagram,  
essentially the outer border of the following diagram  (\ref{L:trans}.1).

(Note: The map 
$\psi\colon\< g_{1\<*}^{}\iGp\W \>w^*\<\<f_{\<2*}^{}f_{\<\<2}^{\<\times}\!
f_{\<\<1}^{\<\times}\< 
\to\iGp\U\>  g_{1\<*}^{} w^*\<\<f_{\<2*}^{}f_{\<\<2}^{\<\times}\! f_{\<\<1}^{\<\times} $\vspace{1pt} in the
middle of  diagram \ref{L:trans}.1 is defined because 
$f_{\<2*}^{}f_{\<\<2}^{\<\times}\! f_{\<\<1}^{\<\times} \set 
\R f_{2*}^{}(f_2^{})\<_{\textup t}^{\<\times}
\<(f_1^{})\<_{\textup t}^{\<\times}
\<\F\in\Dqc(\Z)$, by \mbox{\Pref{Rf-*(qct)}.})

For diagram \ref{L:trans}.1, commutativity of subrectangle (B)  (resp.~(D)) is given
by the definition of~$\beta_2$ (resp.~$\beta_1$.) Commutativity of~(C)
follows from that of the second diagram in \Lref{L:transtheta}.
Commutativity of~(A) is left as an exercise. (It is helpful to
compose with the natural map 
$\iGp\U\>g_{1\<*}^{}g_{2*}^{}v^*\< 
\<f_{\<\<2}^{\<\times}\! f_{\<\<1}^{\<\times} 
\to g_{1\<*}^{}g_{2*}^{}v^*\<\<f_{\<\<2}^{\<\times}\! f_{\<\<1}^{\<\times} $ and to use the characterization of
$\psi$ in~\Pref{uf=gv}(a).) The rest is straightforward.
\end{proof}
$$
\minCDarrowwidth=18pt
\begin{CD}
g_*\iGp\V\> v^*\<\<f_{\<\<2}^{\<\times}\! f_{\<\<1}^{\<\times} 
 @>\mkern30mu\Iso\mkern30mu>> 
   g_{1\<*}^{}g_{2*}^{}\iGp\V v^*\<\<f_{\<\<2}^{\<\times}\! f_{\<\<1}^{\<\times} 
    @. \kern-17pt \overset{\text{via }\beta_2}{\Rarrow{77pt}} @.
     \underset{\UnderElement{}{\downarrow}{7.3ex}{}} 
         {\hz{\hss\hskip-10pt$g_{1\<*}^{}g_{2*}^{}g_2^{\<\times}\! \iGp\W\>
          w^*\!f_{\<\<1}^{\<\<\times} \!$\hss}}
            \\
\vspace{-20pt}
             \\
@V\simeq VV @VV\psi V
             \\
g_*\iGp\V\> v^*\<\<\ftt 
 @. 
    g_{1\<*}^{}\iGp\W g_{2*}^{}v^*\<\<f_{\<\<2}^{\<\times}\! f_{\<\<1}^{\<\times} 
              \hbox to0pt{\hskip46pt{\footnotesize(B)}\hss}
              \\
@V\psi V \hbox{\hskip54pt{\footnotesize{(A)}}} V 
 @V \simeq V g_{1\<*}^{}(\theta_{\<2}'{}^{-1}) V
               \\
\iGp\U g_*v^*\<\<\ftt  
 @.
    g_{1\<*}^{}\iGp\W \>w^*\<f_{\<2*^{}}f_{\<\<2}^{\<\times}\! f_{\<\<1}^{\<\times}  
     @.\Rarrow{93.5pt} @.
      g_{1\<*}^{}\iGp\W \>w^*\<\<f_{\<\<1}^{\<\times} 
               \\
@V\simeq VV @VV\psi V @. @|
                \\
\vspace{-21pt}
             \\
\underset{\UnderElement{\simeq}{\downarrow}{7.3ex}
                                            {\hbox{\hskip51.7pt{\footnotesize{(C)}}}}}
               {\iGp\U\>g_{1\<*}^{}g_{2*}^{}v^*\<\<f_{\<\<2}^{\<\times}\! f_{\<\<1}^{\<\times} }
 @< \!\!\!\text{via }\theta_{\<\<f_2,w}(\ftt)\!\!\! <<
   \iGp\U\>  g_{1\<*}^{} w^*\<\<f_{\<2*}^{}f_{\<\<2}^{\<\times}\! f_{\<\<1}^{\<\times}  
    @>>>
      \iGp\U\>  g_{1\<*}^{} w^*\<\<f_{\<\<1}^{\<\times}  
      @<\psi<<
         \underset{\UnderElement{\hbox{\footnotesize{(D)}\hskip21pt}}
                                        {\downarrow}{7.3ex}{\!\<\<g_{1\<*}^{}\<(\beta_1\<)}}
                         {g_{1\<*}^{}\iGp\W \>w^*\<\<f_{\<\<1}^{\<\times} } 
              \\
\vspace{-20pt}
             \\
@.  @VV\hbox to0pt{$\scriptstyle \theta'_{f_1\<,u}(f_{\<2*}^{}\<\ftt\<)^{-1}$\hss}
V    @VV
\theta_{\<\<1}'{}^{\<-1} V
              \\ 
@.
  \iGp\U\> u^*\<\<f_{\<1*}^{}f_{\<2*}^{}f_{\<\<2}^{\<\times}\! f_{\<\<1}^{\<\times}  
    @>>>
       \iGp\U\> u^*\<\<f_{\<1*}^{}f_{\<\<1}^{\<\times}    
              \\
@.  @VV \simeq V @VVV 
              \\
\iGp\U\>g_*v^*\<\<\ftt
 @>> \mkern27mu\theta'{}^{-1}\mkern27mu>
   \iGp\U\> u^*\<\<f_{\!*}\ftt
     @>>>
       \iGp\U u^*
      @<<<
         g_{1\<*}^{}g_1^{\<\times}\<\< \iGp\U u^*
\end{CD}
$$
\bigskip
\centerline{\bf(\ref{L:trans}.1)}

\end{parag}

\pagebreak[2]

\begin{parag}
This subsection, proving \Lref{L:immbc}, is independent of the preceding one.

\begin{slem}
\label{L:immbc}
\Tref{T:basechange} holds when\/ $f$ is a closed immersion.
\end{slem}

\begin{proof}
The natural isomorphisms
$\R\iGp\V\>v^*\<\<f^!\R\iGp\Y\F\iso\R\iGp\V\>v^*\<\<f^!\F$ and 
$$
g^!u^*\>\R\iGp\Y\F 
\iso 
g^!\R\iGp\U\> u^*\>\R\iGp\Y\F 
\underset{\textup{\ref{P:f* and Gamma}(c)}}\iso
g^!\R\iGp\U u^*\F 
\iso 
g^! u^*\F
$$ 
(see \Cref{C:identities}(b)) 
let us replace~$\F$ by $\R\iGp\Y\F\<$, i.e., we may assume 
$\F\in\Dqct^+(\Y)$.

Recall from \Eref{ft-example}(4) that 
$\Rfs=f_{\!*}\colon\D(\X)\to\D(\Y)$
has a right adjoint~$f^\natural$ such that
$f^\natural(\Dqct^+(\Y))\subset\Dqct^+(\X)$; and that there is a natural
isomorphism
$$
j^{\>\X}_\G\colon \R\iGp\X f^\natural\G\iso\mathbf1^{\<!}\<f^\natural\G\cong f^!\G
\qquad\bigl(\G\in\Dqc^+(\Y)\bigr).
$$
The canonical map $f_{\!*} f^!\to\mathbf 1$ is the natural composition
$$
f_{\!*} f^!\underset{(j^\X)^{-1}}\iso f_{\!*}\R\iGp\X f^\natural\to
f_{\!*} f^\natural \to \mathbf 1.
$$
Similar remarks hold for $g$---also a closed immersion 
\cite[p.\,442, (10.14.5)(ii)]{GD}.

As in the proof of \Lref{L:closed}, the map
$
\theta_{\<\E}\colon u^*\<\<f_{\!*}\>\E\iso g_*v^*\E
$
of \Pref{uf=gv} is an isomorphism for all $\E\in\Dqct(\X)$. 
(Recall \Lref{C:vec-c is qc}.)  This being so, the base-change map~
$\beta_\F$ is easily seen to factor naturally as
$$
\R\iGp\V v^*\<\<f^!\F
\to g^!g_*\R\iGp\V v^*\<\< f^!\<\F
\to g^! g_* v^*\<\<f^!\<\F
\underset{\theta^{-1}}\iso g^! u^*\<f_{\!*}f^!\<\F
\to   g^! u^*\<\F.
$$
Also,  we can define the functorial map
$\beta_{\<\C}^\natural$
to be  the natural composition
$$
v^*\!f^\natural\C\to
g^\natural g_*v^*\!f^\natural\C\underset{\theta^{-1}}\iso
g^\natural u^*\<f_{\!*}f^\natural\C\to 
g^\natural u^*\<\C
\qquad\bigl(\C\in\Dqc(\Y)\bigr).
$$

The maps $\beta_{\<\F}^\natural$ and~$\beta_\F$ are related 
by commutativity of the following diagram, in which $\J$ is an ideal of
definition of $\Y$ (so that $\J\cO_\X$ is an ideal of definition of~$\X$):
$$
\begin{CD}
\R\iGp\V v^*\R\iGp\X f^\natural
@>\Iso>{\textup{\ref{P:f* and Gamma}(b)}}>
\R\iGp\V\>\R\iG{\J\cO_\V} v^*\<\< f^\natural
@>\Iso>>
   \R\iGp\V v^*\<\<f^\natural
     @>\R\iGp\V(\beta^\natural)>>
      \R\iGp\V \>g^\natural\> u^* \\
\vspace{-23pt}\\
@V\R\iGp\V v^*(j^{\>\X}) V\simeq V @. @. @V\simeq V j^\V V 
         \\
\R\iGp\V v^*\<\< f^! 
  @.{}
   @. \mkern-189mu\underset{\beta}{\Rarrow{226pt}} @.
     g^!u^*
\end{CD}
$$
(For the unlabeled isomorphism, see the beginning of the proof of
\Pref{uf=gv}.)  Since $\R\iGp\V$ is right-adjoint to the inclusion
$\Dqct(\V)\hookrightarrow\D(\V)$ (\Pref{Gamma'(qc)}), we can verify
this commutativity after composing with 
$g^!u^*\iso \R\iGp\V \>g^\natural u^*\to g^\natural u^*\<\<$,
at which point the verification is straightforward.

Thus to prove \Lref{L:immbc} we need only show that $\beta_\F^\natural$ is an
isomorphism,\vspace{-.6pt} i.e.\ (since $g$ is a closed immersion), that
$g_*(\beta_\F^\natural)$ is an isomorphism.\vspace{1.6pt} 

For that purpose, consider the unique  functorial map
$$
\sigma=\sigma(\E,\>\G)\colon u^*\R\sHomb_\Y(\E,\>\G)\to
\R\sHomb_\U(u^*\E,u^*\<\G)
\quad\ \: \bigl(\E\in\Dc^-(\Y),\ \G\in\D^+(\Y)\bigr)
$$
which for bounded-below injective complexes~$\G$ is the natural composition
$$
u^*\R\sHomb_\Y(\E,\>\G)\cong u^*\sHomb_\Y(\E,\>\G)\to
\sHomb_\U(u^*\E,u^*\G)\to\R\sHomb_\U(u^*\E,u^*\G).
$$
This map is an \emph{isomorphism}. Indeed, it commutes with localization, so we
need only check for affine~$\Y$, and then, since every coherent $\cO_\Y$-module
is a homomorphic image of a finite-rank free one (\cite[p.\,427, (10.10.2)]{GD}), a
standard way-out argument reduces the problem to the trivial case $\E=\cO_\Y$.

Take $\E\set f_{\!*}\cO_\X={}$(say)\:$\cO_\Y/\I$. The source and target
of~$\sigma(\cO_\Y/\I,\>\F\>)$ are
\begin{gather*}
u^*\R\sHomb(\cO_\Y/\I,\>\F\>)=
u^*\<\<f_{\!*}f^\natural\F \cong g_*v^*\!f^\natural\<\F, \\
\R\sHomb(u^*(\cO_\Y/\I),\>u^*\F\>)= g_*g^\natural u^*\<\<\F .
\end{gather*}
Let $\cK$ be a K-injective $\cO_\U$-complex quasi-isomorphic to $u^*\<\F\<$.
Since the complexes $u^*\sHomb_\Y(\cO_\Y/\I,\>\F\>)$ and  
$\sHomb_\U(u^*\cO_\Y/\I,\>\cK\>)\cong\R\sHomb_\U(u^*\cO_\Y/\I,u^*\F\>)$
are both annihilated\- by~$\I\cO_\U$,  we see that the isomorphism
$\sigma(\cO_\Y/\I,\>\F\>)$ is isomorphic to a map of the form $g_*(\varsigma)$
where
$
\varsigma\colon v^*\<\<f^\natural\F\to g^\natural u^*\<\<\F
$
is a map in~$\D(\V)$. It suffices then to show that
$\varsigma=\beta_\F^\natural\>$, i.e.~(by definition of~$\beta_\F^\natural$),
that the natural composition
$$
u^*\<\<f_{\!*}f^\natural\F\iso
g_*v^*\!f^\natural\F\xrightarrow{g_*(\varsigma)} g_*g^\natural u^*\<\<\F
\xrightarrow{\tau_{\<\<u^{\mkern-1.5mu*}\!\F}^\natural} u^*\<\F
$$
is induced by the natural map\index{ {}$\tau$ (trace map)!$\tau^\natural$}
$$
\tau_{\<\F}^\natural\colon
f_{\!*}f^\natural\F=\R\sHomb(\cO_\Y/\I,\>\F\>)\to
\R\sHomb(\cO_\Y,\>\F\>)=\F.
$$

From \Eref{ft-example}(4) one sees, for injective~$\F\<$,\vspace{-1pt}
that 
$\tau_{\<\F}^\natural$ takes any homomorphism $\varphi\colon \cO_\Y/\I\to\F$ over an
open subset of~$\Y$ to $\varphi(1)$; and similarly for
$\tau_{\<\<u^{\mkern-1.5mu*}\!\F}^\natural\>$. The conclusion follows 
from the above definition of $\sigma(\cO_\Y/\I,\>\F\>)=g_*(\varsigma)$.
\end{proof}

\end{parag}

\begin{parag}\label{reduction}
In this subsection we prove \Tref{T:basechange}
in case $f\colon\X\to\Y$ is a proper map of
ordinary noetherian schemes,
by reduction to the case where $\X$, $\Y$, $\U$ and $\V$ are
\emph{all} ordinary schemes---a case which we take for granted (see
the introductory remarks for section~\ref{sec-basechange}).
Of course when $u$ is \emph{adic} then $\U$ is
already an ordinary scheme, and no reduction is needed at all.

\begin{slem}\label{L:ordbc}  
Let\/ $f\colon\X\to\Y$ be a proper map of
ordinary noetherian schemes. For 
\Tref{T:basechange} to hold with this\/~$f$ it suffices that it hold
whenever\/ $\U$ and\/ $\V$ are ordinary schemes as well.
\end{slem}

\begin{proof}

Without yet assuming that $\X$ and $\Y$ are ordinary schemes, we can reduce
\Tref{T:basechange} to the special case where the formal scheme~$\U$ is
\emph{affine} and $u(\U)$ is contained in an affine open subset of~$\Y$. Indeed,
for the base-change map
$\beta_\F=\beta_{\<\<f,\>u}(\F\>)$  of \Tref{T:basechange} to be an
isomorphism, it clearly suffices that for any composition of fiber squares
$$
\begin{CD}
\V_0@>v_0^{}>>\V@>v>>\X \\
@Vg_0^{}VV @VgVV  @VVfV \\
\U_0@>>\vbox to 0pt{\vskip-1ex\hbox{$\scriptstyle u_0^{}$}\vss}>\U
  @>>\vbox to 0pt{\vskip-1ex\hbox{$\scriptstyle u$}\vss}>\Y
\end{CD}
$$
with $u_0$ the inclusion of an affine open $\U_0\subset\U$ such that $u(\U_0)$ is
contained in an affine open subset of~$\Y$, the map
$$
v_0^*(\beta_\F\>)\colon v_0^*\R\iGp\V v^*\<\<f^!\<\F\to v_0^*g^!u^*\<\F
$$
be an isomorphism. \Rref{R:Dtilde}(6) yields that $\F\in\wDqcp(\Y)\Rightarrow
u^*\F\in\wDqcp(\U)$.\vspace{.5pt} So if  we assume the above-specified special
case, then $\beta_{\<\<f\<,\>uu_0^{}\<}(\F\>)$\vspace{.5pt}
and~$\beta_{\<g,\>u_0^{}\<}(u^*\<\F\>)$ are both isomorphisms. 
From \Pref{Gamma'(qc)}(a) we have a natural
isomorphism  
$$
v_0^*(\beta_\F\>)\cong \R\iGp{\V_0}\<
v_0^*(\beta_{\<\<f,\>u}(\F\>)),
$$
 so \Lref{L:trans}(a) shows that
$v_0^*(\beta_\F\>)$ is in fact an isomorphism.

With reference to the remarks just preceding \Sref{S:trans},  (a) and (b) having
already been proved, only (c) remains, i.e., we need only prove \Tref{T:basechange}
for the rear face of diagram~\eqref{cube}. 

In other words, with the notation of diagram~\eqref{cube},
we may assume in proving \Tref{T:basechange} that $f=f_n$ (a proper map of
ordinary  schemes), and that $u=u_n$.
Moreover $Y_n$ is a closed subscheme of $\Y$, and so if $\U$ is
affine and $u(\U)$ is contained in an affine open subset of~$\Y$, then
$Y_n\times_\Y\U$ is
affine and $u_n(Y_n\times_\Y\U)$ is contained in an affine open subset of~$Y_n$.
It follows that $Y_n\times_\Y\U$ is the completion of an ordinary
affine $Y_n$-scheme. (That can be seen via the one-one correspondence from
maps between affine formal schemes to continuous homomorphisms between their
associated rings \cite[p.\,407, (10.4.6)]{GD}). 
\Tref{T:basechange} is thus reduced to the case depicted in the following diagram,
where
$f\colon\X\to \Y$ is now a proper map of ordinary noetherian schemes, $U$ is
an ordinary affine $\Y$-scheme, $\kappa\colon\U\to U$ is a completion map, and
$u\colon\U\to\Y$ factors as shown.

$$
\begin{CD}
\X\times_\Y\U@>>>\X\times_\Y U@>>>\X \\
@VgV\mkern63mu\text{\footnotesize(1)}V
@VVV@V\text{\footnotesize(2)}\mkern45mu VfV
\\
\U@>>\vbox to 0pt{\vskip-1.1ex\hbox{$\scriptstyle \kappa$}\vss}>U@>>>\Y
\end{CD}
$$

We will show that \Tref{T:basechange} holds for subdiagram (1)
by identifying the base-change map
associated to  $\kappa$
with the \emph{isomorphism}~$\zeta$ in~\Cref{C:compln+basechange}.
As subdiagram (2) is a fiber square of ordinary schemes, \Lref{L:ordbc}
will then result from the preceding reduction and the transitivity
\Lref{L:trans}(a).

\smallskip

\def\fot{f_{\mkern-1.5mu0}^{\<\times}}
It is convenient to re-represent subdiagram~(1) in the 
notation of~\Cref{C:compln+basechange}. 
Consider then a diagram
$$
\begin{CD}
\X@.:=X_{\</Z} @>\kappa_\X^\pd>> X  \\
@V f VV @. @VV f_0^{} V \\
\Y@.:=Y_{/W}@>>\vbox to
0pt{\vskip-1.1ex\hbox{$\scriptstyle\kappa_\Y^\pd$}\vss}> Y
\end{CD}
$$
as in \Cref{C:kappa-f^times-tors}, with $Z=f_0^{-1}W\<$.
That $\zeta$ \emph{is} the base-change map means
 that $\zeta$ is adjoint to the natural composition
$$
\Rfs\R\iGp\X\kappa_\X^*\fot\underset{\psi}\iso 
\R\iGp\Y\Rfs\kappa_\X^*\fot\underset{\theta'{}^{-1}}\iso
\R\iGp\Y\kappa_\Y^*\R f_{\<0*}^{}\fot\lra
\R\iGp\Y\kappa_\Y^*\lra \kappa_\Y^*.
$$
But by definition, $\zeta$ is adjoint to the natural composition
$$
\Rfs\R\iGp\X\kappa_\X^*\fot\underset{\textup{~\ref{Gammas'+kappas}(c)}}\iso
\Rfs\kappa_\X^*\R\iG Z\fot\lra
\Rfs\kappa_\X^*\R\iG Z\fot\<\kappa_{\Y^*}^{}\kappa_\Y^*\\
\underset{\tau_{\textup t}'(\kappa_\Y^*)}\lra\kappa_\Y^*
$$
with $\tau_{\textup t}'$ as in~\Cref{C:kappa-f^times-tors}---so that 
$\tau_{\textup t}'(\kappa_\Y^*)$ factors naturally as
\begin{align*}
\Rfs\kappa_\X^*\R\iG Z\fot\<\kappa_{\Y^*}^{}\kappa_\Y^*
&\underset{\textup{~\ref{C:kappa-f*t}}}\iso
\kappa_\Y^*\R f_{\<0*}^{}\R\iG Z\fot\<\kappa_{\Y^*}^{}\kappa_\Y^*\\
&\lra
\kappa_\Y^*\R f_{\<0*}^{}\fot\<\kappa_{\Y^*}^{}\kappa_\Y^*\\
&\lra
\kappa_\Y^*\kappa_{\Y^*}^{}\kappa_\Y^*\\
&\xrightarrow{\,\>\,\pi\,\,\>} \kappa_\Y^*\<.
\end{align*}
It will suffice then to verify that the following natural diagram commutes
(where, again, we omit all occurrences of $\R$):
$$
\minCDarrowwidth=25pt
\def\1{f_{\!*}\iGp\X\kappa_\X^*\fot}
\def\2{\iGp\Y\< f_{\!*}\kappa_\X^*\fot}
\def\3{\iGp\Y\kappa_\Y^* f_{\<0*}^{}\fot}
\def\4{\iGp\Y\kappa_\Y^*}
\def\5{\kappa_\Y^* f_{\<0*}^{}\fot}
\def\6{\kappa_\Y^* f_{\<0*}^{}\iG Z\fot\<\kappa_{\Y^*}^{}\kappa_\Y^*}
\def\7{\kappa_\Y^* f_{\<0*}^{}\fot\<\kappa_{\Y^*}^{}\kappa_\Y^*}
\def\8{\kappa_\Y^*\kappa_{\Y^*}^{}\kappa_\Y^*}
\def\9{\kappa_\Y^*}
\def\ten{f_{\!*}\kappa_\X^*\iG Z\fot\<\kappa_{\Y^*}^{}\kappa_\Y^*}
\def\lvn{f_{\!*}\kappa_\X^*\iG Z\fot}
\def\twv{\kappa_\Y^* f_{\<0*}^{}\iG Z\fot}
\def\thrn{f_{\!*}\kappa_\X^*\fot}
\begin{CD}
\1 @>\psi>> \2 @>\theta'{}^{-1}>> \3 @>>>\4  \\
\vspace{-19pt}\\
@V\textup{\ref{Gammas'+kappas}(c)}VV
 \text{\footnotesize(A)} @. @VVV @VVV \\
\vspace{-22pt}\\
\lvn 
 @>\Iso>\textup{~\ref{C:kappa-f*t}}> 
  \twv @>>> 
    \5 @>>>
\underset{\UnderElement{}{\downarrow}{2.8ex}
{\!\!\<\raisebox{.4ex}{$\scriptstyle\iota$}}}{\mkern1mu\strut}
\mkern-20mu\hz{$\mkern-8mu\9$\hss} \\
\vspace{-18.7pt}\\
 @VVV @VVV @VVV @A\pi AA\\
\vspace{-20pt}\\
\ten @>\Iso>\textup{~\ref{C:kappa-f*t}}> \6 @>>>\7 @>>>\8\\
\vspace{-12pt}
\end{CD}
$$
Given that $\pi\iota=1$, the verification of commutativity is
straightforward,
except for subrectangle (A).

\enlargethispage{-.7\baselineskip}
Now there is a functorial isomorphism 
$\alpha\colon \R f_{0*}\R\iG Z\iso \R\iG W \R f_{0*}$
which arises 
in the obvious way, via ``K-flabby'' resolutions, from the equality 
$f_{0*}\iG Z=\iG W f_{0*}$ 
(see the last paragraph in the Remark following (3.2.5) 
in~\cite[p.\,25]{AJL}), and whose composition with the natural map
$\R\iG W \R f_{0*}\to\R f_{0*}$ is the natural map 
$\R f_{0*}\R\iG Z\to\R f_{0*}$. And, again, we have the isomorphism 
$\R\iGp\Y\kappa_\Y^*\iso\kappa_\Y^*\R\iG W$ of
\Pref{Gammas'+kappas}(c), whose composition with the natural map
$\kappa_\Y^*\R\iG W\to\kappa_\Y^*$ is the natural map 
$\R\iGp\Y\kappa_\Y^*\to\kappa_\Y^*$. Hence commutativity of (A) follows
from that of the outer border---consisting
entirely of isomorphisms---of the following diagram:

$$
\def\1{f_{\!*}\iGp\X\kappa_\X^*}
\def\2{\iGp\Y \<f_{\!*}\kappa_\X^*}
\def\3{\iGp\Y\kappa_\Y^* f_{\<0*}^{}}
\def\4{\kappa_\Y^* \iG W\< f_{\<0*}^{}}
\def\5{\kappa_\Y^* f_{\<0*}^{}}
\def\lvn{f_{\!*}\kappa_\X^*\iG Z}
\def\twv{\kappa_\Y^* f_{\<0*}^{}\iG Z}
\def\thrn{f_{\!*}\kappa_\X^*}
\def\se{\rotatebox{-30}{\hbox to 41pt{\rightarrowfill}}}
\def\ne{\rotatebox{30}{\hbox to 41pt{\rightarrowfill}}}
\def\sw{\rotatebox{30}{\hbox to 41pt{\leftarrowfill}}}
\def\nw{\rotatebox{-30}{\hbox to 41pt{\leftarrowfill}}}
\begin{CD}
\underset{\UnderElement{\textup{~\ref{Gammas'+kappas}(c)}}{\uparrow}{7.2ex}
{\<\<\simeq}}\1
@>\psi>\vbox to0pt{\vskip2pt\hz{\hss\se\hss}\vskip9pt
  \hz{\hss\ne\hss}\vss}> \2
@.\overset{\theta'}{\Larrow{100pt}} @.
\underset{\UnderElement{\simeq}{\uparrow}{7.2ex}
{\!\!\!\textup{~\ref{Gammas'+kappas}(c)}}}
\3
\\
\vspace{-19pt}\\
@. @VVV  \\
\vspace{-22pt}\\
@. \thrn @<<\theta< \5\\
\vspace{-21pt}\\
@. @. @AAA\\
\lvn @. \overset{\vbox to 0pt{\vss
\hbox to 0pt{\hss$\Iso$\hss}\vskip-6.5pt
    \hbox to
0pt{\hss$\scriptstyle\textup{~\ref{C:kappa-f*t}}$\hss}
 \vskip-12pt\vss}
  }
                {\Rarrow{100pt}} @.
  \twv @>
\vbox to0pt{\vss\hz{\hss\sw\hss}\vskip-6pt
  \hz{\hss\nw\hss}\vskip 2pt\hbox to 0pt{\hss$\Iso$\hss}}
      >\alpha> \4\text{\large\strut}
\end{CD}
$$

\medskip
Since $\R\iGp\Y$ is right-adjoint to the inclusion 
$\Dt(\Y)\hookrightarrow \D(\Y)$ (\Pref{Gamma'(qc)}), we can check
commutativity \emph{after} composing the outer border with the natural
map $\R\iGp\Y f_{\!*}\kappa_\X^*\to f_{\!*}\kappa_\X^*\>$, so that it
suffices to check commutativity of all the subdiagrams of the
preceding one. This is easy to do, as, with $\E\set\fot\<\F$, the maps denoted by~
$\theta_{\<\E}\; (=\theta_{\!f_{\<\halfsize{\footnotesize0}}^{},
  \mkern1.5mu\kappa_\Y^{}\!\<}(\E)$) 
 in \Cref{C:kappa-f*t} and in \Pref{uf=gv} are the same.\vspace{1pt}

This completes the proof of \Lref{L:ordbc}, and of \Tref{T:basechange}.
\end{proof}

\end{parag}

\section{Consequences of the flat base change isomorphism.}
\label{Consequences}\index{base-change isomorphism}

We begin with a flat-base-change theorem for the functor
$\ush f=\BL_\X f^!$ associated to a pseudo\kern.6pt-proper map~$f\colon\X\to Y$
of noetherian formal schemes.
(As before,\vspace{.3pt} $f^!\set\ft$, and $\ush f$ is
right-adjoint to the functor
$\Rfs\R\iGp\X\colon\wDqc(\X)\to\D(\Y)$,\vspace{.5pt} where
$\wDqc(\X)$ is the (full) $\Delta$-subcategory of~$\D(\X)$ such
that\vspace{.25pt}
$$
\F\in\wDqc(\X)\Leftrightarrow\R\iGp\X\F\in\Dqct(\X),
$$ 
see \Cref{C:f*gam-duality}.)

We deduce a sheafified version \Tref{T:sheafify} of \Tref{Th2} of the
Introduction (=\:\Tref{T:qct-duality} + \Cref{C:f*gam-duality}). This
is readily seen equivalent to the case of flat base change where
$u\colon\U\to\Y$ is an open immersion; in other words, it expresses
the local nature, over~$\Y$, of $f^!$ and $\ush f\<$.

\Sref{S:coherent} establishes the local nature of $f^!$ and $\ush f$ over~$\X$.
From this we obtain that  
$\ush f\bigl(\Dc^+(\Y)\bigr)\subset \Dc^+(\X)$ (\Pref{P:coherence}). 
This leads further to an improved base-change theorem for bounded-below
complexes with coherent homology, and to \Tref{T:properdual}, a duality
theorem for such complexes under proper maps.

\medskip
We consider as in \Tref{T:basechange} a fiber square of noetherian
formal schemes 
$$
\begin{CD}
\V@>v>>\X \\
@VgVV @VVfV \\
\U@>>\vbox to 0pt{\vskip-1ex\hbox{$\scriptstyle u$}\vss}>\Y
\end{CD}
$$
with $f$ and $g$ pseudo\kern.6pt-proper,  $u$ and $v$ flat. 

For any $\F\in\wDqcp(\Y)$ we have the composed isomorphism
$$
\vartheta\colon\R\iGp\V\>v^*\<\<\ush f\<\F
\underset{\textup{\ref{P:f* and Gamma}(c)}}\iso
\R\iGp\V\>v^*\>\R\iGp\X\ush f\<\F
\underset{\textup{\ref{C:identities}(a)}}\iso
\R\iGp\V\>v^*\<\<f^!\F
\underset{\textup{\ref{T:basechange}}}\iso
g^!u^*\<\F.
$$
In particular, $v^*\<\<\ush f\<\F\in\wDqc(\V)$.

\pagebreak[3]
\stepcounter{numb}
\renewcommand{\theequation}{\theparag.\arabic{numb}}
\begin{thm}\label{T:sharp-basechange}
Under the preceding conditions, let 
$$
\ush{\beta_{\<\<\F}}\colon v^*\<\<\ush f\<\F\to\ush g
u^*\<\F\qquad\bigl(\F\in\wDqcp(\Y)\bigr)
$$
be the map adjoint to the natural composition
\begin{equation}\label{adjointto}
\R g_*\R\iGp\V\>v^*\<\<\ush f\<\F
\underset{\R g_*\vartheta}\iso\R g_*g^!u^*\<\F\to u^*\F.
\end{equation}
Then the map $\BL_{\<\V}(\ush{\beta_{\<\<\F}})$ is an \emph{isomorphism}
$$
\BL_{\<\V}(\ush{\beta_{\<\<\F}})\colon\BL_{\<\V}\> v^*\<\<\ush f\<\F\iso
\BL_{\<\V}\>\ush g u^*\<\F
\underset{\textup{\ref{C:identities}(a)}}\cong
\ush g u^*\<\F.
$$
Moreover,
if\/ $u$ is an open immersion then\/ $\ush{\beta_{\<\<\F}}$ itself is an isomorphism.
\end{thm}
\stepcounter{numb}
\renewcommand{\theequation}{\theparag.\arabic{numb}}
\begin{proof}
The map $\ush{\beta}$ factors naturally as
\begin{equation}\label{factors}
v^*\<\<\ush f \to \BL_{\<\V}\>v^*\<\<\ush f
\underset{\textup{\ref{R:Gamma-Lambda}(c)}}\iso
\BL_{\<\V}\R\iGp\V\>v^*\<\<\ush f 
\underset{\BL_{\<\V}\vartheta}\iso
\BL_{\<\V} g^!u^*=\ush g u^*\<.
\end{equation}
To see this, one needs to check that \eqref{factors} is adjoint to
\eqref{adjointto}. The natural map $\mathbf 1\to\BL_\V$ factors naturally as 
$\mathbf1\to\BL_\V\R\iGp\V\to\BL_\V$ (easy
check),  and hence the adjointness in question amounts to the readily-verified
commutativity of the outer border of the following diagram (with all occurrences
of $\R$ left out):
$$
\def\1{g_*\iGp\V\>v^*\<\<\ush f}
\def\2{g_*g^!u^*}
\def\5{g_*\iGp\V\BL_{\<\V}\>v^*\<\<\ush f}
\def\6{g_*\iGp\V\BL_{\<\V}\iGp\V\>v^*\<\<\ush f}
\def\7{g_*\iGp\V\BL_{\<\V} g^!u^*}
\begin{CD}
\6@<<<\1@>\Iso >\text{via\;}\vartheta>\2  \\
\vspace{-22pt}\\
@V\simeq VV @AAA @AAA  \\
\5 @>\Iso >> \6 @>\Iso >\text{via\;}\vartheta> \7
\end{CD}
$$

That $\BL_{\<\V}(\ush{\beta_{\<\<\F}})$ is an isomorphism results then from the
idempotence of~$\BL_{\<\V}$ \mbox{(\Rref{R:Gamma-Lambda}(b)).}

When $u$---hence $v$---is an open immersion, we have isomorphisms (the first of
which is obvious): 
$$
\BL_{\<\V}\>v^*\<\<\ush f\iso v^*\<\BL_{\<\X}\ush f
\underset{\textup{\ref{C:identities}(a)}}\iso
 v^*\<\<\ush f,
$$
and the last assertion follows.
\end{proof}

\medskip

Next comes the sheafification of \Tref{Th2}.  Let $f\colon\X\to\Y$ be a map of
locally noetherian formal schemes.  For $\G$ and $\E\in\D(\X)$ we have natural
compositions
$$
\Rfs\R\sHomb_\X\<(\G\<,\E)
\to
\Rfs\R\sHomb_\X\<(\bL f^*\Rfs\G\<,\E)\,
\xrightarrow[\mkern-15mu\textup{\cite[\kern-1pt p.\kern1pt147,
\kern-1pt 6.7]{Sp}}\mkern-15mu]{\Iso}
\,\R\sHomb_{\>\Y}\<(\Rfs\G\<,\Rfs\E)
$$
and
$$
\R\sHomb_\X\<(\G\<,\E)
\to
\R\sHomb_\X\<(\R\iGp\X\>\G\<,\E) 
\xrightarrow[\textup{\ref{C:Hom-Rgamma}}]{\Iso}
\R\sHomb_\X\<(\R\iGp\X\>\G\<,\R\iGp\X\E).
$$
\begin{thm}\label{T:sheafify}\index{Grothendieck Duality!Torsion
(sheafified)}
Let\/ $\X$ and\/ $\Y$ be noetherian formal schemes and let\/
$f\colon\X\to\Y$ be a pseudo\kern.6pt-proper map. Then the following
natural compositions are
\emph{isomorphisms:}
\begin{align*}
\ \ush\delta\<\<\colon\<\Rfs\R\sHomb_{\>\X}\<(\G\<,\>\ush f\<\F\>)
&\to\<
\R\sHomb_{\>\Y}\<(\Rfs\R\iGp\X\>\G\<,\>\Rfs\R\iGp\X\ush f\<\F\>) \\
&\to
\R\sHomb_{\>\Y}\<(\Rfs\R\iGp\X\>\G\<,\>\F\>)\vspace{-4pt}
\quad\ 
\bigl(\G\in\wDqc(\X),\;\F\>\in\wDqcp(\Y)\bigr);
\end{align*}
\vspace{-6pt}

\noindent
$
\ \delta^!\<\<\colon\Rfs\R\sHomb_{\>\X}\<(\G\<,\>f^!\F\>)
\to
\R\sHomb_{\>\Y}\<(\Rfs\G\<,\>\Rfs f^!\F\>)
\to
\R\sHomb_{\>\Y}\<(\Rfs\G\<,\>\F\>)\vspace{3pt}
$

\rightline{$\bigl(\G\in\Dqct(\X),\;\F\>\in\wDqcp(\Y)\bigr).$}
\end{thm}

\begin{proof}
The map $\ush\delta$ is an isomorphism iff the same is true of
$\R\Gamma(\U,\ush\delta)$ for all open $\U\subset\Y$. (For if $\E$---which
may be assumed K-injective---is the vertex of a triangle based
on~$\ush\delta\<$,  then $\ush\delta$ is an isomorphism 
$\Leftrightarrow\E\cong 0\Leftrightarrow H^i(\E)=0$ for all 
$i\in\mathbb Z\Leftrightarrow$ the sheaf associated to the presheaf
$\U\mapsto \textup H^i\Gamma(\U,\E)=\textup H^i\R\Gamma(\U,\E)$ vanishes for
all~$i$.) Set $\V\set f^{-1}\U$, and let $u\colon\U\hookrightarrow\Y$ and
$v\colon\V\hookrightarrow\X$ be the respective inclusions.
We have then the fiber square
$$
\begin{CD}
\V@>v>>\X \\
@VgVV @VVfV \\
\U@>>\vbox to 0pt{\vskip-1ex\hbox{$\scriptstyle u$}\vss}>\Y\hbox to 0pt{,\hss}
\end{CD}
$$
and need only verify that $\R\Gamma(\U,\ush\delta)$ is the
composition of the following sequence of isomorphisms:
\begin{flalign*}
\R\Gamma\bigl(\U,\>\Rfs\R\sHomb_{\>\X}\<(\G\<,\>\ush f\<\F\>)\bigr)
&\<\iso\< \R\Gamma\bigl(\V,\>\R\sHomb_{\>\X}\<(\G\<,\>\ush f\<\F\>)\bigr)
\quad&&\!\!\!\!
\textup{\cite[\!6.4,\:6.7,\:5.15]{Sp}}\\
&\<\iso\< \R\Homb_{\>\V}(v^*\<\G\<,v^*\<\<\ush f\<\F\>)
\quad&&\!\!\!\!
\textup{\cite[\!5.14,\:5.12,\:6.4]{Sp}}\\
&\<\iso\< \R\Homb_{\>\V}(v^*\<\G\<,\>\ush g u^*\F\>)
\quad&&\!\!\!\!(\textup{\Tref{T:sharp-basechange}})\\
&\<\iso\< \R\Homb_{\>\U}\<(\R g_*\R\iGp\V\> v^*\<\G\<,\>u^*\F\>)
\  &&\ \ \!\textup{\bigl(\kern-1pt\ref{C:f*gam-duality},\!
             \ref{R:Dtilde}(6)\kern-1pt\bigr)}\\ 
&\<\iso\< \R\Homb_{\>\U}\<(\R
g_*v^*\R\iGp\X\>\G\<,\>u^*\F\>)
\ &&\quad\ \ \textup{(elementary)}\\
&\<\iso\< \R\Homb_{\>\U}\<(u^*\Rfs\R\iGp\X\>\G\<,\>u^*\F\>)
\ &&\quad\ \ \textup{(elementary)}\\
&\<\iso\< 
\R\Gamma\bigl(\U,\>\R\sHomb_{\>\Y}\<(\Rfs\R\iGp\X\>\G\<,\>\F\>)\bigr)
\hbox to0pt{\qquad\;\textup{(as above).}\hss}
\hskip-14pt
\end{flalign*}
This somewhat tedious verification is left to the reader (who may e.g., refer
to the proof of $(4.3)^{\textup o}\Rightarrow (4.2)$ near the end of \cite{Non
noetherian}).

That $\delta^!$ is an isomorphism can be shown similarly---or 
be deduced  via the natural map $f^!\cong\R\iGp\X \ush f\to \ush f$
(\Cref{C:identities}), which for $\G\in\Dqct(\X)$ induces an isomorphism 
$\R\sHomb_{\>\X}\<(\G\<,\>f^!\F)\iso
\R\sHomb_{\>\X}\<(\G\<,\>\ush f\F)$ (\Cref{C:Hom-Rgamma}).
\end{proof}

\pagebreak

\begin{parag}\label{S:coherent}
For pseudo\kern.6pt-proper~$f\colon\X\to\Y$, the functors $f^!\set \ft$ and
$\ush f$ are \emph{local on~$\X$}, in the following sense.

\begin{sprop}\label{P:local}
Let there be given a commutative diagram
$$
\begin{CD}
\U@>i_1>>\X_1 \\
@V i_2 VV @VV f_1 V \\
\X_2 @>> f_2 > \Y
\end{CD}
$$
of noetherian formal schemes, with\/ $f_1$ and\/~$f_2$ pseudo\kern.6pt-proper and 
$i_1$ and\/~$i_2$ open immersions. Then there are functorial isomorphisms
$$
i_1^*f_1^! \iso i_2^*f_2^!\>,\qquad\quad
i_1^*\ush{f_1} \iso i_2^*\ush{f_2}.
$$
\end{sprop}

\begin{proof}
The second isomorphism results from the first, since for any
$\F\>\in\D(\Y)$ and for $j=1,2$,
\begin{align*}
i_j^*\ush{f_j}\F\>\overset{\text{\ref{C:f*gam-duality}}}=
i_j^*\R\sHomb_{\X_j}\<\<(\R\iGp{\X_j}\<\cO_{\X_j}, f_j^!\F\>)
&\cong\R\sHomb_{\>\U}\<(i_j^*\R\iGp{\X_j}\<\cO_{\X_j}, i_j^*f_j^!\F\>)\\
&\cong
\R\sHomb_{\>\U}\<(\R\iGp{\U}\cO_\U, i_j^*f_j^!\F\>).
\end{align*}

For the first isomorphism, Verdier's\index{Verdier, Jean-Louis} proof of 
\cite[p.\,395, Corollary 1]{f!}---a special case of
\Pref{P:local}---applies verbatim, modulo the following
extensions (a), (b) and~(c) of some elementary properties of schemes
to formal schemes.\looseness=-1

(a) Since pseudo\kern.6pt-proper maps are separated, the graph of~$i_j$ is a
\emph{closed immersion} 
$\gamma\colon \U\hookrightarrow \X_j\<\<\times_{\<\Y}\<\U$
(see \cite[p.\,445, (10.15.4)]{GD}, where the ``finite-type''
hypothesis is used only to ensure that $\X_j\<\times_{\<\Y}\U$ is locally
noetherian, a condition which holds here by the first paragraph in
\Sref{sec-basechange}. And if $\>\U\to\Y$ is an open immersion, then so
is~$\gamma$ (since then both $\pi_j\colon\X_j\!\times_{\<\Y}\<\<\U\to\X_j$ and
$i_j=\pi_j\gamma$ are open immersions).

(b) If $s\colon\U\to\V$ is an open and closed immersion, then the
exact functors~$s_*$ and~$s^*$ are adjoint, and by
\Eref{ft-example}(4) there is a functorial isomorphism
$$
s^!\F\cong s^\natural\F\cong s^*\F\qquad\bigl(\F\in\Dqct(\V)\bigr).
$$

(c) (Formal extension of \cite[p.\,325, (6.10.6)]{GD}.) Let
$\U\overset{\gamma}\hookrightarrow
\W
\overset{w}\hookrightarrow\Z$ be maps of locally noetherian formal
schemes such that $\gamma$~is a closed immersion and $w$ is an open
immersion. (We are interested specifically in the case
$\W\set\X_2\<\times_{\<\Y}\<\U$ and $\Z\set\X_2\<\times_{\<\Y}\<\X_1$,
see (a).) Set $u\set w\gamma$.  Then \emph{the closure\/~$\overline\U$
of\/~$u(\U)$ is a formal subscheme of\/~$\Z$, and the map\/
$\U\to\overline\U$ induced by\/~$u$ is an open immersion.}\vspace{1pt}

Indeed, $\overline\U$ is the support of $\cO_\Z/\I$ where
$\I$ is the kernel of the natural map $\cO_\Z\to u_*\cO_\U\>$; and it
follows from \cite[p.\,441, (10.14.1)]{GD} that we need only show that
$\I$ is \emph{coherent}.
The question being local, we may assume that $\Z$ is affine, say
$\Z=\Spf(A)$. Cover $\U$ by a finite number of affine open 
subschemes~$\U_i\ (1\le i\le n)$, with inclusions
$u_i\colon\U_i\hookrightarrow\U$. Then there is a natural injection
$$
u_*\cO_\U\hookrightarrow
u_*\<\bigl(\!\oplus_{i=1}^nu_{i*}\cO_{\U_i}\bigr)\cong
\oplus_{i=1}^n(uu_i)_*\cO_{\U_i}\>,
$$
so that $\I$ is the intersection of the kernels of the natural maps
$\cO_\Z\to(uu_i)_*\cO_{\U_i}$, giving us a reduction to the case where
$\U$ itself is affine, say $\U=\Spf(B)$. Now if $I$ is the kernel of
the ring-homomorphism $\rho\colon A\to B$ corresponding to~$u$, then for any
$f\in A$ the kernel of the induced map 
$\rho_{\{f\}}\colon A_{\{f\}} \to B_{\{f\}} $ is
$I_{\{f\}}$;\vspace{1pt}
and one deduces that $\I$ is the coherent $\cO_\Z$-module denoted by
$I^\Delta$ in \cite[p.\,427, (10.10.2)]{GD}.
\end{proof}

\begin{sprop}\label{P:coherence}
If\/ $f\colon\X\to\Y$ is a pseudo\kern.6pt-proper map of noetherian formal schemes then
$$
\ush f\bigl(\Dc^+(\Y)\bigr)\subset \Dc^+(\X).
$$
\end{sprop}

 \begin{proof}
Since $\ush f$ commutes with open base change (\Tref{T:sharp-basechange})
we may assume $\Y$ to be affine, say $\Y=\Spf(A)$. Since $f$
is of pseudo-finite type,  every point  of
$\X$ has an open neighborhood~$\U$ such that $f|_\U$ factors as 
$$
\U\overset{i}\hookrightarrow\Spf(B)\xrightarrow{\,p}\Spf(A)=\Y
$$
where $B$ is the completion of a polynomial ring
$P\set A[T_0, T_1,\dots,T_n]$ with respect\- to an ideal~$I$ whose intersection
with
$A$ is open, $i$~is a closed immersion, and $p$~corresponds to the obvious
continuous ring homomorphism $A\to B$ (see footnote in \Sref{maptypes}). This
$\Spf(B)$ is an open subscheme of the completion~$\mathscr P$ of the projective
space ~$\mathbf P_{\!\!\!A}^n$ along the closure of its subscheme
$\Spec(P/I)$. Thus by \Pref{P:local} and
item~(c) in its proof, we can replace $\X$ by a closed formal subscheme of~
$\mathscr P$ having $\U$ as an open subscheme. In other words, we may assume
that $f$ factors as
$
\X\overset{i_1^{}\>}\hookrightarrow\mathscr P\xrightarrow{p_1^{}\>}\Spf(A)=\Y
$
with $i_1$ a closed immersion and $p_1$ the natural map. Then $\ush
f=\ush{i_1}\ush{p_1}$, and we need only consider the two cases (a)  $f=p_1$
and (b) $f=i_1$.

\penalty-1000
Case (a) is given by \Cref{C:completion-proper}.  In case (b) we see
as in example~\ref{ft-example}(4) that for $\F\in\Dc^+(\Y)$ we have
$f^\natural\F\in\Dc^+(\X)$ and
$$
\ush f\<\F=\BL_\X\R\iGp\X f^\natural\<\F\underset
{\text{or
}\ref{C:Hom-Rgamma}}{\overset{\ref{R:Gamma-Lambda}\text{(c)}}{=\!=\!=}}\BL_\X
f^\natural\<\F\overset{\ref{formal-GM}}{=\!=}f^\natural\<\F\in\Dc^+(\X).
\vspace{-4.3ex}
$$
\end{proof}

\begin{scor}\label{C:coh-basechange}\index{base-change isomorphism}
\smallskip
For all\/ $\F\in\Dc^+(\Y)$ the base-change map\/ $\ush{\beta_{\<\<\F}}$ of
\Tref{T:sharp-basechange}  is an \emph{isomorphism}
$$
\ush{\beta_{\<\<\F}}\colon v^*\<\<\ush f\<\F
\iso
\ush g u^*\<\F.
$$
\end{scor}

\begin{proof}
\Pref{formal-GM} gives an isomorphism
$v^*\<\<\ush f\<\F\iso\BL_{\<\V}\> v^*\<\<\ush f\<\F$ .
\end{proof}
\end{parag}
\medskip

We have now the following duality theorem for proper
maps and bounded-below complexes with coherent homology.

\begin{thm}\label{T:properdual}\index{Grothendieck Duality!coherent}
Let\/ $f\colon\X\to\Y$ be a proper map of noetherian formal schemes, so that\/
$\Rfs\(\Dc^+(\X))\subset\Dc^+(\Y)$  and\/ 
$\ush f\bigl(\Dc^+(\Y)\bigr)\subset \Dc^+(\X)$
$($see Propositions~\textup{\ref{P:proper f*}} and~\textup{\ref{P:coherence}).}
Then for\/ $\G\in\Dc^+(\X)$ and $\F\in\Dc^+(\Y)$ there are functorial
\emph{isomorphisms}
\begin{align*}
\Rfs\R\sHomb(\G\<,\ush f\<\F\>)
&\underset{\textup{\ref{T:sheafify}}}\iso
\R\sHomb(\Rfs\R\iGp\X\>\G\<,\>\F\>) \\
&\underset{\textup{\ref{C:f* and Gamma}(d)}}\iso
\R\sHomb(\R\iGp\Y\>\Rfs\G\<,\>\F\>)
\underset{\textup{\ref{formal-GM}}}\iso
\R\sHomb(\Rfs\G\<,\>\F\>).
\end{align*}
In particular, 
$\ush f\colon\Dc^+(\Y)\to \Dc^+(\X)$ is right-adjoint
to\/ $\Rfs\colon\Dc^+(\X)\to\Dc^+(\Y)$.\vspace{1.5pt}

 If\/ $\X$ is properly algebraic we can replace\/
$\ush f$ by the functor\/ $f^{\<\times}$ 
of \Cref{cor-prop-duality}.
\end{thm}

{\sc Proof}
Left to reader.  (For the last assertion see
Corol\-laries~\ref{C:coh-dual} and~\ref{corollary}.)

\enlargethispage*{\baselineskip}


\def\seename{see}
\index{Duality!Torsion|see{Grothendieck Duality}}
\index{Duality!Grothendieck|see{Grothendieck Duality}}
\index{Duality!Greenlees-May|see{Greenlees-May Duality}}
\index{way out|see {boundedness}}
\index{flat base change|see {base-change map, base-change isomorphism}}
\index{Duality!Coherent|see{Grothendieck Duality}}



\begin{thebibliography}{ABC}
\bibitem[AJL$'$]{AJL$'$} \bysame, Greenlees-May duality on formal schemes, this
volume.
\end{thebibliography}
\end{document}